\documentclass[11pt]{article}
\usepackage{amssymb}
\usepackage{amsmath}
\usepackage{times}
\usepackage{graphicx}
\usepackage[all]{xy}
\usepackage{xcolor}

\title{Continuous 2-colorings and topological dynamics\indent}
\author{to appear in Dissertationes Mathematicae\\ \\ Dominique LECOMTE}
\date{\today}

\def\ufootnote#1{\let\savedthfn\thefootnote\let\thefootnote\relax
\footnote{#1}\let\thefootnote\savedthfn\addtocounter{footnote}{-1}}

\newcommand{\ca}{{\bf\Pi}^{1}_{1}}

\newcommand{\boraone}{{\bf\Sigma}^{0}_{1}}
\newcommand{\boratwo}{{\bf\Sigma}^{0}_{2}}

\newcommand{\borone}{{\bf\Delta}^{0}_{1}}
\newcommand{\bortwo}{{\bf\Delta}^{0}_{2}}

\newcommand{\bormone}{{\bf\Pi}^{0}_{1}}

\newcommand{\bormtwo}{{\bf\Pi}^{0}_{2}}

\newcommand{\borxi}{{\bf\Delta}^{0}_{\xi}}

\newtheorem{thm} {Theorem} [section]
\newtheorem{defi} [thm] {Definition}
\newtheorem{cor} [thm] {Corollary}
\newtheorem{lem} [thm] {Lemma}
\newtheorem{prop} [thm] {Proposition}

\newtheorem{them} {Theorem} [subsection]

\newtheorem{coro} [them] {Corollary}
\newtheorem{lemm} [them] {Lemma}
\newtheorem{propo} [them] {Proposition}

\begin{document}

\maketitle

\centerline{$\bullet$ 1) Sorbonne Universit\' e, CNRS, Institut de Math\'ematiques de Jussieu-Paris Rive Gauche,}

\centerline{IMJ-PRG, F-75005 Paris, France}

\centerline{2) Universit\'e de Paris, IMJ-PRG, F-75013 Paris, France}

\centerline{dominique.lecomte@upmc.fr}\medskip

\centerline{$\bullet$ Universit\'e de Picardie, I.U.T. de l'Oise, site de Creil,}

\centerline{13, all\'ee de la fa\"\i encerie, 60100 Creil, France}\medskip\medskip\medskip\medskip\medskip\medskip

\ufootnote{{\it 2020 Mathematics Subject Classification.}~Primary: 03E15, Secondary: 54H05, 37B05, 37B10}

\ufootnote{{\it Keywords and phrases.}~antichain, basis, coloring, compact, conjugacy, continuous chromatic number, continuous homomorphism, graph, minimal dynamical system, odometer, subshift, well-founded, zero-dimensional}

\noindent\emph{Abstract.} We first consider the class $\mathfrak{K}$ of graphs on a zero-dimensional metrizable compact space with continuous chromatic number at least three. We provide a concrete basis of size continuum for $\mathfrak{K}$ made up of countable graphs, comparing them with the quasi-order $\preceq^i_c$ associated with injective continuous homomorphisms. We prove that the size of such a basis is sharp, using odometers. However, using odometers again, we prove that there is no antichain basis in $\mathfrak{K}$, and provide infinite descending chains in $\mathfrak{K}$. Our method implies that the equivalence relation of flip conjugacy of minimal homeomorphisms of 
$2^\omega$ is Borel reducible to the equivalence relation associated with $\preceq^i_c$. We also prove that there is no antichain basis in the class of graphs on a zero-dimensional Polish space with continuous chromatic number at least three. We study the graphs induced by a continuous function, and show that any $\preceq^i_c$-basis for the class of graphs induced by a homeomorphism of a zero-dimensional metrizable compact space with continuous chromatic number at least three must have size continuum, using odometers or subshifts. 

\vfill\eject\baselineskip=13.2pt

\section{$\!\!\!\!\!\!$ Introduction}\indent

 The present article is the continuation of the study of definable colorings initiated in [K-S-T], and continued in [L-Z1] and [C-M-Sc-V1]. All our relations will be binary. The motivation for this work goes back to the following so called $\mathcal{G}_0$-dichotomy, essentially proved in [K-S-T]. 
  
\begin{thm} (Kechris, Solecki, Todor\v cevi\' c) \label{G0} There is a Borel relation $\mathcal{G}_0$ on $2^\omega$ such that, for any Polish space $X$, and for any analytic relation $R$ on $X$, exactly one of the following holds:\smallskip

\noindent (1) there is $c\! :\! X\!\rightarrow\!\omega$ Borel such that $c(x)\!\not=\! c(y)$ if $(x,y)\!\in\! R$ (a \emph{countable Borel coloring} of $R$),\smallskip

\noindent (2) there is $\varphi\! :\! 2^\omega\!\rightarrow\! X$ continuous such that 
$\mathcal{G}_0\!\subseteq\! (\varphi\!\times\!\varphi )^{-1}(R)$.\end{thm}

 If (1) holds, then we say that $R$ has countable Borel chromatic number (a relation $R$ on a set $X$ is a \emph{digraph} if it does not meet the \emph{diagonal} $\Delta (X)\! :=\!\{ (x,x)\mid x\!\in\! X\}$ of $X$; the \emph{Borel chromatic number} $\chi_B(X,R)$ of a digraph $(X,R)$ is the minimum cardinal $\kappa\!\leq\!\aleph_0$ for which there is a Borel coloring of $R$ taking values in 
$\kappa$ (equipped with the discrete topology) if it exists, $2^{\aleph_0}$ otherwise). If (2) holds, then we say that $\varphi$ is a \emph{continuous homomorphism} from $(2^\omega ,\mathcal{G}_0)$ into $(X,R)$, and denote this by 
$(2^\omega ,\mathcal{G}_0)\preceq_c(X,R)$. This result had a lot of developments since. We refer to [K-Ma] for a survey, and to [B0], [B1] and [G-J-Kr-Se] for recent work in continuous combinatorics, which is the topic of the present work. It is natural to ask for a level by level version of Theorem \ref{G0}, with respect to the Borel hierarchy (see the introduction in [K]). This work was initiated in [L-Z1], where the authors prove the following.

\begin{thm} \label{LZ} (Lecomte, Zelen\'y) Let $\xi\!\in\!\{ 1,2,3\}$. Then we can find a zero-dimensional Polish space 
$\mathbb{X}_\xi$, and a Borel relation $\mathbb{R}_\xi$ on $\mathbb{X}_\xi$ such that for any (zero-dimensional if $\xi\! =\! 1$) Polish space $X$, and for any analytic relation $R$ on $X$, exactly one of the following holds:\smallskip  

\noindent (1) there is a countable $\borxi$-measurable coloring of $R$,\smallskip  

\noindent (2) there is $\varphi\! :\!\mathbb{X}_\xi\!\rightarrow\! X$ continuous such that 
$\mathbb{R}_\xi\!\subseteq\! (\varphi\!\times\!\varphi )^{-1}(R)$.\end{thm}

  [C-M-So, Theorem 4.4] gives a version of this for analytic spaces when $\xi\! =\! 2$, and this is also possible when $\xi\! =\! 1$. More recently, the existence of versions of Theorem \ref{G0} for finite Borel colorings was decided. In [T-V], the authors rule out the most straightforward analogs of the $\mathcal{G}_0$-dichotomy for graphs of Borel chromatic number at least $\kappa$, where $4\!\leq\!\kappa\!\leq\!\omega$ (recall that a \emph{graph} is a symmetric digraph). The difficult remaining case has been solved in [C-M-Sc-V1], where the authors introduce a Borel graph $\mathbb{L}_0$ on a zero-dimensional Polish space 
$\mathbb{X}_0$ satisfying the following.

\begin{thm} (Carroy, Miller, Schrittesser, Vidny\'anszky) \label{L0} Let $X$ be a Polish space, and $G$ be an analytic graph on $X$. Exactly one of the following holds:\smallskip

\noindent (1) $G$ has Borel chromatic number at most two,\smallskip

\noindent (2) there is $\varphi\! :\!\mathbb{X}_0\!\rightarrow\! X$ continuous such that 
$\mathbb{L}_0\!\subseteq\! (\varphi\!\times\!\varphi )^{-1}(G)$.\end{thm}

 All this leads to the following question.\medskip
  
\noindent\emph{Question 1.}\ Fix a countable ordinal $\xi\!\geq\! 1$. Is there a Borel graph $\mathbb{G}_\xi$ on a zero-dimensional analytic space $\mathbb{A}_\xi$ which is $\preceq_c$-minimum among analytic graphs on a (zero-dimensional if 
$\xi\! =\! 1$) analytic space with $\borxi$ chromatic number at least three? We will also consider metrizable separable, Polish, compact, and finite spaces.

\vfill\eject\baselineskip=12.9pt

 Since the very beginning of the study of definable chromatic numbers in [K-S-T], injective definable homomorphisms were considered (see also [K-Ma, Sections 4 and 8], [L, Theorem 10], [L-Za, Theorem 1.13], and [L-Z2]). So it is natural to ask the same question with injective continuous homomorphisms instead of continuous homomorphisms (with the notation 
$\preceq^i_c$ instead of $\preceq_c$). In [C-M-Sc-V2], the authors announce the existence of a continuum sized family of closed graphs on a Polish space with Borel chromatic number at least three which are pairwise 
$\preceq^i_B$-incompatible in the class of analytic graphs on a Hausdorff space with Borel chromatic number at least three.\medskip

 We consider the quasi-orders $\preceq^i_c$ and $\preceq_c$ on various classes (a \emph{quasi-order} is a reflexive transitive relation). Let $\Gamma$ be a class, and $\leq$ be a quasi-order on $\Gamma$. A subclass $\mathfrak{B}$ of 
$\Gamma$ is a \emph{basis} for $\Gamma$ if any element $\Gamma$ is $\leq$-above an element of $\mathfrak{B}$. We are looking for basis as small as possible, for the inclusion. In other words, we want the elements of $\mathfrak{B}$ to be pairwise $\leq$-incomparable. A subclass $\mathfrak{B}$ satisfying this property is called an \emph{antichain}. So we are looking for antichain basis, when they exist. In the best case, the antichain basis is a singleton $\{ b\}$, and we say that $b$ is 
\emph{minimum} among elements of $\Gamma$. This is the case in Theorems \ref{G0}, \ref{LZ} and \ref{L0}, but it is not always possible. We are interested in the following questions, very natural when we study a quasi-order.\medskip
 
\noindent (1) Is there an antichain basis?\smallskip

\noindent (2) If there is no antichain basis, is there a reasonably simple basis?\smallskip

\noindent (3) What is the minimal size of a basis?\smallskip

\noindent (4) Are there big antichains?\smallskip

\noindent (5) Are there infinite descending chains?\smallskip

\noindent (6) Can we find minimal elements?\smallskip

\noindent (7) Can we embed a complex quasi-order?\medskip

 In this article, our spaces will be metrizable separable, except in Theorem \ref{eq1++} and its two lemmas. As above, the 
\emph{continuous chromatic number} (\emph{CCN} for short) $\chi_c(X,R)$ of a digraph $(X,R)$ is the minimum cardinal 
$\kappa\!\leq\!\aleph_0$ for which there is a continuous coloring of $R$ taking values in $\kappa$ 
(equipped with the discrete topology) if it exists, $2^{\aleph_0}$ otherwise. We mainly focus on continuous $2$-colorings, even if some other cardinalities will be considered. The case of continuous $2$-colorings is much more complex than in Theorem \ref{LZ} for $\xi\! =\! 1$ and $\omega$ colors, the latter case corresponding directly to the definition of the product topology. In Section 2, we will see that the odd cycles $(2p\! +\! 3,C_{2p+3})$ are witnesses for the fact that any $\preceq^i_c$-basis for the class of graphs on a zero-dimensional metrizable separable (or Polish, or metrizable compact, or finite) space (\emph{0DMS}, 
\emph{0DP}, \emph{0DMC} for short) with CCN at least three must be infinite.  In the compact case, our main results are as follows. Let $\mathfrak{K}$ be the class of graphs on a 0DMC space with CCN at least three.
 
\begin{thm} \label{corcomp''''''} We can find a concrete family $\big( (\mathbb{K}_\alpha ,\mathbb{G}_\alpha )\big)_{\alpha\in 2^\omega}$, where $\mathbb{K}_\alpha$ is a compact subset of $2^\omega$ and $\mathbb{G}_\alpha$ is a countable graph on 
$\mathbb{K}_\alpha$, such that, for any 0DMC space $X$ and any graph $G$ on $X$, exactly one of the following holds:\smallskip

\noindent (1) $G$ has CCN at most two,\smallskip

\noindent (2) we can find $\alpha\!\in\! 2^\omega$ and $\varphi\! :\!\mathbb{K}_\alpha\!\rightarrow\! X$ injective continuous such that $\mathbb{G}_\alpha\!\subseteq\! (\varphi\!\times\!\varphi )^{-1}(G)$.\smallskip

 In other words, $\big( (\mathbb{K}_\alpha ,\mathbb{G}_\alpha )\big)_{\alpha\in 2^\omega}$ is a $\preceq_c^i$-basis (and thus a $\preceq_c$-basis) for $\mathfrak{K}$.\end{thm}

 Recall that an \emph{oriented graph} is an antisymmetric digraph. Theorem \ref{corcomp''''''} and most of our results admit versions for digraphs and oriented graphs. We will come back to this in the last section. It is simpler to work with graphs in Theorem \ref{corcomp''''''}. Note also that in [C-M-Sc-V1], the authors prove that there is no version of Theorem \ref{L0} for oriented graphs.
 
\vfill\eject\medskip\baselineskip=13.2pt

 Recall that $D_2(\bormone )$ is the class of differences of two closed sets, while
$$\boraone\oplus\bormone\! :=\!\{ (O\cap C)\cup (F\!\setminus\! C)\mid C\!\in\!\borone\wedge O\!\in\!\boraone\wedge F\!\in\!\bormone\}$$ 
is the self dual class just after $\borone$ in the Wadge hierarchy of Borel sets (see [K, 22.B and 22.E]). We provide another concrete $\preceq_c$-basis for $\mathfrak{K}$, which is not a $\preceq^i_c$-basis, but is made up of countable $D_2(\bormone )$ graphs with $\boraone\oplus\bormone$ (and thus $\bortwo$ and Borel) chromatic number two, and whose vertices have degree at most one. We will see in Section \ref{Polish} that our basis are not $\preceq_c$-basis for the class of countable graphs on a 0DP space with CCN at least three.\medskip
 
  We next prove that the size of such $\preceq^i_c$-basis is sharp. In order to prove this, we use minimal Cantor dynamical systems. These systems have been widely studied (see, for example, [I-Me], [Ka], [Ku], [Lo], [Me], [P], [Sa-T\"o]). A 
\emph{dynamical system} $(X,f)$ is given by a homeomorphism $f$ of a metrizable compact space $X$. If $X$ is homeomorphic to $2^\omega$, then we say that $(X,f)$ is a \emph{Cantor dynamical system}. A dynamical system (or $f$) is \emph{minimal} if $\mbox{Orb}_f(x)\! :=\!\{ f^i(x)\mid i\!\in\!\mathbb{Z}\}$ is dense in $X$ for each $x\!\in\! X$. If $(Y,g)$ is another dynamical system, we say that these systems (or $f,g$) are 
\emph{orbit}-\emph{equivalent} if there is a homeomorphism $\varphi\! :\! X\!\rightarrow\! Y$ such that 
$\varphi [\mbox{Orb}_f(x)]\! =\!\mbox{Orb}_g\big(\varphi (x)\big)$ for any $x\!\in\! X$. It was known that there is a family of size continuum made up of minimal Cantor dynamical systems which are pairwise not orbit equivalent (see [I-Me]). We consider a property stronger than orbit equivalence, namely flip-conjugacy. We say that two dynamical systems $(X,f),(Y,g)$ (or $f,g$) are \emph{conjugate} (resp., 
\emph{flip}-\emph{conjugate}) if there is a homeomorphism $\varphi\! :\! X\!\rightarrow\! Y$ such that 
$\varphi\!\circ\! f\! =\! g\!\circ\!\varphi$ (resp., $\varphi\!\circ\! f\! =\! g\!\circ\!\varphi$ or 
$\varphi\!\circ\! f\! =\! g^{-1}\!\circ\!\varphi$). We provide a family of size continuum made up of minimal Cantor dynamical systems (in fact odometers) which are pairwise not flip conjugate, and associate to each homeomorphism of this family a graph on a 0DMC space, ensuring the following properties. 

\begin{thm} \label{eantichmin} There is a $\preceq_c$-antichain (and thus $\preceq_c^i$-antichain) 
$\big( (\mathbb{K}_\alpha ,\mathbb{G}_\alpha )\big)_{\alpha\in 2^\omega}$, where\smallskip
 
\noindent (a) $\mathbb{K}_\alpha$ is a 0DMC space,\smallskip

\noindent (b) $\mathbb{G}_\alpha$ is a countable $D_2(\bormone )$ graph on $\mathbb{K}_\alpha$ with CCN three and 
$\boraone\oplus\bormone$ chromatic number two, and whose vertices have degree at most one,\smallskip

\noindent (c) $(\mathbb{K}_\alpha ,\mathbb{G}_\alpha )$ is $\preceq_c^i$-minimal in $\mathfrak{K}$.\smallskip

 In particular, any $\preceq_c^i$-basis for $\mathfrak{K}$ must have size at least continuum.\end{thm}
 
 The minimal examples are particularily important, since they have to be part of any basis, up to equivalence. Note that Theorem \ref{eantichmin} shows that if 
$(\mathbb{A}_1,\mathbb{G}_1)$ exists, then we must have $\chi_c(\mathbb{A}_1,\mathbb{G}_1)\! =\! 3$, 
$\chi_B(\mathbb{A}_1,\mathbb{G}_1)\! =\! 2$, and $(\mathbb{A}_1,\mathbb{G}_1)$ must be strictly 
$\preceq_c$-below $(\mathbb{X}_0,\mathbb{L}_0)$. We will also see in Section \ref{Polish} that $\mathbb{A}_1$ cannot be compact. Theorem \ref{eantichmin} shows that our quasi-orders have large antichains. Moreover, they are not well-founded.

\begin{thm} \label{infdecrcompact} There is a $\preceq_c$ and $\preceq^i_c$-descending chain  
$\big( (\mathbb{K},\mathbb{G}_p)\big)_{p\in\omega}$, where\smallskip
 
\noindent (a) $\mathbb{K}$ is a 0DMC space,\smallskip

\noindent (b) $\mathbb{G}_p$ is a countable $D_2(\bormone )$ graph on $\mathbb{K}$ with CCN three and 
$\boraone\oplus\bormone$ chromatic number two, and whose vertices have degree at most one.\end{thm}
 
 Theorems \ref{eantichmin} and \ref{infdecrcompact} contrast with [C] where it is proved that the closed subsets of a zero-dimensional Polish space are well-quasi-ordered by bi-continuous embeddability (so this quasi-order has finite antichains and descending chains).
 
\vfill\eject
 
 We now give a countable $D_2(\bormone )$ graph on a 0DP space which is not compact.\medskip
  
\noindent\emph{Notation.}\ Recall that if $R$ is a relation on a set $X$ and $l\!\in\!\omega$, then
$$R^l\! :=\!\{ (x,y)\!\in\! X^2\mid\exists (x_i)_{i\leq l}\!\in\! X^{l+1}~~\forall i\! <\! l~~
(x_i,x_{i+1})\!\in\! R\wedge (x,y)\! =\! (x_0,x_l)\}\mbox{,}$$
$R^{-1}\! :=\!\{ (x,y)\!\in\! X^2\mid (y,x)\!\in\! R\}$, and $s(R)\! :=\! R\cup R^{-1}$ is the \emph{symmetrization} of $R$. We now define our graph, on the copy $\mathcal{N}\! :=\! (\{ c,a,\overline{a}\}\cup\omega )^\omega$ of the Baire space 
$\omega^\omega$. We set\medskip

\leftline{$\mathbb{O}_m\! :=\! 
\{ (c^{k+1}a^{j+1}\overline{a}^\infty ,k0^{j+1}\overline{a}^\infty )\mid j,k\!\in\!\omega\}\cup
\{ (ki^{j+1}a^\infty ,k(i\! +\! 1)^{j+1}\overline{a}^\infty )\mid j,k\!\in\!\omega\wedge i\!\leq\! 2k\} ~\cup$}\smallskip

\rightline{$\{ (k(2k\! +\! 1)^{j+1}a^\infty ,c^{k+1}\overline{a}^{j+1}a^\infty )\mid j,k\!\in\!\omega\} .$}\medskip

\noindent and $\mathbb{G}_m\! :=\! s(\mathbb{O}_m)$. The idea of this example is to decompose the graph in levels indexed by $k$, and that the level $k$ is an approximation of the odd cycle on $2k\! +\! 3$ points, the approximation being improved when $k$ increases. Note that the vertices of $\mathbb{G}_m$ have degree at most one.
 
\begin{thm} \label{eq1++} Let $X$ be a first countable topological space, and $G$ be a graph on $X$. The following are equivalent:\smallskip

\noindent (1) $\Delta (X)\cap\overline{\bigcup_{p\in\omega}~\overline{G}^{2p+1}}$ is not empty,\smallskip

\noindent (2) $(\mathcal{N},\mathbb{G}_m)\preceq_c(X,G)$.\smallskip

 Moreover, these conditions imply that $G$ has no continuous 2-coloring.\end{thm}

 The countable $G_\delta$ subset 
$\mathbb{P}\! :=\!\overline{\mbox{proj}[\mathbb{G}_m]}\!\setminus\!\big\{ c^{k+1}\varepsilon^\infty\mid 
k\!\in\!\omega\wedge\varepsilon\!\in\!\{ a,\overline{a}\}\big\}$ of $\mathcal{N}$ has the properties that 
$(\mathbb{P},\mathbb{G}_m)$ has CCN three and 
$(\mathbb{P},\mathbb{G}_m)\prec_c(\mathcal{N},\mathbb{G}_m)$. The next picture describes the level of $\mathbb{G}_m$ corresponding to $k\! =\! 0$, seen in $\mathbb{P}$ (so that the sequences $(c^{k+1}\varepsilon^{j+1}\overline{\varepsilon}^\infty )_{j\in\omega}$ become discrete).\medskip
 
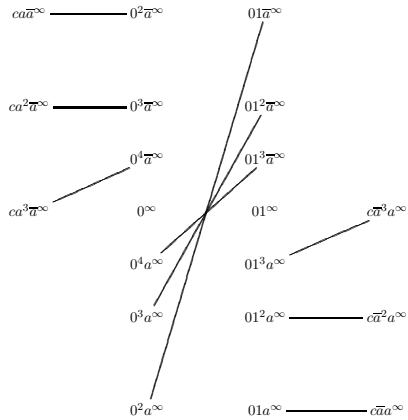
\centerline{\scalebox{0.51}{$$~~~~~\xymatrix{
& & & & 
ca\overline{a}^\infty \ar@{-}[rr] \ar@<0.1mm>@{-}[rr] \ar@<-0.1mm>@{-}[rr] & & 
0^2\overline{a}^\infty & & 
01\overline{a}^\infty & & & & & & & \\ \\ 
& & & & 
ca^2\overline{a}^\infty \ar@{-}[rr] \ar@<0.1mm>@{-}[rr] \ar@<-0.1mm>@{-}[rr] & & 
0^3\overline{a}^\infty & & 
01^2\overline{a}^\infty & & & & & & & \\ 
& & & & & & 
0^4\overline{a}^\infty & & 
01^3\overline{a}^\infty & & & & & & & \\ 
& & & & 
ca^3\overline{a}^\infty \ar@{-}[urr] \ar@<0.1mm>@{-}[urr] \ar@<-0.1mm>@{-}[urr] & & 
\mbox{$0^\infty$} & & 
\mbox{$01^\infty$} & & 
c\overline{a}^3a^\infty & & & & & \\ 
& & & & & & 
0^4a^\infty \ar@{-}[uurr] \ar@<0.1mm>@{-}[uurr] \ar@<-0.1mm>@{-}[uurr] & & 
01^3a^\infty \ar@{-}[urr] \ar@<0.1mm>@{-}[urr] \ar@<-0.1mm>@{-}[urr] & & & & & & & \\ 
& & & & & & 
0^3a^\infty \ar@{-}[uuuurr] \ar@<0.1mm>@{-}[uuuurr] \ar@<-0.1mm>@{-}[uuuurr] & & 
01^2a^\infty \ar@{-}[rr] \ar@<0.1mm>@{-}[rr] \ar@<-0.1mm>@{-}[rr] & & 
c\overline{a}^2a^\infty & & & & & \\ \\ 
& & & & & & 
0^2a^\infty \ar@{-}[uuuuuuuurr] \ar@<0.1mm>@{-}[uuuuuuuurr] \ar@<-0.1mm>@{-}[uuuuuuuurr] & & 
01a^\infty \ar@{-}[rr] \ar@<0.1mm>@{-}[rr] \ar@<-0.1mm>@{-}[rr] & & 
c\overline{a}a^\infty & & & & & 
}$$}}\medskip

 Using subgraphs of $(\mathbb{P},\mathbb{G}_m)$, we prove the following.

\begin{thm} \label{absmin} Let $G$ be a graph on a 0DMS space $Z$, satisfying 
$(Z,G)\preceq^i_c(\mathbb{P},\mathbb{G}_m)$ and having CCN at least three. Then there is a family   
$\big( (P_\alpha ,G_\alpha )\big)_{\alpha\in 2^\omega}$ of graphs on a 0DP space with CCN three, 
$\preceq^i_c$-below $(Z,G)$, and pairwise $\preceq^i_c$-incompatible in the class of graphs on a 0DMS space with CCN at least three.\smallskip

 In particular, there is no $\preceq^i_c$-antichain basis in the class of graphs on a 0DMS (or 0DP) space with CCN at least three, and any $\preceq^i_c$-basis for one of these classes must have size at least continuum.\end{thm}
 
 Theorem \ref{absmin} shows that $(\mathbb{P},\mathbb{G}_m)$ is not $\preceq_c^i$-minimal among graphs on a 0DP (or 0DMS) space with CCN at least three. One can prove that this still holds for $\preceq_c$, but we will not do it here. Theorem \ref{absmin} also gives a negative answer to the version of Question 1 for $\preceq_c^i$. One can prove that no subgraph of one of the examples of Theorem \ref{eantichmin} is $\preceq^i_c$-minimal among graphs on a 0DP space with CCN at least three, but we will not do it here. We also prove a version of Theorem \ref{absmin} in the compact case.
 
\begin{thm} \label{absmincomp} We can find a countable graph $(3^\omega ,\mathbb{G})$ in $\mathfrak{K}$ such that, for each $(K,G)$ in $\mathfrak{K}$ satisfying $(K,G)\preceq^i_c(3^\omega ,\mathbb{G})$, there is a $\preceq^i_c$-antichain 
$\big( (3^\omega ,G_\alpha )\big)_{\alpha\in 2^\omega}$ of graphs with CCN three and 
$\preceq^i_c$-below $(K,G)$. In particular, there is no $\preceq^i_c$-antichain basis in $\mathfrak{K}$.\end{thm}

 We now stated our main results concerning general graphs. The case of graphs induced by a function has been particularily considered in [K-S-T], and also in [Co-M], [Pe] and [TV] for instance. Also, we give at the end of this article a table summarizing the properties of our two quasi-orders for graphs. It is remarkable that the same properties hold for graphs induced by a partial homeomorphism with countable domain, up to, possibly, the existence of the $\preceq^i_c$-concrete basis. If 
$$f\! :\!\mbox{Domain}(f)\!\subseteq\! X\!\rightarrow\!\mbox{Range}(f)\!\subseteq\! X$$ 
is a partial function, then the graph \emph{induced} by $f$ is $G_f\! :=\! s\big(\textup{Graph}(f)\big)\!\setminus\!\Delta (X)$. In 
[K-S-T], it is proved that if $X$ is a standard Borel space and $f$ is a Borel function on $X$ (i.e., has a Borel graph), then the Borel chromatic number of $G_f$ is in $\{ 0,1,2,3,\aleph_0\}$. So it is natural to ask the following.\medskip

\noindent\emph{Question 2.}~Let $X$ be a 0DMS space, and $f\! :\! X\!\rightarrow\! X$ be a partial continuous function with analytic domain. What are the possible values for $\chi_c(X,G_f)$?\medskip
 
 We prove the following.
 
\begin{thm} \label{main} Let $X$ be a 0DMC space, and $f\! :\! X\!\rightarrow\! X$ be a partial continuous injection.\smallskip

\noindent (a) If the domain of $f$ is  open, then $\chi_c(X,G_f)\!\in\!\omega\cup\{ 2^{\aleph_0}\}$, and all these values are possible with fixed point free partial homeomorphisms on a countable space.\smallskip

\noindent (b) If the domain of $f$ is  closed, then $\chi_c(X,G_f)\!\in\!\{ 0,1,2,3,2^{\aleph_0}\}$, and all these values are possible with (total) homeomorphisms of a countable space.\smallskip

 Moreover, we can find a countable Polish space $X$ and a fixed point free partial homeomorphism $f\! :\! X\!\rightarrow\! X$ with open domain and $\chi_c(X,G_f)\! =\!\aleph_0$.\end{thm}
 
  Our method in the proof of Theorem \ref{absmin} shows that there is no $\preceq^i_c$-antichain basis for the class of graphs induced by a partial homeomorphism on a 0DP space with CCN at least three. Also, the method used to prove Theorem \ref{eantichmin} shows that any $\preceq^i_c$-basis for the class of graphs induced by a partial fixed point free continuous involution with countable open domain on a 0DMC space with CCN at least three must have size continuum.\medskip
 
 Theorem \ref{main}(b) leads to consider, for $\kappa\!\leq\! 3$, the class $\mathfrak{G}_\kappa$ of graphs induced by a (total)  homeomorphism of a 0DMC space with CCN strictly bigger than $\kappa$. We will see that in Theorem \ref{LZ} we can take $\mathbb{X}_1\! :=\!\{ 0^\infty\}\cup\{ 0^n1^\infty\mid n\!\in\!\omega\}$ and 
$\mathbb{R}_1\! :=\!\{ (0^{2n}1^\infty ,0^{2n+1}1^\infty )\mid n\!\in\!\omega\}$. Note that the graph $s(\mathbb{R}_1)$ is $G_{f_1}$, where $f_1\! :\!\mathbb{X}_1\!\rightarrow\!\mathbb{X}_1$ is the total homeomorphism defined by $f_1(0^\infty )\! :=\! 0^\infty$ and 
${f_1(0^{2n+\varepsilon}1^\infty )\! :=\! 0^{2n+1-\varepsilon}1^\infty}$. We prove the following. 

\begin{thm} \label{Ckappa} (a) $(1,\emptyset )$ is $\preceq^i_c$-minimum in $\mathfrak{G}_0$.\smallskip

\noindent (b) $(2,G_{\varepsilon\mapsto 1-\varepsilon})$ is $\preceq^i_c$-minimum in $\mathfrak{G}_1$.\smallskip

\noindent (c) Any $\preceq_c^i$-basis for $\mathfrak{G}_2$ must have size continuum.\smallskip

\noindent (d) $(\mathbb{X}_1,G_{f_1})$ is $\preceq^i_c$-minimum in $\mathfrak{G}_3$.\smallskip

Moreover, the $(\mathfrak{G}_\kappa ,\preceq^i_c)$'s and the $(\mathfrak{G}_\kappa ,\preceq_c)$'s are not well-founded. They also  contain antichains of size continuum in the case of $\preceq^i_c$ or when $\kappa\!\not=\! 3$.\end{thm}

 We can also evaluate the descriptive complexity of the $\mathfrak{G}_\kappa$'s. In order to do that, we code the class 
$\mathfrak{G}_\kappa$. By [K, 7.8], any 0DMC space is homeomorphic to a subspace of $2^\omega$, so we can restrict our attention to compact subspaces of $2^\omega$. The Ryll-Nardzewski theorem (see [Kn-R, Corollary 2 and Remark 3]) shows that any homeomorphism on such a subspace can be extended to a homeomorphism of $2^\omega$. The extension map is injective and, conversely, the restriction map is not. But the chromatic number of the graph on the subspace does not depend on the extension, so the fact that the restriction map is not injective creates no problem. The space $\mathcal{K}(X)$ of compact subsets of a metrizable compact space $X$, equipped with the Vietoris topology, is a metrizable compact space. The set 
$\mathcal{H}(2^\omega )$ of homeomorphisms of $2^\omega$ can be equipped with a topology in such a way that it is a Polish group. We set $\mathcal{P}\! :=\!\{ (X,f)\!\in\!\mathcal{K}(2^\omega )\!\times\!\mathcal{H}(2^\omega )\mid f[X]\! =\! X\}$ and code $\mathfrak{G}_\kappa$ with $\mathcal{O}_\kappa\! :=\!\{ (X,f)\!\in\!\mathcal{P}\mid\chi_c(X,G_{f_{\vert X}})\! >\!\kappa\}$.
 
\begin{thm} \label{desc} $\mathcal{P}$ is a Polish space. $\mathcal{O}_0$ is a $\borone$ subset of $\mathcal{P}$, 
$\mathcal{O}_1$ is $\boraone$-complete, and $\mathcal{O}_2,\mathcal{O}_3$ are $\bormtwo$-complete. Moreover, the set 
$\mathcal{O}_2^{\aleph_0}\! :=\!\{ (X,f)\!\in\!\mathcal{O}_2\mid X\mbox{ is countable}\}$ is $\ca$-complete.\end{thm}
 
 Another motivation for studying graphs induced by a function is related to Cantor dynamical systems. We will see that if $f,g$ are minimal homeomorphisms of a Cantor space $X,Y$ respectively, then $f,g$ are flip-conjugate exactly when $(X,G_f)\preceq^i_c(Y,G_g)$. Similar considerations also motivate our study of oriented graphs: in this case, $\textup{Graph}(f),\textup{Graph}(g)$ are oriented graphs, and $f,g$ are conjugate exactly when $\big( X,\textup{Graph}(f)\big)\preceq^i_c\big( Y,\textup{Graph}(g)\big)$. This also leads to study graphs induced by a total homeomorphism. The next result is a version of Theorem \ref{eantichmin} for graphs induced by a total homeomorphism.
 
\begin{thm} \label{eantichmino} There is a $\preceq_c$-antichain (and thus $\preceq_c^i$-antichain) 
$\big( (\mathcal{C}_\alpha ,G_{f_\alpha})\big)_{\alpha\in 2^\omega}$, where\smallskip
 
\noindent (a) $\mathcal{C}_\alpha$ is homeomorphic to $2^\omega$,\smallskip

\noindent (b) $f_\alpha$ is a minimal homeomorphism of $\mathcal{C}_\alpha$ (in fact an odometer), and 
$(\mathcal{C}_\alpha ,G_{f_\alpha})$ has CCN three,\smallskip

\noindent (c) $(\mathcal{C}_\alpha ,G_{f_\alpha})$ is $\preceq_c^i$-minimal in $\mathfrak{G}_2$ and in the class of closed graphs on a 0DMC space with CCN at least three.\smallskip

 In particular, any $\preceq_c^i$-basis for one of these classes must have size continuum.\end{thm}
  
 We also provide a concrete $\preceq_c^i$-basis, made up of graphs induced by an odometer, for the class of elements of 
$\mathfrak{G}_2$ induced by a minimal equicontinuous Cantor dynamical system, and, under the axiom of choice, a 
$\preceq_c^i$-antichain basis for this class. However, we will see that such a basis is far from being a basis for 
$\mathfrak{G}_2$, because of the subshifts associated with irrational rotations, proving a version of Theorem \ref{eantichmino} for them and $\preceq^i_c$. Thanks to subshifts, we also prove a version of Theorem \ref{infdecrcompact} for graphs induced by a total homeomorphism.\medskip

 The next result shows that the situation in the compact case is different from that in the case of spaces which are not compact. The next picture  describes a countable compact subset $K_0$ of $2^\mathbb{Z}$, a two-sided subshift, as well as a  homeomorphism $h_0\! :=\!\sigma_{\vert K_0}\! :\! K_0\!\rightarrow\! K_0$ which is not minimal.
 
\vfill\eject
 
\centerline{\scalebox{0.47}{$$~~~~~\xymatrix{
& & & & & & & \mbox{$(01)^\infty\!\cdot\! 1(01)^\infty$} \ar@{->}[lddddd] \ar@<0.1mm>@{->}[lddddd] \ar@<-0.1mm>@{->}[lddddd] 
& & & & & & & & \\ \\ 
& & & & & & 
(01)^\infty\!\cdot\! (01)1(01)^\infty \ar@{->}[rr] \ar@<0.1mm>@{->}[rr] \ar@<-0.1mm>@{->}[rr] 
& & \ar@{->}[luu] \ar@<0.1mm>@{->}[luu] \ar@<-0.1mm>@{->}[luu]
(10)^\infty\!\cdot\! 11(01)^\infty & & & & & & & \\ \\ \\ 
& & & & & & 
(01)^\infty 1\!\cdot\! (01)^\infty \ar@{->}[rr] \ar@<0.1mm>@{->}[rr] \ar@<-0.1mm>@{->}[rr] & & 
(01)^\infty 10\!\cdot\! (10)^\infty \ar@{->}[llddd] \ar@<0.1mm>@{->}[llddd] \ar@<-0.1mm>@{->}[llddd] & & & & & & & \\ \\ 
& & & & & & 
(01)^\infty\!\cdot\! (01)^21(01)^\infty \ar@{->}[rr] \ar@<0.1mm>@{->}[rr] \ar@<-0.1mm>@{->}[rr] & & 
(10)^\infty\!\cdot\! 1(01)1(01)^\infty \ar@{->}[lluuuuu] \ar@<0.1mm>@{->}[lluuuuu] \ar@<-0.1mm>@{->}[lluuuuu] & & & & & & & \\ 
& & & & & & 
(01)^\infty 1(01)\!\cdot\! (01)^\infty \ar@{->}[rr] \ar@<0.1mm>@{->}[rr] \ar@<-0.1mm>@{->}[rr] & & 
(01)^\infty 1(01)0\!\cdot\! (10)^\infty & & & & & & & \\ 
& & & & & & & \cdots & & & & & & & & \\ 
& & & & & & 
\mbox{$(01)^\infty\!\cdot\! (01)^\infty$} \ar@{<->}[rr] \ar@<0.1mm>@{<->}[rr] \ar@<-0.1mm>@{<->}[rr] & & 
\mbox{$(10)^\infty\!\cdot\! (10)^\infty$} & & & & & & & 
}$$}}\bigskip

 The sequence $(01)^\infty\!\cdot\! (01)^\infty$ is the element $\alpha$ of $2^\mathbb{Z}$ satisfying 
$\alpha (2n\! +\!\varepsilon )\! :=\!\varepsilon$ if $n\!\in\!\omega$ and $\varepsilon\!\in\! 2$ on the positive side, and 
$\alpha (-2n\! -\! 1\! -\!\varepsilon )\! :=\! 1\! -\!\varepsilon$ on the negative side. Here, our space $K_0$ is 
$\mbox{Orb}_\sigma\big( (01)^\infty\!\cdot\! (01)^\infty\big)\bigcup\mbox{Orb}_\sigma\big( (01)^\infty\!\cdot\! 1(01)^\infty\big)$, where $\sigma\! :\! 2^\mathbb{Z}\!\rightarrow\! 2^\mathbb{Z}$ is the shift map. 

\begin{thm} \label{absmincompmain} We can find a countable (0D)MC space $K_0$ and a homeomorphism $h_0$ of $K_0$ such that\smallskip

\noindent (a) $(K_0,G_{h_0})$ has CCN three, and is $\preceq^i_c$-minimal in $\mathfrak{G}_2$ and in the class of closed graphs on a 0DMC space with CCN at least three,\smallskip

\noindent (b) if $S$ is a 0DMS (resp., 0DP) space, $f$ is a homeomorphism of $S$ with the properties that 
$(S,G_f)$ has CCN at least three and $(S,G_f)\preceq^i_c(K_0,G_{h_0})$, then there is a finer 0DMS (resp., 0DP) topology 
$\tau$ on $K_0$ with the properties that $h_0$ is a homeomorphism of $(K_0,\tau )$, $\big( (K_0,\tau ),G_{h_0}\big)$ has CCN  three, and $\big( (K_0,\tau ),G_{h_0}\big)$ is strictly $\preceq^i_c$-below $(S,G_f)$.\smallskip

 In particular, there is no $\preceq^i_c$-antichain basis for the class of graphs induced by a (total) homeomorphism of a 0DMS (or 0DP) space with CCN at least three. Also, any $\preceq^i_c$-basis for one of these classes must have size continuum.\end{thm}

 We saw a number of results describing classes with complex structures. In order to get simpler structures, we can try to study smaller classes, even if a big class may have a minimum element and not a subclass. Using graphs in the style of $(K_0,G_{h_0})$, one can try to study the class of graphs induced by a homeomorphism of a countable MC space with CCN at least three, in this direction. We provide examples of arbitrarily high Cantor-Bendixson rank. 

\begin{thm} \label{CB1intro} (a) Let $\xi\!\geq\! 1$ be a countable ordinal, finite or of the form $\eta\! +\! 3$. Then there is a countable two-sided subshift $\Sigma$ with Cantor-Bendixson rank $\xi$, such that $(\Sigma ,G_{\sigma_{\vert\Sigma}})$ has CCN three, and is $\preceq^i_c$-minimal in $\mathfrak{G}_2$ and in the class of closed graphs on a 0DMC space with CCN at least three.\smallskip

\noindent (b) There is a family $(\Sigma_\alpha )_{\alpha\in 2^\omega}$ of countable two-sided subshifts with Cantor-Bendixson rank three sharing these properties, and such that the family 
$\big( (\Sigma_\alpha ,G_{\sigma_{\vert\Sigma_\alpha}})\big)_{\alpha\in 2^\omega}$ is a $\preceq^i_c$-antichain. In particular, any $\preceq^i_c$-basis for $\mathfrak{G}_2$ or the class of graphs induced by a homeomorphism of a countable (0D)MC space with CCN at least three must have size $2^{\aleph_0}$.\end{thm}

 At this moment, it is still possible to have a $\preceq^i_c$-antichain basis for the class of graphs induced by a homeomorphism of a countable MC space with CCN at least three. Note that the graphs given by Theorem \ref{CB1intro}(a) form a 
$\preceq^i_c$-antichain. The situation for the other values of $\xi$ is not clear.
 
\vfill\eject 
 
 Indeed, recall that if $(X,f)$ is a dynamical system where a compatible metric $d$ on $X$ is fixed, then $(X,f)$ is \emph{expansive} if 
$\exists\varepsilon\! >\! 0~~\forall x\!\not=\! y\!\in\! X~~\exists n\!\in\!\mathbb{Z}~~d\big( f^n(x),f^n(y)\big)\!\geq\!\varepsilon$. The $\sigma_\Sigma$'s are expansive, and there is no expansive homeomorphism of a countable MC space with Cantor-Bendixson rank $\lambda\! +\! 1$ if $\lambda$ is a limit ordinal (see [Ki-Kat-Pa, Theorem 3.2]). We leave this open for future work.\medskip 

 For Question (7), we prove the following.

\begin{thm} \label{embed} We can embed the quasi-order of inclusion on the power set of $\omega$ into\smallskip

\noindent (a) the quasi-order $\preceq^i_c$ on the class of graphs induced by a (total) homeomorphism of a countable 0DMC (and thus 0DP, 0DMS) space with CCN three,\smallskip

\noindent (b) the quasi-order $\preceq_c$ on the class of countable graphs on a 0DMC (and thus 0DP, 0DMS) space with CCN three.\end{thm}

 We can say more about the association between homeomorphisms and graphs mentioned above. The space $\mathbb{M}$ of minimal homeomorphisms of $2^\omega$ is a Polish space. The map associating $(2^\omega ,G_f)$ to $f\!\in\!\mathbb{M}$ is continuous. Moreover, the graph $(2^\omega ,G_f)$ has CCN two or three. The equivalence relations of flip-conjugacy and conjugacy on $\mathbb{M}$ are denoted by $FCO$ and $CO$ respectively. The equivalence relation 
$\preceq^i_c\cap~(\preceq^i_c)^{-1}$ associated with the quasi-order $\preceq^i_c$ on the space 
$$\mathcal{S}_m\! :=\!
\{ (2^\omega ,K)\!\in\!\{ 2^\omega\}\!\times\!\mathcal{K}(2^\omega\!\times\! 2^\omega )\mid 
K\cap\Delta (2^\omega )\! =\!\emptyset\wedge 2\!\leq\!\chi_c(2^\omega ,K)\!\leq\! 3\}$$ 
is denoted by $\equiv^i_c$ (we will check that $\mathcal{S}_m$ is a Polish space). The standard way to compare analytic equivalence relations on standard Borel spaces is the Borel reducibility $\leq_B$ (see, for instance, [G]). Recall that if $X,Y$ are standard Borel spaces and $E,F$ are analytic equivalence relations on $X,Y$ respectively, then 
$(X,E)\leq_B(Y,F)\Leftrightarrow\exists\varphi\! :\! X\!\rightarrow\! Y\mbox{ Borel with }E\! =\! (\varphi\!\times\!\varphi )^{-1}(F)$. 

\begin{thm} \label{redbor} The relations $FCO$, $CO$ and $\equiv^i_c$ are analytic, and $FCO$ is Borel reducible to 
$\equiv^i_c$.\end{thm}
 
 We can also use our countable graphs to prove a version of Theorem \ref{redbor} for graphs of CCN at least three. Using oriented graphs instead of graphs, one can prove that $CO$ is Borel reducible to $\equiv^i_c$. Note that the relation $=^+$ on $\mathbb{R}^\omega$ defined by $x\! =^+\! y\Leftrightarrow\{ x_i\mid i\!\in\!\omega\}\! =\!\{ y_i\mid i\!\in\!\omega\}$ is Borel reducible to $CO$ (this is proved in [Ka]).\medskip

 The present work is organized as follows. In Section 2, we briefly discuss the case of graphs on a finite set. We then study general graphs. Section 3 is about our positive basis results: we prove Theorems \ref{corcomp''''''} and \ref{eq1++}, give our second basis, and start to prepare the proof of Theorem \ref{absmin}. In Section 4, we work in 0DMS spaces and prove Theorem  \ref{absmin}. In Section 5, we study the relation between graphs and dynamical systems, and prove the main part of the version of Theorem \ref{redbor} for graphs of CCN at least three. In Section 6, we start to use odometers and prove Theorems \ref{eantichmin} and \ref{infdecrcompact}. In Section 7, we begin our study with the graphs induced by a function and prove general facts. We then study in Section 8 the graphs induced by a partial function and prove Theorem \ref{main}. In Section 9, we keep on using odometers and prove Theorems \ref{eantichmino}, \ref{embed}(b) and \ref{absmincomp}. Section 10 is devoted to the study of graphs induced by a subshift. In particular, we study the homeomorphisms of a countable compact space and prove Theorem \ref{CB1intro}. In Section 11, we work in 0DMS spaces and prove Theorems \ref{absmincompmain} and \ref{embed}(a). In Section 12, we prove Theorems \ref{Ckappa} and \ref{desc}. In Section 13, we study equivalence relations and prove two versions of Theorem \ref{redbor}. Section 14 is devoted to the versions of our results for digraphs and oriented graphs. Finally, we summarize our work about general graphs in a table in Section 15, which leaves some other open questions for the future.

\section{$\!\!\!\!\!\!$ General graphs on a finite set}\indent

 We briefly discuss the finite case, already showing that the quasi-orders $\preceq_c$ and $\preceq^i_c$ are quite different. In this finite case, we put the discrete topology on the space, so that continuity is automatic. It is known that a graph has chromatic number at most two exactly when it is bipartite, and when it has no odd cycle (see [A-D-H, 2.1]). Thus a graph $G$ on a set $X$ has chromatic number at least three exactly when $\Delta (X)\cap (\bigcup_{p\in\omega}~G^{2p+3})$ is not empty. Recall that a \emph{walk} in a relation $(X,R)$ is a sequence $(x_i)_{i\leq n}\!\in\! X^{n+1}$ such that $(x_i,x_{i+1})\!\in\! R$ for each 
$i\! <\! n$. A walk $(x_i)_{i\leq n}$ is \emph{odd} if $n$ is odd, \emph{closed} if $x_0\! =\! x_n$, and a \emph{cycle} if it is closed, 
$n\!\geq\! 3$ and $(x_i)_{i<n}$ is injective.  We denote, for any natural number $p$, the symmetric cycle on $2p\! +\! 3$ by $C_{2p+3}$. 

\begin{thm} \label{eqfin} Let $X$ be a finite set, and $G$ be a graph on $X$. The following are equivalent:\smallskip

\noindent (1) $(X,G)$ has chromatic number at least three,\smallskip

\noindent (2) $\Delta (X)\cap (\bigcup_{p\in\omega}~G^{2p+3})\!\not=\!\emptyset$,\smallskip

\noindent (3) there is $p\!\in\!\omega$ with $(2p\! +\! 3,C_{2p+3})\preceq^i(X,G)$.\end{thm}

\begin{cor} \label{basisfin} Let $\mathfrak{F}\! :=\!\big( (2p\! +\! 3,C_{2p+3})\big)_{p\in\omega}$. 

\noindent (a) $\mathfrak{F}$ is a $\preceq^i$-antichain basis for the class of graphs on a finite set with chromatic number at least three. In particular, the elements of $\mathfrak{F}$ are $\preceq^i_c$-minimal among graphs on a 0DMS (or 0DP, or 0DMC) space with CCN at least three, and any $\preceq^i_c$-basis for these classes must be infinite.\smallskip

\noindent (b) $\mathfrak{F}$ is a $\preceq$-basis for the class of graphs on a finite set with chromatic number at least three, and is strictly 
$\preceq$-decreasing. In particular, there is neither $\preceq$-antichain basis, nor finite basis for this class. Also, no graph is $\preceq$-minimal in this class.\end{cor}

\noindent\emph{Proof.}\ By [He-N, Corollary 1.4], $(2p\! +\! 3,C_{2p+3})\not\preceq (2q\! +\! 3,C_{2q+3})$ if $p\! <\! q$.\medskip

\noindent (a) Theorem \ref{eqfin} gives the basis. This is an antichain by the argument just above and by injectivity.\medskip
 
\noindent (b) We apply (a), the argument above again, and [He-N, Corollary 1.4].\hfill{$\square$}\medskip

\noindent\emph{Remarks.}\ (1) Let $X$ be a finite set, $G$ be a graph on $X$, $Y$ be a set, $H$ be a graph on $Y$ with the property that $(Y,H)\preceq^i(X,G)$, with witness $\varphi$. We set $V\! :=\!\varphi [Y]$ and 
$E\! :=\! (\varphi\!\times\!\varphi )[H]$. Note that $V$ is a subset of $X$, $E$ is a graph on $V$ contained in $G$ with 
$(V,E)\preceq^i(Y,H)$ with witness $\varphi^{-1}$, and also $(Y,H)\preceq^i(V,E)$. This and the finiteness of $X$ implies that there is no infinite $\preceq^i$-descending chain in the class of graphs on a finite set.\medskip

\noindent (2) There are infinite $\preceq$-antichains in the class of graphs on a finite set with chromatic number at least three. Indeed, this comes from [He-N, Theorem 2.23 and Proposition 3.4]. Following their notation, $S(i,i)$ and $S(j,j)$ are $\preceq$-incomparable if $i\!\not=\! j\!\geq\! 3$ are odd, and have chromatic number at least three.

\section{$\!\!\!\!\!\!$ General graphs on a 0DMC space: the basis}\label{compact}

\label{bascomp}

\subsection{$\!\!\!\!\!\!$ A first basis}\indent

 We now define the concrete family announced in Theorem \ref{corcomp''''''}.\medskip
 
\noindent\emph{Notation.}\ We denote the set of increasing unbounded sequences of natural numbers by
$$\mathcal{S}\! :=\!\{\delta\!\in\!\omega^\omega\mid\forall k\!\in\!\omega ~~
\delta (k)\!\leq\!\delta (k\! +\! 1)\wedge\forall N\!\in\!\omega ~~\exists k\!\in\!\omega ~~\delta (k)\!\geq\! N\} .$$
- In the proof of Theorem \ref{corcomp''''''}, it will be convenient to replace the index set $2^\omega$ with the set 
$\mathcal{I}$ that we now define. We denote a typical element of $\prod_{k\in\omega}~(2^\omega )^{2k+1}\!\!\times\! (2^\omega )^{2k+1}$ by 
$$\gamma\! :=\!\bigg(\Big( \big(\gamma^0_k(i)\big)_{i\leq 2k} ,\big(\gamma^1_k(i)\big)_{i\leq 2k}\Big)\bigg)_{k\in\omega}.$$
We also set\medskip

\leftline{$\mathcal{I}\! :=\!\big\{\gamma\!\in\!\prod_{k\in\omega}~(2^\omega )^{2k+1}\!\!\times\! (2^\omega )^{2k+1}\mid
\gamma^0_k(0),\gamma^1_k(2k)\!\rightarrow\! 0^\infty\wedge
\forall k\!\in\!\omega ~~\forall i\!\leq\! 2k~~\gamma^0_k(i)\!\not=\!\gamma^1_k(i)~\wedge$}\smallskip

\rightline{$\exists\delta\!\in\!\mathcal{S}~~\forall k\!\in\!\omega ~~\forall i\! <\! 2k~~
\gamma^1_k(i)\vert\delta (k)\! =\!\gamma^0_k(i\! +\! 1)\vert\delta (k)\big\}$.}\medskip

\noindent - We then define, for $\gamma\!\in\!\mathcal{I}$, a countable graph $\mathbb{G}_\gamma$ on $2^\omega$ by 
$\mathbb{G}_\gamma\! :=\! s\big(\big\{\big(\gamma^0_k(i),\gamma^1_k(i)\big)\mid k\!\in\!\omega\wedge i\!\leq 2k\big\}\big)$ and set 
$\mathbb{K}_\gamma\! :=\!\overline{\mbox{proj}[\mathbb{G}_\gamma]}^{2^\omega}$, so that $\mathbb{G}_\gamma$ is a graph on the compact set $\mathbb{K}_\gamma$. The next picture represents $\mathbb{G}_\gamma$.\medskip\medskip

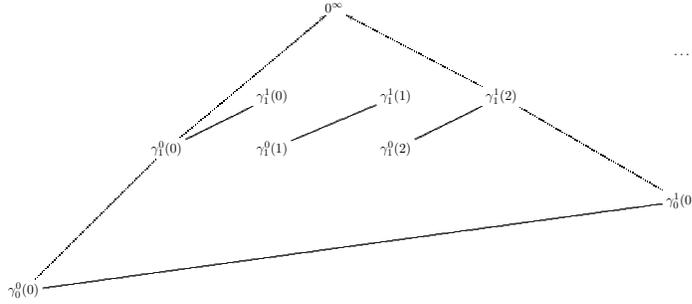
\centerline{\scalebox{0.47}{$$~~~~~\xymatrix{
& & & & & & & 0^\infty & & & & & & & & \\ 
& & & & & & & & & & & & & & \ldots\\ 
& & & & & & \gamma^1_1(0) & & \gamma^1_1(1) & & \gamma^1_1(2) \ar@{..>}[uulll] & & & & & \\ 
& & & & \gamma^0_1(0) \ar@{..>}[uuurrr] \ar@{-}[urr] \ar@<0.1mm>@{-}[urr] \ar@<-0.1mm>@{-}[urr] & & 
\gamma^0_1(1) \ar@{-}[urr] \ar@<0.1mm>@{-}[urr] \ar@<-0.1mm>@{-}[urr] & & 
\gamma^0_1(2) \ar@{-}[urr] \ar@<0.1mm>@{-}[urr] \ar@<-0.1mm>@{-}[urr] & & & & & & & \\ 
& & & & & & & & & & & & & & \gamma^1_0(0) \ar@{..}[uullll] & \\ 
& & & & & & & & & & & & & & \\ 
& \gamma^0_0(0) \ar@{..}[uuurrr] \ar@{-}[uurrrrrrrrrrrrr] \ar@<0.1mm>@{-}[uurrrrrrrrrrrrr] \ar@<-0.1mm>@{-}[uurrrrrrrrrrrrr] & & & & & & & & & & & & & & \\ 
}$$}}\medskip\medskip

 We first prove the exactly part of Theorem \ref{corcomp''''''}.

\begin{propo} \label{K1''''''} Let $\gamma\!\in\!\mathcal{I}$. Then $(\mathbb{K}_\gamma ,\mathbb{G}_\gamma )$ has CCN at least three.\end{propo}

\noindent\emph{Proof.}\ Note that $0^\infty\!\in\!\mathbb{K}_\gamma$. If $(C,\neg C)$ is a coloring of $\mathbb{G}_\gamma$ into clopen subsets of $\mathbb{K}_\gamma$ with $0^\infty\!\in\! C$, then the compactness of $\mathbb{K}_\gamma$ gives 
$l\!\in\!\omega$ with $\alpha ,\beta\!\in\! C$ or $\alpha ,\beta\!\notin\! C$ if $\alpha ,\beta\!\in\!\mathbb{K}_\gamma$ and 
$\alpha\vert l\! =\!\beta\vert l$. Note that $\gamma^0_k(0)\!\in\! C$ if $k\!\geq\! k_0$, where $k_0\!\in\!\omega$ is also large  enough with 
$\delta (k_0)\!\geq\! l$. Assume that $k\!\geq\! k_0$. An induction on $i\! <\! 2k$ shows that 
$\gamma^1_k(i),\gamma^0_k(i\! +\! 1)\!\notin\! C$  if $i$ is even, $\gamma^1_k(i),\gamma^0_k(i\! +\! 1)\!\in\! C$ if 
$i$ is odd, and $\gamma^1_k(2k)\!\notin\! C$. This implies that $0^\infty\!\notin\! C$, which is the desired contradiction. Thus 
$\chi_c(\mathbb{K}_\gamma ,\mathbb{G}_\gamma )\!\geq\! 3$.\hfill{$\square$}\medskip

\noindent\emph{Notation.}\ Let $\mathfrak{C}\! :=\! (\omega\!\setminus\! 2)^\omega$. Fix 
${\bf d}\! =\! (d_j)_{j\in\omega}\!\in\!\mathfrak{C}$.\medskip

\noindent - In the sequel, we denote, for $S\!\subseteq\!\omega$ finite, by $\pi_{j\in S}~d_j$ the natural number, and by 
$\prod_{j\in S}~d_j$ the set of finite sequences of natural numbers. In particular, we set, for $l\!\in\!\omega$, 
$\prod_l\! :=\!\prod_{j<l}~d_j$.\medskip

\noindent - We set $\mathcal{C}\! :=\!\mathcal{C}_{\bf d}\! :=\!\prod_{j\in\omega}~d_j$, so that $\mathcal{C}$ is homeomorphic to 
$2^\omega$. As usual, $N_s\! :=\!\{\alpha\!\in\!\mathcal{C}\mid s\!\subseteq\!\alpha\}$ is a basic clopen set if $s\!\in\!\bigcup_{l\in\omega}~\prod_l$. We extend this notation to other sequential spaces of this kind.\medskip

\noindent - If $R$ is a relation on $\mathcal{C}$, and $n\!\in\!\omega$, then we set 
$R_n\! :=\!\{ (s,t)\!\in\!\prod_n^2\mid (N_s\!\times\! N_t)\cap R\!\not=\!\emptyset\}$ and 
$$^nR\! :=\!\{ (\alpha ,\beta )\!\in\!\mathcal{C}^2\mid\exists (\alpha',\beta')\!\in\! R~~
\alpha\vert n\! =\!\alpha'\vert n\wedge\beta\vert n\! =\!\beta'\vert n\} .$$

\begin{them} \label{eq2''''''} Let ${\bf d}\!\in\!\mathfrak{C}$, and $G$ be a graph on $\mathcal{C}$. The following are equivalent:\smallskip

\noindent (1) $(\mathcal{C},G)$ has CCN at least three,\smallskip

\noindent (2) $\Delta (\mathcal{C})\cap\bigcap_{n\in\omega}~(\bigcup_{p\in\omega}~(^nG)^{2p+1})$ is not empty,\smallskip

\noindent (3) the relations $(\prod_n,G_n)$ have an odd closed walk,\smallskip

\noindent (4) there is $\gamma\!\in\!\mathcal{I}$ with $(\mathbb{K}_\gamma ,\mathbb{G}_\gamma )\preceq^i_c(\mathcal{C},G)$.\end{them}

\noindent\emph{Proof.}\ (4) $\Rightarrow$ (1) This comes from Proposition \ref{K1''''''}.\medskip

\noindent (1) $\Rightarrow$ (2) We first prove that 
$\bigcap_{n\in\omega}\big(\overline{^n\big(\Delta (\mathcal{C})\big)\cup\bigcup_{p\in\omega}~(^nG)^{2p+2}}\cap
\overline{\bigcup_{p\in\omega}~(^nG)^{2p+1}}\big)$ is not empty. We argue by contradiction, which by compactness of 
$\mathcal{C}$ gives $n\!\in\!\omega$ and a clopen relation $O$ on $\mathcal{C}$ separating 
$\bigcup_{p\in\omega}~(^nG)^{2p+1}$ from $^n\big(\Delta (\mathcal{C})\big)\cup\bigcup_{p\in\omega}~(^nG)^{2p+2}$. The compactness of $\mathcal{C}$ gives sequences $(s_j)_{j\leq m}$ and $(t_j)_{j\leq m}$ of finite sequences with the property that 
$O\! =\!\bigcup_{j\leq m}~N_{s_j}\!\times\! N_{t_j}$, and we may assume that all these finite sequences have the same length $l$, and that $l\!\geq\! n$.\medskip

 We define a subset of $\prod_l$ by $V\! :=\!\{ s_j\mid j\!\leq\! m\}\cup\{ t_j\mid j\!\leq\! m\}$. Note that $G_l$ is a graph on 
$\prod_l$ since, for each $s\!\in\!\prod_l$, $(N_s\!\times\! N_s)\cap G\!\subseteq\! ^n\big(\Delta (\mathcal{C})\big)\cap O$. Let 
$(\mathcal{C}_i)_{i\in I}$ be the family of the connected components of $G_l$ restricted to $V$. Fix $i\!\in\! I$. As 
$(\mathcal{C}_i,G_l\cap\mathcal{C}_i^2)$ is a connected graph, we can find an acyclic connected graph 
$\mathcal{G}_i$ on $\mathcal{C}_i$ with $\mathcal{G}_i\!\subseteq\! G_l\cap\mathcal{C}_i^2$. This gives a coloring 
$c_i\! :\!\mathcal{C}_i\!\rightarrow\! 2$ of $\mathcal{G}_i$.\medskip

 We set ${C\! :=\!\bigcup_{i\in I}~\{ N_s\mid c_i(s)\! =\! 0\}}$, so that $C$ is a clopen subset of $\mathcal{C}$. It remains to prove that $G\cap (C^2\cup (\neg C)^2)\! =\!\emptyset$, since this contradicts (a). Towards a contradiction, suppose that there is  
$(\alpha ,\beta )$ in $G\cap C^2$, for example (the other case is similar). As $G\!\subseteq\! O$, we can find $j\!\leq\! m$ with $(\alpha ,\beta )\!\in\! N_{s_j}\!\times\! N_{t_j}$. In particular, $(s_j,t_j)$ is in $G_l\cap V^2$, which gives 
$i\!\in\! I$ with $s_j,t_j\!\in\!\mathcal{C}_i$. As $\alpha ,\beta\!\in\! C$, $c_i(s_j)\! =\! c_i(t_j)$. Let $L\!\in\!\omega$, 
$(u_k)_{k\leq L}$ be the $\mathcal{G}_i$-path from $s_j$ to $t_j$, and, for $k\! <\! L$, 
$(\alpha_k,\beta_k)\!\in\! (N_{u_k}\!\times\! N_{u_{k+1}})\cap G$. Note that $L$ is even since $c_i(s_j)\! =\! c_i(t_j)$. Also, 
$(\alpha ,\beta_0),(\beta_0,\beta_1),(\beta_1,\beta_2),\ldots ,(\beta_{L-3},\beta_{L-2}),(\beta_{L-2},\beta )\!\in\! ^nG$, so that 
$(\alpha ,\beta )$ is in $^n\big(\Delta (\mathcal{C})\big)\cup\bigcup_{p\in\omega}~(^nG)^{2p+2}\!\subseteq\!\neg O$, which is absurd.\medskip

 Pick $(\alpha ,\beta )\!\in\!\bigcap_{n\in\omega}\big(\overline{^n\big(\Delta (\mathcal{C})\big)\cup
\bigcup_{p\in\omega}~(^nG)^{2p+2}}\cap\overline{\bigcup_{p\in\omega}~(^nG)^{2p+1}}\big)$. If $\alpha\! =\!\beta$, then for each $n$ there is $p$ such that $N_{\alpha\vert n}^2$ meets $(^nG)^{2p+1}$. Let $(\gamma_i)_{i\leq 2p+1}\!\in\!\mathcal{C}^{2p+2}$ such that 
$\gamma_0,\gamma_{2p+1}\!\in\! N_{\alpha\vert n}$, and, for each $i\!\leq\! 2p$, $(\gamma_i,\gamma_{i+1})\!\in\! ^nG$. Note that 
$(\alpha ,\gamma_1),(\gamma_{2p},\alpha )\!\in\! ^nG$, so that we may assume that 
$\gamma_0\! =\!\gamma_{2p+1}\! =\!\alpha$ and $(\alpha ,\alpha )\!\in\! (^nG)^{2p+1}$. Thus 
$(\alpha ,\alpha )\!\in\!\Delta (\mathcal{C})\cap\bigcap_{n\in\omega}~(\bigcup_{p\in\omega}~(^nG)^{2p+1})$. So we may assume that $\alpha\!\not=\!\beta$. Note that $(\alpha ,\beta )\!\notin\!\overline{^n\big(\Delta (\mathcal{C})\big)}$ if $n$ is large enough, and that the intersection above is decreasing with respect to $n$. This implies that $(\alpha ,\beta )\!\in\!
\bigcap_{n\in\omega}~\Big(\overline{\bigcup_{p\in\omega}~(^nG)^{2p+2}}\cap\overline{\bigcup_{p\in\omega}~(^nG)^{2p+1}}\Big)$. So we can pick, for each $n\!\in\!\omega$ and each $\varepsilon\!\in\! 2$, a natural number $p^\varepsilon_n$, and 
$$(\alpha_{2n+1+\varepsilon},\beta_{2n+1+\varepsilon})\!\in\! 
(N_{\alpha\vert n}\!\times\! N_{\beta\vert n})\cap (^nG)^{2p^\varepsilon_n+1+\varepsilon} .$$
Let $(\gamma^{2n+1+\varepsilon}_i)_{i\leq 2p^\varepsilon_n+1+\varepsilon}\!\in\!
\mathcal{C}^{2p^\varepsilon_n+2+\varepsilon}$ such that $\gamma^{2n+1+\varepsilon}_0\! =\!\alpha_{2n+1+\varepsilon}$, 
$\gamma^{2n+1+\varepsilon}_{2p^\varepsilon_n+1+\varepsilon}\! =\!\beta_{2n+1+\varepsilon}$, and, for each 
$i\!\leq\! 2p^\varepsilon_n\! +\!\varepsilon$, $(\gamma^{2n+1+\varepsilon}_i,\gamma^{2n+1+\varepsilon}_{i+1})\!\in\! ^nG$. Fix $n\!\in\!\omega$. Note that 
$(\gamma^{2n+1}_0,\gamma^{2n+2}_1),(\gamma^{2n+2}_{2p^1_n+1},\gamma^{2n+1}_{2p^0_n+1})$ are in $^nG$. This implies that $(\alpha ,\alpha )\!\in\! (^nG)^{2p^0_n+1+2p^1_n}\!\subseteq\!\bigcup_{p\in\omega}~(^nG)^{2p+1}$.\medskip

\noindent (2) $\Rightarrow$ (4) We choose, for each $j\!\in\!\omega$, $l_j\!\geq\! 1$ with $d_j\!\leq\! 2^{l_j}$. This defines an injection 
$i_j\! :\! d_j\!\rightarrow\! 2^{l_j}$. We define, for $x\!\in\!\mathcal{C}$, $\psi (x)\!\in\! 2^\omega$ by 
$\psi (x)\! :=\! {^\frown}_{j\in\omega}~i_j\big( x(j)\big)$. Note that $\psi$ is a continuous injection, and thus a homeomorphism onto its range $R$. We set $H\! :=\! (\psi\!\times\!\psi )[G]$, so that $H$ is a graph on $2^\omega$. Moreover, if $n\!\in\!\omega$ and $(x,y)\!\in\! ^nG$, then $\big(\psi (x),\psi (y)\big)\!\in\! ^nH$. Pick 
$(x,x)\!\in\!\bigcap_{n\in\omega}(\bigcup_{p\in\omega}(^nG)^{2p+1})$, and set $\alpha\! :=\!\psi (x)$.

\vfill\eject

 Then $(\alpha ,\alpha )\!\in\!\bigcap_{n\in\omega}(\bigcup_{p\in\omega}(^nH)^{2p+1})$. This gives, for each $n\!\in\!\omega$, a natural number $p_n$ with $(\alpha ,\alpha )\!\in\! (^nH)^{2p_n+1}$. Note that, extracting a subsequence if necessary, we may assume that $(p_n)_{n\in\omega}$ is constant or strictly increasing. Let $(\gamma^n_j)_{j\leq 2p_n+1}$ in 
$(2^\omega )^{2p_n+2}$ such that $\gamma^n_0\! =\!\alpha$, $\gamma^n_{2p_n+1}\! =\!\alpha$, and, for each $i\!\leq\! 2p_n$, 
$(\gamma^n_i,\gamma^n_{i+1})$ in $^nH$. This gives, when $i\!\leq\! 2p_n$, 
$(\alpha^n_i,\beta^n_i)\!\in\! (N_{\gamma^n_i\vert n}\!\times\! N_{\gamma^n_{i+1}\vert n})\cap H$.\medskip

 We now define $\upsilon\!\in\!\prod_{k\in\omega}~(2^\omega )^{2k+1}\!\!\times\! (2^\omega )^{2k+1}$ as follows. If 
$k\!\leq\! p_0$ and $i\!\leq\! 2k$, then we set $\upsilon^0_k(i)\! :=\!\alpha^0_i$ and $\upsilon^1_k(i)\! :=\!\beta^0_i$.\medskip

 Assume first that $(p_n)_{n\in\omega}$ is constant. If $k\! >\! p_0$ and $i\!\leq\! 2p_0$, then we set $\upsilon^0_k(i)\! :=\!\alpha^k_i$ and 
$\upsilon^1_k(i)\! :=\!\beta^k_i$. If $k\! >\! p_0$ and $2p_0\! <\! i\!\leq\! 2k$, then we set 
$$\upsilon^0_k(i)\! :=\!\left\{\!\!\!\!\!\!\!
\begin{array}{ll}
& \alpha^k_{2p_0}\mbox{ if }i\mbox{ is even,}\cr
& \beta^k_{2p_0}\mbox{ if }i\mbox{ is odd,}\cr
\end{array}
\right.
~~~~~
\upsilon^1_k(i)\! :=\!\left\{\!\!\!\!\!\!\!
\begin{array}{ll}
& \beta^k_{2p_0}\mbox{ if }i\mbox{ is even,}\cr
& \alpha^k_{2p_0}\mbox{ if }i\mbox{ is odd.}
\end{array}
\right.$$

 Assume now that $(p_n)_{n\in\omega}$ is strictly increasing. If $p_j\! <\! k\!\leq\! p_{j+1}$ and $i\!\leq\! 2p_j$, then we set 
$\upsilon^0_k(i)\! :=\!\alpha^{p_j}_i$ and $\upsilon^1_k(i)\! :=\!\beta^{p_j}_i$. If $p_j\! <\! k\!\leq\! p_{j+1}$ and $2p_j\! <\! i\!\leq\! 2k$, then we set 
$$\upsilon^0_k(i)\! :=\!\left\{\!\!\!\!\!\!\!
\begin{array}{ll}
& \alpha^{p_j}_{2p_j}\mbox{ if }i\mbox{ is even,}\cr
& \beta^{p_j}_{2p_j}\mbox{ if }i\mbox{ is odd,}\cr
\end{array}
\right.
~~~~~
\upsilon^1_k(i)\! :=\!\left\{\!\!\!\!\!\!\!
\begin{array}{ll}
& \beta^{p_j}_{2p_j}\mbox{ if }i\mbox{ is even,}\cr
& \alpha^{p_j}_{2p_j}\mbox{ if }i\mbox{ is odd.}
\end{array}
\right.$$
Note that $\upsilon^0_k(0)\! =\!\alpha^{m_k}_0$, with $\mbox{lim}_{k\rightarrow\infty}~m_k\! =\!\infty$. As 
$\alpha^{m_k}_0\!\in\! N_{\gamma^{m_k}_0\vert m_k}\! =\! N_{\alpha\vert m_k}$, $\big(\upsilon^0_k(0)\big)_{k\in\omega}$ converges to $\alpha$. Similarly, if $(p_n)_{n\in\omega}$ is constant, then $\upsilon^1_k(2k)\! =\!\beta^k_{2p_0}$ if $k$ is large  enough. As $\beta^k_{2p_0}$ is in 
$N_{\gamma^k_{2p_0+1}\vert k}\! =\! N_{\alpha\vert k}$, $\upsilon^1_k(2k)\vert k\! =\!\alpha\vert k$. If the sequence 
$(p_n)_{n\in\omega}$ is strictly increasing, then $\upsilon^1_k(2k)\! =\!\beta^{p_j}_{2p_j}$ if $k$ is large enough. As 
$\beta^{p_j}_{2p_j}\!\in\! N_{\gamma^{p_j}_{2p_j+1}\vert p_j}\! =\! N_{\alpha\vert p_j}$, 
$\upsilon^1_k(2k)\vert p_j\! =\!\alpha\vert p_j$. Thus $\big(\upsilon^1_k(2k)\big)_{k\in\omega}$ converges to $\alpha$.\medskip

 Note that we chose $\upsilon$ in such a way that $\big(\upsilon^0_k(i),\upsilon^1_k(i)\big)$ is in the graph $H$, so that 
$\upsilon^0_k(i)\!\not=\!\upsilon^1_k(i)$. If $i\! <\! 2k$, then we also ensured that 
$$\left\{\!\!\!\!\!\!\!
\begin{array}{ll}
& \upsilon^1_k(i)\vert 0\! =\!\upsilon^0_k(i\! +\! 1)\vert 0\mbox{ if }k\!\leq\! p_0\mbox{,}\cr
& \upsilon^1_k(i)\vert k\! =\!\upsilon^0_k(i\! +\! 1)\vert k\mbox{ if }k\! >\! p_0\wedge i\! <\! 2p_0\mbox{,}\cr
& \upsilon^1_k(i)\! =\!\upsilon^0_k(i\! +\! 1)\mbox{ if }k\! >\! p_0\wedge i\!\geq\! 2p_0\mbox{,}\cr
& \upsilon^1_k(i)\vert p_j\! =\!\upsilon^0_k(i\! +\! 1)\vert p_j\mbox{ if }p_j\! <\! k\!\leq\! p_{j+1}\wedge i\! <\! 2p_j\mbox{,}\cr
& \upsilon^1_k(i)\! =\!\upsilon^0_k(i\! +\! 1)\mbox{ if }p_j\! <\! k\!\leq\! p_{j+1}\wedge i\!\geq\! 2p_j.
\end{array} 
\right.$$
This defines an element $\zeta$ of $\mathcal{S}$ as desired. Note that the map $h\! :\! 2^\omega\!\rightarrow\! 2^\omega$ defined by
$$h(\beta )(n)\! :=\!\left\{\!\!\!\!\!\!\!
\begin{array}{ll}
& \beta (n)\mbox{ if }\alpha (n)\! =\! 0\mbox{,}\cr
& 1\! -\!\beta (n)\mbox{ if }\alpha (n)\! =\! 1\mbox{,}
\end{array}
\right.$$ 
is a homeomorphism sending $\alpha$ to $0^\infty$. We set 
$\gamma^\varepsilon_k(i)\! :=\! h\big(\upsilon^\varepsilon_k(i)\big)$, which defines 
$\gamma$ in the set $\prod_{k\in\omega}~(2^\omega )^{k+1}\!\!\times\! (2^\omega )^{k+1}$. The sequence $\delta\! :=\!\zeta$ is in $\mathcal{S}$ and is a witness for the fact that $\gamma\!\in\!\mathcal{I}$. Moreover, the map $h^{-1}$ is a witness for the fact that $(2^\omega ,\mathbb{G}_\gamma )\preceq^i_c(2^\omega ,H)$. We set $K\! :=\! h[R]$. As $R$ is compact, $K$ is too. Note that $H\!\subseteq\! R^2$, so that $\mathbb{G}_\gamma\!\subseteq\! K^2$. Thus 
$\mbox{proj}[\mathbb{G}_\gamma ]\!\subseteq\! K$ and 
${\mathbb{K}_\gamma\! =\!\overline{\mbox{proj}[\mathbb{G}_\gamma ]}\!\subseteq\! K}$. We are done since 
$(\mathbb{K}_\gamma ,\mathbb{G}_\gamma )\preceq^i_c(K,\mathbb{G}_\gamma )\preceq^i_c(R,H)\preceq^i_c
(\mathcal{C},G)$.

\vfill\eject

\noindent (2) $\Rightarrow$ (3) Assume that 
$(\alpha ,\alpha )\!\in\!\Delta (\mathcal{C})\cap\bigcap_{n\in\omega}~(\bigcup_{p\in\omega}~(^nG)^{2p+1})$. Fix 
$n\!\in\!\omega$. We can find $p$ and $(\gamma_i)_{i\leq 2p+1}\!\in\!\mathcal{C}^{2p+2}$ with 
$\gamma_0\! =\!\gamma_{2p+1}\! =\!\alpha$ and $(\gamma_i,\gamma_{i+1})\!\in\! {}^nG$ if $i\!\leq\! 2p$. This gives, for each 
$i\!\leq\! 2p$, $(\alpha_i,\beta_i)\!\in\! G$ with $\alpha_i\vert n\! =\!\gamma_i\vert n$ and 
$\beta_i\vert n\! =\!\gamma_{i+1}\vert n$. We set, for each $i\!\leq\! 2p\! +\! 1$, $s_i\! :=\!\gamma_i\vert n$. Then 
$(s_i)_{i\leq 2p+1}$ is an odd closed walk in $(\prod_n,G_n)$.\medskip

\noindent (3) $\Rightarrow$ (2) Let $(s^k_i)_{i\leq 2p_k+1}$ be an odd closed walk in $(\prod_k,G_k)$. As 
$(s^k_i,s^k_{i+1})\!\in\! G_k$, we can find $(\alpha^k_i,\beta^k_i)\!\in\! (N_{s^k_i}\!\times\! N_{s^k_{i+1}})\cap G$ if 
$i\!\leq\! 2p_k$. The compactness of $\mathcal{C}$ provides $\alpha\!\in\!\mathcal{C}$ and $(k_n)_{n\in\omega}$ strictly increasing such that $\alpha_0^{k_n}\vert n\! =\!\alpha\vert n$ for each $n$. Note that\medskip

\noindent - $(\alpha^{k_n}_0,\beta^{k_n}_0)\!\in\! G$, $\alpha\vert n\! =\!\alpha^{k_n}_0\vert n$ and 
$\alpha^{k_n}_1\vert n\! =\! s^{k_n}_1\vert n\! =\!\beta^{k_n}_0\vert n$,\smallskip

\noindent - $(\alpha^{k_n}_i,\beta^{k_n}_i)\!\in\! G$ and 
$\alpha^{k_n}_{i+1}\vert n\! =\! s^{k_n}_{i+1}\vert n\! =\!\beta^{k_n}_i\vert n$ if $1\!\leq\! i\! <\! 2p_{k_n}$,\smallskip

\noindent - $(\alpha^{k_n}_{2p_{k_n}},\beta^{k_n}_{2p_{k_n}})\!\in\! G$ and $\alpha\vert n\! =\!\alpha^{k_n}_0\vert n\! =\! 
s^{k_n}_0\vert n\! =\! s^{k_n}_{2p_{k_n}+1}\vert n\! =\!\beta^{k_n}_{2p_{k_n}}\vert n$,\medskip

 This implies that $(\alpha,\alpha^{k_n}_1,\ldots ,\alpha^{k_n}_{2p_{k_n}},\alpha )$ is in $({}^nG)^{2p_{k_n}+1}$.
\hfill{$\square$}\medskip

\noindent\emph{Remark.}\ The cycle $C_3$ on 3 (pairwise different) points $p,q,r$ is a graph with (continuous) chromatic number $3$, and we may assume that $p\! =\! 0^\infty$, $q\! =\! 1^\infty$, $r\! =\! (01)^\infty$. In this case, the element 
$\gamma$ of $\mathcal{I}$ given by Theorem \ref{eq2''''''}(3) highly lacks of ``injectivity". For instance, we can take 
$\gamma\!\in\!\mathcal{I}$ given by the equalities $\gamma^0_k(0)\! :=\! 0^\infty\! =:\!\gamma^1_k(2k)$, 
$\gamma^0_k(2l\! +\! 1)\! :=\! 1^\infty\! =:\!\gamma^1_k(2l)$, and $\gamma^0_k(2l\! +\! 2)\! :=\! (01)^\infty\! =:\!\gamma^1_k(2l\! +\! 1)$.\medskip

\noindent\emph{Proof of Theorem \ref{corcomp''''''}.}\ By Proposition \ref{K1''''''}, (1) and (2) cannot hold simultaneously. By [K, 4.2], $X$ is Polish, and by [K, 7.8], $X$ is homeomorphic to a subspace of $2^\omega$, which has to be compact and therefore closed. So we may assume that $X$ is a closed subset of $2^\omega$. Assume that the problem is solved for 
$X\! =\! 2^\omega$, and that (1) does not hold. Note that (1) does not hold in $2^\omega$ since it does not hold in $X$. This gives $\gamma\!\in\!\mathcal{I}$ and ${\varphi\! :\!\mathbb{K}_\gamma\!\rightarrow\! 2^\omega}$ injective continuous with 
$\mathbb{G}_\gamma\!\subseteq\! (\varphi\!\times\!\varphi )^{-1}(G)$. In particular, $\mathbb{G}_\gamma$ is contained in the closed set $\big(\varphi^{-1}(X)\big)^2$. Thus $\mbox{proj}[\mathbb{G}_\gamma ]\!\subseteq\!\varphi^{-1}(X)$ and 
$\mathbb{K}_\gamma\! =\!\overline{\mbox{proj}[\mathbb{G}_\gamma ]}\!\subseteq\!\varphi^{-1}(X)$, so that $\varphi$ is a witness for the fact that $(\mathbb{K}_\gamma ,\mathbb{G}_\gamma)\preceq^i_c(X,G)$. So we may assume that 
$X\! =\! 2^\omega$. It remains to apply Theorem \ref{eq2''''''}.\hfill{$\square$}

\subsection{$\!\!\!\!\!\!$ A second basis}\indent \label{sb}

 We now provide another $\preceq_c$-basis, closer to the examples used later.\medskip
 
\noindent\emph{Notation.}\ It will be convenient to use the index set $\mathcal{J}$ that we now define. Fix 
${\bf d}\! =\! (d_j)_{j\in\omega}\!\in\!\mathfrak{C}$. We denote a typical element of 
$\big( (\bigcup_{m\in\omega}~\prod_{m+1})^{<\omega}\big)^\omega$ by 
$\beta\! :=\!\Big(\big(s_l(i)\big)_{i<\lambda_l}\Big)_{l\in\omega}$. We then set
$$\mathcal{J}\! :=\!\big\{\beta\!\in\!\big( (\bigcup_{m\in\omega}~
\prod_{m+1})^{<\omega}\big)^\omega\mid (\lambda_l)_{l\in\omega}\!\in\!\mathcal{S}\wedge\forall l\!\in\!\omega ~~
\lambda_l\! >\! 0\mbox{ is even }\wedge\forall i\! <\lambda_l~~\vert s_l(i)\vert\! =\! l\! +\! 1\big\}\mbox{,}$$
$$\mathcal{J}^c\! :=\!\big\{\beta\!\in\!\mathcal{J}\mid\forall i\!\in\!\omega ~~
\big( s_l(i)0^\infty\big)_{l\in\omega ,\lambda_l>i}\mbox{ converges to some }\gamma_i\!\in\!\mathcal{C}\big\} .$$
Let $c,a,\overline{a}$ be pairwise different not in $\omega$, and $\overline{\overline{a}}\! :=\! a$. We define, for 
$\beta\!\in\!\mathcal{J}$, a countable digraph 
$\mathbb{O}_\beta$ on $\mathcal{K}_{\bf d}\! :=\!\prod_{j\in\omega}~(d_j\cup\{ c,a,\overline{a}\} )$ by\medskip

\leftline{$\mathbb{O}_\beta\! :=\!\{ (c^{l+1}a\overline{a}^\infty ,s_l(0)\overline{a}a^\infty )\mid l\!\in\!\omega\}\cup
\{ (s_l(i)a^{i+1}\overline{a}^\infty ,s_l(i\! +\! 1)\overline{a}^{i+2}a^\infty )\mid l\!\in\!\omega\wedge i\!\leq\!\lambda_l\! -\! 2\} ~\cup$}\smallskip

\rightline{$\{ (s_l(\lambda_l\! -\! 1)a^{\lambda_l}\overline{a}^\infty ,c^{l+1}\overline{a}a^\infty )\mid l\!\in\!\omega\} .$}\medskip

 This allows us to define the graph $\mathbb{G}_\beta\! :=\! s(\mathbb{O}_\beta )$. We then set 
$\mathbb{K}_\beta\! :=\!\overline{\mbox{proj}[\mathbb{G}_\beta ]}^{\mathcal{K}_{\bf d}}$, so that $\mathbb{K}_\beta$ is a 0DMC space and $\mathbb{G}_\beta$ is a graph on $\mathbb{K}_\beta$ whose vertices have degree at most one.

\begin{lemm} \label{chromgen} Let ${\bf d}\!\in\!\mathfrak{C}$, and $\beta\!\in\!\mathcal{J}$. Then 
$(\mathbb{K}_\beta ,\mathbb{G}_\beta )$ has CCN at least three and $\boraone\oplus\bormone$ chromatic number two. If moreover $s_l(2i\! +\!\varepsilon )(0)\! =\!\varepsilon$ for each $l\!\in\!\omega$, each 
$i\!\leq\!\frac{\lambda_l-2}{2}$ and each $\varepsilon\!\in\! 2$, then $(\mathbb{K}_\beta ,\mathbb{G}_\beta )$ has CCN three.\end{lemm}

\noindent\emph{Proof.}\ If $(C,\neg C)$ is a coloring of $\mathbb{G}_\beta$ into clopen subsets of 
$\mathbb{K}_\beta$ which are not empty and satisfy $c^\infty\!\in\! C$, then $N_{c^{j_0}}\!\subseteq\! C$. The compactness of 
$\mathbb{K}_\beta$ gives $l_0\!\geq\! j_0$ with $x,y\!\in\! C$ or $x,y\!\notin\! C$ if $x,y\!\in\!\mathbb{K}_\beta$ and 
$x\vert l_0\! =\! y\vert l_0$. Assume that $l\!\geq\! l_0$. An induction on $i\! <\!\lambda_l$ shows that 
$s_l(i)\overline{a}^{i+1}a^\infty ,s_l(i)a^{i+1}\overline{a}^\infty\!\notin\! C$ if $i$ is even, 
$s_l(i)\overline{a}^{i+1}a^\infty ,s_l(i)a^{i+1}\overline{a}^\infty\!\!\in\! C$ if $i$ is odd, and $c^{l+1}\overline{a}a^\infty\!\notin\! C$, which is the desired contradiction. Thus $\chi_c(\mathbb{K}_\beta ,\mathbb{G}_\beta )\!\geq\! 3$.\medskip
 
 We define an open subset of $\mathbb{K}_\beta$ by 
$O\! :=\!\big\{ x\!\in\!\mathbb{K}_\beta\mid\exists n\!\in\!\omega ~~x(n)\! =\! a\wedge x\vert n\!\in\!\prod_{j<n}~(d_j\cup\{ c\} )\big\}$. The 
$\boraone\oplus\bormone$ partition $(O,\neg O)$ of $\mathbb{K}_\beta$ is a witness for the fact that 
$2\!\leq\!\chi_{\boraone\oplus\bormone}(\mathbb{K}_\beta ,\mathbb{G}_\beta )\!\leq\! 2$.\medskip
 
 If moreover $s_l(2i\! +\!\varepsilon )(0)\! =\!\varepsilon$ for each $l\!\in\!\omega$, each $i\!\leq\!\frac{\lambda_l-2}{2}$ and each 
$\varepsilon\!\in\! 2$, then the clopen partition $(N_c,N_0,N_1)$ of $\mathbb{K}_\beta$ is a witness for the fact that 
$\chi_c(\mathbb{K}_\beta ,\mathbb{G}_\beta )\!\leq\! 3$.\hfill{$\square$}

\begin{lemm} \label{D2gen} Let ${\bf d}\!\in\!\mathfrak{C}$, and $\beta\!\in\!\mathcal{J}$. Then the graph $\mathbb{G}_\beta$ is  $D_2(\bormone )$.\end{lemm}

\noindent\emph{Proof.}\ We check that 
$\overline{\mathbb{G}_\beta }\! =\!\mathbb{G}_\beta\cup\big(\overline{\mathbb{G}_\beta}\cap (\{ c^\infty\}\cup\mathcal{C})^2\big)$. For the left to right inclusion, assume that $(x,y)\!\in\!\overline{\mathbb{G}_\beta}\!\setminus\!\mathbb{G}_\beta$. We may assume that 
$(x,y)$ is the limit of a sequence $\big( (x_m,y_m)\big)_{m\in\omega}$ such that 
$d\big( (x_m,y_m),(\{ c^\infty\}\cup\mathcal{C})^2\big)\!\leq\! 2^{-m}$. We are done since $\{ c^\infty\}\cup\mathcal{C}$ is closed in 
$\mathbb{K}_\beta$. As the first union is the disjoint union of $\mathbb{G}_\beta$ and a closed relation on 
$\mathbb{K}_\beta$, $\mathbb{G}_\beta$ is $D_2(\bormone )$.\hfill{$\square$}\medskip

 The point $c^\infty$ will often be crucial to ensure a big CCN.

\begin{lemm} \label{c} Let ${\bf d}\!\in\!\mathfrak{C}$, $\beta\!\in\!\mathcal{J}$ with $s_l(2i\! +\!\varepsilon )(0)\! =\!\varepsilon$ for each $l\!\in\!\omega$, each $i\!\leq\!\frac{\lambda_l-2}{2}$ and each $\varepsilon\!\in\! 2$, $X$ be a topological space, and $G$ be a digraph on $X$ having CCN at least three such that $(X,G)\preceq_c(\mathbb{K}_\beta ,\mathbb{G}_\beta )$, with witness $\varphi$. Then $c^\infty\!\in\!\varphi [X]$.\end{lemm}

\noindent\emph{Proof.}\ We argue by contradiction. Let 
$C\! :=\!\varphi^{-1}(N_0\cup\bigcup_{l\in\omega}~N_{c^{l+1}\overline{a}})$. Then $(C,\neg C)$ is a coloring of $G$ into clopen sets since $C\! =\!\varphi^{-1}(N_0\cup\bigcup_{l\in\omega}~N_{c^{l+1}\overline{a}}\cup\{ c^\infty\} )$, which is absurd.\hfill{$\square$}\medskip

 We now prove that, for ${\bf d}\! =\! 2^\infty$, $\big( (\mathbb{K}_\beta ,\mathbb{G}_\beta )\big)_{\beta\in\mathcal{J}^c}$ is a 
$\preceq_c$-basis for $\mathfrak{K}$.
  
\begin{them} \label{corcomp'''''''} Let $X$ be a 0DMC space, and $G$ be a graph on $X$. Then exactly one of the following holds:\smallskip

\noindent (1) $(X,G)$ has CCN at most two,\smallskip

\noindent (2) there is $\beta\!\in\!\mathcal{J}^c$ (for ${\bf d}\! =\! 2^\infty$) such that 
$(\mathbb{K}_\beta ,\mathbb{G}_\beta )\preceq_c(X,G)$.\end{them}

\noindent\emph{Proof.}\ By Theorem \ref{corcomp''''''} and Lemma \ref{chromgen}, it is enough to prove that if $\gamma\!\in\!\mathcal{I}$, then we can find $\beta\!\in\!\mathcal{J}^c$ (for ${\bf d}\! =\! 2^\infty$) such that 
$(\mathbb{K}_\beta ,\mathbb{G}_\beta )\preceq_c(\mathbb{K}_\gamma ,\mathbb{G}_\gamma )$. As 
$\mbox{lim}_{k\rightarrow\infty}~\gamma^\varepsilon_k\big(\varepsilon\!\cdot\! (2k)\big)\! =\! 0^\infty$, we can find 
$\delta_\varepsilon\!\in\!\mathcal{S}$ such that 
$0^{\delta_\varepsilon (k)}\!\subseteq\!\gamma^\varepsilon_k\big(\varepsilon\!\cdot\! (2k)\big)$ for each $k\!\in\!\omega$. We define 
$\Delta\!\in\!\mathcal{S}$ by setting $\Delta (k)\! :=\!\mbox{min}\big(\delta (k),\delta_0(k),\delta_1(k)\big)$. Let 
$(k_q)_{q\in\omega}\!\in\!\mathcal{S}$ such that $\Delta (k_q)\! >\! q$ and $k_q\! >\! 0$. We set $\lambda'_q\! :=\! 2k_q$, so that 
$(\lambda'_q)_{q\in\omega}\!\in\!\mathcal{S}$. We then set, for $i\! <\! 2k_q$, $s'_q(i)\! :=\!\gamma^1_{k_q}(i)\vert (q\! +\! 1)$. This allows us to set $\beta'\! :=\!\big(\big( s'_q(i)\big)_{i<2k_q}\big)_{q\in\omega}$. Note that $\beta'\!\in\!\mathcal{J}$. This defines $(\mathbb{K}_{\beta'},\mathbb{G}_{\beta'})$, and we will define $\beta$ later.\medskip

 We now check that $(\mathbb{K}_{\beta'},\mathbb{G}_{\beta'})\preceq_c(\mathbb{K}_\gamma ,\mathbb{G}_\gamma )$. We have to define $\varphi\! :\!\mathbb{K}_{\beta'}\!\rightarrow\!\mathbb{K}_\gamma$.

\vfill\eject
 
 We first define a function 
$\varphi_0\! :\!\mbox{proj}[\mathbb{G}_{\beta'}]\!\rightarrow\!\mbox{proj}[\mathbb{G}_\gamma ]\!\subseteq\!
\mathbb{K}_\gamma$, by setting
$$\left\{\!\!\!\!\!\!\!
\begin{array}{ll}
& \varphi_0(c^{q+1}a\overline{a}^\infty )\! :=\!\gamma^0_{k_q}(0)\mbox{,}\cr
& \varphi_0(c^{q+1}\overline{a}a^\infty )\! :=\!\gamma^1_{k_q}(2k_q)\mbox{,}\cr
& \varphi_0(s'_q(i)a^{i+1}\overline{a}^\infty )\! :=\!\gamma^0_{k_q}(i\! +\! 1)\mbox{ if }i\! <\! 2k_q\mbox{,}\cr
& \varphi_0(s'_q(i)\overline{a}^{i+1}a^\infty )\! :=\!\gamma^1_{k_q}(i)\mbox{ if }i\! <\! 2k_q.
\end{array}\right.$$
Note that $\mathbb{G}_{\beta'}\!\subseteq\! (\varphi_0\!\times\!\varphi_0)^{-1}(\mathbb{G}_\gamma )$. Let us prove that $\varphi_0$ is uniformly continuous on the projection $\mbox{proj}[\mathbb{G}_{\beta'}]$. We set, for 
$x\!\in\!\mbox{proj}[\mathbb{G}_{\beta'}]$, $\delta (x)\! :=\!\mbox{min}\big\{ n\!\in\!\omega\mid x(n)\!\in\!\{ a,\overline{a}\}\big\}$. Note that if 
$\varepsilon\!\in\!\{ a,\overline{a}\}$, then $\varphi_0(c^{q+1}\varepsilon\overline{\varepsilon}^\infty )\vert (q\! +\! 1)\! =\! 0^{q+1}$ and 
$\varphi_0(s'_q(i)\varepsilon^{i+1}\overline{\varepsilon}^\infty )\vert (q\! +\! 1)\! =\! s'_q(i)$, by the choice of $(k_q)_{q\in\omega}$. Let 
$q_0\!\in\!\omega$. We want to find $n\!\in\!\omega$ such that $\varphi_0(x)\vert (q_0\! +\! 1)\! =\!\varphi_0(y)\vert (q_0\! +\! 1)$ if 
$x\vert (n\! +\! 1)\! =\! y\vert (n\! +\! 1)$ and $x,y\!\in\!\mbox{proj}[\mathbb{G}_{\beta'}]$. If $\delta (x),n\! >\! q_0$, then 
$\delta (y)\! >\! q_0$ and $\varphi_0(x)\vert (q_0\! +\! 1)\! =\!\varphi_0(y)\vert (q_0\! +\! 1)$ by the previous facts. Note that there are finitely many $z\!\in\!\mbox{proj}[\mathbb{G}_{\beta'}]$ with $\delta (z)\!\leq\! q_0$. We choose $n\! >\! q_0$ large enough so that $z\! =\! t$ if $z,t\!\in\!\mbox{proj}[\mathbb{G}_{\beta'}]$, $\delta (z),\delta (t)\!\leq\! q_0$ and $z\vert n\! =\! t\vert n$, so that $n$ is as desired. The theorem of extension of uniformly continuous maps (see [Bo, chapter II, \textsection 3, Section 6, Theorem 2]) provides 
$\varphi\! :\!\overline{\mbox{proj}[\mathbb{G}_{\beta'}]}^{\mathcal{K}_{2^\infty}}\!\rightarrow\!\mathbb{K}_\gamma$ continuous extending $\varphi_0$. As $\overline{\mbox{proj}[\mathbb{G}_{\beta'}]}^{\mathcal{K}_{2^\infty}}\! =\!\mathbb{K}_{\beta'}$, the map $\varphi$ is as desired.\smallskip

 It remains to find $\beta\!\in\!\mathcal{J}^c$ with $(\mathbb{K}_\beta ,\mathbb{G}_\beta )\preceq_c(\mathbb{K}_{\beta'},\mathbb{G}_{\beta'})$. By compactness of 
$\{ 0,1,c,a,\overline{a}\}^\omega$, we can find $(\gamma_i)_{i<\lambda'_0}\!\in\! (2^\omega )^{\lambda'_0}$ and 
$(q^0_j)_{j\in\omega}$ strictly increasing such that $q^0_0\! =\! 0$ and, for each $i\! <\!\lambda'_{q^0_0}$, 
$\big( s'_{q^0_j}(i)0^\infty\big)_{j\in\omega}$ converges to $\gamma_i$. Extracting a further subsequence if necessary, we may assume that $s'_{q^0_j}(i)\vert j\! =\!\gamma_i\vert j$ if $i\! <\!\lambda'_{q^0_0}$. We can find 
$(\gamma'_i)_{i<\lambda'_{q^0_1}}$ and $(q^1_j)_{j\in\omega}$ strictly increasing such that ${q^1_0\! =\! q^0_1}$, 
$\{ q^1_j\mid j\!\in\!\omega\}\!\subseteq\!\{ q^0_m\mid m\! >\! 0\}$, and, for each $i\! <\!\lambda'_{q_0^1}$, 
$(s'_{q^1_j}(i)0^\infty )_{j\in\omega}$ converges to $\gamma'_i$. As $(q^1_j)_{j\in\omega}$ is a subsequence of 
$(q^0_m)_{m\in\omega}$, $\gamma'_i\! =\!\gamma_i$ if $i\! <\!\lambda'_{q^0_0}$. For this reason, we may set, for 
$i\! <\!\lambda'_{q_0^1}$, ${\gamma_i\! :=\!\gamma'_i}$. Note that, extracting a further subsequence if necessary, we may assume that $s'_{q^1_j}(i)\vert j\! =\!\gamma_i\vert j$ if ${i\! <\!\lambda'_{q_0^1}}$. Then, inductively, we can find 
$(\gamma_i)_{i<2\lambda'_{q^k_1}}$ and $(q^{k+1}_j)_{j\in\omega}$ strictly increasing with the properties that 
$q^{k+1}_0\! =\! q^k_1$, ${\{ q^{k+1}_j\mid j\!\in\!\omega\}\!\subseteq\!\{ q^k_m\mid m\! >\! 0\}}$, and, for each 
$i\! <\!\lambda'_{q^{k+1}_0}$, $\big( s'_{q^{k+1}_j}(i)0^\infty\big)_{j\in\omega}$ converges to $\gamma_i$. We can also ensure that $s'_{q^{k+1}_j}(i)\vert j\! =\!\gamma_i\vert j$ if $i\! <\!\lambda'_{q^{k+1}_0}$. Note that $q_0^{k+1}\! =\! q^k_1\! >\! q^k_0$. We then set ${\lambda_l\! :=\!\lambda'_{q^{2l+2}_0}}$ and, for 
$i\! <\!\lambda_{l}$, $s_l(i)\! :=\! s'_{q^{2l+2}_0}(i)\vert (l\! +\! 1)$. Fix $i\!\in\!\omega$. Let $l$ minimal such that 
$i\! <\!\lambda'_{q^{l+1}_0}$. Note that $q^{2l+2}_0\! =\! q^{2l+1}_1\! >\! q^{2l+1}_0\! =\! q^{2l}_1\! >\! q^{2l}_0\! =\! 
q^{2l-1}_1\! >\!\ldots\! =\! q^{l+1}_1\! >\! q^{l+1}_0$, 
which shows that $q^{2l+2}_0\! =\! q^{l+1}_{j_{l+1}}$ for some $j_{l+1}\!\geq\! l\! +\! 1$. This gives 
$s'_{q^{l+1}_{j_{l+1}}}(i)\vert j_{l+1}\! =\!\gamma_i\vert j_{l+1}$, and 
${\gamma_i\vert (l\! +\! 1)\! =\! s'_{q_0^{2l+2}}(i)\vert (l\! +\! 1)\! =\! s_l(i)}$. In other words, 
$\big( s_l(i)0^\infty\big)_{\lambda_l>i}$ converges to $\gamma_i$ for each $i\!\in\!\omega$, so that 
$\beta\! :=\!\big(\big( s_l(i)\big)_{i<\lambda_l}\big)_{l\in\omega}\!\in\!\mathcal{J}^c$.\smallskip

 Let us prove that $(\mathbb{K}_\beta ,\mathbb{G}_\beta )\preceq^i_c(\mathbb{K}_{\beta'},\mathbb{G}_{\beta'})$. We define $\psi\! :\!\mathbb{K}_\beta\!\rightarrow\!\mathbb{K}_{\beta'}$, which will be the identity on 
$\{ c^\infty\}\cup (\mathbb{K}_\beta\cap 2^\omega )$. Let $\varepsilon\!\in\!\{ a,\overline{a}\}$, $l\!\in\!\omega$, and 
$i\!\in\!\omega$ with either $i\! =\! 0$ if $s\! =\! c^{l+1}$, or $s\! =\! s_l(i)$. We define $\psi (s\varepsilon^{i+1}\overline{\varepsilon}^\infty )$ in such a way that 
$s\!\subseteq\!\psi (s\varepsilon^{i+1}\overline{\varepsilon}^\infty )$ and 
$\psi (s\varepsilon^{i+1}\overline{\varepsilon}^\infty )$ ends with $\varepsilon^{i+1}\overline{\varepsilon}^\infty$. We set
$$\left\{\!\!\!\!\!\!\!
\begin{array}{ll}
& \psi (c^{l+1}a\overline{a}^\infty )\! :=\! c^{q^{2l+2}_0+1}a\overline{a}^\infty\mbox{,}\cr
& \psi (c^{l+1}\overline{a}a^\infty )\! :=\! c^{q^{2l+2}_0+1}\overline{a}a^\infty\mbox{,}\cr
& \psi (s_l(i)a^{i+1}\overline{a}^\infty )\! :=\! s'_{q^{2l+2}_0}(i)a^{i+1}\overline{a}^\infty\mbox{ if }
i\! <\!\lambda_l\! =\!\lambda'_{q^{2l+2}_0}\mbox{,}\cr
& \psi (s_l(i)\overline{a}^{i+1}a^\infty )\! :=\! s'_{q^{2l+2}_0}(i)\overline{a}^{i+1}a^\infty\mbox{ if }i\! <\!\lambda'_{q^{2l+2}_0}.
\end{array}\right.$$ 
The map $\psi$ is injective continuous as desired.\hfill{$\square$}

\vfill\eject

\noindent\emph{Remark.}\ The $\mathbb{K}_\beta$'s are infinite, so that 
$\big( (\mathbb{K}_\beta ,\mathbb{G}_\beta )\big)_{\beta\in\mathcal{J}^c}$ is not a $\preceq^i_c$-basis, because of the (finite) odd cycles. On the other hand, the $(\mathbb{K}_\beta ,\mathbb{G}_\beta )$'s have $\boraone\oplus\bormone$ chromatic number two, by Lemma \ref{chromgen}. In particular, they have Borel chromatic number two.\medskip

 The next two results will help us to compare the subgraphs of the $\mathbb{G}_\beta$'s.

\begin{lemm} \label{flip'} Let ${\bf d}\!\in\!\mathfrak{C}$, $\beta\!\in\!\mathcal{J}$, $V_{\bf d}\!\subseteq\!\mathbb{K}_\beta$, 
$E_{\bf d}\!\subseteq\!\mathbb{G}_\beta\cap V_{\bf d}^2$, ${\bf d}'$, $\beta'$, $V_{{\bf d}'}$, $E_{{\bf d}'}$ having the corresponding properties, satisfying $(V_{\bf d},E_{\bf d})\preceq_c(V_{{\bf d}'},E_{{\bf d}'})$ with witness $\varphi$, and 
$x\!\in\!\overline{\mbox{proj}[s(E_{\bf d})]}^{\mathbb{K}_\beta}\cap V_{\bf d}\cap\mathcal{C}_{\bf d}$ with 
$\varphi (x)\!\notin\!\mathcal{C}_{{\bf d}'}$. Then we can find $t\!\in\!\bigcup_{l\in\omega}~\prod_l$ such that\smallskip

- either $\varphi [\overline{\mbox{proj}[s(E_{\bf d})]}^{\mathbb{K}_\beta}\cap V_{\bf d}\cap\mathcal{C}_{\bf d}\cap N_t]
\!\subseteq\!\{ c^\infty\}$,\smallskip

- or we can find $\varepsilon\!\in\!\{ a,\overline{a}\}$, ${l,m\!\in\!\omega}$ and 
$s\!\in\!\{ c^{l+1}\}\cup\{ s_l(i)\mid i\! <\!\lambda_l\}$ with the property that 
$\varphi [\overline{\mbox{proj}[s(E_{\bf d})]}^{\mathbb{K}_\beta}\cap V_{\bf d}\cap\mathcal{C}_{\bf d}\cap N_t]\!\subseteq\!
\{ s\varepsilon^{m+1}\overline{\varepsilon}^\infty\}$.\end{lemm}

\noindent\emph{Proof.}\ If $\varphi (x)$ is of the form $s\varepsilon^{m+1}\overline{\varepsilon}^\infty$, then there is 
$m_0\!\in\!\omega$ such that $\varphi (z)\!\supseteq\! s\varepsilon^{m+1}\overline{\varepsilon}$ if $z$ is in 
${V_{\bf d}\cap N_{x\vert m_0}}$, by continuity of $\varphi$. Assume that 
$z\!\in\!\overline{\mbox{proj}[s(E_{\bf d})]}^{\mathbb{K}_\beta}\cap V_{\bf d}\cap\mathcal{C}_{\bf d}\cap N_{x\vert m_0}$. Then 
$z\! =\!\mbox{lim}_{j\rightarrow\infty}~z_j$, where $z_j$ is in $\mbox{proj}[s(E_{\bf d})]\cap N_{x\vert m_0}$, which gives 
$u_j\!\in\! V_{\bf d}$ with $(z_j,u_j)\!\in\! s(E_{\bf d})\!\subseteq\! V_{\bf d}^2$. This implies that 
$\big(\varphi (z_j),\varphi (u_j)\big)\!\in\! s(E_{{\bf d}'})$, $\varphi (z_j)\! =\! s\varepsilon^{m+1}\overline{\varepsilon}^\infty$ and $\varphi (z)\! =\! s\varepsilon^{m+1}\overline{\varepsilon}^\infty$.\medskip

 If now $\varphi (x)\! =\! c^\infty$, then there is $m_0\!\in\!\omega$ such that $\varphi (z)(0)\! =\! c$ if $z\!\in\! V_{\bf d}\cap N_{x\vert m_0}$, by continuity of $\varphi$. As in the previous case, we get $(z_j,u_j)$. This time, $\varphi (z_j)(0)\! =\! c$, so that $\varphi (z_j)$ is of the form $c^{k_j+1}\varepsilon\overline{\varepsilon}^\infty$ and 
$\varphi (z)\!\in\!\{ c^\infty\}\cup\big\{ c^{k+1}\varepsilon\overline{\varepsilon}^\infty\mid 
k\!\in\!\omega\wedge\varepsilon\!\in\!\{ a,\overline{a}\}\big\}$. By the previous point, we may assume that $\varphi (z)\! =\! c^\infty$.
\hfill{$\square$}\medskip

 We get a condition sufficient to send $\mathcal{C}_{\bf d}$ into $\mathcal{C}_{{\bf d}'}$.

\begin{lemm} \label{incl} Let ${\bf d}\!\in\!\mathfrak{C}$, $\beta\!\in\!\mathcal{J}$, $V_{\bf d}\!\subseteq\!\mathbb{K}_\beta$, 
$E_{\bf d}\!\subseteq\!\mathbb{G}_\beta\cap V_{\bf d}^2$, ${\bf d}'$, $\beta'$, $V_{{\bf d}'}$, $E_{{\bf d}'}$ having the corresponding properties, satisfying $\mathcal{C}_{\bf d}\!\subseteq\!
\overline{V_{\bf d}\cap\mathcal{C}_{\bf d}}^{\mathcal{C}_{\bf d}}\cap\overline{\mbox{proj}[s(E_{\bf d})]}^{\mathbb{K}_\beta}$, and also $(V_{\bf d},E_{\bf d})\preceq^i_c(V_{{\bf d}'},E_{{\bf d}'})$ with witness $\varphi$. Then 
$\varphi [V_{\bf d}\cap\mathcal{C}_{\bf d}]\!\subseteq\!\mathcal{C}_{{\bf d}'}$.\end{lemm}

\noindent\emph{Proof.}\ Towards a contradiction, suppose that there is $x\!\in\! V_{\bf d}\cap\mathcal{C}_{\bf d}$ with 
$\varphi (x)\!\notin\!\mathcal{C}_{{\bf d}'}$. Note that $x$ is in $\overline{\mbox{proj}[s(E_{\bf d})]}^{\mathbb{K}_\beta}$. Lemma \ref{flip'} provides $t\!\in\!\bigcup_{l\in\omega}~\prod_l$ such that $\varphi [V_{\bf d}\cap\mathcal{C}_{\bf d}\cap N_t]$ has at most one element, which contradicts the injectivity of $\varphi$ since 
$\mathcal{C}_{\bf d}\!\subseteq\!\overline{V_{\bf d}\cap\mathcal{C}_{\bf d}}^{\mathcal{C}_{\bf d}}$.\hfill{$\square$}\medskip

We give a condition sufficient to send $c^\infty$ to itself.

\begin{lemm} \label{fcc} Let ${\bf d}\!\in\!\mathfrak{C}$, $\beta\!\in\!\mathcal{J}$ with 
$s_l(2i\! +\!\varepsilon )(0)\! =\!\varepsilon$ for each $l\!\in\!\omega$, each $i\!\leq\!\frac{\lambda_l -2}{2}$ and each $\varepsilon\!\in\! 2$, 
$V_{\bf d}\!\subseteq\!\mathbb{K}_\beta$, $E_{\bf d}\!\subseteq\!\mathbb{G}_\beta\cap V_{\bf d}^2$ satisfying  
$\chi_c(V_{\bf d},E_{\bf d})\!\geq\! 3$, ${\bf d}'$, $\beta'$, $V_{{\bf d}'}$, $E_{{\bf d}'}$ having the corresponding properties, satisfying $V_{\bf d}\!\setminus\! (\{ c^\infty\}\cup\mathcal{C}_{\bf d})\!\subseteq\!\mbox{proj}[s(E_{\bf d})]$, and 
$(V_{\bf d},E_{\bf d})\preceq_c(V_{{\bf d}'},E_{{\bf d}'})$ with witness $\varphi$ satisfying 
$\varphi [V_{\bf d}\cap\mathcal{C}_{\bf d}]\!\subseteq\!\mathcal{C}_{{\bf d}'}$. Then 
$\varphi [V_{\bf d}\!\setminus\! (\{ c^\infty\}\cup\mathcal{C}_{\bf d})]\!\subseteq\! 
V_{{\bf d}'}\!\setminus\! (\{ c^\infty\}\cup\mathcal{C}_{{\bf d}'})$ and $\varphi (c^\infty )\! =\! c^\infty$.\end{lemm}

\noindent\emph{Proof.}\ As $V_{\bf d}\!\setminus\! (\{ c^\infty\}\cup\mathcal{C}_{\bf d})\!\subseteq\!\mbox{proj}[s(E_{\bf d})]$, 
$\varphi [V_{\bf d}\!\setminus\! (\{ c^\infty\}\cup\mathcal{C}_{\bf d})]\!\subseteq\!
\mbox{proj}[s(E_{{\bf d}'})]\!\subseteq\! V_{{\bf d}'}\!\setminus\! (\{ c^\infty\}\cup\mathcal{C}_{{\bf d}'})$. It remains to apply Lemma \ref{c}.\hfill{$\square$}

\subsection{$\!\!\!\!\!\!$ Lower bounds}\indent

 In this subsection, we prove Theorem \ref{eq1++}, among other things. We first recover an implication in the style of 
(2) $\Rightarrow$ (1) in Theorem \ref{eq2''''''}.

\vfill\eject

\begin{lemm} \label{chrominter+} Let $X$ be a first countable topological space, and $G$ be a graph on $X$ with the property  that $\Delta (X)\cap\overline{\bigcup_{p\in\omega}~\overline{G}^{2p+1}}$ is not empty. Then $(X,G)$ has no continuous 2-coloring.\end{lemm}

\noindent\emph{Proof.}\ Let $x\!\in\! X$ with $(x,x)\!\in\!\overline{\bigcup_{p\in\omega}~\overline{G}^{2p+1}}$. As $X$ is first countable, we can find $\big( (x_n,y_n)\big)_{n\in\omega}$ in $(X^2)^\omega$ converging to $(x,x)$ with 
$(x_n,y_n)\!\in\!\bigcup_{p\in\omega}~\overline{G}^{2p+1}$. Let $(p_n)_{n\in\omega}\!\in\!\omega^\omega$ with 
$(x_n,y_n)\!\in\!\overline{G}^{2p_n+1}$. Let $(z^n_i)_{i\leq 2p_n+1}\!\in\! X^{2p_n+2}$ with $z^n_0\! =\! x_n$, 
$z^n_{2p_n+1}\! =\! y_n$, and $(z^n_i,z^n_{i+1})\!\in\!\overline{G}$ if $i\!\leq\! 2p_n$. This gives 
$\big( (x^{n,i}_j,y^{n,i}_j)\big)_{j\in\omega}\!\in\! G^\omega$ converging to $(z^n_i,z^n_{i+1})$. If $(C,X\!\setminus\! C)$ is a coloring of $G$ into clopen subsets of $X$ which are not empty and satisfy $x\!\in\! C$, then 
$x_n,y_n\!\in\! C$ for some $n$ large enough. In particular, $x^{n,0}_j\!\in\! C$ if $j$ is large enough. An induction on 
$i\!\leq\! 2p_n\! +\! 1$ shows that $y^{n,i}_j\!\notin\! C$ if $j$ is large enough, $z^n_{i+1}\!\notin\! C$ and 
$x^{n,i+1}_j\!\notin\! C$ if $j$ is large enough when if $i$ is even, and $y^{n,i}_j\!\in\! C$ if $j$ is large enough, 
$z^n_{i+1}\!\in\! C$, $x^{n,i+1}_j\!\in\! C$ if $j$ is large enough when $i$ is odd, which is the desired contradiction. Thus 
$\chi_c(X,G)\!\geq\! 3$.\hfill{$\square$}

\begin{lemm} \label{chromgenm} $(\mathcal{N},\mathbb{G}_m)$ has CCN three and $\boraone\oplus\bormone$ chromatic number two.\end{lemm}

\noindent\emph{Proof.}\ Note  that 
$(c^\infty ,c^\infty )\!\in\!\Delta (\mathcal{N})\cap\overline{\bigcup_{p\in\omega}~\overline{\mathbb{O}_m}^{2p+1}}$. Indeed, if 
$k\!\in\!\omega$, then $(c^{k+1}a^\infty ,c^{k+1}\overline{a}^\infty )$ is in $\overline{\mathbb{O}_m}^{2k+3}$, with witness 
$\big( c^{k+1}a^\infty ,k0^\infty ,k1^\infty ,\ldots ,k(2k\! +\! 1)^\infty ,c^{k+1}\overline{a}^\infty\big)$. Lemma \ref{chrominter+} implies that  
$\chi_c(\mathcal{N},\mathbb{G}_m)\!\geq\! 3$. The clopen partition 
$(N_a\cup N_{\overline{a}}\cup N_c,\bigcup_{k,i\in\omega}~N_{k(2i)},\bigcup_{k,i\in\omega}~N_{k(2i+1)})$ shows that 
$\chi_c(\mathcal{N},\mathbb{G}_m)\! =\! 3$.\medskip
 
 We define an open subset of $\mathcal{N}$ by 
$O\! :=\!\big\{ x\!\in\!\mathcal{N}\mid\exists n\!\in\!\omega ~~x(n)\! =\! a\wedge x\vert n\!\in\! (\omega\cup\{ c\} )^n\big\}$. The 
$\boraone\oplus\bormone$ partition $(O,\neg O)$ of $\mathcal{N}$ is a witness for the fact that 
$\chi_{\boraone\oplus\bormone}(\mathcal{N},\mathbb{G}_m)\! =\! 2$.\hfill{$\square$}

\begin{propo} \label{D2genm} The graph $\mathbb{G}_m$ is $D_2(\bormone )$.\end{propo}

\noindent\emph{Proof.}\ Note that\medskip

\leftline{$\overline{\mathbb{O}_m}\! =\!\mathbb{O}_m\cup\{ (c^{k+1}a^\infty ,k0^\infty )\mid k\!\in\!\omega\}\cup
\{ (ki^\infty ,k(i\! +\! 1)^\infty )\mid k\!\in\!\omega\wedge i\!\leq\! 2k\} ~\cup$}\smallskip

\rightline{$\{ (k(2k\! +\! 1)^\infty ,c^{k+1}\overline{a}^\infty )\mid k\!\in\!\omega\} .$}\medskip

\noindent Thus $\overline{\mathbb{O}_m}$ is the disjoint union of $\mathbb{O}_m$ and a closed relation on 
$\mathcal{N}$, so that $\mathbb{O}_m$ is $D_2(\bormone )$. The proof for $\mathbb{G}_m$ is similar.\hfill{$\square$}

\begin{lemm} \label{eq1+} Let $X$ be a first countable topological space, and $G$ be a graph on $X$. The following are equivalent:\smallskip

\noindent (1) $\Delta (X)\cap\overline{\bigcup_{p\in\omega}~\overline{G}^{2p+1}}\!\not=\!\emptyset$,\smallskip

\noindent (2) $(\mathcal{N},\mathbb{G}_m)\preceq_c(X,G)$.\end{lemm}

\noindent\emph{Proof.}\ (2) $\Rightarrow$ (1) Let $\varphi$ be a witness for the fact that (a) holds. We set 
$x\! :=\!\varphi (c^\infty )$. Let $U$ be an open neighborhood of $x$, so that $\varphi^{-1}(U)\!\times\!\varphi^{-1}(U)$ is an open neighborhood of $(c^\infty ,c^\infty )$ and contains $N_{c^k}\!\times\! N_{c^k}$ for some $k\!\in\!\omega$. As 
$(c^{k+1}a^\infty ,c^{k+1}\overline{a}^\infty )\!\in\!\overline{\mathbb{O}_m}^{2k+3}$, 
$\big(\varphi (c^{k+1}a^\infty ),\varphi (c^{k+1}\overline{a}^\infty )\big)$ is in $\overline{G}^{2k+3}\cap (U\!\times\! U)$. This shows that $(x,x)\!\in\!\overline{\bigcup_{p\in\omega}~\overline{G}^{2p+1}}$.\medskip

\noindent (1) $\Rightarrow$ (2) Let $(x,x)\!\in\!\Delta (X)\cap\overline{\bigcup_{p\in\omega}~\overline{G}^{2p+1}}$, and 
$(U_n)$ be a decreasing countable basis of open neighborhoods of $x$. There is, for each 
$n\!\in\!\omega$, $p_n\!\geq\! 1$ such that $\overline{G}^{2p_n+1}\cap U_n^2\!\not=\!\emptyset$. Note that we may assume that the sequence $(p_n)_{n\in\omega}$ is constant or strictly increasing.

\vfill\eject

 Let $(x_n,y_n)\!\in\!\overline{G}^{2p_n+1}\cap U_n^2$, with witness $(z^n_i)_{i\leq 2p_n+1}$, and 
$\big( (x^{n,i}_j,y^{n,i}_j)\big)_{j\in\omega}\!\in\! G^\omega$ converging to $(z^n_i,z^n_{i+1})$.\medskip

 The map $\varphi\! :\!\mathcal{N}\!\rightarrow\! X$ sends $c^\infty$ to $x$. If $k\! <\! p_0$, then $\varphi$ sends\medskip
 
\noindent - $c^{k+1}a^\infty$ to $z^0_0\! =\! x_0$, $\bigcup_{q\not= a}~N_{c^{k+1}a^{j+1}q}$ to $x^{0,0}_j$,\smallskip

\noindent - $ki^\infty$ to $z^0_{i+1}$ if $i\!\leq\! 2k\! +\! 1$, $N_{ki^{j+1}\overline{a}}$ to $y^{0,i}_j$ if $i\!\leq\! 2k\! +\! 1$, 
$\bigcup_{m\not= i,\overline{a}}~N_{ki^{j+1}m}$ to $x^{0,i+1}_j$ if $i\!\leq\! 2k\! +\! 1$,\smallskip

\noindent - $N_{ki}\cup N_\varepsilon$ to $x$ if $i\! >\! 2k\! +\! 1$ or $i\!\notin\!\omega$, and $\varepsilon\!\in\!\{ a,\overline{a}\}$,\smallskip
 
\noindent - $\bigcup_{m\not= c,a}~\{ c^{k+1}m^\infty\}$ to $z^0_{2k+3}$, $\bigcup_{m\not= c,a,\text{ and }q\not= m}~N_{c^{k+1}m^{j+1}q}$ to $y^{0,2k+2}_j$.\medskip

 If $(p_n)_{n\in\omega}$ is constant and $k\!\geq\! p_0$, then $\varphi$ sends\medskip
 
\noindent - $c^{k+1}a^\infty$ to $z^k_0\! =\! x_k$, $\bigcup_{q\not= a}~N_{c^{k+1}a^{j+1}q}$ to $x^{k,0}_j$,\smallskip
 
\noindent - $ki^\infty$ to $z^k_{i+1}$ if $i\!\leq\! 2p_0\! -\! 2$, $N_{ki^{j+1}\overline{a}}$ to $y^{k,i}_j$ if $i\!\leq\! 2p_0\! -\! 2$, 
$\bigcup_{m\not= i,\overline{a}}~N_{ki^{j+1}m}$ to $x^{k,i+1}_j$ if $i\!\leq\! 2p_0\! -\! 2$,\smallskip

\noindent - $k(2l+1)^\infty$ to $z^k_{2p_0}$ if $p_0\! -\! 1\!\leq\! l\! <\! k$, $\bigcup_{m\not= 2l+1}~N_{k(2l+1)^{j+1}m}$ to $y^{k,2p_0-1}_j$ if $p_0\! -\! 1\!\leq\! l\! <\! k$,\smallskip

\noindent - $k(2l)^\infty$ to $z^k_{2p_0-1}$ if $p_0\! -\! 1\! <\! l\!\leq\! k$, $\bigcup_{m\not= 2l}~N_{k(2l)^{j+1}m}$ to $x^{k,2p_0-1}_j$ if 
$p_0\! -\! 1\! <\! l\!\leq\! k$,\smallskip
 
\noindent - $k(2k+1)^\infty$ to $z^k_{2p_0}$, $N_{k(2k+1)^{j+1}\overline{a}}$ to $y^{k,2p_0-1}_j$, 
$\bigcup_{m\not= 2k+1,\overline{a}}~N_{k(2k+1)^{j+1}m}$ to $x^{k,2p_0}_j$.\smallskip

\noindent - $N_{ki}\cup N_\varepsilon$ to $x$ if $i\! >\! 2k\! +\! 1$ or $i\!\notin\!\omega$, and $\varepsilon\!\in\!\{ a,\overline{a}\}$,\smallskip

\noindent - $\bigcup_{m\not= c,a}~\{ c^{k+1}m^\infty\}$ to $z^k_{2p_0+1}$, $\bigcup_{m\not= c,a,\text{ and }q\not= m}~N_{c^{k+1}m^{j+1}q}$ to $y^{k,2p_0}_j$.\medskip

 If $(p_n)_{n\in\omega}$ is strictly increasing and $p_n\!\leq\! k\! <\! p_{n+1}$, then $\varphi$ sends\medskip
 
\noindent - $c^{k+1}a^\infty$ to $z^n_0\! =\! x_n$, $\bigcup_{q\not= a}~N_{c^{k+1}a^{j+1}q}$ to $x^{n,0}_j$,\smallskip
 
\noindent - $ki^\infty$ to $z^n_{i+1}$ if $i\!\leq\! 2p_n\! -\! 2$, $N_{ki^{j+1}\overline{a}}$ to $y^{n,i}_j$ if $i\!\leq\! 2p_n\! -\! 2$, 
$\bigcup_{m\not= i,\overline{a}}~N_{ki^{j+1}m}$ to $x^{n,i+1}_j$ if $i\!\leq\! 2p_n\! -\! 2$,\smallskip
  
\noindent - $k(2l+1)^\infty$ to $z^n_{2p_n}$ if $p_n\! -\! 1\!\leq\! l\! <\! k$, $\bigcup_{m\not= 2l+1}~N_{k(2l+1)^{j+1}m}$ to $y^{n,2p_n-1}_j$ if 
$p_n\! -\! 1\!\leq\! l\! <\! k$.\smallskip
 
\noindent - $k(2l)^\infty$ to $z^n_{2p_n-1}$ if $p_n\! -\! 1\! <\! l\!\leq\! k$, $\bigcup_{m\not= 2l}~N_{k(2l)^{j+1}m}$ to $x^{n,2p_n-1}_j$ if 
$p_n\! -\! 1\! <\! l\!\leq\! k$,\smallskip
 
\noindent - $k(2k+1)^\infty$ to $z^n_{2p_n}$, $N_{k(2k+1)^{j+1}\overline{a}}$ to $y^{n,2p_n-1}_j$, 
$\bigcup_{m\not= 2k+1,\overline{a}}~N_{k(2k+1)^{j+1}m}$ to $x^{n,2p_n}_j$.\smallskip

\noindent - $N_{ki}\cup N_\varepsilon$ to $x$ if $i\! >\! 2k\! +\! 1$ or $i\!\notin\!\omega$, and $\varepsilon\!\in\!\{ a,\overline{a}\}$,\smallskip

\noindent - $\bigcup_{m\not= c,a}~\{ c^{k+1}m^\infty\}$ to $z^n_{2p_n+1}$, $\bigcup_{m\not= c,a,\text{ and }q\not= m}~N_{c^{k+1}m^{j+1}q}$ to $y^{n,2p_n}_j$.\medskip

 Note that $\varphi$ is a witness for (2).\hfill{$\square$}\medskip

 Theorem \ref{eq1++} is now a consequence of Lemmas \ref{eq1+} and \ref{chrominter+}.\medskip

\noindent\emph{Remark.}\ We saw in the proof of Proposition \ref{D2genm} that\medskip

\leftline{$\overline{\mathbb{O}_m}\! =\!\mathbb{O}_m\cup\{ (c^{k+1}a^\infty ,k0^\infty )\mid k\!\in\!\omega\}\cup
\{ (ki^\infty ,k(i\! +\! 1)^\infty )\mid k\!\in\!\omega\wedge i\!\leq\! 2k\} ~\cup$}\smallskip

\rightline{$\{ (k(2k\! +\! 1)^\infty ,c^{k+1}\overline{a}^\infty )\mid k\!\in\!\omega\} .$}\medskip

\noindent Moreover,
$$\overline{\mbox{proj}[\mathbb{G}_m]}\! =\!\mbox{proj}[\mathbb{G}_m]\cup
\big\{ c^{k+1}\varepsilon^\infty\mid k\!\in\!\omega\wedge\varepsilon\!\in\!\{ a,\overline{a}\}\big\}\cup\{ c^\infty\}\cup
\{ ki^\infty\mid k\!\in\!\omega\wedge i\!\leq\! 2k\! +\! 1\}$$
is a closed countable subset of $\mathcal{N}$. As 
$\Delta (\mathbb{P})\cap\overline{\bigcup_{p\in\omega}~\overline{\mathbb{G}_m}^{2p+1}}$ is empty (its only possible element could be $(c^\infty ,c^\infty )$, this is not the case since the $c^{k+1}\varepsilon^\infty$'s are not in $\mathbb{P}$), 
$(\mathbb{P},\mathbb{G}_m)\prec_c(\mathcal{N},\mathbb{G}_m)$ by Lemma \ref{eq1+}, as announced in the introduction.\medskip

 We now characterize the subgraphs of $\mathbb{G}_m$ having a big CCN. In fact, we will need some generalizations of 
$\mathbb{G}_m$ in the sequel, that we now describe.\medskip
  
\noindent\emph{Notation.}\ We define a copy of $(\mathbb{P},\mathbb{G}_m)$ and subgraphs of it. We set, for 
$\delta\!\in\! 2^\omega$, $\mathbb{G}_\delta\! :=\! s(\mathbb{O}_\delta )$, where\medskip

\leftline{$\mathbb{O}_\delta\! :=\! 
\{ (c^{k+1}0ja^\infty ,k0^{j+2}\overline{a}^\infty )\mid\delta (k)\! =\! 1\wedge j\!\in\!\omega\} ~\cup$}\smallskip

\rightline{$\{ (ki0^{j+1}a^\infty ,k(i\! +\! 1)0^{j+1}\overline{a}^\infty )\mid\delta (k)\! =\! 1\wedge i\!\leq\! 2k\wedge j\!\in\!\omega\} ~\cup$}\smallskip

\rightline{$\{ (k(2k\! +\! 1)0^{j+1}a^\infty ,c^{k+1}1j\overline{a}^\infty )\mid\delta (k)\! =\! 1\wedge j\!\in\!\omega\}\mbox{,}$}\medskip

\noindent and $\mathbb{P}_\delta\! =\!\mbox{proj}[\mathbb{G}_\delta ]\cup\{ c^\infty\}\cup
\{ ki0^\infty\mid\delta (k)\! =\! 1\wedge i\!\leq\! 2k\! +\! 1\}$. Note that the vertices of $\mathbb{G}_\delta$ have degree at most one. We also set $\mathbb{P}_\infty\! :=\!\{\alpha\!\in\! 2^\omega\mid\exists^\infty n\!\in\!\omega ~~\alpha (n)\! =\! 1\}$.

\begin{lemm} \label{copy} $(\mathbb{P}_{1^\infty},\mathbb{G}_{1^\infty})\equiv^i_c(\mathbb{P},\mathbb{G}_m)$.\end{lemm}
 
\noindent\emph{Proof.}\ We define $\varphi\! :\!\mathbb{P}_{1^\infty}\!\rightarrow\!\mathbb{P}$ by $\varphi (c^\infty )\! :=\! c^\infty$, 
$\varphi (ki0^\infty )\! :=\! ki^\infty$, $\varphi (c^{k+1}0ja^\infty )\! :=\! c^{k+1}a^{j+1}\overline{a}^\infty$, 
$\varphi (ki0^{j+1}\varepsilon^\infty )\! :=\! ki^{j+1}\varepsilon^\infty$, and 
$\varphi (c^{k+1}1j\overline{a}^\infty )\! :=\! c^{k+1}\overline{a}^{j+1}a^\infty$. The map $\varphi$ is a witness for the fact that 
$(\mathbb{P}_{1^\infty},\mathbb{G}_{1^\infty})\preceq^i_c(\mathbb{P},\mathbb{G}_m)$, and $\varphi^{-1}$ is a witness for the fact that 
$(\mathbb{P},\mathbb{G}_m)\preceq^i_c(\mathbb{P}_{1^\infty},\mathbb{G}_{1^\infty})$.\hfill{$\square$}\medskip

\begin{lemm} \label{chromalph} Let $\delta\!\in\!\mathbb{P}_\infty$. Then $(\mathbb{P}_\delta ,\mathbb{G}_\delta )$ has CCN three and $\boraone\oplus\bormone$ chromatic number two.\end{lemm}

\noindent\emph{Proof.}\ For the upper bounds, we prove that ${\chi_c(\mathbb{P}_{1^\infty},\mathbb{G}_{1^\infty})\!\leq\! 3}$ and 
$\chi_{\boraone\oplus\bormone}(\mathbb{P}_{1^\infty},\mathbb{G}_{1^\infty})\!\leq\! 2$ since 
$\mathbb{P}_\delta\!\subseteq\!\mathbb{P}_{1^\infty}$ and $\mathbb{G}_\delta\!\subseteq\!\mathbb{G}_{1^\infty}$. This comes from Lemmas \ref{chromgenm} and \ref{copy}. For the first lower bound, towards a contradiction, suppose that there is a clopen subset $C$ of 
$\mathbb{P}_\delta$ with $\mathbb{G}_\delta\cap (C^2\cup (\mathbb{P}_\delta\!\setminus\! C)^2)\! =\!\emptyset$. We may assume that $c^\infty\!\in\! C$, which gives $k_0\!\geq\! 1$ such that $N_{c^{k_0+1}}\cap\mathbb{P}_\delta\!\subseteq\! C$. Assume that 
$k\!\geq\! k_0$ and $\delta (k)\! =\! 1$. As $(c^{k+1}0ja^\infty ,k0^{j+2}\overline{a}^\infty )\!\in\!\mathbb{G}_\delta$ and 
$c^{k+1}0ja^\infty\!\in\! C$, $k0^{j+2}\overline{a}^\infty\!\in\!\mathbb{P}_\delta\!\setminus\! C$. Thus 
$k0^\infty\!\in\!\mathbb{P}_\delta\!\setminus\! C$, which gives $j_0\!\in\!\omega$ such that 
$k0^{j+2}a^\infty\!\in\!\mathbb{P}_\delta\!\setminus\! C$ if $j\!\geq\! j_0$. As 
$(k0^{j+2}a^\infty ,k10^{j+1}\overline{a}^\infty )\!\in\!\mathbb{G}_\delta$, 
$k10^{j+1}\overline{a}^\infty\!\in\! C$ if $j\!\geq\! j_0$. Thus $k10^\infty\!\in\! C$. An induction on 
$i\!\leq\! 2k\! +\! 1$ shows that $ki0^\infty\!\in\! C$ if $i$ is odd, and $ki0^\infty\!\in\!\mathbb{P}_\delta\!\setminus\! C$ if $i$ is even. This implies that $c^{k+1}1j\overline{a}^\infty\!\in\!\mathbb{P}_\delta\!\setminus\! C$ if $j$ is large enough, which is absurd.\hfill{$\square$}

\begin{propo} \label{D2alph} Let $\delta\!\in\! 2^\omega$. Then $\mathbb{G}_\delta$ is a $D_2(\bormone )$ graph on the 0DP space $\mathbb{P}_\delta$.\end{propo}

\noindent\emph{Proof.}\ As $\mathbb{P}_\delta$ is a closed subset of $\mathcal{N}$, it is a 0DP space. Note that
$$\overline{\mathbb{O}_\delta}\! =\!\mathbb{O}_\delta\cup
\big\{\big( ki0^\infty ,k(i\! +\! 1)0^\infty\big)\mid\delta (k)\! =\! 1\wedge i\!\leq\! 2k\big\} .$$
Thus $\overline{\mathbb{O}_\delta}$ is the disjoint union of $\mathbb{O}_\delta$ and a closed relation on 
$\mathcal{N}$, so that $\mathbb{O}_\delta$ is $D_2(\bormone )$. The proof for $\mathbb{G}_\delta$ is similar.\hfill{$\square$}\medskip

 We now characterize the subgraphs of $\mathbb{G}_\delta$ having a big CCN.

\begin{lemm} \label{charsubbeta} Let $V\!\subseteq\!\mathbb{P}_\delta$, and 
$E\!\subseteq\!\mathbb{G}_\delta\cap V^2$. The following are equivalent:\smallskip

\noindent (1) the digraph $(V,E)$ has CCN at least three,\smallskip

\noindent (2) $c^\infty\!\in\! V$ and there is $I\!\subseteq\!\{ k\!\in\!\omega\mid\delta (k)\! =\! 1\}$ infinite such that, for each 
$k\!\in\! I$,\smallskip

(a) $\forall i\!\leq\! 2k\! +\! 1~~ki0^\infty\!\in\! V$,\smallskip 

(b) $\exists^\infty j\!\in\!\omega ~~(c^{k+1}0ja^\infty ,k0^{j+2}\overline{a}^\infty )\!\in\! s(E)$,\smallskip 

(c) $\forall i\!\leq\! 2k~~\exists^\infty j\!\in\!\omega ~~(ki0^{j+1}a^\infty ,k(i\! +\! 1)0^{j+1}\overline{a}^\infty )\!\in\! s(E)$,\smallskip 

(d) $\exists^\infty j\!\in\!\omega ~~(k(2k\! +\! 1)0^{j+1}a^\infty ,c^{k+1}1j\overline{a}^\infty )\!\in\! s(E)$.\end{lemm}
 
\noindent\emph{Proof.}\ Note first that 
$\chi_c(V,E)\!\leq\!\chi_c\big(V,s(E)\big)\!\leq\!\chi_c(\mathbb{P}_\delta ,\mathbb{G}_\delta )\! =\! 3$, by Lemma \ref{chromalph}. We may and will assume that $E\! =\! s(E)$ is a graph.\medskip

\noindent (2) $\Rightarrow$ (1) Towards a contradiction, suppose that there is a clopen subset $C$ of $V$ with the property that 
$E\cap (C^2\cup (V\!\setminus\! C)^2)\! =\!\emptyset$. We may assume that $c^\infty\!\in\! C$, which gives 
$k_0\!\in\!\omega$ such that $N_{c^{k_0+1}}\cap V\!\subseteq\! C$. Assume that $k\!\geq\! k_0$ is in $I$. By (b), there are infinitely many $j$'s with $(c^{k+1}0ja^\infty ,k0^{j+2}\overline{a}^\infty )\!\in\! E\!\subseteq\! V^2$, which implies that 
$c^{k+1}0ja^\infty\!\in\! C$ and $k0^{j+2}\overline{a}^\infty\!\notin\! C$. By (a), 
$k0^\infty\!\in\! V$, so that $k0^\infty\!\in\! V\!\setminus\! C$. By (c), there are infinitely many $j$'s with the property that 
$(k0^{j+2}a^\infty ,k10^{j+1}\overline{a}^\infty )\!\in\! E\!\subseteq\! V^2$, so that we may assume that 
$k0^{j+2}a^\infty\!\in\! V\!\setminus\! C$ and thus $k10^{j+1}\overline{a}^\infty\!\in\! C$. By (a) again, $k10^\infty\!\in\! V$, so that $k10^\infty\!\in\! C$. An induction on $i\!\leq\! 2k\! +\! 1$ shows that $ki0^\infty\!\in\! V\!\setminus\! C$ if $i$ is even, and 
$ki0^\infty\!\in\! C$ if $i$ is odd. By (d), this gives infinitely many $j$'s such that 
$c^{k+1}1j\overline{a}^\infty\!\in\! V\!\setminus\! C$, which is the desired contradiction.\medskip

\noindent (1) $\Rightarrow$ (2) If $c^\infty\!\notin\! V$, then we set 
$C\! :=\!\big(\bigcup_{k\in\omega}~(\bigcup_{i\leq 2k+1\text{ even}}~N_{ki}\cup N_{c^{k+1}1})\big)\cap V$. Then 
$(C,V\!\setminus\! C)$ is a coloring of $E$ into clopen sets, which is absurd. If (2) does not hold, then there is 
$k_0\!\in\!\omega$ such that one of the properties (a)-(d) does not hold if $k\!\geq\! k_0$. We will use the following notation.\medskip

\noindent - If (a) does not hold, then $i_k\!\leq\! 2k\! +\! 1$ will be minimal with $ki_k0^\infty\!\notin\! V$,\smallskip

\noindent - If (b) does not hold, then $j_k\!\in\!\omega$ will be minimal such that 
$(c^{k+1}0ja^\infty ,k0^{j+2}\overline{a}^\infty )\!\notin\! E$ if $j\!\geq\! j_k$,\smallskip

\noindent - If (c) does not hold, then $i_k\!\leq\! 2k$ and $j_k\!\in\!\omega$ will be minimal such that 
$(ki_k0^{j+1}a^\infty ,k(i_k\! +\! 1)0^{j+1}\overline{a}^\infty )$ is not in $E$ if $j\!\geq\! j_k$,\smallskip

\noindent - If (d) does not hold, then $j_k\!\in\!\omega$ will be minimal such that 
$(k(2k\! +\! 1)0^{j+1}a^\infty ,c^{k+1}1j\overline{a}^\infty )\!\notin\! E$ if $j\!\geq\! j_k$.\medskip

 We then set\medskip

\leftline{$C'\! :=\!\big(\bigcup_{k<k_0}~(\bigcup_{i\leq 2k+1\text{ even}}~N_{ki}\cup N_{c^{k+1}1})~\cup$}\smallskip

\rightline{$\bigcup_{k\geq k_0,\neg (a),i_k\text{ even}}~(\bigcup_{i<i_k\text{ even}}~
N_{ki}\cup\bigcup_{j\in\omega}~N_{ki_k0^{j+1}\overline{a}}\cup\bigcup_{i_k<i\leq 2k+1\text{ odd}}~N_{ki})~\cup$}\smallskip

\rightline{$\bigcup_{k\geq k_0,\neg (a),i_k\text{ odd}}~(\bigcup_{i<i_k\text{ even}}~N_{ki}\cup
\bigcup_{j\in\omega}~N_{ki_k0^{j+1}a}\cup\bigcup_{i_k<i\leq 2k+1\text{ odd}}~N_{ki})~\cup$}\smallskip

\rightline{$\bigcup_{k\geq k_0,(a),\neg (b)}~(\bigcup_{j<j_k}~N_{k0^{j+2}\overline{a}}\cup
\bigcup_{i\leq 2k+1\text{ odd}}~N_{ki})~\cup$}\smallskip

\rightline{$\bigcup_{k\geq k_0,(a),(b),\neg (c),i_k\text{ even}}~\big(\bigcup_{i<i_k\text{ even}}~N_{ki}\cup
\big( N_{ki_k}\!\setminus\! (\bigcup_{j<j_k}~N_{ki_k0^{j+1}a})\big)\cup\bigcup_{i_k<i\leq 2k+1\text{ odd}}~N_{ki}\big) ~\cup$}\smallskip

\rightline{$\bigcup_{k\geq k_0,(a),(b),\neg (c),i_k\text{ odd}}~\big(\bigcup_{i<i_k\text{ even}}~N_{ki}\cup
\bigcup_{j<j_k}~N_{ki_k0^{j+1}a}\cup\bigcup_{i_k<i\leq 2k+1\text{ odd}}~N_{ki}\big) ~\cup$}\smallskip

\rightline{$\bigcup_{k\geq k_0,(a),(b),(c),\neg (d)}~(\bigcup_{i\leq 2k+1\text{ even}}~N_{ki}\cup
\bigcup_{j<j_k}~N_{k(2k+1)0^{j+1}a})\big)\cap V.$}\medskip

\noindent Then $(C',V\!\setminus\! C')$ is a coloring of $E$ into clopen sets, which is absurd.\hfill{$\square$}

\section{$\!\!\!\!\!\!$ General graphs on a 0DMS space}\label{Polish}

\emph{Remark.}\ We study the limits of Theorem \ref{eq2''''''}. In its proof, we used the compactness of 
$\mathcal{C}_{\bf d}$. This is essential. Indeed, if we replace $\mathcal{C}_{\bf d}$ with $\omega^\omega$ or $\mathcal{N}$, then the notation $^nG$ still makes sense and the following hold. The implications $\mbox{(4) }\Rightarrow\mbox{ (1)}$ and 
$\mbox{(2) }\Rightarrow\mbox{ (3)}$ still hold, with the same proof. Also, the implication $\mbox{(4) }\Rightarrow\mbox{ (2)}$ still holds, using uniform continuity.\medskip

\noindent (a) The implication (1) $\Rightarrow$ (4) does not hold. Indeed, if $X$ is a 0DMC space and $G$ is a graph on $X$ with CCN at least three, then $(X,G)$ is not $\preceq_c$-below $(\mathcal{N},\mathbb{G}_m)$. Indeed, we argue by contradiction to see that, which gives a continuous map $\varphi\! :\! X\!\rightarrow\!\mathcal{N}$. We set 
$V\! :=\!\varphi [X]$ and $E\! :=\! (\varphi\!\times\!\varphi )[G]$, so that the graph $(V,E)$ has CCN three.

\vfill\eject

 The compactness of $X$ implies that the first coordinate of the elements of 
$V\!\setminus\! (N_c\cup N_a\cup N_{\overline{a}})$ is bounded by some natural number $k_0$. We set 
$C\! :=\!\big(\bigcup_{k\leq k_0}~(\bigcup_{i\leq 2k+1\text{ even}}~N_{ki}\cup N_{c^{k+1}\overline{a}})\big)\cap V$, so that 
$C$ is a clopen subset of $V$ and $E\cap(C^2\cup (V\!\setminus\! C)^2)\! =\!\emptyset$, which is the desired contradiction. In particular, if $(\mathbb{A}_1,\mathbb{G}_1)$ exists, then $\mathbb{A}_1$ cannot be compact, and 
$\big( (\mathbb{K}_\alpha ,\mathbb{G}_\alpha )\big)_{\alpha\in 2^\omega}$ given by Theorem \ref{corcomp''''''} is not a 
$\preceq_c$-basis for the class of countable graphs on a 0DP space with CCN at least three.\medskip

\noindent (b) The implication (2) $\Rightarrow$ (4) does not hold. Indeed, note that 
$(c^\infty ,c^\infty )\!\in\!\bigcap_{n\in\omega}~(\bigcup_{p\in\omega}~(^n\mathbb{G}_m)^{2p+1})$.\medskip

\noindent (c) The implication (2) $\Rightarrow$ (1) does not hold. Indeed, consider the following countable graph on 
$\omega^\omega$:\medskip
 
\leftline{$\mathbb{T}\! :=\! s\big(\{ (0^{2k+1}1^\infty ,(2k\! +\! 2)0^\infty )\mid k\!\in\!\omega\}\cup
\big\{\big( (2k\! +\! 2)i0^k1^\infty ,(2k\! +\! 2)(i\! +\! 1)0^\infty\big)\mid k\!\in\!\omega\wedge i\!\leq\! 2k\big\} ~\cup$}\smallskip

\rightline{$\big\{\big( (2k\! +\! 2)(2k\! +\! 1)0^k1^\infty ,0^{2k+2}1^\infty\big)\mid k\!\in\!\omega\big\}\big) .$}\medskip

\noindent Then $(0^\infty ,0^\infty )\!\in\!\Delta (\omega^\omega )\cap\bigcap_{n\in\omega}~(^n\mathbb{T})^{2n+1}$. We set 
$C\! :=\! N_0\cup\bigcup_{k\in\omega ,j\leq 2k}~N_{(2k+2)(j+1)0^{k+1}}$. Then 
$\mathbb{T}\cap (C^2\cup (\neg C)^2)\! =\!\emptyset$ and $C$ is a clopen subset of $\omega^\omega$, so that 
$\chi_c(\omega^\omega ,\mathbb{T})\! =\! 2$.\medskip

 We now turn to the proof of Theorem \ref{absmin}.
  
\begin{lem} \label{generalab} Let $(Q,\leq)$ be quasi-order for which there is $q_0\!\in\! Q$ such that, for any $q\!\leq\! q_0$, there are 
$\leq$-incomparable $p_0,p_1\!\leq\! q$. Then $(Q,\leq)$ has no antichain basis.\end{lem}

\noindent\emph{Proof.}\ Towards a contradiction, suppose that there is an antichain basis $B$ for $(Q,\leq)$. As $B$ is a basis, there is 
$q\!\in\! B$ with $q\!\leq\! q_0$. Our assumption gives $p_0,p_1\!\leq\! q$ with $p_0\!\perp\! p_1$. As $B$ is a basis, there is, for each $\varepsilon\!\in\! 2$, $q_\varepsilon\!\in\! B$ with $q_\varepsilon\!\leq\! p_\varepsilon$. As $q_\varepsilon\!\leq\! q$ are in the antichain $B$, $q_0\! =\! q\! =\! q_1$. Thus $p_0\!\leq\! q\! =\! q_1\!\leq\! p_1$ and $p_0\!\leq\! p_1$, which contradicts the 
$\leq$-incomparability of  $p_0,p_1$.\hfill{$\square$}\medskip

 Recall the graph $(\mathbb{P}_\delta ,\mathbb{G}_\delta )$ defined before Lemma \ref{copy}.

\begin{lem} \label{below0} Let $\delta\!\in\! 2^\omega$, and $G$ be a graph on a 0DMS space $Z$, with CCN at least three and satisfying 
$(Z,G)\preceq^i_c(\mathbb{P}_\delta ,\mathbb{G}_\delta )$. Then there is $\delta'\!\in\!\mathbb{P}_\infty$ such that 
$\{ k\!\in\!\omega\mid\delta'(k)\! =\! 1\}\!\subseteq\!\{ k\!\in\!\omega\mid\delta (k)\! =\! 1\}$ and 
$(\mathbb{P}_{\delta'},\mathbb{G}_{\delta'})\preceq^i_c(Z,G)$.\end{lem}

\noindent\emph{Proof.}\ Assume that $(Z,G)\preceq^i_c(\mathbb{P}_\delta ,\mathbb{G}_\delta )$, with witness $\varphi$. We set 
$V\! :=\!\varphi [Z]$ and $E\! :=\! (\varphi\!\times\!\varphi )[G]$, so that, by Lemma \ref{charsubbeta}, $c^\infty\!\in\! V$ and the set 
$I\!\subseteq\!\{ k\!\in\!\omega\mid\delta (k)\! =\! 1\}$ of $k$'s satisfying (a)-(d) is infinite. We set $(\eta_a,\eta_{\overline{a}})\! :=\! (0,1)$, and define\medskip

\noindent - a singleton $\{\nu\}\! :=\!\varphi^{-1}(\{ c^\infty\} )$,\smallskip

\noindent - singletons $N^{k,i}\! :=\!\{ n^{k,i}\}\! :=\!\varphi^{-1}(\{ ki0^\infty\} )$ (for $k\!\in\! I$ and $i\!\leq\! 2k\! +\! 1$),\smallskip

\noindent - infinite sets $J^a_k\! :=\! J^{\overline{a}}_{k,0}\! :=\!
\{ j\!\in\!\omega\mid (c^{k+1}0ja^\infty ,k0^{j+2}\overline{a}^\infty )\!\in\! E\}$ and
$$J^{\overline{a}}_k\! :=\! J^a_{k,2k+1}\! :=\!
\{ j\!\in\!\omega\mid (k(2k\! +\! 1)0^{j+1}a^\infty ,c^{k+1}1j\overline{a}^\infty )\!\in\! E\}$$ 
(for $k\!\in\! I$),\smallskip

\noindent - infinite sets $J^a_{k,i}\! :=\! J^{\overline{a}}_{k,i+1}\! :=\!
\{ j\!\in\!\omega\mid (ki0^{j+1}a^\infty ,k(i\! +\! 1)0^{j+1}\overline{a}^\infty )\!\in\! E\}$ (for $k\!\in\! I$ and $i\!\leq\! 2k$),\smallskip

\noindent - singletons 
$Z^{k,\varepsilon ,j}\! :=\!\{ z^{k,\varepsilon ,j}\}\! :=\!\varphi^{-1}(\{ c^{k+1}\eta_\varepsilon j\varepsilon^\infty\} )$ (for 
$k\!\in\! I$, $\varepsilon\!\in\!\{ a,\overline{a}\}$ and $j\!\in\! J^\varepsilon_k$),\smallskip

\noindent - singletons $Z^{k,i,j,\varepsilon}\! :=\!\{ z^{k,i,j,\varepsilon}\}\! :=\!\varphi^{-1}(\{ ki0^{j+1}\varepsilon^\infty\} )$ (for 
$k\!\in\! I$, $i\!\leq\! 2k\! +\! 1$, $j\!\in\! J^\varepsilon_{k,i}$, and $\varepsilon\!\in\!\{ a,\overline{a}\}$).\medskip

 By [K, 7.8], we may assume that $Z\!\subseteq\!\omega^\omega$. We set
$$I'\! :=\!\big\{ k\!\in\! I\mid\forall i\!\leq\! 2k~~
(n^{k,i},n^{k,i+1})\!\in\!\overline{\{ (z^{k,i,j,a},z^{k,i+1,j,\overline{a}})\mid j\!\in\! J^a_{k,i}\}}\big\} .$$ 

 If $k\!\in\! I\!\setminus\! I'$, then $i_k\!\leq\! 2k$ will be minimal with 
$(n^{k,i_k},n^{k,i_k+1})\!\notin\!\overline{\{ (z^{k,i_k,j,a},z^{k,i_k+1,j,\overline{a}})\mid j\!\in\! J^a_{k,i_k}\}}$. This gives $l_k\!\in\!\omega$ with $(N_{n^{k,i_k}\vert l_k}\!\times\! N_{n^{k,i_k+1}\vert l_k})\cap\{ (z^{k,i_k,j,a},z^{k,i_k+1,j,\overline{a}})\mid 
j\!\in\! J^a_{k,i_k}\}\! =\!\emptyset$. If we set $J^{a,-}_{k,i_k}\! :=\!\{ j\!\in\! J^a_{k,i_k}\mid z^{k,i_k,j,a}\!\notin\! N_{n^{k,i_k}\vert l_k}\}$ and 
$J^{a,+}_{k,i_k}\! :=\!\{ j\!\in\! J^a_{k,i_k}\mid z^{k,i_k+1,j,\overline{a}}\!\notin\! N_{n^{k,i_k+1}\vert l_k}\}$, then 
$J^a_{k,i_k}\! =\! J^{a,-}_{k,i_k}\cup J^{a,+}_{k,i_k}$.\medskip

 Let us prove that\medskip
 
\leftline{$\forall l\!\in\!\omega ~~\exists^\infty k\!\in\! I'~~\forall p\!\in\!\omega ~~\exists j\!\in\! J^a_k~~
Z^{k,a,j}\cap N_{\nu\vert l},Z^{k,0,j,\overline{a}}\cap N_{n^{k,0}\vert p}\!\not=\!\emptyset ~\wedge$}\smallskip

\rightline{$\forall p\!\in\!\omega ~~\exists j\!\in\! J^{\overline{a}}_k~~
Z^{k,2k+1,j,a}\cap N_{n^{k,2k+1}\vert p},Z^{k,\overline{a},j}\cap N_{\nu\vert l}\!\not=\!\emptyset .$}\medskip
 
\noindent Towards a contradiction, suppose that we can find $l_0,k_0\!\in\!\omega$ such that, for each $k\!\geq\! k_0$ in $I'$, either there is 
$p^a_k\!\in\!\omega$ such that, for $j\!\in\! J^a_k$, $Z^{k,a,j}\cap N_{\nu\vert l_0}\! =\!\emptyset$ or 
$Z^{k,0,j,\overline{a}}\cap N_{n^{k,0}\vert p^a_k}\! =\!\emptyset$, or there is $p^{\overline{a}}_k\!\in\!\omega$ such that, for 
$j\!\in\! J^{\overline{a}}_k$, $Z^{k,2k+1,j,a}\cap N_{n^{k,2k+1}\vert p^{\overline{a}}_k}\! =\!\emptyset$ or 
$Z^{k,\overline{a},j}\cap N_{\nu\vert l_0}\! =\!\emptyset$. If $k\!\geq\! k_0$ and (b) from Lemma \ref{charsubbeta} does not hold, then $j_k$ will be minimal such that $(c^{k+1}0ja^\infty ,k0^{j+2}\overline{a}^\infty )\!\notin\! s(E)$ if $j\!\geq\! j_k$.\medskip

 We set, for $\varepsilon\!\in\!\{ a,\overline{a}\}$, $S^\varepsilon\! :=\!\{ k\!\in\! I'\mid k\!\geq\! k_0\wedge p^\varepsilon_k\mbox{ exists}\}$. We also set, for $k\!\in\! S^a$, 
$$S^{k,a}\! :=\!\{ j\!\in\! J^a_k\mid Z^{k,0,j,\overline{a}}\cap N_{n^{k,0}\vert p^a_k}\! =\!\emptyset\}$$ 
and, for $k\!\in\! S^{\overline{a}}$,  
$S^{k,\overline{a}}\! :=\!\{ j\!\in\! J^{\overline{a}}_k\mid Z^{k,2k+1,j,a}\cap N_{n^{k,2k+1}\vert p^{\overline{a}}_k}\! =\!\emptyset\}$. We then set\medskip

\leftline{$C\! :=\!\bigcup_{k<k_0}~\varphi^{-1}(\bigcup_{i\leq 2k+1\text{ even}}~N_{ki}\cup N_{c^{k+1}1})~\cup$}\smallskip

\rightline{$\bigcup_{k\geq k_0,\neg (a),i_k\text{ even}}~\varphi^{-1}(\bigcup_{i<i_k\text{ even}}~
N_{ki}\cup\bigcup_{j\in\omega}~N_{ki_k0^{j+1}\overline{a}}\cup\bigcup_{i_k<i\leq 2k+1\text{ odd}}~N_{ki})~\cup$}\smallskip

\rightline{$\bigcup_{k\geq k_0,\neg (a),i_k\text{ odd}}~\varphi^{-1}(\bigcup_{i<i_k\text{ even}}~N_{ki}\cup
\bigcup_{j\in\omega}~N_{ki_k0^{j+1}a}\cup\bigcup_{i_k<i\leq 2k+1\text{ odd}}~N_{ki})~\cup$}\smallskip

\rightline{$\bigcup_{k\geq k_0,(a),\neg (b)}~\varphi^{-1}(\bigcup_{j<j_k}~N_{k0^{j+2}\overline{a}}\cup
\bigcup_{i\leq 2k+1\text{ odd}}~N_{ki})~\cup$}\smallskip

\rightline{$\bigcup_{k\geq k_0,(a),(b),\neg (c),i_k\text{ even}}~\varphi^{-1}\big(\bigcup_{i<i_k\text{ even}}~N_{ki}\cup
\big( N_{ki_k}\!\setminus\! (\bigcup_{j<j_k}~N_{ki_k0^{j+1}a})\big)\cup\bigcup_{i_k<i\leq 2k+1\text{ odd}}~N_{ki}\big) ~\cup$}\smallskip

\rightline{$\bigcup_{k\geq k_0,(a),(b),\neg (c),i_k\text{ odd}}~\varphi^{-1}\big(\bigcup_{i<i_k\text{ even}}~N_{ki}\cup
\bigcup_{j<j_k}~N_{ki_k0^{j+1}a}\cup\bigcup_{i_k<i\leq 2k+1\text{ odd}}~N_{ki}\big) ~\cup$}\smallskip

\rightline{$\bigcup_{k\geq k_0,(a),(b),(c),\neg (d)}~\varphi^{-1}(\bigcup_{i\leq 2k+1\text{ even}}~N_{ki}\cup
\bigcup_{j<j_k}~N_{k(2k+1)0^{j+1}a})~\cup$}\smallskip

\rightline{$\bigcup_{k\geq k_0,k\in I\setminus I',i_k\text{ even}}~\big(\varphi^{-1}(
\bigcup_{i<i_k\text{ even}}~N_{ki}\cup\bigcup_{j\in\omega}~N_{ki_k0^{j+1}\overline{a}}\cup
\bigcup_{j\notin J^{a,-}_{k,i_k}}~N_{ki_k0^{j+1}a}~\cup$}\smallskip

\rightline{$\bigcup_{j\in\omega}~N_{k(i_k+1)0^{j+1}a}\cup
\bigcup_{j\notin J^{a,+}_{k,i_k}}~N_{k(i_k+1)0^{j+1}\overline{a}}\cup
\bigcup_{i_k+1<i\leq 2k+1\text{ odd}}~N_{ki})\cup N_{n^{k,i_k}\vert l_k}\cup N_{n^{k,i_k+1}\vert l_k}\big) ~\cup$}\smallskip

\rightline{$\bigcup_{k\geq k_0,k\in I\setminus I',i_k\text{ odd}}~\big(\varphi^{-1}(
\bigcup_{i<i_k\text{ even}}~N_{ki}\cup\bigcup_{j\notin J^{a,-}_{k,i_k}}~N_{ki_k0^{j+1}a}~\cup$}\smallskip

\rightline{$\bigcup_{j\notin J^{a,+}_{k,i_k}}~N_{k(i_k+1)0^{j+1}\overline{a}}\cup
\bigcup_{i_k+1<i\leq 2k+1\text{ odd}}~N_{ki})\cup N_{n^{k,i_k}\vert l_k}\cup N_{n^{k,i_k+1}\vert l_k}\big) ~\cup$}\smallskip

\rightline{$\bigcup_{k\in S^{\overline{a}}}~\varphi^{-1}(\bigcup_{i\leq 2k+1\text{ even}}~N_{ki}\cup
\bigcup_{j\in S^{k,\overline{a}}}~N_{k(2k+1)0^{j+1}a}\cup\bigcup_{j\in J^{\overline{a}}_k\setminus S^{k,\overline{a}}}~N_{c^{k+1}1j\overline{a}})~\cup$}\smallskip

\rightline{$\bigcup_{k\in S^a\setminus S^{\overline{a}}}~\varphi^{-1}(\bigcup_{j\in J^a_k\setminus S^{k,a}}~N_{c^{k+1}0ja}\cup
\bigcup_{j\in S^{k,a}}~N_{k0^{j+2}\overline{a}}\cup\bigcup_{i\leq 2k+1\text{ odd}}~N_{ki}) .$}\medskip

\noindent As $\varphi^{-1}(N_k\cap\mathbb{P}_\delta )$ is a clopen subset of $Z$, the only possible limit point of 
$(z^{k,i,j,\varepsilon})_{j\in\omega}$ is $n^{k,i}$. Also, the only possible limit point of $(z^{k,\varepsilon ,j})_{j,k\in\omega}$ is $\nu$. This implies that $(C,Z\!\setminus\! C)$ is a coloring of $G$ into clopen sets, which is absurd.\medskip

 We then set, for $l\!\in\!\omega$,\medskip
 
\leftline{$S_l\! :=\!\{ k\!\in\! I'\mid\forall p\!\in\!\omega ~~\exists j\!\in\! J^a_k~~Z^{k,a,j}\cap N_{\nu\vert l},
Z^{k,0,j,\overline{a}}\cap N_{n^{k,0}\vert p}\!\not=\!\emptyset ~\wedge$}\smallskip

\rightline{$\forall p\!\in\!\omega ~~\exists j\!\in\! J^{\overline{a}}_k~~
Z^{k,2k+1,j,a}\cap N_{n^{k,2k+1}\vert p},Z^{k,\overline{a},j}\cap N_{\nu\vert l}\!\not=\!\emptyset\}\mbox{,}$}\medskip
 
\noindent so that $(S_l)_{l\in\omega}$ is decreasing.

\vfill\eject

 We inductively define $k_0\! :=\!\mbox{min}~S_0$, and 
$k_{l+1}\! :=\!\mbox{min}~S_{l+1}\cap (k_l,\infty)$, so that $(k_l)_{l\in\omega}$ is stricly increasing. We define 
$\delta'\!\in\! 2^\omega$ by $\delta'(k)\! =\! 1\Leftrightarrow\exists l\!\in\!\omega ~~k\! =\! k_l$, so that 
$\delta'\!\in\!\mathbb{P}_\infty$. We pick, for $l,p\!\in\!\omega$ and $i\!\leq\! 2k_l$,\smallskip

\noindent - $j_{l,p,0}\!\in\! J^a_{k_l}$, $z_{k_l,a,j_{l,p,0}}\!\in\! N_{\nu\vert l}$ and 
$z_{k_l,0,j_{l,p,0},\overline{a}}\!\in\! N_{n^{k_l,0}\vert p}$, ensuring the injectivity of $(j_{l,p,0})_{p\in\omega}$,\smallskip 

\noindent - $j_{l,p,i+1}\!\in\! J^a_{k_l,i}$, $z_{k_l,i,j_{l,p,i+1},a}\!\in\! N_{n^{k_l,i}\vert p}$, 
$z_{k_l,i+1,j_{l,p,i+1},\overline{a}}\!\in\! N_{n^{k_l,i+1}\vert p}$, ensuring the injectivity of $(j_{l,p,i+1})_{p\in\omega}$,\smallskip 

\noindent - $j_{l,p,2k_l+2}\!\in\! J^{\overline{a}}_{k_l}$, $z_{k_l,2k_l+1,j_{l,p,2k_l+2},a}\!\in\! N_{n^{k_l,2k_l+1}\vert p}$ and 
$z_{k_l,\overline{a},j_{l,p,2k_l+2}}\!\in\! N_{\nu\vert l}$, ensuring the injectivity of $(j_{l,p,2k_l+2})_{p\in\omega}$.\medskip

 We are now ready to construct $\psi\! :\!\mathbb{P}_{\delta'}\!\rightarrow\! Z$. Note that\medskip
 
\leftline{$\mathbb{P}_{\delta'}\! =\!
\big\{ c^{k_l+1}\eta_\varepsilon j\varepsilon^\infty\mid j,l\!\in\!\omega\wedge\varepsilon\!\in\!\{ a,\overline{a}\}\big\}\cup
\big\{ k_li0^{j+1}\varepsilon^\infty\mid j,l\!\in\!\omega\wedge i\!\leq\! 2k_l\! +\! 1\wedge\varepsilon\!\in\!\{ a,\overline{a}\}\big\} ~\cup$}\smallskip

\rightline{$\{ c^\infty\}\cup\{ k_li0^\infty\mid l\!\in\!\omega\wedge i\!\leq\! 2k_l\! +\! 1\} .$}\medskip

\noindent We first set $\psi (c^\infty )\! :=\!\nu$. If $l\!\in\!\omega$, then $\psi$ sends $c^{k_l+1}0pa^\infty$ to 
$z_{k_l,a,j_{l,p,0}}$, $k_li0^\infty$ to $n^{k_l,i}$ if $i\!\leq\! 2k_l\! +\! 1$, $k_li0^{p+1}a^\infty$ to $z_{k_l,i,j_{l,p,i+1},a}$ and 
$k_li0^{p+1}\overline{a}^\infty$ to $z_{k_l,i,j_{l,p,i},\overline{a}}$ if $i\!\leq\! 2k_l\! +\! 1$, and $c^{k_l+1}1p\overline{a}^\infty$ to $z_{k_l,\overline{a},j_{l,p,2k_l+2}}$. Note that $\psi$ is as desired.\hfill{$\square$}

\begin{lem} \label{less} Let $\delta ,\delta'\!\in\!\mathbb{P}_\infty$ with 
$(\mathbb{P}_\delta ,\mathbb{G}_\delta )\preceq^i_c(\mathbb{P}_{\delta'},\mathbb{G}_{\delta'})$. Then there is 
$k_0\!\in\!\omega$ with the property that 
$\{ k\!\geq\! k_0\mid\delta (k)\! =\! 1\}\!\subseteq\!\{ k\!\in\!\omega\mid\delta'(k)\! =\! 1\}$.\end{lem}

\noindent\emph{Proof.}\ Let $\varphi\! :\!\mathbb{P}_\delta\!\rightarrow\!\mathbb{P}_{\delta'}$ be a witness for the fact that 
$(\mathbb{P}_\delta ,\mathbb{G}_\delta )\preceq^i_c(\mathbb{P}_{\delta'},\mathbb{G}_{\delta'})$. We set 
$V\! :=\!\varphi [\mathbb{P}_\delta ]$ and 
$E\! :=\! (\varphi\!\times\!\varphi )[\mathbb{G}_\delta ]$, so that $\chi_c(V,E)\! =\! 3$. By Lemma \ref{charsubbeta}, $c^\infty\!\in\! V$. Moreover, $\varphi\big[\mbox{proj}[\mathbb{G}_\delta ]\big]\!\subseteq\!\mbox{proj}[\mathbb{G}_{\delta'}]$.\medskip

 Let us prove that $\varphi [\mathbb{P}_\delta\!\setminus\! (\mbox{proj}[\mathbb{G}_\delta ]\cup\{ c^\infty\} )]\!\subseteq\!
\mathbb{P}_{\delta'}\!\setminus\! (\mbox{proj}[\mathbb{G}_{\delta'}]\cup\{ c^\infty\} )$. Towards a contradiction, suppose that we can find $k_1$ such that $\delta (k_1)\! =\! 1$ and $i\!\leq\! 2k_1\! +\! 1$ with 
$\varphi (k_1i0^\infty )\!\in\!\mbox{proj}[\mathbb{G}_{\delta'}]\cup\{ c^\infty\}$. We set $(\eta_a,\eta_{\overline{a}})\! :=\! (0,1)$. If 
${\varphi (k_1i0^\infty )\!\in\!\mbox{proj}[\mathbb{G}_{\delta'}]}$, then there are $i',j',k',\varepsilon$ with the property that either 
$\varphi (k_1i0^\infty )\! =\! c^{k'+1}\eta_\varepsilon j'\varepsilon^\infty$, or 
$\varphi (k_1i0^\infty )\! =\! k'i'0^{j'+1}\varepsilon^\infty$. The continuity of $\varphi$ provides a natural number $j_0$ with 
$\varphi [N_{k_1i0^{j_0+1}}]\!\subseteq\! N_{c^{k'+1}\eta_\varepsilon j'\varepsilon}$ or 
$\varphi [N_{k_1i0^{j_0+1}}]\!\subseteq\! N_{k'i'0^{j'+1}\varepsilon}$. This implies that the sequence 
$\big(\varphi (k_1i0^{j+1}a^\infty )\big)_{j\geq j_0}$ is constant, which contradicts the injectivity of $\varphi$. If 
$\varphi (k_1i0^\infty )\! =\! c^\infty$, then $\delta''\!\in\! 2^\omega$ defined by 
$\delta''(k)\! =\! 1\Leftrightarrow\delta (k)\! =\! 1\wedge k\!\not=\! k_1$ is in $\mathbb{P}_\infty$ and 
$\varphi_{\vert\mathbb{P}_{\delta''}}$ is a witness for the fact that 
$(\mathbb{P}_{\delta''},\mathbb{G}_{\delta''})\preceq^i_c(\mathbb{P}_{\delta'},\mathbb{G}_{\delta'})$. But the injectivity of 
$\varphi$ implies that $c^\infty\!\notin\!\varphi [\mathbb{P}_{\delta''}]$, which implies that 
$\chi_c(\mathbb{P}_{\delta''},\mathbb{G}_{\delta''})\!\leq\! 2$ by Lemma \ref{charsubbeta}, and contradicts Lemma \ref{chromalph}.\medskip

 This implies that $\varphi (c^\infty )\! =\! c^\infty$, and gives $k_0$ with $\varphi [N_{c^{k_0+1}}]\!\subseteq\! N_c$. Pick $k\!\geq\! k_0$ with $\delta (k)\! =\! 1$. This gives, for each $p\!\in\!\omega$, $j_p,k_p,\varepsilon_p$ with $\delta'(k_p)\! =\! 1$ and 
$\varphi (c^{k+1}0pa^\infty )\! =\! c^{k_p+1}\eta_{\varepsilon_p}j_p\varepsilon_p^\infty$. Extracting a subsequence if necessary, we may assume that the $\varepsilon_p$'s are equal to $\varepsilon$. Thus  
$$\varphi (k0^{p+2}\overline{a}^\infty )\! =\!\left\{\!\!\!\!\!\!\!
\begin{array}{ll}
& k_p0^{j_p+2}\overline{a}^\infty\mbox{ if }\varepsilon\! =\! a\mbox{,}\cr 
& k_p(2k_p\! +\! 1)0^{j_p+1}a^\infty\mbox{ if }\varepsilon\! =\!\overline{a}.
\end{array}
\right.$$
The continuity of $\varphi$ implies that $\big(\varphi (k0^{p+2}\overline{a}^\infty )\big)_{p\in\omega}$ converges to $\varphi (k0^\infty )$. Extracting a subsequence if necessary, we may assume that the $k_p$'s are equal to $k'$ and $(j_p)_{p\in\omega}$ is injective. Thus 
$\varphi (k0^\infty )$ is $k'0^\infty$ or $k'(2k'\! +\! 1)0^\infty$. Now note that the sequence 
$\big(\varphi (k0^{p+2}a^\infty )\big)_{p\in\omega}$ also converges to $\varphi (k0^\infty )$. This gives, for each $p\!\in\!\omega$, 
$j'_p,\varepsilon'_p$ with $(j'_p)_{p\in\omega}$ is injective (up to an extraction) and 
$$\varphi (k0^{p+2}a^\infty )\! =\!\left\{\!\!\!\!\!\!\!
\begin{array}{ll}
& k'0^{j'_p+2}{\varepsilon'_p}^\infty\mbox{ if }\varepsilon\! =\! a\mbox{,}\cr 
& k'(2k'\! +\! 1)0^{j'_p+1}{\varepsilon'_p}^\infty\mbox{ if }\varepsilon\! =\!\overline{a}.
\end{array}
\right.$$

 Extracting a subsequence if necessary, we may assume that the $\varepsilon'_p$'s are equal to $\varepsilon'$. Thus  
$$\varphi (k10^{p+1}\overline{a}^\infty )\! =\!\left\{\!\!\!\!\!\!\!
\begin{array}{ll}
& \left\{\!\!\!\!\!\!\!
\begin{array}{ll}
& k'10^{j'_p+1}\overline{a}^\infty\mbox{ if }\varepsilon'\! =\! a\mbox{,}\cr 
& c^{k'+1}0(j'_p\! +\! 1)a^\infty\mbox{ if }\varepsilon'\! =\!\overline{a}\mbox{,}
\end{array}
\right.
\mbox{ if }\varepsilon\! =\! a\mbox{,}\cr\cr
& \left\{\!\!\!\!\!\!\!
\begin{array}{ll}
& c^{k'+1}1(j'_p\! +\! 1)\overline{a}^\infty\mbox{ if }\varepsilon'\! =\! a\mbox{,}\cr 
& k'(2k')0^{j'_p+1}a^\infty\mbox{ if }\varepsilon'\! =\!\overline{a}\mbox{,}
\end{array}
\right.
\mbox{ if }\varepsilon\! =\!\overline{a}.
\end{array}
\right.$$
As $\big(\varphi (k10^{p+1}\overline{a}^\infty )\big)_{p\in\omega}$ converges to 
$\varphi (k10^\infty )\!\in\!\mathbb{P}_{\delta'}\!\setminus\! (\mbox{proj}[\mathbb{G}_{\delta'}]\cup\{ c^\infty\} )$, 
$\varphi (k10^{p+1}\overline{a}^\infty )\!\notin\! N_c$, and the second and third cases are not possible if $p$ is large enough. Thus $\varepsilon'\! =\!\varepsilon$, and $\varphi (k10^\infty )$ is $k'10^\infty$ or $k'(2k')0^\infty$. Now note that the sequence 
$\big(\varphi (k10^{p+1}a^\infty )\big)_{p\in\omega}$ also converges to $\varphi (k10^\infty )$. This gives, for each 
$p\!\in\!\omega$, $j''_p,\varepsilon''_p$ with $(j''_p)_{p\in\omega}$ is injective (up to an extraction) and 
$$\varphi (k10^{p+1}a^\infty )\! =\!\left\{\!\!\!\!\!\!\!
\begin{array}{ll}
& k'10^{j''_p+1}{\varepsilon''_p}^\infty\mbox{ if }\varepsilon\! =\! a\mbox{,}\cr 
& k'(2k')0^{j''_p+1}{\varepsilon''_p}^\infty\mbox{ if }\varepsilon\! =\!\overline{a}.
\end{array}
\right.$$
Extracting a subsequence if necessary, we may assume that the $\varepsilon''_p$'s are equal to $\varepsilon''$. If $k\! >\! 0$, then 
the continuity of $\varphi$ implies that $\big(\varphi (k20^{p+1}\overline{a}^\infty )\big)_{p\in\omega}$ converges to 
$\varphi (k20^\infty )\!\in\!\mathbb{P}_{\delta'}\!\setminus\! (\mbox{proj}[\mathbb{G}_{\delta'}]\cup\{ c^\infty\} )$. This implies that 
$\varphi (k20^\infty )\!\notin\! N_c$, and $\varphi (k20^{p+1}\overline{a}^\infty )\!\notin\! N_c$ if $p$ is large enough. So we may assume that
$$\varphi (k20^{p+1}\overline{a}^\infty )\! =\!\left\{\!\!\!\!\!\!\!
\begin{array}{ll}
& \left\{\!\!\!\!\!\!\!
\begin{array}{ll}
& k'20^{j''_p+1}\overline{a}^\infty\mbox{ if }\varepsilon''\! =\! a\mbox{,}\cr 
& k'0^{j''_p+2}a^\infty\mbox{ if }\varepsilon''\! =\!\overline{a}\mbox{,}
\end{array}
\right.
\mbox{ if }\varepsilon\! =\! a\mbox{,}\cr\cr
& \left\{\!\!\!\!\!\!\!
\begin{array}{ll}
& k'(2k'\! +\! 1)0^{j''_p+1}\overline{a}^\infty\mbox{ if }\varepsilon''\! =\! a\mbox{,}\cr 
& k'(2k'\! -\! 1)0^{j''_p+1}a^\infty\mbox{ if }\varepsilon''\! =\!\overline{a}\mbox{,}
\end{array}
\right.
\mbox{ if }\varepsilon\! =\!\overline{a}.
\end{array}
\right.$$
The injectivity of $\varphi$ and the value of $\varphi (k0^\infty )$ imply that second and third cases are not possible if $p$ is large enough. Thus $\varepsilon''\! =\!\varepsilon$, and $\varphi (k20^\infty )$ is $k'20^\infty$ or $k'(2k'\! -\! 1)0^\infty$. This implies that $k'\! >\! 0$. If now $k\! =\! 0$, then 
$$\varphi (c^{k+1}1p\overline{a}^\infty )\! =\!\left\{\!\!\!\!\!\!\!
\begin{array}{ll}
& c^{k'+1}1j''_p\overline{a}^\infty\mbox{ if }\varepsilon\! =\! a\mbox{,}\cr 
& c^{k'+1}0(j''_p\! +\! 1)a^\infty\mbox{ if }\varepsilon\! =\!\overline{a}.
\end{array}
\right.$$
since $\varphi (c^{k+1}1p\overline{a}^\infty )\!\in\! N_c$. This implies that $k'\! =\! 0$. Iterating this argument, we see that $k'\! =\! k$ and $\delta'(k)\! =\! 1$. Thus $\delta'(k)\! =\! 1$ if $\delta (k)\! =\! 1$ and $k\!\geq\! k_0$.\hfill{$\square$}\medskip
  
\noindent\emph{Proof of Theorem \ref{absmin}.}\ Lemmas \ref{copy} and \ref{below0} provide $\delta\!\in\!\mathbb{P}_\infty$ such that $(\mathbb{P}_\delta ,\mathbb{G}_\delta )\preceq^i_c\big(Z,G\big)$. We enumerate injectively 
$\{ n\!\in\!\omega\mid\delta (n)\! =\! 1\}\! =:\!\{ n_p\mid p\!\in\!\omega\}$. Let $(p_n)_{n\in\omega}$ be the sequence of prime numbers. We define, for each $\alpha\!\in\! 2^\omega$, 
$S_\alpha\!\subseteq\!\omega$ by $S_\alpha\! :=\!\{ p_0^{\alpha (0)+1}\ldots p_n^{\alpha (n)+1}\mid n\!\in\!\omega\}$. Note that 
$S_\alpha$ is infinite, and $S_\alpha\cap S_\beta$ is finite if $\alpha\!\not=\!\beta$. We define $\delta_\alpha\!\in\! 2^\omega$ by 
$\delta_\alpha (n)\! =\! 1\Leftrightarrow\exists p\!\in\! S_\alpha ~~n\! =\! n_p$. Note that $\delta_\alpha\!\in\!\mathbb{P}_\infty$ and 
$\delta_\alpha (n)\!\leq\!\delta (n)$ for each $n$. We set 
$(P_\alpha ,G_\alpha )\! :=\! (\mathbb{P}_{\delta_\alpha},\mathbb{G}_{\delta_\alpha})$, so that 
$P_\alpha\!\subseteq\!\mathbb{P}_\delta$, $G_\alpha\!\subseteq\!\mathbb{G}_\delta$, 
$(P_\alpha ,G_\alpha )\preceq^i_c(Z,G)$, and $\chi_c(P_\alpha ,G_\alpha )\! =\! 3$ by Lemma \ref{chromalph}.\medskip

 Let us prove that $(P_\alpha ,G_\alpha ),(P_{\alpha'},G_{\alpha'})$ are $\preceq^i_c$-incompatible among graphs on a 0DMS space with CCN at least three if $\alpha\!\not=\!\alpha'$. We argue by contradiction, which provides a 0DMS space $Z$ and a graph $G$ on $Z$ with CCN at least three and $(Z,G)\preceq^i_c(P_\alpha ,G_\alpha ),(P_{\alpha'},G_{\alpha'})$. Lemma \ref{below0} gives $\delta'\!\in\!\mathbb{P}_\infty$ such that 
$\{ k\!\in\!\omega\mid\delta'(k)\! =\! 1\}\!\subseteq\!\{ k\!\in\!\omega\mid\delta_\alpha (k)\! =\! 1\}$ and 
$(\mathbb{P}_{\delta'},\mathbb{G}_{\delta'})\preceq^i_c(Z,G)$. Lemma \ref{less} gives $k_0\!\in\!\omega$ such that 
$\{ k\!\geq\! k_0\mid\delta'(k)\! =\! 1\}\!\subseteq\!\{ k\!\in\!\omega\mid\delta_{\alpha'}(k)\! =\! 1\}$. This implies that 
$\{ k\!\in\!\omega\mid\delta'(k)\! =\! 1\}$ is finite, contradicting $\delta'\!\in\!\mathbb{P}_\infty$.

\vfill\eject

 Lemma \ref{generalab} then implies that there is no $\preceq^i_c$-antichain basis in the class of graphs on a 0DMS (or 0DP) space with CCN at least three. If $B$ is a basis for the class of graphs on a 0DP space with CCN at least three and  
$\alpha\!\in\! 2^\omega$, then there is $b_\alpha\!\in\! B$ with $b_\alpha\preceq^i_c(P_\alpha ,G_\alpha )$. The previous point shows that the sequence $(b_\alpha )_{\alpha\in 2^\omega}$ is injective, so that $B$ has size at least 
continuum.\hfill{$\square$}

\section{$\!\!\!\!\!\!$ General graphs and dynamical systems} \label{gdn} \indent

 We now prepare the proof of Theorem \ref{redbor}. We establish preliminary results holding not only for our examples, and clarify the relation between Cantor dynamical systems and our graphs. We first introduce examples in the style of the 
$\mathbb{G}_\gamma$'s defined at the beginning of Section \ref{compact}.\medskip
 
\noindent\emph{Notation.}\ Fix ${\bf d}\!\in\!\mathfrak{C}$ (defined before Theorem \ref{eq2''''''}). We associate to $\bf d$ the following objects:\medskip

\noindent - an increasing unbounded sequence $(n_l)_{l\in\omega}$ of natural numbers, sometimes denoted by 
$(n^{\bf d}_l)_{l\in\omega}$,\smallskip

\noindent - sequences $(L_l)_{l\in\omega},(R_l)_{l\in\omega}$ of integers with $R_l\! -\! L_l\! =\! 2n_l\! +\! 1$.\medskip

 Let $f\! :\!\mathcal{C}\!\rightarrow\!\mathcal{C}$ be a homeomorphism, so that $(\mathcal{C},f)$ is a Cantor dynamical system. We will associate a graph to $(\mathcal{C},f)$, as follows. Recall the definition of $\mathcal{J}$ at the beginning of Subsection \ref{sb}. We define, for $l\!\in\!\omega$ and $L_l\!\leq\! i\!\leq\! R_l$, $f_{l,i}\! :=\! f^{\bf d}_{l,i}\! :=\! f^i(0^\infty )\vert (l\! +\! 1)$ in $\prod_{l+1}$. This defines $\beta\!\in\!\mathcal{J}$ by setting $\lambda_l\! :=\! 2n_l\! +\! 2$, and $s_l(i)\! :=\! f_{l,L_l+i}$. We set 
$\mathcal{C}^+\! :=\!\mathbb{K}_\beta$, $\mathbb{O}_f\! :=\!\mathbb{O}_\beta$, and $\mathbb{G}_f\! :=\! s(\mathbb{O}_f)$, so that $\mathbb{G}_f\! =\!\mathbb{G}_\beta$. 

\begin{lem} \label{chrom} $(\mathcal{C}^+,\mathbb{G}_f)$ has CCN at least three and 
$\boraone\oplus\bormone$ chromatic number two. If moreover $d_0\! =\! 2$ and $f(x)(0)\!\not=\! x(0)$ for each 
$x\!\in\!\mathcal{C}$, then $(\mathcal{C}^+,\mathbb{G}_f)$ has CCN three.\end{lem}

\noindent\emph{Proof.}\ We first apply Lemma \ref{chromgen}. For the end, we argue as in the proof of Lemma \ref{chromgen}.\hfill{$\square$}\medskip
    
 We now want to compare the subgraphs of the $\mathbb{G}_f$'s. 

\begin{lem} \label{Z} Let $X$ be a topological space, $f\! :\! X\!\rightarrow\! X$ be a homeomorphism, $Y,g$ having the corresponding properties, $x\!\in\! X$, and $\varphi\! :\!\mbox{Orb}_f(x)\!\rightarrow\! Y$ such that 
$\varphi\big( f(z)\big)\! =\! g\big(\varphi (z)\big)$ for each $z\!\in\!\mbox{Orb}_f(x)$. Then 
${\varphi\big( f^i(x)\big)\! =\! g^i\big(\varphi (x)\big)}$ for each $i\!\in\!\mathbb{Z}$. Similarly, if 
${\varphi\big( f(z)\big)\! =\! g^{-1}\big(\varphi (z)\big)}$ for each $z\!\in\!\mbox{Orb}_f(x)$, then 
${\varphi\big( f^i(x)\big)\! =\! g^{-i}\big(\varphi (x)\big)}$ for each $i\!\in\!\mathbb{Z}$. In particular, 
$\varphi [\mbox{Orb}_f(x)]\! =\!\mbox{Orb}_g\big(\varphi (x)\big)$ in both cases.\end{lem}

\noindent\emph{Proof.}\ Inductively, we see that $\varphi\big( f^i(z)\big)\! =\! g^i\big(\varphi (z)\big)$ for each $i\!\in\!\omega$ and each $z\!\in\!\mbox{Orb}_f(x)$. In particular, 
$\varphi (x)\! =\!\varphi\Big( f^i\big( f^{-i}(x)\big)\Big)\! =\! g^i\Big(\varphi\big ( f^{-i}(x)\big)\Big)$, so that 
$\varphi\big( f^{-i}(x)\big)\! =\! g^{-i}\big(\varphi (x)\big)$. This implies that $\varphi\big( f^i(x)\big)\! =\! g^i\big(\varphi (x)\big)$ if 
$i\!\in\!\mathbb{Z}$. The other case is similar.\hfill{$\square$}

\begin{defi} We say that the tuple $({\bf d},{\bf d}',f_{\bf d},f_{{\bf d}'},V_{\bf d},V_{{\bf d}'},E_{\bf d},E_{{\bf d}'},\varphi )$ is a 
\emph{continuous tuple} if ${\bf d}\!\in\!\mathfrak{C}$, ${f_{\bf d}\! :\!\mathcal{C}_{\bf d}\!\rightarrow\!\mathcal{C}_{\bf d}}$ is a homeomorphism, $V_{\bf d}\!\subseteq\!\mathcal{C}^+_{\bf d}$, $E_{\bf d}\!\subseteq\!\mathbb{G}_{f_{\bf d}}\cap V_{\bf d}^2$ is a graph, $\chi_c(V_{\bf d},E_{\bf d})\!\geq\! 3$, ${\bf d}',f_{{\bf d}'},V_{{\bf d}'}$, $E_{{\bf d}'}$ have the corresponding properties, and $(V_{\bf d},E_{\bf d})\preceq_c(V_{{\bf d}'},E_{{\bf d}'})$ with witness $\varphi$.\end{defi}

 The next results are steps towards flip-conjugacy. 

\begin{lem} \label{or} Let $({\bf d},{\bf d}',f_{\bf d},f_{{\bf d}'},V_{\bf d},V_{{\bf d}'},E_{\bf d},E_{{\bf d}'},\varphi )$ be a continuous tuple, and 
$$\big( x,f_{\bf d}(x)\big)\!\in\!\overline{E_{\bf d}}^{(\mathcal{C}_{\bf d}^+)^2}\cap V_{\bf d}^2$$ 
with $\varphi (x),\varphi\big( f_{\bf d}(x)\big)\!\in\!\mathcal{C}_{{\bf d}'}$. Then 
$\varphi\big( f_{\bf d}(x)\big)\! =\! f_{{\bf d}'}\big(\varphi (x)\big)$ or 
$\varphi\big( f_{\bf d}(x)\big)\! =\! f_{{\bf d}'}^{-1}\big(\varphi (x)\big)$.\end{lem}

\noindent\emph{ Proof.}\ Note that $\big( x,f_{\bf d}(x)\big)\! =\!\mbox{lim}_{l\rightarrow\infty}~(x_l,u_l)$, where 
$(x_l,u_l)\!\in\! E_{\bf d}$. This successively implies that 
$\big(\varphi (x_l),\varphi (u_l)\big)\!\in\! E_{{\bf d}'}\!\subseteq\!\mathbb{G}_{f_{{\bf d}'}}$ and 
$\Big(\varphi (x),\varphi\big( f_{\bf d}(x)\big)\Big)\! =\!\! =\!\mbox{lim}_{l\rightarrow\infty}~\big(\varphi (x_l),\varphi (u_l)\big)$. As 
$\varphi (x),\varphi\big( f_{\bf d}(x)\big)$ are in $\mathcal{C}_{{\bf d}'}$, the first coordinate of $\varphi (x_l),\varphi (u_l)$ is in $d_0$ if $l$ is large enough. The definition of $\mathbb{G}_{f_{{\bf d}'}}$ provides $k_l\!\in\!\omega$, $i_l\!\in\!\mathbb{Z}$, and 
$(\rho_l,\varepsilon_l)\!\not=\! (\theta_l,\eta_l)$ in $\{ (0,a),(1,\overline{a})\}$ with the properties that $i_l\!\leq\! 2n_{k_l}$, 
$$\big(\varphi (x_l),\varphi (u_l)\big)\! =\! 
(f_{{\bf d}'}^{L_{k_l}+i_l+\rho_l}(0^\infty )\vert (k_l\! +\! 1)\varepsilon_l^{i_l+\rho_l+1}\overline{\varepsilon_l}^\infty ,
f_{{\bf d}'}^{L_{k_l}+i_l+\theta_l}(0^\infty )\vert (k_l\! +\! 1)\eta_l^{i_l+\theta_l+1}\overline{\eta_l}^\infty )$$ 
for such a $l$. Extracting a subsequence if necessary, we may assume that the sequence $(\rho_l)$ is constant. Moreover, the fact that $\varphi (x)\!\in\!\mathcal{C}_{{\bf d}'}$ implies that we may also assume that the sequence $(k_l)$ is strictly increasing. This implies that $\Big(\varphi (x),\varphi\big( f_{\bf d}(x)\big)\Big)$ is at distance zero from the closed set 
$\textup{Graph}(f_{{\bf d}'})^{1-2\rho_0}\cap V^2_{{\bf d}'}$. Thus 
$\varphi\big( f_{\bf d}(x)\big)\! =\! f_{{\bf d}'}\big(\varphi (x)\big)$ or 
$\varphi\big( f_{\bf d}(x)\big)\! =\! f_{{\bf d}'}^{-1}\big(\varphi (x)\big)$.\hfill{$\square$}

\begin{cor} \label{or+} Let $({\bf d},{\bf d}',f_{\bf d},f_{{\bf d}'},V_{\bf d},V_{{\bf d}'},E_{\bf d},E_{{\bf d}'},\varphi )$ be a continuous tuple such that $\varphi$ is injective, 
$\mathcal{C}_{\bf d}\!\subseteq\!\overline{V_{\bf d}\cap\mathcal{C}_{\bf d}}^{\mathcal{C}_{\bf d}}$ and 
$\textup{Graph}(f_{\bf d})\!\subseteq\!\overline{E_{\bf d}}^{(\mathcal{C}_{\bf d}^+)^2}$, and assume that 
$x,f_{\bf d}(x)\!\in\! V_{\bf d}\cap\mathcal{C}_{\bf d}$. Then either 
$\varphi\big( f_{\bf d}(x)\big)\! =\! f_{{\bf d}'}\big(\varphi (x)\big)$, or 
$\varphi\big( f_{\bf d}(x)\big)\! =\! f_{{\bf d}'}^{-1}\big(\varphi (x)\big)$.\end{cor}

\noindent\emph{Proof.}\ By Lemma \ref{incl}, $\varphi (x),\varphi\big( f_{\bf d}(x)\big)\!\in\!\mathcal{C}_{{\bf d}'}$. It remains to apply Lemma \ref{or}.\hfill{$\square$}

\begin{lem} \label{commonlemma} Let $X$ be a 0DMS space, ${f\! :\! X\!\rightarrow\! X}$ be a homeomorphism, 
$V\!\subseteq\! X$, $I$ be a subset of $\{ x\!\in\! V\mid f(x)\!\in\! V\}$, $Y,g,W$ having the corresponding properties, and 
$\varphi\! :\! V\!\rightarrow\! W$ be a continuous injection. We assume that 
$\varphi\big( f(x)\big)\! =\! g\big(\varphi (x)\big)\mbox{ or }\varphi\big( f(x)\big)\! =\! g^{-1}\big(\varphi (x)\big)\mbox{ if }
x\!\in\! I$.\smallskip 

\noindent (a)  Assume that $g^2$ is fixed point free, $x,f(x)\!\in\! I$, and $f^2(x)\!\not=\! x$. Then 
$\varphi\big( f(x)\big)\! =\! g\big(\varphi (x)\big)$ and $\varphi\big( f^2(x)\big)\! =\! g^2\big(\varphi (x)\big)$, or 
$\varphi\big( f(x)\big)\! =\! g^{-1}\big(\varphi (x)\big)$ and $\varphi\big( f^2(x)\big)\! =\! g^{-2}\big(\varphi (x)\big)$.\smallskip

\noindent (b) Assume that $g^2$ is fixed point free, $\mbox{Orb}_f(x)$ is a dense subset of $I$, and $f^2(x)\!\not=\! x$. Then either $\varphi\big( f(z)\big)\! =\! g\big(\varphi (z)\big)$ for each $z\!\in\! I$, or 
${\varphi\big( f(z)\big)\! =\! g^{-1}\big(\varphi (z)\big)}\mbox{ for each }z\!\in\! I$.\smallskip

\noindent (c) Assume that $g^2$ is fixed point free, $\mbox{Orb}_f(x)\!\subseteq\! I$, and $f^2(x)\!\not=\! x$. Then 
$\varphi [\mbox{Orb}_f(x)]\! =\!\mbox{Orb}_g\big(\varphi (x)\big)$.\end{lem}

\noindent\emph{Proof.}\ (a) We set $P\! :=\!\{ z\!\in\! I\mid\varphi\big( f(z)\big)\! =\! g\big(\varphi (z)\big)\}$ and 
$M\! :=\!\{ z\!\in\! I\mid\varphi\big( f(z)\big)\! =\! g^{-1}\big(\varphi (z)\big)\}$, so that $(P,M)$ is a covering of $I$ into closed sets. As $g^2$ is fixed point free, $P,M$ are disjoint. Thus $P$ is clopen in $I$. If $x\!\in\! P$ and $f(x)\!\in\! M$, then 
$\varphi\big( f^2(x)\big)\! =\! g^{-1}\Big(\varphi\big( f(x)\big)\Big)\! =\!\varphi (x)$, which contradicts the fact that 
$f^2(x)\!\not=\! x$ since $\varphi$ is injective. The argument is similar if we exchange $P$ and $M$.\medskip

\noindent (b) By (a), either $\mbox{Orb}_f(x)\!\subseteq\! P$, or $\mbox{Orb}_f(x)\!\subseteq\! M$. By density, 
$P\!\in\!\{ I,\emptyset\}$.\medskip

\noindent (c) If $\mbox{Orb}_f(x)\!\subseteq\! P$, then, by Lemma \ref{Z}, $\varphi\big( f^i(x)\big)\! =\! g^i\big(\varphi (x)\big)$ if 
$i\!\in\!\mathbb{Z}$ and we are done. Otherwise, by (a), we are in the similar case $\mbox{Orb}_f(x)\!\subseteq\! M$, so that 
$\varphi [\mbox{Orb}_f(x)]\! =\!\mbox{Orb}_g\big(\varphi (x)\big)$ in both cases.\hfill{$\square$}

\begin{lem} \label{orbeq} Let $({\bf d},{\bf d}',f_{\bf d},f_{{\bf d}'},V_{\bf d},V_{{\bf d}'},E_{\bf d},E_{{\bf d}'},\varphi )$ be a continuous tuple such that $f_{\bf d}^2,f_{{\bf d}'}^2$ are fixed point free, $\varphi$ is injective, 
$\mathcal{C}_{\bf d}\!\subseteq\!\overline{V_{\bf d}\cap\mathcal{C}_{\bf d}}^{\mathcal{C}_{\bf d}}$ and 
$\textup{Graph}(f_{\bf d})\!\subseteq\!\overline{E_{\bf d}}^{(\mathcal{C}_{\bf d}^+)^2}$.\smallskip

\noindent (a) Assume that $x,f_{\bf d}(x),f^2_{\bf d}(x)\!\in\! V_{\bf d}\cap\mathcal{C}_{\bf d}$. Then 
$\varphi\big( f_{\bf d}(x)\big)\! =\! f_{{\bf d}'}\big(\varphi (x)\big)\mbox{ and }
\varphi\big( f_{\bf d}^2(x)\big)\! =\! f_{{\bf d}'}^2\big(\varphi (x)\big)$, or 
$\varphi\big( f_{\bf d}(x)\big)\! =\! f_{{\bf d}'}^{-1}\big(\varphi (x)\big)$ and 
$\varphi\big( f_{\bf d}^2(x)\big)\! =\! f_{{\bf d}'}^{-2}\big(\varphi (x)\big)$.\smallskip

\noindent (b) If $\mbox{Orb}_{f_{\bf d}}(x)$ is a dense subset of 
$I\! :=\!\{ z\!\in\! V_{\bf d}\cap\mathcal{C}_{\bf d}\mid f_{\bf d}(z)\!\in\! V_{\bf d}\}$ for some $x$, then either 
$\varphi\big( f_{\bf d}(z)\big)\! =\! f_{{\bf d}'}\big(\varphi (z)\big)$ for each $z\!\in\! I$, or 
${\varphi\big( f_{\bf d}(z)\big)\! =\! f_{{\bf d}'}^{-1}\big(\varphi (z)\big)}\mbox{ for each }z\!\in\! I$.\medskip

\noindent (c) Assume that $\mbox{Orb}_{f_{\bf d}}(x)\!\subseteq\! V_{\bf d}\cap\mathcal{C}_{\bf d}$. Then 
$\varphi [\mbox{Orb}_{f_{\bf d}}(x)]\! =\!\mbox{Orb}_{f_{{\bf d}'}}\big(\varphi (x)\big)$.\end{lem}

\noindent\emph{Proof.}\ We apply Lemma \ref{commonlemma} to $X\! :=\! \mathcal{C}_{\bf d}$, $f\! :=\! f_{\bf d}$, 
$V\! :=\! V_{\bf d}\cap\mathcal{C}_{\bf d}$, $I$ defined in (b), $Y\! :=\! \mathcal{C}_{{\bf d}'}$, $g\! :=\! f_{{\bf d}'}$, 
$W\! :=\! V_{{\bf d}'}\cap\mathcal{C}_{{\bf d}'}$, and $\varphi\! :=\!\varphi_{\vert V}$, which is possible by Lemma \ref{incl} and Corollary \ref{or+}.\hfill{$\square$}

\begin{lem} \label{flip} Let ${\bf d},{\bf d}'\!\in\!\mathfrak{C}$, and 
$f_{\bf d}\! :\!\mathcal{C}_{\bf d}\!\rightarrow\!\mathcal{C}_{\bf d}$, 
$f_{{\bf d}'}\! :\!\mathcal{C}_{{\bf d}'}\!\rightarrow\!\mathcal{C}_{{\bf d}'}$ be minimal homeomorphisms such that 
$\textup{Graph}(f_{\bf d})\!\subseteq\!\overline{\mathbb{G}_{f_{\bf d}}}^{(\mathcal{C}_{\bf d}^+)^2}$. If 
$(\mathcal{C}^+_{\bf d},\mathbb{G}_{f_{\bf d}})\preceq_c^i(\mathcal{C}^+_{{\bf d}'},\mathbb{G}_{f_{{\bf d}'}})$, then 
$f_{\bf d},f_{{\bf d}'}$ are flip-conjugate.\end{lem}

\noindent\emph{Proof.}\ Let $\varphi\! :\!\mathcal{C}^+_{\bf d}\!\rightarrow\!\mathcal{C}^+_{{\bf d}'}$ be a witness for the fact that 
$(\mathcal{C}^+_{\bf d},\mathbb{G}_{f_{\bf d}})\preceq_c^i(\mathcal{C}^+_{{\bf d}'},\mathbb{G}_{f_{{\bf d}'}})$. By Lemma \ref{chrom}, $({\bf d},{\bf d}',f_{\bf d},f_{{\bf d}'},\mathcal{C}^+_{\bf d},\mathcal{C}^+_{{\bf d}'},\mathbb{G}_{f_{\bf d}},
\mathbb{G}_{f_{{\bf d}'}},\varphi )$ is a continuous tuple satisfying the assumptions of Lemmas \ref{incl} and \ref{orbeq}. In particular, $\varphi [\mathcal{C}_{\bf d}]\!\subseteq\!\mathcal{C}_{{\bf d}'}$ and the map 
$\psi\! :=\!\varphi_{\vert\mathcal{C}_{\bf d}}\! :\!\mathcal{C}_{\bf d}\!\rightarrow\!\mathcal{C}_{{\bf d}'}$ is a witness for the fact that $f_{\bf d},f_{{\bf d}'}$ are flip-conjugate. Indeed, Lemma \ref{orbeq} implies that 
$\varphi [\mbox{Orb}_{f_{\bf d}}(x)]\! =\!\mbox{Orb}_{f_{{\bf d}'}}\big(\varphi (x)\big)$ if $x\!\in\!\mathcal{C}_{\bf d}$. As 
$f_{{\bf d}'}$ is minimal, the compact set $\psi [\mathcal{C}_{\bf d}]$ is dense in $\mathcal{C}_{{\bf d}'}$, showing that $\psi$ is onto, and thus a homeomorphism by compactness of $\mathcal{C}_{\bf d}$.\hfill{$\square$}\medskip

\noindent\emph{Notation.}\ For the converse, we give a definition of the sequences $(n_l)_{l\in\omega}$, $(L_l)_{l\in\omega}$ and $(R_l)_{l\in\omega}$. We define $(n_l)_{l\in\omega}$ by $n_l\! :=\! l$, so that $(n_l)_{l\in\omega}$ is increasing unbounded.\medskip

\noindent - We define a map $\zeta\! :\!\omega\!\rightarrow\!\mathbb{Z}$ having the property that each integer appears infinitely many times in the range of $\zeta$.\medskip

\noindent - We define sequences $(L_l)_{l\in\omega},(R_l)_{l\in\omega}$ of integers by $L_{2m}\! :=\! R_{2m+1}\! :=\!\zeta (m)$ and $R_l\! -\! L_l\! :=\! 2l\! +\! 1$, so that the sequences 
$\big( f^{L_l}(0^\infty )\big)_{l\in\omega},\big( f^{R_l}(0^\infty )\big)_{l\in\omega}$ are dense if $f$ is minimal.\medskip

 Note also that $\{ f^i(0^\infty )\mid\exists^\infty l\!\in\!\omega ~~L_l\!\leq\! i\! <\! R_l\}$ is dense in $\mathcal{C}$ if $f$ is minimal, which implies that 
$\textup{Graph}(f_{\bf d})\!\subseteq\!\overline{\mathbb{G}_{f_{\bf d}}}^{(\mathcal{C}_{\bf d}^+)^2}$.

\begin{lem} \label{conv} Let ${\bf d},{\bf d}'\!\in\!\mathfrak{C}$, $f_{\bf d}\! :\!\mathcal{C}_{\bf d}\!\rightarrow\!\mathcal{C}_{\bf d}$, 
$f_{{\bf d}'}\! :\!\mathcal{C}_{{\bf d}'}\!\rightarrow\!\mathcal{C}_{{\bf d}'}$ be minimal homeomorphisms, and 
$(n_l)_{l\in\omega},(L_l)_{l\in\omega},(R_l)_{l\in\omega}$ just defined. If $f_{\bf d},f_{{\bf d}'}$ are flip-conjugate, then 
$(\mathcal{C}^+_{\bf d},\mathbb{G}_{f_{\bf d}})\preceq_c^i(\mathcal{C}^+_{{\bf d}'},\mathbb{G}_{f_{{\bf d}'}})$.\end{lem}

\noindent\emph{Proof.}\ As $f_{\bf d},f_{{\bf d}'}$ are flip-conjugate, we get a homeomorphism 
$\psi\! :\!\mathcal{C}_{\bf d}\!\rightarrow\!\mathcal{C}_{{\bf d}'}$. We have to define a function 
$\varphi\! :\!\mathcal{C}^+_{\bf d}\!\rightarrow\!\mathcal{C}^+_{{\bf d}'}$. We first set $\varphi (c^\infty )\! :=\! c^\infty$, and 
$\varphi (x)\! :=\!\psi (x)$ if $x\!\in\! \mathcal{C}_{\bf d}$. As $f_{{\bf d}'}$ is uniformly continuous, for any $l\!\in\!\omega$, there is $U\!\geq\! l$ such that, for any $y,z\!\in\!\mathcal{C}_{{\bf d}'}$, $f_{{\bf d}'}(y)\vert (l\! +\! 1)\! =\! f_{{\bf d}'}(z)\vert (l\! +\! 1)$ if 
$y\vert (U\! +\! 1)\! =\! z\vert (U\! +\! 1)$, which defines $U\! :\!\omega\!\rightarrow\!\omega$.\medskip

 Assume first that $\psi\!\circ\! f_{\bf d}\! =\! f_{{\bf d}'}\!\circ\!\psi$. By Lemma \ref{Z}, 
$\psi\big( f_{\bf d}^i(0^\infty )\big)\! =\! f_{{\bf d}'}^i\big(\psi (0^\infty )\big)$ for each $i\!\in\!\mathbb{Z}$. We define 
$\big(\varphi (f^{\bf d}_{l,L_l+i}\varepsilon^{i+1}\overline{\varepsilon}^\infty )\big)_{i\leq 2l+1}$ by induction on $l$, ensuring that 
$\varphi (f^{\bf d}_{l,L_l+i}\varepsilon^{i+1}\overline{\varepsilon}^\infty )$ is of the form 
$f^{{\bf d}'}_{l'(l,i,\varepsilon ),L_{l'(l,i,\varepsilon )}+i'(l,i,\varepsilon )}\varepsilon^{i'(l,i,\varepsilon )+1}
\overline{\varepsilon}^\infty$ with $l'(l,i\! +\! 1,\overline{a})\! =\! l'(l,i,a)$ and $i'(l,i\! +\! 1,\overline{a})\! =\! i'(l,i,a)$ if $i\!\leq\! 2l$. Fix $l\!\in\!\omega$. Let 
$M\! :=\!\mbox{sup}\big\{ l'(k,j,\eta )\mid k\! <\! l\wedge j\!\leq\! 2k\! +\! 1\wedge\eta\!\in\!\{ a,\overline{a}\}\big\}$.\medskip

\noindent - We choose $m\! :=\! l'(l,0,\overline{a})\! >\!\mbox{max}(l,M)$ such that   
${f_{{\bf d}'}^{L_m}(0^\infty )\vert (l\! +\! 1)\! =\! f_{{\bf d}'}^{L_l}\big(\psi (0^\infty )\big)\vert (l\! +\! 1)}$, which is possible since the sequence  
$\big( f_{{\bf d}'}^{L_m}(0^\infty )\big)_{m\in\omega}$ is dense. We set 
$\varphi (c^{l+1}a\overline{a}^\infty )\! :=\! c^{m+1}a\overline{a}^\infty$ and 
$\varphi (f^{\bf d}_{l,L_l}\overline{a}a^\infty )\! :=\! f^{{\bf d}'}_{m,L_m}\overline{a}a^\infty$. Note that 
$\big(\varphi (c^{l+1}a\overline{a}^\infty ),\varphi (f^{\bf d}_{l,L_l}\overline{a}a^\infty )\big)\!\in\!\mathbb{G}_{f_{{\bf d}'}}$, as desired.\medskip

\noindent - Now fix $i\!\leq\! 2n_l$. As $f_{{\bf d}'}$ is minimal, we can find $k_i\!\in\!\mathbb{Z}$ with 
$$f_{{\bf d}'}^{k_i}(0^\infty )\vert\big( U(l)\! +\! 1\big)\! =\! f_{{\bf d}'}^{L_l+i}\big(\psi (0^\infty )\big)\vert\big( U(l)\! +\! 1\big) .$$
Note that $f_{{\bf d}'}^{k_i+1}(0^\infty )\vert (l\! +\! 1)\! =\! f_{{\bf d}'}^{L_l+i+1}\big(\psi (0^\infty )\big)\vert (l\! +\! 1)$. We choose 
$m'\! :=\! l'(l,i,a)\! >\! l'(l,i,\overline{a})$ with the property that $L_{m'}\!\leq\! k_i\! =\! L_{m'}\! +\! i'_i\! <\! R_{m'}$. We set 
$\varphi (f^{\bf d}_{l,L_l+i}a^{i+1}\overline{a}^\infty )\! :=\! f^{{\bf d}'}_{m',k_i}a^{i'_i+1}\overline{a}^\infty$ and 
$\varphi (f^{\bf d}_{l,L_l+i+1}\overline{a}^{i+2}a^\infty )\! :=\! f^{{\bf d}'}_{m',k_i+1}\overline{a}^{i'_i+2}a^\infty$. Note that 
$\big(\varphi (f^{\bf d}_{l,L_l+i}a^{i+1}\overline{a}^\infty ),\varphi (f^{\bf d}_{l,L_l+i+1}\overline{a}^{i+2}a^\infty )\big)$ is in 
$\mathbb{G}_{f_{{\bf d}'}}$, as desired.

\vfill\eject

\noindent - We then choose $m''\! :=\! l'(l,2l\! +\! 1,a)\! >\! l'(l,2l\! +\! 1,\overline{a})$ with 
$f_{{\bf d}'}^{R_{m''}}(0^\infty )\vert (l\! +\! 1)\! =\! f_{{\bf d}'}^{R_l}\big(\psi (0^\infty )\big)\vert (l\! +\! 1)$, which is possible since 
$\big( f_{{\bf d}'}^{R_m}(0^\infty )\big)_{m\in\omega}$ is dense. We set 
$\varphi (f^{\bf d}_{l,R_l}a^{2l+2}\overline{a}^\infty )\! :=\! f^{{\bf d}'}_{m'',R_{m''}}a^{2m''+2}\overline{a}^\infty$ and 
$\varphi (c^{l+1}\overline{a}a^\infty )\! :=\! c^{m''+1}\overline{a}a^\infty$. Note that 
$\big(\varphi (f^{\bf d}_{l,R_l}a^{2l+2}\overline{a}^\infty ),\varphi (c^{l+1}\overline{a}a^\infty )\big)\!\in\!\mathbb{G}_{f_{{\bf d}'}}$, as desired.\medskip

 This completes the definition of $\varphi$. Our construction implies that $\varphi$ is a homomorphism from 
$(\mathcal{C}^+_{\bf d},\mathbb{G}_{f_{\bf d}})$ into $(\mathcal{C}^+_{{\bf d}'},\mathbb{G}_{f_{{\bf d}'}})$. Our choice of the 
$l'(l,i,\varepsilon )$'s implies the injectivity of $\varphi$. For the continuity, note first that the sequence 
$\big(\varphi (c^{l+1}\varepsilon\overline{\varepsilon}^\infty )\big)_{l\in\omega}$ converges to $c^\infty$ since $m,m''\!\geq\! l$. If now 
$(f^{\bf d}_{l_k,L_{l_k}+i_k}\varepsilon_k^{i_k+1}\overline{\varepsilon_k}^\infty )_{k\in\omega}$ converges to 
$x\!\in\!\mathcal{C}_{\bf d}^+$, then we may assume that $x\!\in\!\mathcal{C}_{\bf d}$. Note that the sequence  
$\big( f_{\bf d}^{L_{l_k}+i_k}(0^\infty )\big)_{k\in\omega}$ converges to $x$, and 
$\Big( f_{{\bf d}'}^{L_{l_k}+i_k}\big(\psi (0^\infty )\big)\Big)_{k\in\omega}$ converges to $\psi (x)\! =\!\varphi (x)$. Our construction ensures that $f_{{\bf d}'}^{L_{l'(l,i,\varepsilon )}+i'(l,i,\varepsilon )}(0^\infty )\vert (l\! +\! 1)\! =\! 
f_{{\bf d}'}^{L_l+i}\big(\psi (0^\infty )\big)\vert (l\! +\! 1)$, and 
$$\begin{array}{ll}
\varphi (f^{\bf d}_{l_k,L_{l_k}+i_k}\varepsilon_k^{i_k+1}\overline{\varepsilon_k}^\infty )\!\!\!\!
& \! =\! f^{{\bf d}'}_{l'(l_k,i_k,\varepsilon_k),L_{l'(l_k,i_k,\varepsilon_k)}+i'(l_k,i_k,\varepsilon_k)}\varepsilon_k^{i'(l_k,i_k,\varepsilon_k)+1}\overline{\varepsilon_k}^\infty\cr\cr
& \! =\! f_{{\bf d}'}^{L_{l'(l_k,i_k,\varepsilon_k)}+i'(l_k,i_k,\varepsilon_k)}(0^\infty )\vert\big( l'(l_k,i_k,\varepsilon_k)\! +\! 1\big)\varepsilon_k^{i'(l_k,i_k,\varepsilon_k)+1}\overline{\varepsilon_k}^\infty .
\end{array}$$ 

 As $l'(l,i,\varepsilon)\!\geq\! l$, we get 
$$\varphi (f^{\bf d}_{l_k,L_{l_k}+i_k}\varepsilon_k^{i_k+1}\overline{\varepsilon_k}^\infty )\vert (l_k\! +\! 1)
\! =\! f_{{\bf d}'}^{L_{l'(l_k,i_k,\varepsilon_k)}+i'(l_k,i_k,\varepsilon_k)}(0^\infty )\vert (l_k\! +\! 1)
\! =\! f_{{\bf d}'}^{L_{l_k}+i_k}\big(\psi (0^\infty )\big)\vert (l_k\! +\! 1).$$ 
Thus $\big(\varphi (f^{\bf d}_{l_k,L_{l_k}+i_k}\varepsilon_k^{i_k+1}\overline{\varepsilon_k}^\infty )\big)_{k\in\omega}$ converges to 
$\varphi (x)$, proving the continuity of $\varphi$.\medskip

 The case where $\psi\!\circ\! f_{\bf d}\! =\! f_{{\bf d}'}^{-1}\!\circ\!\psi$ is similar.\hfill{$\square$}\medskip
     
\noindent\emph{Notation.}\ We set $(X,R)\equiv^i_c(Y,S)\Leftrightarrow (X,R)\preceq^i_c(Y,S)\wedge (Y,S)\preceq^i_c(X,R)$, so that 
$\equiv^i_c$ is the equivalence relation associated with the quasi-order $\preceq^i_c$.
 
\begin{cor} \label{corflip} Let ${\bf d},{\bf d}'\!\in\!\mathfrak{C}$, 
$f_{\bf d}\! :\!\mathcal{C}_{\bf d}\!\rightarrow\!\mathcal{C}_{\bf d}$, 
$f_{{\bf d}'}\! :\!\mathcal{C}_{{\bf d}'}\!\rightarrow\!\mathcal{C}_{{\bf d}'}$ be minimal homeomorphisms, and 
$(n_l)_{l\in\omega},(L_l)_{l\in\omega},(R_l)_{l\in\omega}$ just defined. Then 
$(\mathcal{C}^+_{\bf d},\mathbb{G}_{f_{\bf d}})\equiv_c^i(\mathcal{C}^+_{{\bf d}'},\mathbb{G}_{f_{{\bf d}'}})$  if and only if 
$f_{\bf d},f_{{\bf d}'}$ are flip-conjugate.\end{cor}

\section{$\!\!\!\!\!\!$ General graphs and odometers}\indent

 We now provide a countable graph on a(n infinite) 0DMC space with CCN three which is strictly $\preceq_c$-below the odd cycles. Our example is based on odometers. We give some notation useful for the sequel.\medskip
 
\noindent\emph{Notation.}\ Fix ${\bf d}\!\in\!\mathfrak{C}$ (defined before Theorem \ref{eq2''''''}). The \emph{odometer} 
$o\! :=\! o_{\bf d}\! :\!\mathcal{C}\!\rightarrow\!\mathcal{C}$ is defined by
$$o(\alpha )\! :=\!\left\{\!\!\!\!\!\!\!
\begin{array}{ll}
& 0^\infty\mbox{ if }\forall j\!\in\!\omega ~~\alpha (j)\! =\! d_j\! -\! 1\mbox{,}\cr
& 0^n\big(\alpha (n)\! +\! 1\big)\alpha (n\! +\! 1)\alpha (n\! +\! 2)\ldots\mbox{ if }\alpha (n)\! <\! d_n\! -\! 1\wedge
\forall j\! <\! n~~\alpha (j)\! =\! d_j\! -\! 1.
\end{array}
\right.$$ 
As $\mbox{Orb}^+_o(\alpha )\! :=\!\{ o^i(\alpha )\mid i\! >\! 0\}$ sees all the words of length $n$ in the first $n$ coordinates for any $x$, $o$ is a minimal homeomorphism. We sometimes extend the definition of $o$ to finite sequences.

\vfill\eject

 We set 
$\mathfrak{D}\! :=\!\{ {\bf d}\!\in\!\mathfrak{C}\mid d_0\! =\! 2\wedge\forall j\!\geq\! 1~~d_j\mbox{ is odd}\}$.\medskip

\noindent - We define, for ${\bf d}\!\in\!\mathfrak{D}$, $(n_l)_{l\in\omega}\!\in\!\mathcal{S}$ by $n_0\! :=\! 0$ and 
$n_{l+1}\! :=\!\frac{(\pi_{1\leq j\leq l+1}~d_j)-1}{2}$. Note that
$$o^{2n_l+1}(0^\infty )\! =\! 1{^\frown}_{1\leq j\leq l}~\frac{d_j-1}{2}0^\infty$$
converges to $\mu\! :=\! 1{^\frown}_{j\geq 1}~\frac{d_j-1}{2}$ as $l$ goes to infinity, and 
$o^{4n_l+1}(0^\infty )\! =\! {^\frown}_{j\leq l}~(d_j\! -\! 1)0^\infty$.\medskip

\noindent - We define, for $l\!\in\!\omega$ and $i\!\leq\! 2n_l\! +\! 1$, $o_{l,i}\! :=\! o^{\bf d}_{l,i}\! :=\! o^i(0^\infty )\vert (l\! +\! 1)\!\in\!\prod_{l+1}$. This defines $\beta\!\in\!\mathcal{J}^c$, by setting $\lambda_l\! :=\! 2n_l\! +\! 2$ and $s_l(i)\! :=\! o_{l,i}$. We set $\mathcal{C}^+\! :=\!\mathbb{K}_\beta$, $\mathbb{O}_o\! :=\!\mathbb{O}_\beta$, and 
$\mathbb{G}_o\! :=\! s(\mathbb{O}_o)$, so that $\mathbb{G}_o\! =\!\mathbb{G}_\beta$. 

\begin{prop} \label{belowfin} Let ${\bf d}\!\in\!\mathfrak{D}$. Then $(\mathcal{C}^+,\mathbb{G}_o)$ has CCN three, 
$\boraone\oplus\bormone$ chromatic number two, and is strictly $\preceq_c$-below the examples of Corollary \ref{basisfin}.\end{prop}

\noindent\emph{Proof.}\ Lemma \ref{chromgen} proves the assertions about chromatic numbers. Let $p\!\in\!\omega$. We define a function $\varphi\! :\!\mathcal{C}^+\!\rightarrow\! 2p\! +\! 3$ as follows. Fix $l_0\!\in\!\omega$ minimal such that $n_{l_0}\!\geq\! p$. We set $\varphi (x)\! :=\! 0$ if $x(0)\! =\! c$. If $i\!\leq\! 4n_{l_0}\! +\! 1$ and $o^i(0^{l_0+1})\!\subseteq\! x$ (recall that $o$ can be extended to finite sequences), then we set
$$\varphi (x)\! :=\!\left\{\!\!\!\!\!\!\!
\begin{array}{ll}
& i\! +\! 1\mbox{ if }i\!\leq\! 2p\! +\! 1\mbox{,}\cr
& 2p\! +\! 2\mbox{ if }2p\! +\! 1\! <\! i\!\leq\! 2n_{l_0}\! +\! 1\wedge i\mbox{ is odd,}\cr
& 2p\! +\! 1\mbox{ if }2p\! +\! 1\! <\! i\!\leq\! 2n_{l_0}\! +\! 1\wedge i\mbox{ is even,}\cr
& 2p\! +\! 2n_{l_0}\! +\! 3\! -\! i\mbox{ if }2n_{l_0}\! +\! 1\! <\! i\!\leq\! 2p\! +\! 2n_{l_0}\! +\! 1\mbox{,}\cr
& 2\mbox{ if }2p\! +\! 2n_{l_0}\! +\! 1\! <\! i\!\leq\! 4n_{l_0}\! +\! 1\wedge i\mbox{ is odd,}\cr
& 1\mbox{ if }2p\! +\! 2n_{l_0}\! +\! 1\! <\! i\!\leq\! 4n_{l_0}\! +\! 1\wedge i\mbox{ is even.}
\end{array}
\right.$$ 
It remains to define $\varphi (o_{l,i}\varepsilon^{i+1}\overline{\varepsilon}^\infty )$ if $\varepsilon\!\in\!\{ a,\overline{a}\}$, 
$l\! <\! l_0$ and $i\!\leq\! 2n_l\! +\! 1$. We set $\varphi (o_{l,i}\varepsilon^{i+1}\overline{\varepsilon}^\infty )\! :=\! i\! +\! 1$. This defines a continuous homomorphism $\varphi$ from $(\mathcal{C}^+,\mathbb{G}_o)$ into $(2p\! +\! 3,C_{2p+3})$. The inequality 
$(\mathcal{C}^+,\mathbb{G}_o)\prec_c(2p\! +\! 3,C_{2p+3})$ is strict because of Corollary \ref{basisfin}(b).\hfill{$\square$}\medskip

\noindent\emph{Remark.}\ We clarify the limits of Theorem \ref{eqfin}. In its proof, we used the finiteness of 
$X$. This is essential. Indeed, if we replace $\chi$ with $\chi_c$, $\preceq^i$ with $\preceq^i_c$ and $X$ with 
$\mathcal{C}^+$, then the following hold. The implications $\mbox{(2) }\Rightarrow\mbox{ (1)}$, 
$\mbox{(3) }\Rightarrow\mbox{ (1)}$ and $\mbox{(3) }\Leftrightarrow\mbox{ (2)}$ still hold. The implications 
$\mbox{(1) }\Rightarrow\mbox{ (2)}$ and $\mbox{(1) }\Rightarrow\mbox{ (3)}$ do not hold, because of Proposition \ref{belowfin}.\medskip

 We now characterize the subgraphs of $\mathbb{G}_o$ having a big CCN.

\begin{lem} \label{charsub++} Let ${\bf d}\!\in\!\mathfrak{D}$, $V\!\subseteq\!\mathcal{C}^+$, and 
$E\!\subseteq\!\mathbb{G}_o\cap V^2$. The following are equivalent:\smallskip

\noindent (1) the digraph $(V,E)$ has CCN at least three,\smallskip

\noindent (2) the following hold:\smallskip

(a) $0^\infty ,\mu, c^\infty\!\in\! V$ and $\mathcal{C}\!\subseteq\!\overline{V\cap\mathcal{C}}^\mathcal{C}$,\smallskip 

(b) $\{ (c^\infty ,0^\infty ),(\mu, c^\infty )\}\cup\textup{Graph}(o)\!\subseteq\!\overline{s(E)}^{(\mathcal{C}^+)^2}$.\end{lem}
 
\noindent\emph{Proof.}\ Note first that 
$\chi_c(V,E)\!\leq\!\chi_c\big(V,s(E)\big)\!\leq\!\chi_c(\mathcal{C}^+,\mathbb{G}_o)\! =\! 3$, by Lemma \ref{chromgen}. We may and will assume that $E\! =\! s(E)$ is a graph.\medskip

\noindent (1) $\Rightarrow$ (2).(a) For $0^\infty$, we argue by contradiction. Let 
$C\! :=\! (N_1\cup\bigcup_{l\in\omega}~N_{0^{l+1}\overline{a}})\cap V$. Then $(C,V\!\setminus\! C)$ is a coloring of $E$ into clopen sets since $C\! =\! (N_1\cup\bigcup_{l\in\omega}~N_{0^{l+1}\overline{a}}\cup\{ 0^\infty\} )\cap V$, which is absurd. For $\mu$, we argue similarly, with $(N_0\cup\bigcup_{l\in\omega}~N_{\mu\vert(l+1)a})\cap V$. For $c^\infty$, we apply Lemma \ref{c}.\medskip

\noindent\emph{Claim.}\it\ Let $l\!\in\!\omega$. Then either for each $0\! <\! i\! <\! 2n_l\! +\! 1$ there is $x\!\in\! V\cap\mathcal{C}$ with $o_{l,i}\!\subseteq\! x$, or for each $2n_l\! +\! 1\! <\! j\! <\! 4n_l\! +\! 2$ there is $y\!\in\! V\cap\mathcal{C}$ with 
$o_{l,j}\!\subseteq\! y$.\rm

\vfill\eject

 Indeed, towards a contradiction, suppose that $l_0,i,j$ exist. Assume, for example, that $i,j$ are even, the other cases being similar. We extend the notation $o_{l,i}$ for $i\!\leq\! 4n_l\! +\! 1$ (not only $i\!\leq\! 2n_l\! +\! 1$), and set\medskip
 
\leftline{$C'\! :=\!\big(\bigcup_{1\leq l\leq l_0}~(N_{c^l\overline{a}}\cup\bigcup_{s\in\prod_l}~N_{s\overline{a}})\cup
\bigcup_{k<i,k\text{ even}}~N_{o_{l_0,k}}\cup\bigcup_{l\geq l_0,s\in\prod_{l_0<j\leq l}~d_j}~N_{o_{l_0,i}s\overline{a}}~\cup$}\smallskip

\rightline{$\bigcup_{i<k<j,k\text{ odd}}~N_{o_{l_0,k}}\cup\bigcup_{l\geq l_0,t\in\prod_{l_0<j\leq l}~d_j}~N_{o_{l_0,j}ta}\cup\bigcup_{j<k<4n_{l_0}+2,k\text{ even}}~N_{o_{l_0,k}}\big)\cap V$.}\medskip

\noindent Then $(C',V\!\setminus\! C')$ is a coloring of $E$ into clopen sets, which is absurd.\hfill{$\diamond$}\medskip

 Now let $l\!\in\!\omega$, $t\!\in\!\prod_l$, $i\! <\! 2n_l\! +\! 1$ with $t0\! =\! o_{l,i}$, and $2n_l\! +\! 1\!\leq\! j\! <\! 4n_l\! +\! 2$ with 
$t(d_l\! -\! 1)\! =\! o_{l,j}$. We may assume that $i\!\not=\! 0$ and $j\!\not=\! 2n_l\! +\! 1$ since $0^\infty ,\mu\!\in\! V$. Then there is $x\!\in\! V\cap\mathcal{C}$ with $t\!\subseteq\! x$ by the claim, so that $x\!\in\! V\cap\mathcal{C}\cap N_t$.\medskip

\noindent (b) For $(c^\infty ,0^\infty )$, we argue by contradiction. If $l\!\geq\! l_0$ is large enough, then 
$(c^{l+1}a\overline{a}^\infty ,0^{l+1}\overline{a}a^\infty )$ is not in $E$. Let 
$C''\! :=\! (\bigcup_{1\leq l\leq l_0}~(N_{c^l\overline{a}}\cup\bigcup_{s\in\prod_l}~N_{s\overline{a}})\cup
\bigcup_{1s\in\prod_{j\leq l_0}~d_j}~N_{1s})\cap V$. Then 
$(C'',V\!\setminus\! C'')$ is a coloring of $E$ into clopen sets, which is absurd. For $(\mu ,c^\infty )$, we argue similarly, with the clopen set $(\bigcup_{1\leq l\leq l_0}~(N_{c^l\overline{a}}\cup\bigcup_{s\in\prod_l}~N_{s\overline{a}})\cup \bigcup_{0s\in\prod_{j\leq l_0}~d_j}~N_{0s})\cap V$.\medskip

 In order to prove that $\textup{Graph}(o)\!\subseteq\!\overline{E}^{(\mathcal{C}^+)^2}$, towards a contradiction, suppose that there is $i_0\!\in\!\omega$ such that $\big( o^{i_0}(0^\infty ),o^{i_0+1}(0^\infty )\big)\!\notin\!\overline{E}^{(\mathcal{C}^+)^2}$. If $l\!\geq\! l_0$ is large enough, then $i_0\!\leq\! 2n_l$ and
$$(o_{l,i_0}a^{i_0+1}\overline{a}^\infty ,o_{l,i_0+1}\overline{a}^{i_0+2}a^\infty )\!\notin\! E.$$
We set $P_-\! :=\!\{ o^i(0^{l_0+1})\mid i\! <\! i_0\}$, $P_+\! :=\!\{ o^i(0^{l_0+1})\mid i_0\! <\! i\! <\!\pi_{j\leq l_0}~d_j\}$, which defines a partition $(P_-,\{ o_{l_0,i_0}\},P_+)$ of $\prod_{l_0+1}$. Assume first that $i_0$ is even. Let
$$C'''\! :=\! (\bigcup_{1\leq l\leq l_0}~(N_{c^l\overline{a}}\cup\bigcup_{s\in\prod_l}~N_{s\overline{a}})\cup
\bigcup_{s\in P_-\cup\{ o_{l_0,i_0}\},s(0)=0}~N_s\cup\bigcup_{s\in P_+,s(0)=1}~N_s)\cap V.$$
Then $(C''',V\!\setminus\! C''')$ is a coloring of $E$ into clopen sets, which is absurd. If $i_0$ is odd, then we consider 
$$C'''\! :=\! (\bigcup_{1\leq l\leq l_0}~(N_{c^l\overline{a}}\cup\bigcup_{s\in\prod_l}~N_{s\overline{a}})\cup N_c\cup
\bigcup_{s\in P_-\cup\{ o_{l_0,i_0}\},s(0)=1}~N_s\cup\bigcup_{s\in P_+,s(0)=0}~N_s)\cap V$$
(for instance $c^{l+1}a\overline{a}^\infty\!\in\! C'''$ and $o_{l,0}\overline{a}a^\infty\!\notin\! C'''$).\medskip
 
\noindent (2) $\Rightarrow$ (1) Towards a contradiction, suppose that there is a clopen subset $C$ of $V$ with the property that 
$E\cap\big( C^2\cup (V\!\setminus\! C)^2\big)\! =\!\emptyset$, and by (a) we may assume that $c^\infty\!\in\! C$. (b) gives infinitely many $l$'s such that $(c^{l+1}a\overline{a}^\infty ,0^{l+1}\overline{a}a^\infty )$ is in $E$, and infinitely many $l$'s with 
$(\mu\vert (l\! +\! 1)a^{2n_l+2}\overline{a}^\infty ,c^{l+1}\overline{a}a^\infty )\!\in\! E$. For these large enough $l$'s, 
$0^{l+1}\overline{a}a^\infty\!\notin\! C$ and $\mu\vert (l\! +\! 1)a^{2n_l+2}\overline{a}^\infty\!\notin\! C$. By (a), 
$0^\infty ,\mu\!\in\! V\!\setminus\! C$. This gives $l_0\!\in\!\omega$ such that 
$V\cap (N_{0^{l_0+1}}\cup N_{\mu\vert (l_0+1)})\!\subseteq\! V\!\setminus\! C$.\medskip

 (a) provides $x_1\!\in\! V\cap\mathcal{C}$ with $o_{l_0,1}\!\subseteq\! x_1$, which gives $l_1\! >\! l_0$ such that 
$V\cap N_{x_1\vert l_1}\!\subseteq\! C$ or $V\cap N_{x_1\vert l_1}\!\subseteq\! V\!\setminus\! C$. Let $i_1\!\leq\! 2n_{l_1}$ with $o_{l_1,i_1}\! =\! (x_1\vert l_1)0$. (a) provides $x_2\!\in\! V\cap\mathcal{C}$ with $o_{l_1,i_1+1}\!\subseteq\! x_2$, which gives 
$l_2\! >\! l_1$ such that $V\cap N_{x_2\vert l_2}\!\subseteq\! C$ or $V\cap N_{x_2\vert l_2}\!\subseteq\! V\!\setminus\! C$. Continuing like this we get, for each $j\!\leq\! 2n_{l_0}$, $x_{j+1}\!\in\! V\cap\mathcal{C}$ with $o_{l_j,i_j+1}\!\subseteq\! x_{j+1}$, and $l_{j+1}\! >\! l_j$ such that $V\cap N_{x_{j+1}\vert l_{j+1}}\!\subseteq\! C$ or $V\cap N_{x_{j+1}\vert l_{j+1}}\!\subseteq\! V\!\setminus\! C$, with 
$i_0\! :=\! 0$.\medskip

 By (b), $\big( o^{i_1-1}(0^\infty ),o^{i_1}(0^\infty )\big)\!\in\!\overline{E}^{(\mathcal{C}^+)^2}$, so that $E$ meets 
$N_{0^{l_0+1}}\!\times\! N_{o_{l_1,i_1}}$. As $V\cap N_{0^{l_0+1}}$ is contained in $V\!\setminus\! C$, this implies that 
$V\cap N_{o_{l_1,i_1}}\!\subseteq\! C$. By (b), $\big( o^{i_2-1}(0^\infty ),o^{i_2}(0^\infty )\big)\!\in\!\overline{E}^{(\mathcal{C}^+)^2}$, so that $E$ meets $N_{o_{l_1,i_1}}\!\times\! N_{o_{l_2,i_2}}$. As $V\cap N_{o_{l_1,i_1}}\!\subseteq\! C$, this implies that $V\cap N_{o_{l_2,i_2}}\!\subseteq\! V\!\setminus\! C$.

\vfill\eject

 More generally, if $1\!\leq\! j\!\leq\! 2n_{l_0}\! +\! 1$, then 
$V\cap N_{o_{l_j,i_j}}\!\subseteq\! C$ when $j$ is odd, and $V\cap N_{o_{l_j,i_j}}\!\subseteq\! V\!\setminus\! C$ when $j$ is even. As 
$\mu\vert (l_0\! +\! 1)$ is an initial segment of $o_{l_{2n_{l_0}+1},i_{2n_{l_0}+1}}$, 
$$V\cap N_{o_{l_{2n_{l_0}+1},i_{2n_{l_0}+1}}}\!\subseteq\! C\cap N_{\mu\vert (l_0+1)}\!\subseteq\! C\!\setminus\! C\mbox{,}$$ which is the desired contradiction.\hfill{$\square$}\medskip

 The compactness ensures some surjectivity.

\begin{lem} \label{ont} Let ${\bf d}\!\in\!\mathfrak{D}$, $V\!\subseteq\!\mathcal{C}^+$, $X$ be a 0DMC space, and $G$ be a digraph on $X$ having CCN at least three such that $(X,G)\preceq_c(V,\mathbb{G}_o)$, with witness $\varphi$. Then $\varphi$ is onto $\mathcal{C}\!\subseteq\! V$.\end{lem}

\noindent\emph{Proof.}\ Note that $(\varphi [X],(\varphi\!\times\!\varphi )[G])$ has CCN three. By Lemma \ref{charsub++}, 
$\mathcal{C}$ is contained in $\overline{\varphi [X]\cap\mathcal{C}}^\mathcal{C}$. As $\varphi [X]$ is compact, 
$\overline{\varphi [X]\cap\mathcal{C}}^\mathcal{C}\! =\!\varphi [X]\cap\mathcal{C}$, and thus 
$\mathcal{C}\!\subseteq\!\varphi [X]$.\hfill{$\square$}\medskip

 We now prove some minimality of the $\mathbb{G}_o$'s. 
 
\begin{thm} \label{Gomin} Let ${\bf d}\!\in\!\mathfrak{D}$, $V$ be a compact subspace of $\mathcal{C}^+$, and 
$E\!\subseteq\!\mathbb{G}_o\cap V^2$ be a graph with CCN at least three. Then $(\mathcal{C}^+,\mathbb{G}_o)\preceq_c^i(V,E)$.\end{thm}
 
\noindent\emph{Proof.}\ Note that $\mathcal{C}^+\! =\!\mbox{proj}[\mathbb{G}_o]\cup\{ c^\infty\}\cup\mathcal{C}$. By Lemma \ref{charsub++} and compactness of $V$, $\{ c^\infty\}\cup\mathcal{C}\!\subseteq\! V$. We have to define 
$\varphi\! :\!\mathcal{C}^+\!\rightarrow\! V$. The map $\varphi$ will be the identity on $\{ c^\infty\}\cup\mathcal{C}$. Let 
$\varepsilon\!\in\!\{ a,\overline{a}\}$, $l\!\in\!\omega$, and either $s\! =\! c^{l+1}$ and $i\! =\! 0$, or $s\! =\! o_{l,i}$. We define 
$\varphi (s\varepsilon^{i+1}\overline{\varepsilon}^\infty )$ by induction on $l$, in such a way that 
$s\!\subseteq\!\varphi (s\varepsilon^{i+1}\overline{\varepsilon}^\infty )$ and 
$\varphi (s\varepsilon^{i+1}\overline{\varepsilon}^\infty )$ finishes with some $\varepsilon^{j+1}\overline{\varepsilon}^\infty$. As $(c^\infty ,0^\infty )\!\in\!\overline{E}\! :=\!\overline{E}^{(\mathcal{C}^+)^2}$, we can find 
$\big(\varphi (ca\overline{a}^\infty ),\varphi (0\overline{a}a^\infty )\big)$ in $E$ as desired. As 
$(0^\infty ,10^\infty )\!\in\!\overline{E}$, we can find $\big(\varphi (0a\overline{a}^\infty ),\varphi (1\overline{a}^2a^\infty )\big)$ in $E$ as desired. As $(\mu ,c^\infty )\!\in\!\overline{E}$, we can find 
$\big(\varphi (1a^2\overline{a}^\infty ),\varphi (c\overline{a}a^\infty )\big)\!\in\! E$ as desired. Note that we are done for 
$l\! =\! 0$. The general case is similar, we ensure the injectivity of $\varphi$ by avoiding the finitely many previously chosen sequences.\hfill{$\square$}\medskip
 
 The next lemma will provide several $\preceq_c$-antichains.

\begin{lem} \label{prepanti} Let ${\bf d},{\bf d}'\!\in\!\mathfrak{C}$ such that the $d_j,d'_l$'s are prime, $d'_l$ is not in 
$\{ d_j\mid j\!\in\!\omega\}$ if $d'_l\!\not=\! 3$ and $l$ is large enough, ${\bf d}'$ is unbounded, and 
$(\mathcal{C}_{\bf d},G_{o_{\bf d}})\preceq_c(\mathcal{C}_{{\bf d}'},G_{o_{{\bf d}'}})$ with witness $\varphi$. Then $\varphi$ is not onto.\end{lem}

\noindent\emph{Proof.}\ We argue by contradiction. As $\varphi$ is uniformly continuous, there is, for each $l\!\in\!\omega$, 
$L\! :=\! L(l)\!\geq\! l$ with the property that $\varphi (x)\vert (l\! +\! 1)\! =\!\varphi (y)\vert (l\! +\! 1)$ if 
$x,y\!\in\!\mathcal{C}_{\bf d}$ and $x\vert (L\! +\! 1)\! =\! y\vert (L\! +\! 1)$. We choose $l\!\geq\! 2$ with the property that 
$d'_{l+1}\!\notin\!\{ d_j\mid j\!\in\!\omega\}$, which is possible for infinitely many $l$'s. As $\varphi$ is onto, we can find a surjection ${\Gamma_l\! :\!\prod_{L+1}\!\rightarrow\!\prod_{j\leq l}~d'_j}$ such that 
$\varphi [\mathcal{C}_{\bf d}\cap N_s]\!\subseteq\! N_{\Gamma_l(s)}$. We enumerate 
${\prod_{L+1}\! :=\!\{ s^L_i\mid i\! <\!\pi_{j\leq L}~d_j\}}$ in the order defined by $o_{\bf d}$, starting with $0^{L+1}$. Note that, respecting this order, $\{ s^{L(l+1)}_i\mid i\! <\!\pi_{j\leq L(l+1)}~d_j\}$ is 
$$\{ s^L_i0^{L(l+1)-L}\mid i\! <\!\pi_{j\leq L}~d_j\}\cup\ldots\cup\{ s^L_i{^\frown}_{L<j\leq L(l+1)}~(d_j\! -\! 1)\mid i\! <\!\pi_{j\leq L}~d_j\} .$$
This implies that $\{\Gamma_{l+1}(s^{L(l+1)}_i)\mid i\! <\!\pi_{j\leq L(l+1)}~d_j\}$ is 
$$\{\Gamma_l(s^L_i)\varepsilon_i^0\mid i\! <\!\pi_{j\leq L}~d_j\}\cup\ldots\cup\{\Gamma_l(s^L_i)\varepsilon_i^{(\pi_{L<j\leq L(l+1)}~d_j)-1}\mid i\! <\!\pi_{j\leq L}~d_j\} .$$ 

 As $o_{{\bf d}'}(x)\vert 2\!\not=\! o_{{\bf d}'}^{-1}(x)\vert 2$ for each $x\!\in\!\mathcal{C}_{{\bf d}'}$, 
$o_{{\bf d}'}\big(\Gamma_l(s^L_i)\big)\!\not=\! o^{-1}_{{\bf d}'}\big(\Gamma_l(s^L_i)\big)$ if $i\! <\!\pi_{j\leq L}~d_j$. As $\varphi$ is a homomorphism, $\Gamma_l(s^L_{i+1})$ can only be the image or the inverse image of $\Gamma_l(s^L_i)$ by the map $o_{{\bf d}'}$ if 
$i\! +\! 1\! <\!\pi_{j\leq L}~d_j$. Similarly, $\Gamma_{l+1}(s^{L(l+1)}_{k+1})$ can only be the image or the inverse image of 
$\Gamma_{l+1}(s^{L(l+1)}_k)$ by the map $o_{{\bf d}'}$ if $k\! +\! 1\! <\!\pi_{j\leq L(l+1)}~d_j$. If 
$s^{L(l+1)}_k$ extends $s^L_i$, then $\Gamma_{l+1}(s^{L(l+1)}_k)$ extends $\Gamma_l(s^L_i)$, and $\Gamma_l(s^L_{i+1})$ is the image of $\Gamma_l(s^L_i)$ if and only if $\Gamma_{l+1}(s^{L(l+1)}_{k+1})$ is the image of $\Gamma_{l+1}(s^{L(l+1)}_k)$ since 
$o_{{\bf d}'}\big(\Gamma_l(s^L_i)\big)\!\not=\! o^{-1}_{{\bf d}'}\big(\Gamma_l(s^L_i)\big)$ if $i\! <\!\pi_{j\leq L}~d_j$. This implies that 
${m\! :=\!\varepsilon^k_{\pi_{j\leq L}~d_j-1}\! -\!\varepsilon^k_0\ (\mbox{mod }d'_{l+1})}$ does not depend on $k$. Similarly, either 
$\varepsilon^{k+1}_0\! =\!\varepsilon^k_{\pi_{j\leq L}~d_j-1}$ for each $k$, or 
${\varepsilon^{k+1}_0\!\equiv\!\varepsilon^k_{\pi_{j\leq L}~d_j-1}\! +\! 1\ (\mbox{mod }d'_{l+1})}$ for each $k$, or 
${\varepsilon^{k+1}_0\!\equiv\!\varepsilon^k_{\pi_{j\leq L}~d_j-1}\! -\! 1\ (\mbox{mod }d'_{l+1})}$ for each $k$.\medskip

 Assume first that $\varepsilon^{k+1}_0\! =\!\varepsilon^k_{\pi_{j\leq L}~d_j-1}$ for each $k$. An induction shows that 
$$\varepsilon^k_{\pi_{j\leq L}~d_j-1}\!\equiv\!\varepsilon^0_0\! +\! (k\! +\! 1)m\ (\mbox{mod }d'_{l+1}) .$$ 
Thus $\varepsilon^0_0\!\equiv\!\varepsilon^0_0\! +\! (\pi_{L<j\leq L(l+1)}~d_j)m\ (\mbox{mod }d'_{l+1})$ and $d'_{l+1}$ divides $m$ since $d'_{l+1}$ is not in $\{ d_j\mid j\!\in\!\omega\}$.\medskip

 Assume now that $\varepsilon^{k+1}_0\!\equiv\!\varepsilon^k_{\pi_{j\leq L}~d_j-1}\! +\! 1\ (\mbox{mod }d'_{l+1})$ for each $k$. An induction shows that 
$$\varepsilon^k_{\pi_{j\leq L}~d_j-1}\!\equiv\!\varepsilon^0_0\! +\! (k\! +\! 1)(m\! +\! 1)\! -\! 1\ (\mbox{mod }d'_{l+1}).$$ 
Thus $\varepsilon^0_0\!\equiv\!\varepsilon^0_0\! +\! (\pi_{L<j\leq L(l+1)}~d_j)(m\! +\! 1)\ (\mbox{mod }d'_{l+1})$ and 
$d'_{l+1}$ divides $m\! +\! 1$.\medskip

 Assume now that $\varepsilon^{k+1}_0\!\equiv\!\varepsilon^k_{\pi_{j\leq L}~d_j-1}\! -\! 1\ (\mbox{mod }d'_{l+1})$ for each $k$. An induction shows that 
$$\varepsilon^k_{\pi_{j\leq L}~d_j-1}\!\equiv\!\varepsilon^0_0\! +\! (k\! +\! 1)(m\! -\! 1)\! +\! 1\ (\mbox{mod }d'_{l+1}).$$ 
Thus $\varepsilon^0_0\!\equiv\!\varepsilon^0_0\! +\! (\pi_{L<j\leq L(l+1)}~d_j)(m\! -\! 1)\ (\mbox{mod }d'_{l+1})$ and $d'_{l+1}$ divides $m\! -\! 1$.\medskip

 In all cases, this shows that $\varepsilon^k_i$ does not depend on $k$. This argument can be extended to any length strictly greater than $l$. This cannot always hold since the sequence ${\bf d}'$ is unbounded.\hfill{$\square$}\medskip
 
 In our applications of Lemma \ref{prepanti}, $d_l,d'_l$ can be $3$ for infinitely many $l$'s.
 
\begin{thm} \label{anticha} There is a map $\Phi\! :\! 2^\omega\!\rightarrow\!\mathfrak{D}$ such that 
$(\mathcal{C}^+_{\Phi (\alpha )},\mathbb{G}_{o_{\Phi (\alpha )}})\not\preceq_c(\mathcal{C}_{\Phi (\beta )}^+,\mathbb{G}_{o_{\Phi (\beta )}})$ if $\alpha\!\not=\!\beta$.\end{thm}
 
\noindent\emph{Proof.}\ Let $(p_n)_{n\in\omega}$ be the sequence of prime numbers. We define, for each 
$\alpha\!\in\! 2^\omega$, $S_\alpha\!\subseteq\!\omega$ by 
$S_\alpha\! :=\!\{ 0\}\cup\{ p_0^{\alpha (0)+1}\ldots p_n^{\alpha (n)+1}\mid n\!\in\!\omega\}$. Note that $S_\alpha$ is infinite,  contains $0$, and $S_\alpha\cap S_\beta$ is finite if $\alpha\!\not=\!\beta$. In this proof, we consider $(d_\alpha )_0\! =\! 2$, 
$(d_\alpha )_j\! =\! 3$ if $j\!\notin\! S_\alpha$, $(d_\alpha )_j\! =\! p_{j+1}$ if $0\! <\! j\!\in\! S_\alpha$, so that 
$\Phi (\alpha )\! :=\! {\bf d}_\alpha\!\in\!\mathfrak{D}$ is unbounded, the $(d_\alpha )_j$'s are prime, $(d_\beta )_l$ is not in 
$\{ (d_\alpha )_j\mid j\!\in\!\omega\}$ if $\alpha\!\not=\!\beta$, $(d_\beta )_l\!\not=\! 3$ and $l$ is large enough.\medskip

 So assume that ${\bf d},{\bf d}'\!\in\!\mathfrak{D}$, the $d_j,d'_l$'s are prime, $d'_l$ is not in $\{ d_j\mid j\!\in\!\omega\}$ if 
$d'_l\!\not=\! 3$ and $l$ is large enough, ${\bf d}'$ is unbounded, and 
$(\mathcal{C}^+_{\bf d},\mathbb{G}_{o_{\bf d}})\preceq_c(\mathcal{C}_{{\bf d}'}^+,\mathbb{G}_{o_{{\bf d}'}})$ with witness 
$\varphi$. By Lemmas \ref{c}, \ref{ont} and Proposition \ref{belowfin}, $\varphi$ is onto $\{ c^\infty\}\cup\mathcal{C}_{{\bf d}'}$. As $\varphi$ is uniformly continuous, there is, for each $l\!\in\!\omega$, $L\! :=\! L(l)\!\geq\! l$ with the property that 
$\varphi (x)\vert (l\! +\! 1)\! =\!\varphi (y)\vert (l\! +\! 1)$ if $x,y\!\in\!\mathcal{C}^+_{\bf d}$ and 
$x\vert (L\! +\! 1)\! =\! y\vert (L\! +\! 1)$.\medskip

\noindent\emph{Claim 1.}\it\ If $l\!\in\!\omega$, then there is $L'_0\!\geq\! L(l)$ such that, for each $L'\!\geq\! L'_0$, each 
$i\!\leq\! 2n_{L'}\! +\! 1$ and each $\varepsilon\!\in\!\{ a,\overline{a}\}$, 
$\varphi (o^{\bf d}_{L',i}\varepsilon^{i+1}\overline{\varepsilon}^\infty )\vert (l\! +\! 1)\!\in\!\prod_{j\leq l}~(d'_j\cup\{ c\} )$.\rm\medskip

 Indeed, towards a contradiction, suppose that there is $l_0$ such that, for each $k\!\geq\! L(l_0)$, we can find 
${L'_k\!\geq\! k}$, $i_k,\varepsilon_k$ and $m_k\!\leq\! l_0$ such that 
$\varphi (o^{\bf d}_{L'_k,i_k}\varepsilon_k^{i_k+1}\overline{\varepsilon_k}^\infty )(m_k)\!\in\!\{ a,\overline{a}\}$. As $\varphi$ is a homomorphism, the uniform continuity of $\varphi$ and an induction on $i\! -\! i_k$  show that 
$\varphi (o^{\bf d}_{L'_k,i}\varepsilon^{i+1}\overline{\varepsilon}^\infty )(m_k),
\varphi (c^{L'_k+1}\varepsilon\overline{\varepsilon}^\infty )(m_k)$ are in $\{ a,\overline{a}\}$ if $i\!\leq\! 2n_{L'_k}\! +\! 1$ and 
$\varepsilon\!\in\!\{ a,\overline{a}\}$. If now $i\!\in\!\omega$, then $\big( o_{\bf d}^i(0^\infty ),o_{\bf d}^{i+1}(0^\infty )\big)$ is the limit of couples of the form $(o_{L'_k,i}^{\bf d}a^{i+1}\overline{a}^\infty ,o_{L'_k,i+1}^{\bf d}\overline{a}^{i+2}a^\infty )$, where we may assume that $(m_k)_k$ is constant. The continuity of $\varphi$ implies that $\varphi\big( o_{\bf d}^i(0^\infty )\big)(m_0)$ is in $\{ a,\overline{a}\}$. This shows that $\mathcal{C}^+_{\bf d}\!\subseteq\!\varphi^{-1}(\mbox{proj}[\mathbb{G}_{o_{{\bf d}'}}])$ and $c^\infty\!\notin\!\varphi [\mathcal{C}^+_{\bf d}]$, which is absurd.\hfill{$\diamond$}\medskip

\noindent\emph{Claim 2.}\it\ $\varphi [\mathcal{C}_{\bf d}]\! =\!\mathcal{C}_{{\bf d}'}$ and $\varphi (c^\infty )\! =\! c^\infty$.\rm\medskip

 Indeed, by Lemma \ref{fcc} and Proposition \ref{belowfin}, it is enough to see that 
$\varphi [\mathcal{C}_{\bf d}]\!\subseteq\!\mathcal{C}_{{\bf d}'}$. By Lemma \ref{flip'} and Proposition \ref{belowfin}, it is enough to see that 
$\varphi [\mathcal{C}_{\bf d}\cap N_t]\!\not\subseteq\!\{ c^\infty\}$ and 
$\varphi [\mathcal{C}_{\bf d}\cap N_t]\!\not\subseteq\!\{ s\varepsilon^{m+1}\overline{\varepsilon}^\infty\}$ if 
$t\!\in\!\bigcup_{l\in\omega}~\prod_l$. We argue by contradiction.\medskip

 If the singleton is of the form $\{ s\varepsilon^{m+1}\overline{\varepsilon}^\infty\}$ with $s\! =\! c^{k+1}$ or $s\!\in\!\prod_{j\leq k}~d'_j$, then we choose $i\!\in\!\omega$ with $t\!\subseteq\! o^i(0^\infty )$, so that 
$\varphi\big( o^i(0^\infty )\big)\! =\! s\varepsilon^{m+1}\overline{\varepsilon}^\infty$. We may assume that 
$s\varepsilon^{m+1}\overline{\varepsilon}\!\subseteq\!\varphi (z)$ if $z\!\in\!\mathcal{C}^+_{\bf d}\cap N_t$. We apply Claim 1 to 
$l\! :=\!\vert s\vert$, which gives $L'_0$. We choose $L'\!\geq\!\mbox{max}(\vert t\vert ,L'_0)$ with the property that 
$i\!\leq\! 2n_{L'}\! +\! 1$. Then 
$\varphi (o^{\bf d}_{L',i}a^{i+1}\overline{a}^\infty )\vert (l\! +\! 1)\! =\! s\varepsilon\!\in\!\prod_{j\leq l}~(d'_j\cup\{ c\} )$, which is the desired contradiction. If the singleton is $\{ c^\infty\}$, then we may assume that $\varphi (z)(0)\! =\! c$ if 
$z\!\in\!\mathcal{C}^+_{\bf d}\cap N_t$. We fix $l\! >\!\vert t\vert$ with $n^{{\bf d}'}_l\! >\! 2n^{\bf d}_{\vert t\vert +1}\! +\! 2$. Let 
$w\!\in\!\mathcal{C}_{{\bf d}'}\!\setminus\!\{\varphi (c^\infty )\}$, and $u\!\in\!\mathcal{C}_{\bf d}$ with $\varphi (u)\! =\! w$. Note that there is $i\!\leq\! 2n^{\bf d}_{L+1}\! +\! 1$ with $o^{\bf d}_{L+1,i}\! =\! u\vert (L\! +\! 1)0$, and 
$\varphi (u\vert (L\! +\! 1)0a^{i+1}\overline{a}^\infty )\vert (l\! +\! 1)\! =\! w\vert (l\! +\! 1)$. Also, there is 
$k\!\leq\! 2n^{\bf d}_{L+1}\! +\! 1$ such that $t\!\subseteq\! o^{\bf d}_{L+1,k}$ and 
$\vert k\! -\! i\vert\!\leq\! 2n^{\bf d}_{\vert t\vert +1}\! +\! 1$. Note that $\varphi (o^{\bf d}_{L+1,k}a^{k+1}\overline{a}^\infty )(0)\! =\! c$, and $\varphi (o^{\bf d}_{L+1,k}a^{k+1}\overline{a}^\infty )\vert (l\! +\! 1)\! =\! c^{l+1}$ since $\varphi$ is uniformly continuous. Also, there is 
${i'\! <\!\pi_{j\leq l}~d'_j}$ with $w\vert (l\! +\! 1)\! =\! o^{{\bf d}'}_{l,i'}$. Note that 
$i'\!\leq\! 2n^{\bf d}_{\vert t\vert +1}\! +\! 1$ or $\vert 2n^{{\bf d}'}_l\! +\! 1\! -\! i'\vert\!\leq\! 2n^{\bf d}_{\vert t\vert +1}\! +\! 1$ since $\varphi$ is a homomorphism. Thus $\{ y\vert (l\! +\! 1)\mid y\!\in\!\mathcal{C}_{{\bf d}'}\}$ cannot be $\prod_{j\leq l}~d'_j$, which contradicts the fact that $\varphi$ is onto $\mathcal{C}_{{\bf d}'}$.\hfill{$\diamond$}\medskip

 We set $\psi\! :=\!\varphi_{\vert\mathcal{C}_{\bf d}}$. By Claim 2, $\psi$ takes values in $\mathcal{C}_{{\bf d}'}$ and is onto 
$\mathcal{C}_{{\bf d}'}$. By Lemma \ref{prepanti}, it remains to see that $\varphi\!\times\!\varphi$ sends $G_{o_{\bf d}}$ into 
$G_{o_{{\bf d}'}}$. As $o_{\bf d}$ is a minimal homeomorphism, it is fixed point free and thus 
$G_{o_{\bf d}}\! =\! s\big(\textup{Graph}(o_{\bf d})\big)$, and similarly with ${\bf d}'$. By Proposition \ref{belowfin} and Lemma \ref{charsub++}, $\textup{Graph}(o_{\bf d})$ is contained in $\overline{\mathbb{G}_{o_{\bf d}}}^{(\mathcal{C}^+_{\bf d})^2}$, so that $\varphi\!\times\!\varphi$ sends $G_{o_{\bf d}}$ into 
$\overline{\mathbb{G}_{o_{{\bf d}'}}}^{(\mathcal{C}^+_{{\bf d}'})^2}\cap\mathcal{C}^2_{{\bf d}'}$ by Claim 2. It remains to note that $\overline{\mathbb{G}_{o_{{\bf d}'}}}^{(\mathcal{C}^+_{{\bf d}'})^2}$ is contained in 
$\mathbb{G}_{o_{{\bf d}'}}\cup s\big(\{ (c^\infty ,0^\infty ),(\mu, c^\infty )\}\cup\textup{Graph}(o_{{\bf d}'})\big)$ to conclude.
\hfill{$\square$}\medskip

 In the compact case, the $\preceq^i_c$-minimality can be seen on subgraphs.

\begin{lem} \label{cpctmin} Let $(X,G)$ in $\mathfrak{K}$. The following are equivalent:\smallskip

\noindent (1) $(X,G)$ is $\preceq^i_c$-minimal in $\mathfrak{K}$,\smallskip

\noindent (2) $(X,G)\preceq^i_c(V,E)$ if $V$ is a compact subset of $X$ and $E\!\subseteq\! G$ is a graph on $V$ with $\chi_c(V,E)\!\geq\! 3$.\end{lem}

\noindent\emph{Proof.}\ It is enough to see that (2) implies (1). So let $X'$ be a 0DMC space, and $G'$ be a graph on $X'$ with 
$\chi_c(X',G')\!\geq\! 3$ and $(X',G')\preceq^i_c(X,G)$ with witness $\varphi$. We set $V\! :=\!\varphi [X']$ and also 
$E\! :=\! (\varphi\!\times\!\varphi )[G']$. As $X'$ is compact, so is $V$, and $E\!\subseteq\! G$ is a graph on $V$. Note that $\varphi$ is a witness for the fact that $(X',G')\preceq^i_c(V,E)$, so that $\chi_c(V,E)\!\geq\! 3$. By (2), $(X,G)\preceq^i_c(V,E)$. By compactness of $X'$ again, $\varphi^{-1}$ is a witness for the fact that $(V,E)\preceq^i_c(X',G')$, which as desired implies that $(X,G)\preceq^i_c(X',G')$.\hfill{$\square$}

\vfill\eject\medskip

\noindent\emph{Proof of Theorem \ref{eantichmin}.}\ Fix ${\bf d}\!\in\!\mathfrak{D}$. Note that $(\mathcal{C}^+,\mathbb{G}_o)$ satisfies the properties (a) and (b), by Lemma \ref{D2gen} and Proposition \ref{belowfin}. For (c), i.e., the minimality of 
$(\mathcal{C}^+,\mathbb{G}_o)$, we apply Proposition \ref{belowfin}, Theorem \ref{Gomin} and Lemma \ref{cpctmin}. It remains to apply Theorem \ref{anticha}.\hfill{$\square$}\medskip

\noindent\emph{Remark.}\ [P, Theorem 11.38] shows that the $o_{\Phi (\alpha )}$'s involved in Theorem \ref{eantichmin} are pairwise not flip-conjugate, as announced in the introduction just before Theorem \ref{eantichmin}.\medskip

\noindent\emph{Proof of Theorem \ref{infdecrcompact}.}\ Fix ${\bf d}\!\in\!\mathfrak{D}$. The idea is to modify $\mathbb{G}_o$. Let $d$ be a letter not in $\omega\cup\{ c,a,\overline{a}\}$. We set, for 
$l\!\in\!\omega$,\medskip

\leftline{$\mathbb{H}_l\! :=\!\{ (c^{l+1}d^{j+1}a\overline{a}^\infty ,o_{l,0}d^{j+1}\overline{a}a^\infty)\mid j\!\in\!\omega\} ~\cup$}\smallskip

\rightline{$\{ (o_{l,i}d^{j+1}a^{i+1}\overline{a}^\infty ,o_{l,i+1}d^{j+1}\overline{a}^{i+2}a^\infty )\mid j\!\in\!\omega\wedge i\!\leq\! 2n_l\} ~\cup$}\smallskip

\rightline{$\{ (o_{l,2n_l+1}d^{j+1}a^{2n_l+2}\overline{a}^\infty ,c^{l+1}d^{j+1}\overline{a}a^\infty )\mid j\!\in\!\omega\} .$}\medskip

\noindent We then set $\mathbb{O}_p\! :=\!\bigcup_{l\geq p}~\mathbb{H}_l$ and $\mathbb{G}_p\! :=\! s(\mathbb{O}_p)$, so that $\mathbb{G}_p$ is a countable graph on the compact space 
$\mathcal{K}\! :=\!\prod_{j\in\omega}~(d_j\cup\{ c,a,\overline{a},d\} )$. We set 
$\mathbb{K}\! :=\!\overline{\mbox{proj}[\mathbb{G}_0]}^{\mathcal{K}}$, so that $\mathbb{K}$ is a 0DMC space and 
$\mathbb{G}_p$ is a graph on $\mathbb{K}$, $(\mathbb{G}_p)_{p\in\omega}$ is $\subseteq$-decreasing and thus 
$\big( (\mathbb{K},\mathbb{G}_p)\big)_{p\in\omega}$ is $\preceq_c$ and $\preceq^i_c$-decreasing. As in the proofs of Lemmas \ref{chromgen} and \ref{D2gen}, we see that $\mathbb{G}_p$ has CCN three, 
$\boraone\oplus\bormone$ chromatic number two, and is $D_2(\bormone )$.\medskip

 It remains to see that $(\mathbb{K},\mathbb{G}_p)$ is not $\preceq_c$-below $(\mathbb{K},\mathbb{G}_{p+1})$. Towards a contradiction, suppose that there is $\varphi\! :\!\mathbb{K}\!\rightarrow\!\mathbb{K}$. The continuity of $\varphi$ implies that 
$\overline{\mathbb{G}_p}\!\subseteq\! (\varphi\!\times\!\varphi )^{-1}(\overline{\mathbb{G}_{p+1}})$. Note that 
$$(c^{p+1}d^\infty ,o_{p,0}d^\infty ,\ldots ,o_{p,2n_p+1}d^\infty )$$ 
is a $\overline{\mathbb{G}_p}$-cycle of length $2n_p\! +\! 3$, and  therefore has to be sent in a 
$\overline{\mathbb{G}_{p+1}}$-cycle of length at most $2n_p\! +\! 3$. But such a $\overline{\mathbb{G}_{p+1}}$-cycle does not exist since\medskip

\leftline{$\overline{\mathbb{G}_p}\! :=\!\mathbb{G}_p\cup s\big(\{ (c^{l+1}d^\infty ,o_{l,0}d^\infty)\mid l\!\geq\! p\}\cup
\{ (o_{l,i}d^\infty ,o_{l,i+1}d^\infty )\mid l\!\geq\! p\wedge i\!\leq\! 2n_l\} ~\cup$}\smallskip

\rightline{$\{ (o_{l,2n_l+1}d^\infty ,c^{l+1}d^\infty )\mid l\!\geq\! p\}\cup\{ (c^\infty ,0^\infty ),(\mu ,c^\infty )\}\cup\textup{Graph}(o)\big) .$}\medskip

\noindent This finishes the proof.\hfill{$\square$}

\section{$\!\!\!\!\!\!$ Graphs induced by a function: general facts}

\emph{Remarks.}\ (1) In the case of finite spaces, the quasi-order $\preceq$ on the class of graphs induced by a partial bijection  with chromatic number at least three is linear. Indeed, such a space can be decomposed in pairwise $f$-unrelated injective walks of the form $\{ x,f(x),\ldots ,f^l(x)\}$. As the chromatic number is at least three, one of these walks gives an odd cycle. The graph induced by the bijection is $\preceq$-equivalent to its odd cycle of minimal length. As the odd cycles are $\preceq$-comparable, so too are all these graphs.\medskip

\noindent (2) Note that the map $A\!\mapsto\!\oplus_{p\in A}~(2p\! +\! 3,C_{2p+3})$ is an embedding of the quasi-order of inclusion on the set of finite subsets of $\omega$ into the quasi-order $\preceq^i$ on the class of graphs induced by a bijection on a finite set with chromatic number at least three.\medskip

 Under some relatively weak assumptions, we can characterize when the CCN of $G_f$ is big. Note that we extend Theorem \ref{eq1++} under these assumptions.

\vfill\eject

\begin{thm} \label{CNG} Let $X$ be a 0DMS space of cardinality at least two, $f\! :\! X\!\rightarrow\! X$ be a homeomorphism, and $x\!\in\! X$ with $\overline{\mbox{Orb}_f(x)}\! =\! X$. The following are equivalent:\smallskip

\noindent (1) $(X,G_f)$ has CCN at least three,\smallskip

\noindent (2) $\overline{\{ f^{2n}(x)\mid n\!\in\!\mathbb{Z}\}}\cap
\overline{\{ f^{2p+1}(x)\mid p\!\in\!\mathbb{Z}\}}$ is not empty,\smallskip

\noindent (3) $\Delta (X)\cap\overline{\bigcup_{p\in\omega}~G_f^{2p+1}}$ is not empty,\smallskip

\noindent (4) $\Delta (X)\cap\overline{\bigcup_{p\in\omega}~\overline{G_f}^{2p+1}}$ is not empty,\smallskip

\noindent (5) $(\mathcal{N},\mathbb{G}_m)\preceq_c(X,G_f)$.\end{thm}

\noindent\emph{Proof.}\ (1) $\Rightarrow$ (2) Note that $f$ sends $\{ f^{2n}(x)\mid n\!\in\!\mathbb{Z}\}$ onto 
$\{ f^{2n+1}(x)\mid n\!\in\!\mathbb{Z}\}$. As $f$ is a homeomorphism, $f$ sends 
$C_e\! :=\!\overline{\{ f^{2n}(x)\mid n\!\in\!\mathbb{Z}\}}$ onto $C_o\! :=\!\overline{\{ f^{2p+1}(x)\mid p\!\in\!\mathbb{Z}\}}$. Note that 
$$X\! =\!\overline{\mbox{Orb}_f(x)}\! =\! C_e\cup C_o.$$ 
If $C_e$ is disjoint from $C_o$, then $(C_e,C_o)$ defines a continuous coloring of $(X,G_f)$, contradicting (1).\medskip

\noindent (2) $\Rightarrow$ (3) By (2), there is $y$ in the intersection. Let $O$ be an open neighborhood of $y$, and $m,n$ be integers with $f^{2n}(x),f^{2m+1}(x)\!\in\! O$. As $X$ has cardinality at least two, $f_{\vert\text{Orb}_f(x)}$ is fixed point free. We put $p\! :=\!\vert m\vert\! +\!\vert n\vert$, so that $f^{2n}(x),\ldots ,x,\ldots ,f^{2m+1}(x)$ is a witness for the fact that $O^2$ meets $G_f^{2p+1}$.\medskip

\noindent (3) $\Rightarrow$ (4) Note that $G_f\!\subseteq\!\overline{G_f}$.\medskip

\noindent (4) $\Rightarrow$ (1) We apply Lemma \ref{chrominter+}.\medskip

\noindent (4) $\Leftrightarrow$ (5) We apply Lemma \ref{eq1+}.\hfill{$\square$}\medskip

\noindent\emph{Notation.}\ The set of fixed points of a map $f$ is very much related to the CCN of 
$G_f$. Let $X$ be a set, and $f\! :\! X\!\rightarrow\! X$ be a partial map. The set 
$F_1\! :=\!\{ x\!\in\!\mbox{Domain}(f)\mid f(x)\! =\! x\}$ of fixed points of $f$ is sometimes also denoted by $F_1^f$. 

\begin{prop} \label{fp} Let $X$ be a 0DMS space, and $f\! :\! X\!\rightarrow\! X$ be a partial continuous function. If $F_1$ is not an open subset of $\mbox{Domain}(f)$, then $\chi_c(X,G_f)\! =\! 2^{\aleph_0}$.\end{prop}

\noindent\emph{Proof.}\ Let $(C_i)_{i\in\omega}$ be a partition of $X$ into clopen sets. As $F_1$ is not open in $\mbox{Domain}(f)$, we can find $x\!\in\! F_1$ and $(x_n)_{n\in\omega}\!\in\! (\mbox{Domain}(f)\!\setminus\! F_1)^\omega$ converging to $x$. Note that $f(x_n)$ is different from $x_n$, and $\big( f(x_n)\big)_{n\in\omega}$ converges to $f(x)\! =\! x$. Let $i$ with $x\!\in\! C_i$. Then we may assume that $x_n,f(x_n)\!\in\! C_i$. This implies that $\big( x_n,f(x_n)\big)\!\in\! G_f\cap C_i^2$.\hfill{$\square$}
 
\begin{cor} \label{corfp} Let $X$ be a 0DMS space, and $f\! :\! X\!\rightarrow\! X$ be a partial continuous function with closed domain.\smallskip

\noindent (a) Exactly one of the following holds:\smallskip

(1) $F_1$ is an open subset of $\mbox{Domain}(f)$,\smallskip

(2) $\chi_c(X,G_f)\! =\! 2^{\aleph_0}$.\smallskip

\noindent (b) If $F_1$ is an open subset of $\mbox{Domain}(f)\!\in\!\borone (X)$ and $f$ is injective, then $\chi_c(X,G_f)\! =\! 0$ if $X\! =\!\emptyset$, $1$ if $F_1\! =\!\mbox{Domain}(f)$ and 
$X\!\not=\!\emptyset$, $\chi_c(X\!\setminus\! F_1,G_f\cap (X\!\setminus\! F_1)^2)$ if $F_1\!\not=\!\mbox{Domain}(f)$.\end{cor}

\noindent\emph{Proof.}\ (a) Assume that (1) holds. Note that 
$s\big(\textup{Graph}(f_{\vert\text{Domain}(f)\setminus F_1})\big)$ and $\Delta (X)$ are disjoint closed relations on the metrizable space $X$. By [K, 22.16], there is a clopen relation $C$ on $X$ separating $\Delta (X)$ from 
$s\big(\textup{Graph}(f_{\vert\text{Domain}(f)\setminus F_1})\big)$. This relation gives a countable continuous coloring of 
$(X,G_f)$ since $X$ is zero-dimensional and second countable. So (2) does not hold.\medskip

 If $F_1$ is not an open subset of $\mbox{Domain}(f)$, then we apply Proposition \ref{fp}.\medskip
  
\noindent (b) If $F_1\!\not=\!\mbox{Domain}(f)$, then we can find $2\!\leq\! n\!\leq\!\omega$ and a continuous coloring 
$c\! :\! X\!\setminus\! F_1\!\rightarrow\! n$ of $(X\!\setminus\! F_1,G_f\cap (X\!\setminus\! F_1)^2)$, by (a). As $f$ is injective, 
$f[\mbox{Domain}(f)\!\setminus\! F_1]\cap F_1\! =\!\emptyset$, so that $F_1$ and $X\!\setminus\! F_1$ are $f$-invariant. The extension of $c$ by $0$ on $F_1$ is a continuous coloring of $(X,G_f)$. Conversely, any continuous coloring of $(X,G_f)$ gives a coloring of $(X\!\setminus\! F_1,G_f\cap (X\!\setminus\! F_1)^2)$, by restriction.\hfill{$\square$}\medskip

 In the introduction, we announced a version of Theorem \ref{LZ} for analytic spaces when $\xi\! =\! 1$. Here are the argument and some precisions. Recall that $\mathbb{X}_1\! :=\!\{ 0^\infty\}\cup\{ 0^n1^\infty\mid n\!\in\!\omega\}$, 
$f_1\! :\!\mathbb{X}_1\!\rightarrow\!\mathbb{X}_1$ is defined by $f_1(0^\infty )\! :=\! 0^\infty$ and 
${f_1(0^{2n+\varepsilon}1^\infty )\! :=\! 0^{2n+1-\varepsilon}1^\infty}$, and 
$\mathbb{R}_1\! :=\!\{ (0^{2n}1^\infty ,0^{2n+1}1^\infty )\mid n\!\in\!\omega\}$. We also define 
$f_0\! :\!\mathbb{X}_1\!\rightarrow\!\mathbb{X}_1$ by $f_0(\alpha )\! :=\! 0^\infty$.

\begin{prop} \label{LZ1} (a) (Lecomte-Zelen\'y) Let $X$ be a zero-dimensional Lindel\"of first countable space, and $R$ be a relation on $X$. Then exactly one of the following holds:\smallskip

\noindent (1) there is a countable continuous coloring of $R$,\smallskip  

\noindent (2) there is $f\! :\!\mathbb{X}_1\!\rightarrow\! X$ continuous such that 
$\mathbb{R}_1\!\subseteq\! (f\!\times\! f)^{-1}(R)$.\smallskip

\noindent In particular, $(\mathbb{X}_1,G_{f_1})$ is $\preceq_c$-minimum in the class of graphs on a 0DMS space with uncountable CCN.\smallskip

\noindent (b) $\{ (\mathbb{X}_1,G_{f_0}),(\mathbb{X}_1,G_{f_1})\}$ is a $\preceq^i_c$-antichain basis for the class of graphs on a 0DMS space with uncountable CCN.\end{prop}

\noindent\emph{Proof.}\ (a) If $\Delta (X)\cap\overline{R}\! =\!\emptyset$, then for each $x\!\in\! X$ there is a clopen neighborhood $C_x$ of $x$ with $R\cap C_x^2\! =\!\emptyset$. As $X$ is Lindel\"of, the covering $(C_x)_{x\in X}$ of $X$ can be replaced with a covering 
$(C_n)_{n\in\omega}$. Replacing $C_n$ with $C_n\!\setminus\! (\bigcup_{m<n}~C_m)$ if necessary, we may assume that the $C_n$'s are pairwise disjoint, which gives a countable continuous coloring of $R$. If there is $(x,x)\!\in\!\overline{R}$, then the fact that $X$ is first countable provides a sequence $(x_n)$ converging to $x$ with $(x_{2n},x_{2n+1})\!\in\! R$. It remains to set $\varphi (0^\infty )\! :=\! x$ and $\varphi (0^n1^\infty )\! :=\! x_n$ to see that (2) holds. If $C$ is a clopen subset of 
$\mathbb{X}_1$ containing $0^\infty$, then we can find $n$ with $0^{2n}1^\infty ,0^{2n+1}1^\infty\!\in\! C$, so that 
$(0^{2n}1^\infty ,0^{2n+1}1^\infty )\!\in\! G_{f_1}\cap C^2$ and $(0^{2n}1^\infty ,0^\infty )\!\in\! G_{f_0}\cap C^2$. This implies that the $(\mathbb{X}_1,G_{f_\varepsilon})$'s have uncountable CCN. It remains to note that $s(\mathbb{R}_1)\! =\! G_{f_1}$ to see that (1) and (2) cannot hold simultaneously.\medskip

\noindent (b) Let $X$ be a 0DMS space, and $R$ be a graph on $X$ with uncountable CCN. We use the proof of (a), which gives $(x_n)$. As $x_{2n}\!\not=\! x_{2n+1}$, we may assume, extracting a subsequence if necessary, that the sequence 
$\big( (x_{2n},x_{2n+1})\big)_{n\in\omega}$ is injective, and that 
$(x_{2n})_{n\in\omega}$ or $(x_{2n+1})_{n\in\omega}$ is injective too. By symmetry, we may assume that 
$(x_{2n})_{n\in\omega}$ is injective and does not take the value $x$. We may also assume that $(x_{2n+1})_{n\in\omega}$ is either injective and does not take the value $x$, or is constant with value $x$. If $(x_{2n+1})_{n\in\omega}$ is injective, then we may assume that $\{ x_{2n}\mid n\!\in\!\omega\}$ and $\{ x_{2n+1}\mid n\!\in\!\omega\}$ are disjoint. In this case, we define $\varphi\! :\!\mathbb{X}_1\!\rightarrow\! X$ by $\varphi (0^\infty )\! :=\! x$, ${\varphi (0^{2n}1^\infty )\! :=\! x_{2n}}$, and $\varphi (0^{2n+1}1^\infty )\! :=\! x_{2n+1}$, so that $\varphi$ is an injective continuous homomorphism from 
$(\mathbb{X}_1,G_{f_1})$ into $(X,R)$. If $(x_{2n+1})_{n\in\omega}$ is constant with value $x$, then we define a function 
${\psi\! :\!\mathbb{X}_1\!\rightarrow\! X}$ by ${\psi (0^\infty )\! :=\! x}$, and $\psi (0^n1^\infty )\! :=\! x_{2n}$, so that $\psi $ is an injective continuous homomorphism from $(\mathbb{X}_1,G_{f_0})$ into $(X,R)$. This shows that 
$\{ (\mathbb{X}_1,G_{f_0}),(\mathbb{X}_1,G_{f_1})\}$ is a $\preceq^i_c$-basis.

\vfill\eject

 If $f\! :\!\mathbb{X}_1\!\rightarrow\!\mathbb{X}_1$ is a witness for 
$(\mathbb{X}_1,G_{f_1})\preceq^i_c(\mathbb{X}_1,G_{f_0})$, then $f$ sends   
$\{ 0^{2n}1^\infty ,0^{2n+1}1^\infty\}$ onto some $\{ 0^{p_n}1^\infty ,0^\infty\}$, which contradicts the injectivity of $f$. If 
$g\! :\!\mathbb{X}_1\!\rightarrow\!\mathbb{X}_1$ is a witness for the inequality  
$(\mathbb{X}_1,G_{f_0})\preceq^i_c(\mathbb{X}_1,G_{f_1})$, then $g$ sends $(0^p1^\infty ,0^\infty )$ to some 
$(0^{2n_p+\varepsilon_p}1^\infty ,0^{2n_p+1-\varepsilon_p}1^\infty )$, and 
$(n_p)_{p\in\omega}$, $(\varepsilon_p)_{p\in\omega}$ have to be constant, which contradicts the injectivity of $g$. This implies that $\{ (\mathbb{X}_1,G_{f_0}),(\mathbb{X}_1,G_{f_1})\}$ is a $\preceq^i_c$-antichain. Note then that 
$(\mathbb{X}_1,G_{f_0}),(\mathbb{X}_1,G_{f_1})$ are in our class.\hfill{$\square$}\medskip

\noindent\emph{Example.}\ Recall $(\mathbb{X}_1,f_1)$ defined before Proposition \ref{LZ1}. Then 
$\mathbb{X}_1\!\not=\!\emptyset$ is a 0DMC space, and $f_1$ is a homeomorphism whose only fixed point is $0^\infty$. By Proposition \ref{fp}, $\chi_c(\mathbb{X}_1,G_{f_1})\! =\! 2^{\aleph_0}$. If we restrict $f_1$ to the open set 
$\mathbb{X}_1\!\setminus\!\{ 0^\infty\}$, then $F_1$ becomes open in the domain of the restriction, and 
$\chi_c(\mathbb{X}_1,G_{{f_1}_{\vert\mathbb{X}_1\setminus\{ 0^\infty\}}})\! =\! 2^{\aleph_0}$. Indeed, if $(C_i)_{i\in\omega}$ be a partition of $\mathbb{X}_1$ into clopen sets, then there is $i$ with $0^\infty\!\in\! C_i$. We may assume that 
$0^{2n}1^\infty ,0^{2n+1}1^\infty \!\in\! C_i$. This implies that $(0^{2n}1^\infty ,0^{2n+1}1^\infty )$ is in 
$G_{{f_1}_{\vert\mathbb{X}_1\setminus\{ 0^\infty\}}}\cap C_i^2$. This shows that we cannot extend Corollary \ref{corfp}(a) when the domain of $f$ is open.\medskip

 We now turn to the study of involutions.
 
\begin{prop} \label{invol} Let $X\!\not=\!\emptyset$ be a 0DMS space, and 
${f\!\! :\! X\!\rightarrow\! X}$ be a fixed point free continuous involution. Then $\chi_c(X,G_f)\! =\! 2$.\end{prop}

\noindent\emph{Proof.}\ Note that $f$ is a continuous bijection with inverse $f$, so that it is a homeomorphism. If $x\!\in\! X$, then 
$f(x)\!\not=\! x$ since $f$ is fixed point free, which gives a clopen neighborhood $N$ of $x$ with $f(x)\!\notin\! N$. As $f$ is a homeomorphism, $C^0\! :=\! N\cap f^{-1}(X\!\setminus\! N)$ and $C^1\! :=\! f[C^0]\! =\! f[N]\!\setminus\! N$ are disjoint clopen subsets of $X$. In particular, $C^0\! =\! f[C^1]$ and $C\! :=\! C^0\cup C^1$ is a $f$-invariant clopen neighborhood of $x$. As $X$ has the Lindel\"of property, we can cover $X$ with countably many such $f$-invariant clopen sets, say $(C_n)_{n\in\omega}$. In particular, 
$\cup_{p<n}~C_n$, $X\!\setminus\! (\cup_{p<n}~C_n)$ and $O_n\! :=\! C_n\!\setminus\! (\cup_{p<n}~C_n)$ are $f$-invariant clopen sets with union $X$. If $U\!\subseteq\! C$ is a $f$-invariant clopen set, then we set $U^\varepsilon\! :=\! U\cap C^\varepsilon$, so that $U$ is the disjoint union of $U^0$ and $U^1$, and $U^{1-\varepsilon}\! =\! f[U^\varepsilon ]$ for each $\varepsilon\!\in\! 2$. We can apply this to $O_n\!\subseteq\! C_n$, so that $X$ is the disjoint union of the family of clopen sets 
$(O_n^\varepsilon )_{n\in\omega ,\varepsilon\in 2}$, and $O_n^{1-\varepsilon}\! =\! f[O_n^\varepsilon ]$ for each 
$(n,\varepsilon )\!\in\!\omega\!\times\! 2$. We then define $c\! :\! X\!\rightarrow\! 2$ by $c(x)\! :=\!\varepsilon$ if 
$x\!\in\! O_n^\varepsilon$ for some $n$, and $c$ is a continuous coloring of $(X,G_f)$.\hfill{$\square$}\medskip

 Proposition \ref{invol} implies some minimality of $G_{f_1}$.

\begin{prop} \label{LZ1+} $(\mathbb{X}_1,G_{f_1})$ is $\preceq^i_c$-minimal, but not $\preceq_c$-minimal, in 
$\mathfrak{K}$.\end{prop}

\noindent\emph{Proof.}\ Let $V$ be compact subset of $\mathbb{X}_1$, and $E\!\subseteq\! G_{f_1}$ be a graph on $V$ with 
${\chi_c(V,E)\!\geq\! 3}$. By Lemma \ref{cpctmin}, it is enough to see that $(\mathbb{X}_1,G_{f_1})\preceq^i_c(V,E)$. As 
$P\! :=\!\mbox{proj}[E]\!\subseteq\! V\cap\mathbb{X}_1\!\setminus\!\{ 0^\infty\}$, we can find, for each $\varepsilon\!\in\! 2$, 
$S_\varepsilon\!\subseteq\!\omega$ with 
$P\! =\!\{ 0^{2n+\varepsilon}1^\infty\mid\varepsilon\!\in\! 2\wedge n\!\in\! S_\varepsilon\}$. As $E$ is a graph, 
$$n\!\in\! S_0\Leftrightarrow 0^{2n}1^\infty\!\in\! P\Leftrightarrow 0^{2n+1}1^\infty\!\in\! P\Leftrightarrow n\!\in\! S_1\mbox{,}$$  
so that $S_0\! =\! S_1$. Thus $E\! =\! G_{{f_1}_{\vert P}}$ is the graph induced by the fixed point free involution 
${f_1}_{\vert P}\! :\! P\!\rightarrow\! P$, which is continuous since $P$ is discrete. As $P\!\subseteq\!\mathbb{X}_1$ is not empty, $\chi_c(P,E)\! =\! 2$ by Proposition \ref{invol}. This gives 
$C\!\in\!\borone (P)$ with $E\cap\big( C^2\cup (P\!\setminus\! C)^2\big)\! =\!\emptyset$. If $0^\infty\!\notin\! V$, then $C$ is a clopen subset of the discrete space $V$, $E\cap\big( C^2\cup (V\!\setminus\! C)^2\big)\! =\!\emptyset$ and 
$\chi_c(V,E)\!\leq\! 2$, which cannot be. Thus $0^\infty\!\in\! V$. Note also that there are infinitely many $p$'s with 
$(0^{2p}1^\infty ,0^{2p+1}1^\infty )\!\in\! E$, otherwise $(0^{2p}1^\infty ,0^{2p+1}1^\infty )\!\notin\! E$ if $p\!\geq\! p_0$. We then set $C'\! :=\! (\bigcup_{p<p_0}~N_{0^{2p}1})\cap V$, so that $C'$ is a clopen subset of $V$, 
$E\cap\big( (C')^2\cup (V\!\setminus\! C')^2)\! =\!\emptyset$, and $\chi_c(V,E)\!\leq\! 2$. This implies that 
$(0^\infty ,0^\infty )\!\in\!\overline{E}^{V^2}$, and gives an injective sequence $(p_n)_{n\in\omega}$ with 
$(0^{2p_n}1^\infty ,0^{2p_n+1}1^\infty )\!\in\! E$. We then define $g\! :\!\mathbb{X}_1\!\rightarrow\! V$ by 
$g(0^\infty )\! :=\! 0^\infty$ and $g(0^{2n+\varepsilon}1^\infty )\! :=\! 0^{2p_n+\varepsilon}1^\infty$, so that $g$ is a witness for 
$(\mathbb{X}_1,G_{f_1})\preceq^i_c(V,E)$, as desired.

\vfill\eject

 We now set $(0^+,1^+,2^+)\! =\! (1,2,0)$, $X\! :=\!\{\varepsilon^\infty\mid\varepsilon\!\in\! 3\}\cup
 \{\varepsilon^{n+1}(\varepsilon^+)^\infty\mid\varepsilon\!\in\! 3\wedge n\!\in\!\omega\}$, and 
$G\! :=\! s(\{ (\varepsilon^{2p+1}(\varepsilon^+)^\infty ,(\varepsilon^+)^{2p+2}\big( (\varepsilon^+)^+\big)^\infty )\mid
\varepsilon\!\in\! 3\wedge p\!\in\!\omega\} )$. Note that $X$ is a 0DMC space, and $G$ is a graph on $X$. As $(\varepsilon^\infty ,(\varepsilon^+)^\infty )\!\in\!\overline{G}$, $(0^\infty ,0^\infty )\!\in\!\Delta (X)\cap\overline{G}^3$, and thus $\chi_c(X,G)\!\geq\! 3$ by Theorem \ref{eq1++}. In particular, $(X,G)\!\in\!\mathfrak{K}$. The map 
$\alpha\!\mapsto\!\alpha (0)$ is a $3$-continuous coloring of $(X,G)$, so that ${\chi_c(X,G)\! =\! 3}$. As 
 $(\mathbb{X}_1,G_{f_1})$ has uncountable CCN, $(\mathbb{X}_1,G_{f_1})\not\preceq_c(X,G)$. However, $(X,G)\preceq_c(\mathbb{X}_1,G_{f_1})$, with witness $\varepsilon^\infty\!\mapsto\! 0^\infty$ and 
$\varepsilon^{n+1}(\varepsilon^+)^\infty\!\mapsto\! 0^{n+1}1^\infty$.\hfill{$\square$}\medskip

 The proof of of Theorem \ref{absmin} implies, as announced in the introduction, the following.

\begin{thm} \label{c2Polish} There is no $\preceq^i_c$-antichain basis for the class of graphs induced by a partial homeomorphism on a 0DMS (or 0DP) space with CCN at least three. In fact, we can even restrict this class to the case where the spaces are countable Polish and the functions are fixed point free involutions with open domain.\end{thm}

\noindent\emph{Proof.}\ By Proposition \ref{D2alph}, $\mathbb{P}_\delta$ is countable Polish. Note that
$$\mbox{proj}[\mathbb{G}_\delta ]\! =\!
\mathbb{P}_\delta\!\setminus\! (\{ c^\infty\}\cup\{ ki0^\infty\mid\delta (k)\! =\! 1\wedge i\!\leq\! 2k\! +\! 1\} )$$
is an open subset of $\mathbb{P}_\delta$. We define 
$f_\delta\! :\!\mbox{proj}[\mathbb{G}_\delta ]\!\rightarrow\!\mbox{proj}[\mathbb{G}_\delta ]$ by 
$f_\delta (\alpha )\! :=\!\mbox{the unique }\beta\!\in\!\mbox{proj}[\mathbb{G}_\delta ]$ with $(\alpha ,\beta )\!\in\!\mathbb{G}_\delta$, so that $f_\delta$ is a fixed point free involution, and $\mathbb{G}_\delta\! =\!\textup{Graph}(f_\delta )\! =\! G_{f_\delta}$. As 
$\mbox{proj}[\mathbb{G}_\delta ]$ is discrete, $f_\delta$ is continuous, and thus a homeomorphism. By Lemma \ref{chromalph}, 
$\chi_c(\mathbb{P}_\delta ,G_{f_\delta})\! =\! 3$ if $\delta$ has infinitely many ones. Let $(X,G_f)$ in our class with 
$(X,G_f)\preceq^i_c(\mathbb{P}_{1^\infty},G_{f_{1^\infty}})$. Theorem \ref{absmin} provides a family 
$(\delta_\gamma )_{\gamma\in 2^\omega}$ in $\mathbb{P}_\infty$ such that 
$(\mathbb{P}_{\delta_\gamma},G_{f_{\delta_\gamma}})\preceq^i_c(X,G_f)$ and the 
$(\mathbb{P}_{\delta_\gamma},G_{f_{\delta_\gamma}})$'s are pairwise $\preceq^i_c$-incompatible in our class. We then apply Lemma \ref{generalab}.\hfill{$\square$}\medskip

 The next results will help us to prove a condition sufficient to get the minimality of some $G_f$'s.
 
\begin{lem} \label{return} Let $X$ be a 0DMC space, $f\! :\! X\!\rightarrow\! X$ be a minimal homeomorphism, and 
$C\!\not=\!\emptyset$ be a clopen subset of $X$. Then there is $L\! >\! 0$ such that, for each $x\!\in\! X$, we can find 
$0\! <\! l\! <\! L$ with $f^l(x)\!\in\! C$, and $0\! <\! l'\! <\! L$ with $f^{-l'}(x)\!\in\! C$.\end{lem}

\noindent\emph{Proof.}\ This is standard. By minimality, 
$X\!\subseteq\!\bigcup_{m\in\mathbb{Z}}~\{ x\!\in\! X\mid f^m(x)\!\in\! C\}$. The compactness of $X$ gives $M\! >\! 0$ with 
$X\!\subseteq\!\bigcup_{-M<m<M}~\{ x\!\in\! X\mid f^m(x)\!\in\! C\}$. We put $L\! :=\! 2M$. If $x\!\in\! X$, then there is 
$-M\! <\! m\! <\! M$ with $f^{m+M}(x)\!\in\! C$, and $0\! <\! l\! :=\! m\! +\! M\! <\! L$. This is similar for $l'$.\hfill{$\square$}\medskip

 Lemma \ref{return} allows us to define $r_C\! :\! X\!\rightarrow\! L$ by $r_C(x)\! :=\!\mbox{min}\{ l\! <\! L\mid f^l(x)\!\in\! C\}$, and $r_C$ is continuous. Similarly, we can define $r'_C\! :\! X\!\rightarrow\! L$ continuous by $r'_C(x)\! :=\!\mbox{min}\{ l'\! <\! L\mid f^{-l'}(x)\!\in\! C\}$.

\begin{lem} \label{gendense} Let $X$ be a 0DMC space, $f\! :\! X\!\rightarrow\! X$ be a minimal  homeomorphism with the property that $\chi_c(X,G_f)\!\geq\! 3$, and $E\!\subseteq\! G_f$ be a graph. The following are equivalent:\smallskip

\noindent (1) $\chi_c(X,E)\!\geq\! 3$,\smallskip

\noindent (2) $E$ is dense in $G_f$.\end{lem}

\noindent\emph{Proof.}\ (1) $\Rightarrow$ (2) Towards a contradiction, suppose that we can find open subsets $U,V$ of $X$ with $G_f\cap (U\!\times\! V)\!\not=\!\emptyset$ and $E\cap (U\!\times\! V)\! =\!\emptyset$. Pick 
$(x,y)\!\in\! G_f\cap (U\!\times\! V)$. If $y\! =\! f(x)$, then $x\!\in\! U\cap f^{-1}(V)$, and there is a clopen neighborhood 
$C\!\subseteq\! U\cap f^{-1}(V)$ of $x$ with $C\cap f[C]\! =\!\emptyset$ since $f$ is fixed point free. By symmetry of $E$, 
$(C\!\times\! f[C])\cap E\! =\! (f[C]\!\times\! C)\cap E\! =\!\emptyset$. We first define 
$s\! :\! X\!\rightarrow\!\omega$ by $s(x)\! :=\! r_C(x)\! +\! r'_{f[C]}(x)$, so that $s$ is continuous.

\vfill\eject

 We then define $c\! :\! X\!\rightarrow\! 2$ by 
$c(x)\! :=\!\mbox{parity}\big( r_C(x)\big)$ if $x\!\notin\! f[C]$, and $c(x)\! :=\!\mbox{parity}\big( s(x)\big)$ if $x\!\in\! f[C]$. Note that $c$ is continuous. It is enough to see that $c$ is a coloring of $(X,E)$. It is enough to see that 
$c\big( f(x)\big)\!\not=\! c(x)$ if $\big( x,f(x)\big)\!\in\! E$. Note that $x\!\notin\! C$. If $x\!\notin\! f[C]$, then $f(x)\!\notin\! f[C]$ and $r_C(x)\! =\! r_C\big( f(x)\big)\! +\! 1$. If $x\!\in\! f[C]$, then $f(x)\!\notin\! f[C]$ and $s(x)\! =\! r_C(x)$, and we conclude similarly.\medskip

\noindent (2) $\Rightarrow$ (1) Towards a contradiction, suppose that there is a clopen subset $C$ of $X$ with the property that 
$E\cap\big( C^2\cup (X\!\setminus\! C)^2\big)\! =\!\emptyset$. As $C$ is clopen and $E$ is dense in $G_f$, 
$$G_f\cap\big( C^2\cup (X\!\setminus\! C)^2\big)\!\subseteq\!\overline{E}\cap\big( C^2\cup (X\!\setminus\! C)^2\big)\! =\!\emptyset\mbox{,}$$ 
which is the desired contradiction.\hfill{$\square$}\medskip

 Lemma \ref{gendense} essentially implies that $G_f$ is minimal if $f$ is.

\begin{lem} \label{mino} Let $X$ be a 0DMC space, $f\! :\! X\!\rightarrow\! X$ be a minimal homeomorphism such that  
$(X,G_f)$ has CCN at least three, $K$ be a 0DMC space, $G$ be a closed graph on $K$ such that $(K,G)\preceq^i_c(X,G_f)$. Then exactly one of the following holds:\smallskip

\noindent (1) $(K,G)$ has CCN at most two,\smallskip

\noindent (2) $(X,G_f)\preceq^i_c(K,G)$.\smallskip

 In other words, $(X,G_f)$ is $\preceq^i_c$-minimal in $\mathfrak{G}_2$ and in the class of closed graphs on a 0DMC space with CCN at least three.\end{lem}

\noindent\emph{Proof.}\ Let $V$ be compact subset of $X$, and $E\!\subseteq\! G_f$ be a compact graph on $V$ with 
$\chi_c(V,E)\!\geq\! 3$ (which implies that $\chi_c(X,E)\!\geq\! 3$). As in the proof of Lemma \ref{cpctmin}, it is enough to see that $(X,G_f)\preceq^i_c(V,E)$ to see that (1) or (2) holds. Let $P\! :=\!\mbox{proj}[E]$, which is compact like $E$. The next two claims give our result.\medskip

\noindent\emph{Claim 1.}\it\ $P\! =\! V\! =\! X$.\rm\medskip

 Indeed, by compactness it is enough to see that $P$ is dense in $X$. Towards a contradiction, suppose that there is a clopen subset $C\!\not=\!\emptyset$ of $X$ disjoint from $P$, so that $r_C$ is defined. We define $c\! :\! P\!\rightarrow\! 2$ by 
$c(x)\! :=\!\mbox{parity}\big( r_C(x)\big)$. Note that $c$ is continuous. It is enough to see that $c$ is a coloring of $(P,E)$, since this implies that $\chi_c(V,E)\!\leq\! 2$ by [E, Theorem 2.1(1)]. It is enough to see that $c\big( f(x)\big)\!\not=\! c(x)$ if $\big( x,f(x)\big)\!\in\! E$ since $f$ is fixed point free and by symmetry. The equality $r_C(x)\! =\! r_C\big( f(x)\big)\! +\! 1$ gives the result.\hfill{$\diamond$}\medskip

\noindent\emph{Claim 2.}\it\ $E\! =\! G_f$.\rm\medskip

 Indeed, Lemma \ref{gendense} implies that $E$ is dense in $G_f$. It remains to note that $E$ is compact.\hfill{$\diamond$}\medskip

 Claims 1 and 2 imply that $(V,E)\! =\! (X,G_f)$.\medskip
 
 For $\mathfrak{G}_2$, assume that $K$ is a 0DMC space, $h$ is a homeomorphism of $K$, $\chi_c(K,G_h)\!\geq\! 3$, and $(K,G_h)\preceq^i_c(X,G_f)$. Then $\chi_c(K,G_h)\!\leq\!\chi_c(X,G_f)\! <\! 2^{\aleph_0}$ by Corollary \ref{corfp}. Corollary \ref{corfp} then implies that $F_1^h$ is an open subset of $K$, $F_1^h\!\not=\! K$ and 
$\chi_c(K,G_h)\! =\!\chi_c(K\!\setminus\! F_1^h,G_h\cap (K\!\setminus\! F_1^h)^2)$. Note that $K\!\setminus\! F_1^h$ is a 0DMC space, $h_{\vert K\setminus F_1^h}$ is a homeomorphism of $K\!\setminus\! F_1^h$, 
$G_{h_{\vert K\setminus F_1^h}}\! =\! G_h\cap (K\!\setminus\! F_1^h)^2$, 
$\chi_c(K\!\setminus\! F_1^h,G_{h_{\vert K\setminus F_1^h}})\!\geq\! 3$ and 
$(K\!\setminus\! F_1^h,G_{h_{\vert K\setminus F_1^h}})\preceq^i_c(X,G_f)$. Moreover, $h_{\vert K\setminus F_1^h}$ is fixed point free, so that $G_{h_{\vert K\setminus F_1^h}}$ is closed. Note that 
$(X,G_f)\preceq^i_c(K\!\setminus\! F_1^h,G_{h_{\vert K\setminus F_1^h}})\preceq^i_c(K,G_h)$, by the first part of the present theorem. In other words, $(X,G_f)$ is $\preceq^i_c$-minimal in $\mathfrak{G}_2$.\hfill{$\square$}\medskip

\noindent\emph{Remark.}\ Proposition \ref{LZ1+} shows that the converse of Lemma \ref{mino} does not hold, since $0^\infty$ is a fixed point of $f_1$.

\vfill\eject\medskip

 The next result is in the style of Corollary \ref{corflip}.
 
\begin{lem} \label{flipo} Let $X,Y$ be 0DMC spaces of cardinality at least three, and 
$f\! :\! X\!\rightarrow\! X$, $g\! :\! Y\!\rightarrow\! Y$ be minimal homeomorphisms. Then $f,g$ are flip-conjugate if and only if 
$(X,G_f)\preceq_c^i(Y,G_g)$, with the same witness. In particular, $(X,G_f)\preceq_c^i(Y,G_g)$ implies that 
$(Y,G_g)\preceq_c^i(X,G_f)$.\end{lem}

\noindent\emph{Proof.}\ Let $\varphi$ be a witness for the fact that $(X,G_f)\preceq_c^i(Y,G_g)$. If $x\!\in\! X$, then 
$\big( x,f(x)\big)\!\in\! G_f$ since $f$ is fixed point free, so that $\Big(\varphi(x),\varphi\big( f(x)\big)\Big)\!\in\! G_g$. So  
$\varphi\big( f(x)\big)\! =\! g\big(\varphi (x)\big)$ or ${\varphi (x)\! =\! g\Big(\varphi\big( f(x)\big)\Big)}$. We apply Lemma \ref{commonlemma} to $f$, $V\! :=\! I\! :=\! X$, $g$, $W\! :=\! Y$, and $\varphi$. In particular, 
$\varphi [\mbox{Orb}_f(x)]\! =\!\mbox{Orb}_g\big(\varphi (x)\big)$. As $g$ is minimal, the compact set 
$\varphi[X]$ is dense in $Y$, showing that $\varphi$ is onto, and thus a homeomorphism by compactness of $X$.\medskip

 Conversely, assume that $f,g$ are flip-conjugate, which gives a homeomorphism $\varphi\! :\! X\!\rightarrow\! Y$ with  
$\varphi\!\circ\! f\! =\! g\!\circ\!\varphi$ or $\varphi\!\circ\! f\! =\! g^{-1}\!\circ\!\varphi$. If $x\!\in\! X$, then 
$g\big(\varphi (x)\big)\! =\!\varphi\big( f(x)\big)$ or $\varphi (x)\! =\! g\Big(\varphi\big( f(x)\big)\Big)$, so that 
$\Big(\varphi (x),\varphi\big( f(x)\big)\Big)\!\in\! G_g$. Thus $\varphi$ is a witness for the fact that 
$(X,G_f)\preceq_c^i(Y,G_g)$.\hfill{$\square$}

\section{$\!\!\!\!\!\!$ Possible chromatic numbers}\indent

 The main goal of this section is to prove Theorem \ref{main}. The next result is essentially [Kra-St, Corollary 2.3].
 
\begin{thm} \label{homeocompact} (Krawczyk-Steprans) Let $X\!\not=\!\emptyset$ be a 0DMC space, and $f\! :\! X\!\rightarrow\! X$ be a fixed point free continuous map. Then $\chi_c(X,G_f)\!\in\!\{ 2,3\}$.\end{thm}

\noindent\emph{Proof.}\ [Kra-St, Corollary 2.3] shows that $\chi_c(X,G_f)\!\leq\! 3$. As $X$ is not empty, 
$\chi_c(X,G_f)\!\not=\! 0$. As $f$ is fixed point free, $\chi_c(X,G_f)\!\not=\! 1$.\hfill{$\square$}

\begin{cor} \label{homeocompactcomplete} Let $X$ be a 0DMC space, and $f\! :\! X\!\rightarrow\! X$ be a continuous injection. Then $\chi_c(X,G_f)$ is in $\{ 0,1,2,3,2^{\aleph_0}\}$, and all these values are possible with homeomorphisms of a countable metrizable compact space.\end{cor}

\noindent\emph{Proof.}\ If $F_1$ is not an open subset of $X$, then $\chi_c(X,G_f)\! =\! 2^{\aleph_0}$ by Proposition \ref{fp}. If $F_1$ is an open subset of $X$, then we may assume that $\chi_c(X,G_f)\! =\!\chi_c(X\!\setminus\! F_1,G_f\cap (X\!\setminus\! F_1)^2)$ by Corollary \ref{corfp}. In other words, we may assume that $X$ is not empty and $f$ is fixed point free. It remains to apply Theorem \ref{homeocompact} for the possible values.\medskip

 If $X\! =\!\emptyset$, then $\chi_c(X,G_f)\! =\! 0$. If $X\! =\! 1$ and $f\! =\!\mbox{Id}$, then $\chi_c(X,G_f)\! =\! 1$. If ${X\! =\! 2}$ and ${f(\varepsilon )\! :=\! 1\! -\!\varepsilon}$, then $c\! :\! 2\!\rightarrow\! 2$ defined by $c(\varepsilon )\! :=\!\varepsilon$ is a continuous coloring of $(X,G_f)$, so that $\chi_c(X,G_f)\! =\! 2$. If $X\! =\! 3$ and 
$f(\varepsilon )\! :=\!\varepsilon^+$ where $(0^+,1^+,2^+)\! =\! (1,2,0)$, then ${c\! :\! 3\!\rightarrow\! 3}$ defined by $c(\varepsilon )\! :=\!\varepsilon$ is a continuous coloring of 
$(X,G_f)$, so that $\chi_c(X,G_f)\! =\! 3$ since $(0,1)$, $(1,2)$ and $(2,0)$ are in $G_f$. We conclude with the example just after Corollary \ref{corfp}.\hfill{$\square$}\medskip

 We will now extend Theorem \ref{homeocompact} to some partial injections. In order to do that, we prove a fixed point free version of the Ryll-Nardzewski theorem (see [Kn-R]). We need to emphasize one point in the [Kn-R] proof, and give the full proof for completeness.   
 
\begin{thm} \label{Ryll} Let $X$ be a Cantor space, $P,Q$ be closed nowhere dense subsets of $X$, and 
$h\! :\! P\!\rightarrow\! Q$ be a fixed point free homeomorphism. Then there is a fixed point free homeomorphism 
$h^*\! :\! X\!\rightarrow\! X$ extending $h$.\end{thm}

\noindent\emph{Proof.}\ We may assume that $X\! =\! 2^\omega$ and, considering $\varepsilon\alpha\!\mapsto\! (1\! -\!\varepsilon )\alpha$, that $P,Q$ are not empty. Let $S\! :=\!\{ s\!\in\! 2^{<\omega}\mid N_s\cap P\! =\!\emptyset\}$, so that 
$2^\omega\!\setminus\! P\! =\!\bigcup_{s\in S}~N_s$ since $P$ is closed. As $P$ is not empty and nowhere dense, $P$ is not clopen and $S$ is infinite. We enumerate $S$ in the increasing order of the lengths of the finite binary sequences, which gives $\{ s_i\mid i\!\in\!\omega\}$. Note that we may assume that the $N_{s_i}$'s are pairwise disjoint, so that the length of $s_i$ goes to infinity.\medskip

\noindent\emph{Claim 1.}\it\ (a) The sequence $\big( d(N_{s_i},P)\big)_{i\in\omega}$ converges to zero.\smallskip

\noindent (b) Let $p\!\in\! P$ and $l\!\in\!\omega$. Then there is one (and thus infinitely many) $i\!\in\!\omega$ with 
$2^{-l}\!\geq\! d(N_{s_i},p)$.\rm\medskip

 Indeed, we argue by contradiction for (a), which gives $l\!\in\!\omega$ such that, for each $j\!\in\!\omega$, there is $i\!\geq\! j$ such that, for each $\alpha\!\in\! 2^\omega$ and each $\gamma\!\in\! P$, $\gamma\vert l\!\not\subseteq\! s_i\alpha$. This provides a strictly increasing sequence $(i_k)_{k\in\omega}$ such that, for each $k\!\in\!\omega$ and each $\gamma\!\in\! P$, 
$\gamma\vert l\!\not\subseteq\! s_{i_k}0^\infty$. Extracting a further subsequence if necessary, we may assume that 
$(s_{i_k}0^\infty )_{k\in\omega}$ converges to some $\delta\!\in\! 2^\omega$, by compactness. As the length of $s_i$ goes to infinity, 
$\delta\!\notin\! P$. This gives $j$ with $\delta\!\in\! N_{s_j}$. Thus $N_{s_j}$ meets infinitely many $N_{s_{i_k}}$'s, which is the desired contradiction.\medskip

 For (b), towards a contradiction, suppose that we can find $p\!\in\! P$ and $l\!\in\!\omega$. As 
$2^\omega\!\setminus\! P\! =\!\bigcup_{i\in\omega}~N_{s_i}$ and $P$ is nowhere dense, we can find $i$ with 
$s_i0^\infty\!\in\! N_{p\vert l}$.\hfill{$\diamond$}\medskip

 Similarly, we can find $\{ t_j\mid j\!\in\!\omega\}\!\subseteq\! 2^{<\omega}$ such that $2^\omega\!\setminus\! Q$ is the disjoint union of the $N_{t_j}$'s, the length of $t_j$ goes to infinity, the sequence $\big( d(N_{t_j},Q)\big)_{j\in\omega}$ converges to zero, and, for each $q\!\in\! Q$ and each $l\!\in\!\omega$, there is $j\!\in\!\omega$ with $2^{-l}\!\geq\! d(N_{t_j},q)$. We fix, for each $i$, $p_i\!\in\! P$ with $d(N_{s_i},p_i)\! =\! d(N_{s_i},P)$, as well as, for each $j$, $q_j\!\in\! Q$ with $d(N_{t_j},q_j)\! =\! d(N_{t_j},Q)$. Note that, by Claim 1,\medskip
 
\noindent (I) there is $f\! :\!\omega\!\rightarrow\!\omega$ injective such that $d(N_{s_i},p_i)\!\geq\! d\big( N_{t_{f(i)}},h(p_i)\big)$,\smallskip

\noindent (II) there is $g\! :\!\omega\!\rightarrow\!\omega$ injective such that $d(N_{t_j},q_j)\!\geq\! d\big( N_{s_{g(j)}},h^{-1}(q_j)\big)$.\medskip

 [Ba, Theorem 1] provides partitions $(I',I''),(J',J'')$ of $\omega$ with $f[I']\! =\! J'$ and $g[J'']\! =\! I''$.\medskip

\noindent\emph{Claim 2.}\it\ Let $s,t\!\in\! 2^{<\omega}$. Then there is a fixed point free homeomorphism $\phi\! :\! N_s\!\rightarrow\! N_t$.\rm\medskip

 Indeed, if $s,t$ are incompatible, we just set $\phi (s\alpha )\! :=\! t\alpha$. If $s\!\subseteq\! t$, then we set 
$$\phi (s\alpha )\! :=\! t\big( 1\! -\!\alpha (\vert t\vert\! -\!\vert s\vert )\big)\alpha (0)\ldots\alpha (\vert t\vert\! -\!\vert s\vert\! -\! 1)
\alpha (\vert t\vert\! -\!\vert s\vert\! +\! 1)\alpha (\vert t\vert\! -\!\vert s\vert\! +\! 2)\ldots$$ 
If $t\!\subsetneqq\! s$, then we set $\phi (s\alpha )\! :=\! t\alpha (1)\ldots\alpha (\vert s\vert\! -\!\vert t\vert )
\big( 1\! -\!\alpha (0)\big)\alpha (\vert s\vert\! -\!\vert t\vert\! +\! 1)\alpha (\vert s\vert\! -\!\vert t\vert\! +\! 2)\ldots$\hfill{$\diamond$}\medskip

 If $i'\!\in\! I'$, Claim 2 provides a fixed point free homeomorphism $\varphi_{i'}\! :\! N_{s_{i'}}\!\rightarrow\! N_{t_{f(i')}}$. The sum of these maps provides a fixed point free homeomorphism 
$\varphi\! :\! U'\! :=\!\bigcup_{i'\in I'}~N_{s_{i'}}\!\rightarrow\! V'\! :=\!\bigcup_{j'\in J'}~N_{t_{j'}}$, since $f$ is injective. Similarly, if $j''\!\in\! J''$, Claim 2 provides a fixed point free homeomorphism $\psi_{j''}\! :\! N_{t_{j''}}\!\rightarrow\! N_{s_{g(j'')}}$. This gives a fixed point free homeomorphism 
$$\psi\! :\! V''\! :=\!\bigcup_{j''\in J''}~N_{t_{j''}}\!\rightarrow\! U''\! :=\!\bigcup_{i''\in I''}~N_{s_{i''}}\mbox{,}$$ 
by injectivity of $g$.

\vfill\eject

 As $2^\omega\!\setminus\! P$ is the sum of $U',U''$ and $2^\omega\!\setminus\! Q$ is the sum of $V',V''$, the function $k\! :\! 2^\omega\!\setminus\! P\!\rightarrow\! 2^\omega\!\setminus\! Q$ defined by 
$$k(\alpha )\! :=\!\left\{\!\!\!\!\!\!\!
\begin{array}{ll}
& \varphi (\alpha )\mbox{ if }\alpha\!\in\! U'\cr
& \psi^{-1}(\alpha )\mbox{ if }\alpha\!\in\! U''
\end{array}
\right.$$ 
is a fixed point free homeomorphism. It remains to prove that the bijection $h^*\! :\! 2^\omega\!\rightarrow\! 2^\omega$ defined by 
$$h^*(\alpha )\! :=\!\left\{\!\!\!\!\!\!\!
\begin{array}{ll}
& h(\alpha )\mbox{ if }\alpha\!\in\! P\cr
& k(\alpha )\mbox{ if }\alpha\!\notin\! P
\end{array}
\right.$$ 
is continuous at each point of $P$. So let $p\!\in\! P$, and 
$(p_n)_{n\in\omega}\!\in\! (2^\omega\!\setminus\! P)^\omega$ converging to $p$. There is $i_n$ with $p_n\!\in\! N_{s_{i_n}}$, and the set $\{ i_n\mid n\!\in\!\omega\}$ is infinite. Note that $d(N_{s_{i_n}},p_{i_n})$ tends to $0$. As the length of $s_i$ goes to infinity, 
$\big( d(p_n,p_{i_n})\big)_{n\in\omega}$ converges to zero, as well as $\big( d(p,p_{i_n})\big)_{n\in\omega}$. Thus 
$(p_{i_n})_{n\in\omega}$ converges to $p$, and $\big( h(p_{i_n})\big)_{n\in\omega}$ converges to $h(p)$. Let $(i'_n)$ be the sequence of $i_n$'s in $I'$. As $d(N_{s_{i'_n}},p_{i'_n})\!\geq\! d\big( N_{t_{f(i'_n)}},h(p_{i'_n})\big)$, 
$\Big(\! d\big( N_{t_{f(i'_n)}},h(p_{i'_n})\big)\!\Big)_{n}$ tends to zero. Call $p'_n$ the fixed element of $N_{s_{i'_n}}$, as in the notation $p_n\!\in\! N_{s_{i_n}}$. Note that the point $\varphi (p'_n)\! =\!\varphi_{i'_n}(p'_n)$ is in $N_{t_{f(i'_n)}}$. As the length of $t_j$ goes to infinity, $\Big( d\big(\varphi (p'_n),h(p_{i'_n})\big)\Big)_{n}$ converges to zero, as well as 
$\Big( d\big(\varphi (p'_n),h(p)\big)\Big)_{n}$. So we proved that $\Big( d\big(\varphi (p'_n),h(p)\big)\Big)_{n}$ converges to zero if $p'_n\!\in\! N_{s_{i'_n}}$ and $\big( d(p'_n,p)\big)_{n\in\omega}$ converges to zero, i.e., 
$\Big( d\big( h^*(p'_n),h^*(p)\big)\Big)_{n}$ converges to zero if $\big( d(p'_n,p)\big)_{n\in\omega}$ does. Similarly, 
$\Big( d\big(\psi (q''_n),h^{-1}(q)\big)\Big)_{n}$ converges to zero if $q\!\in\! Q$, 
${q''_n\!\in\! N_{t_{j''_n}}}$ and $\big( d(q''_n,q)\big)_{n\in\omega}$ converges to zero. Consider now the sequence $(i''_n)$ be the sequence of $i_n$'s in $I''$. Let $j''_n\!\in\! J''$ with $g(j''_n)\! =\! i''_n$. Note that there is $q''_n\!\in\! N_{t_{j''_n}}$ with $p''_n\! =\!\psi (q''_n)\! =\!\psi_{j''_n}(q''_n)$. Let $F\! :=\!\{ q_{j''_n}\mid n\}$. We will check that $\overline{F}\!\setminus\! F\!\subseteq\!\{ h(p)\}$. Let $q\!\in\!\overline{F}\!\setminus\! F$. Assume, for the simplicity of the notation, that $(q_{j''_n})_n$ converges to $q$. As $\big( d(N_{t_j},Q)\big)_{j\in\omega}$ converges to zero, 
$\big( d(N_{t_{j''_n}},q_{j''_n})\big)_n$ converges to zero. As the length of $t_j$ goes to infinity, $\big( d(q''_n,q_{j''_n})\big)_n$ converges to zero, as well as $\big( d(q''_n,q)\big)_n$, $\Big( d\big(\psi (q''_n),h^{-1}(q)\big)\Big)_{n}$ and 
$\Big( d\big( p''_n,h^{-1}(q)\big)\Big)_{n}$. Thus ${p\! =\! h^{-1}(q)}$, as desired. Thus $\big( d(q''_n,h(p)\big)_n$ converges to zero. So $\Big( d\big(\psi^{-1}(p''_n),h(p)\big)\Big)_{n}$ converges to zero if $\big( d(p''_n,p)\big)_n$ does, and 
$\Big( d\big( h^*(p''_n),h^*(p)\big)\Big)_{n}$ converges to zero if $\big( d(p''_n,p)\big)_n$ does.\hfill{$\square$}

\begin{cor} \label{injcontcompact} Let $X$ be a 0DMC space, and $f\! :\! X\!\rightarrow\! X$ be a fixed point free partial continuous injection whose domain is not empty and closed. Then $\chi_c(X,G_f)\!\in\!\{ 2,3\}$.\end{cor}

\noindent\emph{Proof.}\ By [K, 7.8], we may assume that $X\!\subseteq\! 2^\omega$. Note that the map 
$I\! :\! 2^\omega\!\rightarrow\! 2^\omega$ defined by ${I(\alpha )\! :=\! (0,\alpha (0),0,\alpha (1),\ldots )}$ is a homeomorphism onto its nowhere dense range. If moreover we define $g\! :\! I[\mbox{Domain}(f)]\!\rightarrow\! I[\mbox{Range}(f)]$ by 
$g(y)\! :=\! I\Big( f\big( I^{-1}(y)\big)\Big)$, then $g$ is a fixed point free partial continuous injection whose domain is not empty and closed, $(X,G_f)\preceq_c\big( I[X],G_g\big)$ with witness $I$, and thus 
$2\!\leq\!\chi_c(X,G_f)\!\leq\!\chi_c\big( I[X],G_g\big)$. So we may assume that $X$ is closed nowhere dense in $2^\omega$, as well as $\mbox{Domain}(f)$ and $\mbox{Range}(f)$. By compactness, $f$ is a fixed point free homeomorphism from 
$\mbox{Domain}(f)$ onto $\mbox{Range}(f)$. Theorem \ref{Ryll} provides a fixed point free homeomorphism 
$f^*\! :\! 2^\omega\!\rightarrow\! 2^\omega$ extending $f$. By Theorem  \ref{homeocompact}, 
$\chi_c(2^\omega ,G_{f^*})\!\in\!\{ 2,3\}$. Thus $\chi_c(X,G_f)\!\in\!\{ 2,3\}$.\hfill{$\square$}\medskip

\noindent\emph{Remark.}\ The conclusion of this corollary does not hold if the domain of $f$ is open, by the example after Corollary \ref{corfp}.

\vfill\eject

\begin{cor} \label{injcontcompactcomplete} Let $X$ be a 0DMC space, and $f\! :\! X\!\rightarrow\! X$ be a partial continuous injection with closed domain. Then $\chi_c(X,G_f)\!\in\!\{ 0,1,2,3,2^{\aleph_0}\}$, and all these values are possible with  homeomorphisms of a countable metrizable compact space.\end{cor}

\noindent\emph{Proof.}\ If $F_1$ is not an open subset of $\mbox{Domain}(f)$, then $\chi_c(X,G_f)\! =\! 2^{\aleph_0}$ by Proposition \ref{fp}. If $F_1$ is an open subset of $\mbox{Domain}(f)$, then it is a clopen subset of $\mbox{Domain}(f)$. This gives an open subset $O$ of $X$ with $F_1\! =\! O\cap\mbox{Domain}(f)$. Note that $F_1$ and $X\!\setminus\! O$ are disjoint closed subsets of the zero-dimensional metrizable space $X$, which gives a clopen subset $C'$ of $X$ with $F_1\!\subseteq\! C'\!\subseteq\! O$, so that 
$F_1\! =\! C'\cap\mbox{Domain}(f)$. Note then that $f[F_1]\! =\! F_1$ and $f[\mbox{Domain}(f)\!\setminus\! F_1]$ are disjoint compact subsets of $X$, which gives a clopen subset $C''$ of $X$ with 
$F_1\!\subseteq\! C''\!\subseteq\! X\!\setminus\! f[\mbox{Domain}(f)\!\setminus\! F_1]$. We set $C\! :=\! C'\cap C''$, so that $C$ is also a clopen subset of $X$ with $F_1\! =\! C\cap\mbox{Domain}(f)$. Note that $F\! :=\! (X\!\setminus\! C)\cap f^{-1}(X\!\setminus\! C)$ is a clopen subset of $\mbox{Domain}(f)$, and thus a closed subset of $X$. Moreover, 
$G_{f_{\vert F}}\! =\! G_f\cap (X\!\setminus\! C)^2$. Note that $X\!\setminus\! C$ is a 0DMC space, and 
$f_{\vert F}\! :\! X\!\setminus\! C\!\rightarrow\! X\!\setminus\! C$ is a fixed point free partial continuous injection with closed domain. Corollary \ref{injcontcompact} provides a continuous coloring 
$c'\! :\! X\!\setminus\! C\!\rightarrow\! 3$ of $(X\!\setminus\! C,G_f\cap (X\!\setminus\! C)^2)$. We extend $c'$ by $0$ on $C$, which defines $c\! :\! X\!\rightarrow\! 3$ continuous. If $f(x)\!\not=\! x$ is defined, then 
$x\!\in\!\mbox{Domain}(f)\!\setminus\! F_1\!\subseteq\! X\!\setminus\! C$. If $f(x)\!\in\! F_1$, then $f^2(x)\! =\! f(x)$, and 
$f(x)\! =\! x$ by injectivity of $f$, which is absurd, proving that $f(x)\!\notin\! F_1$. If $f(x)\!\in\!\mbox{Domain}(f)$, then 
$f(x)\!\in\! X\!\setminus\! C$ and $c(x)\!\not=\! c\big( f(x)\big)$ since $c'\! :\! X\!\setminus\! C\!\rightarrow\! 3$ is a coloring of 
$(X\!\setminus\! C,G_f\cap (X\!\setminus\! C)^2)$. If $f(x)\!\notin\!\mbox{Domain}(f)$, then either $f(x)\!\in\! X\!\setminus\! C$ and $c(x)\!\not=\! c\big( f(x)\big)$ again, or 
$f(x)\!\in\! C\!\subseteq\! C''\!\subseteq\! X\!\setminus\! f[\mbox{Domain}(f)\!\setminus\! F_1]$, which is absurd. So $c$ is a coloring of $(X,G_f)$.\medskip

 We then apply Corollary \ref{homeocompactcomplete}.\hfill{$\square$}\medskip
 
\noindent\emph{Proof of Theorem \ref{main}.}\ The case (b) of a closed domain comes from Corollary \ref{injcontcompactcomplete}. In fact, if $G$ is an arbitrary graph on $X$, two cases can happen. Either $\Delta (X)$ meets $\overline{G}$ in $(x,x)$, in which case, for any countable partition $(C_i)_{i\in\omega}$ of $X$ into clopen sets, there is $i$ with $x\!\in\! C_i$, and $G$ meets $C_i^2$, so that 
$\chi_c(X,G)\! =\! 2^{\aleph_0}$. Or $\Delta (X)$ does not meet $\overline{G}$, in which case the compactness of $X$ provides a finite continuous coloring of $G$. So $\chi_c(X,G)$ cannot be $\aleph_0$. For the values $0,1,2^{\aleph_0}$ in the open case (a), we use the proof of Corollary \ref{homeocompactcomplete}.\medskip

 So let $1\!\leq\! n\! <\!\omega$. We set 
$K_n\! :=\!\{ p^\infty\mid p\!\leq\! n\}\cup\{ p^{j+1}m^\infty\mid\ p\!\not=\! m\!\leq\! n\wedge j\!\in\!\omega\}$ 
and $\mbox{Domain}(f_n)\! :=\! K_n\!\setminus\!\{ p^\infty\mid p\!\leq\! n\}$, so that $\mbox{Domain}(f_n)$ is an open subset of the countable metrizable compact space $K_n$. We set $(0^+,1^+,\ldots ,n^+)\! :=\! (1,2,\ldots ,n,0)$ and, for $p\!\leq\! n$ and $r\! <\! n$, $p^{+^{r+1}}\! :=\! p^+$ if $r\! =\! 0$, $p^{+^{r+1}}\! :=\! (p^{+^r})^+$ if $r\! >\! 0$. If $\alpha\!\in\!\mbox{Domain}(f_n)$, then we can find 
$p\!\not=\! m\!\leq\! n$, $q\!\in\!\omega$, and $r\! <\! n$ with $\alpha\! =\! p^{nq+r+1}m^\infty$. We then set 
$f_n(\alpha )\! :=\! (p^{+^{r+1}})^{nq+r+1}(m^{+^{r+1}})^\infty$. Note that $f_n$ takes values in $\mbox{Domain}(f_n)$, is continuous, and 
$f_n^{n+1}\! =\!\mbox{Id}$. In particular, $f_n$ is a bijection and $f_n^{-1}\! =\! f_n^n$ is continuous, so that $f_n$ is a homeomorphism. The map $\alpha\!\mapsto\!\alpha (0)$ is a continuous $(n\! +\! 1)$-coloring of $(K_n,G_{f_n})$ (in fact, $f_n$ is fixed point free). If 
$c\! :\! K_n\!\rightarrow\! n$ is continuous, then we can find $p\! <\! l\!\leq\! n$ with $c(p^\infty )\! =\! c(l^\infty )\! =:\! i$. This gives $q\!\in\!\omega$ with $N_{p^q}\cup N_{l^q}\!\subseteq\! c^{-1}(\{ i\} )$. Let $r\! <\! n$ with $p^{+^{r+1}}\! =\! l$.  Then 
$\big( p^{nq+r+1}(p\! +\! 1)^\infty ,(p^{+^{r+1}})^{nq+r+1}\big( (p\! +\! 1)^{+^{r+1}}\big)^\infty\big)\!\in\! 
G_{f_n}\cap\big( c^{-1}(\{ i\} )\big)^2$, so that the function $c$ is not a coloring of $(K_n,G_{f_n})$. Thus 
$\chi_c(K_n,G_{f_n})\! =\! n\! +\! 1$.\medskip

 We then set ${X\! :=\!\oplus_{n\geq 1}~K_n}$, $\mbox{Domain}(f)\! :=\!\oplus_{n\geq 1}~\mbox{Domain}(f_n)$, and 
$f(n,\alpha )\! :=\!\big( n,f_n(\alpha )\big)$. Then $X$ is a countable Polish space $X$, $\mbox{Domain}(f)$ is an open subset of $X$, and  $f$ is a fixed point free partial homeomorphism from $\mbox{Domain}(f)$ onto it. As $(K_n,G_{f_n})\preceq_c(X,G_f)$, 
${\chi_c(X,G_f)\!\geq\!\aleph_0}$. The map $(n,\alpha )\!\mapsto\!\alpha (0)$ is an $\aleph_0$-coloring of $(X,G_f)$, so that 
$\chi_c(X,G_f)\! =\!\aleph_0$.\hfill{$\square$}

\vfill\eject\medskip
 
 In the case of spaces which are not compact, the first space to look at is $\omega$.
 
\begin{cor} \label{injcontomegacomplete} Let $f\! :\!\omega\!\rightarrow\!\omega$ be a partial (continuous) injection. Then 
$\chi_c(\omega ,G_f)\!\in\!\{ 1,2,3\}$, and all these values are possible.\end{cor}

\noindent\emph{Proof.}\ Let $n\!\in\!\omega$. If there is $p\!\in\!\omega$ with $f^p(n)\! =\! n$, we take it minimal, in which case the orbit of $n$ is $\{ f^i(n)\mid 0\!\leq\! i\! <\! p\}$. If $p$ is odd, then we set $c\big( f^i(n)\big)\! :=\! 0$ if $i\! <\! p\! -\! 1$ is even, $c\big( f^i(n)\big)\! :=\! 1$ if $i\! <\! p$ is odd, $c\big( f^{p-1}(n)\big)\! :=\! 2$, so that $c$ is a coloring of $G_f$ on the orbit. If $p$ is even, then we set $c\big( f^i(n)\big)\! :=\! 0$ if $i\! <\! p$ is even, $c\big( f^i(n)\big)\! :=\! 1$ if $i\! <\! p$ is odd, so that $c$ is a coloring of $G_f$ on the orbit. If there is no $p\!\in\!\omega$ with $f^p(n)\! =\! n$, then we can find ordinals 
$\xi ,\eta\!\leq\!\omega$ such that the orbit of $n$ is $\{ f^i(n)\mid -\xi\! <\! i\! <\!\eta\}$. We set $c\big( f^i(n)\big)\! :=\! 0$ if $i$ is even, $c\big( f^i(n)\big)\! :=\! 1$ if $i$ is odd, so that $c$ is a coloring of $G_f$ on the orbit. We defined a continuous coloring $c\! :\!\omega\!\rightarrow\! 3$ of $(\omega ,G_f)$ since $\omega$ is discrete. So $\chi_c(X,G_f)\!\leq\! 3$, and 
$\chi_c(X,G_f)\!\not=\! 0$ since $\omega$ is not empty.\medskip 

 If $f\! =\!\mbox{Id}$, then $\chi_c(X,G_f)\! =\! 1$. If $f(2n\! +\!\varepsilon )\! :=\! 2n\! +\! (1\! -\!\varepsilon )$, then 
$\chi_c(X,G_f)\! =\! 2$ by Proposition \ref{invol}. If $f(3n\! +\!\varepsilon )\! :=\! 3n\! +\!\varepsilon^+$, where 
$(0^+,1^+,2^+)\! =\! (1,2,0)$, then $c\! :\!\omega\!\rightarrow\! 3$ defined by $c(3n\! +\!\varepsilon )\! :=\!\varepsilon$ is a continuous coloring of $(X,G_f)$, so that $\chi_c(X,G_f)\! =\! 3$ since $(0,1),(1,2),(2,0)\!\in\! G_f$.\hfill{$\square$}\medskip

 The next natural space to look at is the Baire space $\omega^\omega$.
 
\section{$\!\!\!\!\!\!$ Graphs induced by a function and odometers}\indent

 We now study graphs induced by a homeomorphism of an uncountable 0DMC space.\medskip

\noindent\emph{Remark.}\ We set, for ${\bf d}\!\in\!\mathfrak{D}$, 
$X_o\! :=\!\overline{\mbox{proj}[\mathbb{G}_o]}^{\mathcal{K}_{\bf d}}$, so that $X_o$ is a 0DMC space. This space is  
${X_o\! =\!\mbox{proj}[\mathbb{G}_o]\cup\{ c^\infty\}\cup\mathcal{C}}$, so that 
$\mbox{proj}[\mathbb{G}_o]\! =\! X_o\!\setminus\! (\{ c^\infty\}\cup\mathcal{C})$ is a countable open subset of $X_o$. We define $f_o\! :\!\mbox{proj}[\mathbb{G}_o]\!\rightarrow\!\mbox{proj}[\mathbb{G}_o]$ by 
$f_o(\alpha )\! :=\!\mbox{the unique }\beta\!\in\!\mbox{proj}[\mathbb{G}_o]$ with $(\alpha ,\beta )\!\in\!\mathbb{G}_o$, so that $f_o$ is a fixed point free involution, and $\mathbb{G}_o\! =\!\textup{Graph}(f_o)\! =\! G_{f_o}$. As $\mbox{proj}[\mathbb{G}_o]$ is discrete, $f_o$ is continuous, and thus a homeomorphism. By Proposition \ref{belowfin}, $\chi_c(X_o,G_{f_o})\! =\! 3$.
\medskip

 The proof of Theorem \ref{eantichmin} provides a $\preceq_c$-antichain made up of $\preceq^i_c$-minimal graphs in the class of graphs induced by a partial fixed point free continuous involution with countable open domain on a 0DMC space with CCN at least three. In particular, any $\preceq^i_c$-basis for this class must have size continuum, as announced in the introduction.
\medskip

 We now turn to the proof of Theorem \ref{eantichmino}, i.e., we study the $G_o$'s instead of the $\mathbb{G}_o$'s.\medskip
 
\noindent\emph{Notation.}\ We set $\mathcal{O}\! :=\!\{ {\bf d}\!\in\!\mathfrak{C}\mid\forall j\!\in\!\omega ~~d_j\mbox{ is odd}\}$.
 
\begin{prop} \label{cno} Let ${\bf d}\! =\! (d_j)_{j\in\omega}\!\in\!\mathfrak{C}$. Then $\chi_c(\mathcal{C},G_o)\! =\! 3$ if 
${\bf d}\!\in\!\mathcal{O}$, $\chi_c(\mathcal{C},G_o)\! =\! 2$ otherwise.\end{prop}

\noindent\emph{Proof.}\ The key remark is that $o^{\pi_{j<l}~d_j}(0^\infty )\!\in\! N_{0^l}$ for each $l\!\in\!\omega$. Assume that there is $j_0\!\in\!\omega$ such that $d_{j_0}$ is even. We define $c\! :\!\mathcal{C}\!\rightarrow\! 2$ by 
$c(\alpha )\! :=\!\mbox{parity}(i)$ if $i\! <\!\pi_{j\leq j_0}~d_j$ and $o^i(0^{j_0+1})\!\subseteq\!\alpha$. Then $c$ is continuous, and a coloring of $(\mathcal{C},G_o)$ by the key remark and the fact that $\pi_{j\leq j_0}~d_j$ is even. Conversely, assume that there is a coloring $c'\! :\!\mathcal{C}\!\rightarrow\! 2$ of $(\mathcal{C},G_o)$. Let $\varepsilon\! :=\! c'(0^\infty )$, and also $C\! :=\! (c')^{-1}(\{\varepsilon\} )$. As $C$ is a clopen subset of $\mathcal{C}$, there is $l_0\!\in\!\omega$ with 
$N_{0^{l_0}}\!\subseteq\! C$. An induction shows that $o^i(0^\infty )\!\in\! C$ if $i$ is even, $o^i(0^\infty )\!\notin\! C$ if $i$ is odd since $c'$ is a coloring of $(\mathcal{C},G_o)$. The key remark shows that $\pi_{j<l}~d_j$ is even if $l\!\geq\! l_0$. This gives $j\!\in\!\omega$ such that $d_j$ is even. This shows that $\chi_c(\mathcal{C},G_o)\!\geq\! 3$ if 
${\bf d}\!\in\!\mathcal{O}$. Let ${\bf d}\!\in\!\mathcal{O}$. We define $c''\! :\!\mathcal{C}\!\rightarrow\! 3$ by 
$c''(\alpha )\! :=\!\mbox{parity}\big(\alpha (0)\big)$ if $\alpha (0)\! <\! d_0\! -\! 1$ and $c''(\alpha )\! :=\! 2$ if 
$\alpha (0)\! =\! d_0\! -\! 1$. Then $c''$ is continuous, and a coloring of $(\mathcal{C},G_o)$.\hfill{$\square$}

\vfill\eject

 Proposition \ref{cno} implies that we can apply Lemma \ref{mino} to $(X,G_f)\! =\! (\mathcal{C},G_o)$ if 
${\bf d}\! =\! (d_j)_{j\in\omega}\!\in\!\mathcal{O}$. It is also important to assume that $(X,G_f)$ has CCN at least three in Lemma \ref{mino}. Indeed, if ${\bf d}\!\in\!\mathfrak{C}\!\setminus\!\mathcal{O}$, then $\mathcal{C}_{\bf d}$ is a 0DMC space, 
$o_{\bf d}$ is a minimal homeomorphism with $\chi_c(\mathcal{C}_{\bf d},G_{o_{\bf d}})\!\geq\! 2$ by Proposition \ref{cno}, and the strict inequality $(2,G_{\varepsilon\mapsto 1-\varepsilon})\prec^i_c(\mathcal{C}_{\bf d},G_{o_{\bf d}})$ holds.\medskip

 A consequence of Lemma \ref{flipo} is the existence, announced in the introduction, of $\preceq_c^i$-basis in the case of equicontinuity. 
 
\begin{defi} We say that a dynamical system $(X,f)$, where a compatible metric $d$ on $X$ is fixed, is \emph{equicontinuous} if 
$\forall x\!\in\! X~~\forall\varepsilon\! >\! 0~~\exists\delta\! >\! 0~~\forall y\!\in\! B_d(x,\delta )~~\forall n\!\in\!\omega ~~
d\big( f^n(x),f^n(y)\big)\! <\!\varepsilon$.\end{defi}
 
 This means that the family $(f^n)_{n\in\omega}$ is equicontinuous. For instance, if 
${\bf d}\!\in\! (\omega\!\setminus\! 2)^\omega$, then $(\mathcal{C},o)$ is equicontinuous (see [Ku, 4.1.2]). We set
$$\mathfrak{G}_2^e\! :=\!\{ (X,G_f)\!\in\!\mathfrak{G}_2\mid (X,f)\mbox{ is equicontinuous}\wedge\exists x\!\in\! X~~
\mbox{Orb}^+_f(x)\mbox{ is dense infinite}\} .$$ 
 
\begin{prop} \label{me} (a) $\{ (\mathcal{C},G_o)\mid {\bf d}\!\in\! (\omega\!\setminus\! 2)^\omega\cap\mathcal{O}\}$ is a 
$\preceq_c^i$-basis for $\mathfrak{G}_2^e$.\smallskip

\noindent (b) Under the axiom of choice, there is a $\preceq_c^i$-antichain basis for $\mathfrak{G}_2^e$.\end{prop}

\noindent\emph{Proof.}\ (a) Let $(X,G_f)\!\in\!\mathfrak{G}_2^e$. [Ku, Theorem 2.9, Corollary 2.34, and Section 4.1 (in particular Theorem 4.4)] provide ${\bf d}\!\in\!\mathfrak{C}$ such that $(\mathcal{C},o)$ is conjugate to $(X,f)$. This gives a homeomorphism $\varphi\! :\!\mathcal{C}\!\rightarrow\! X$ with $\varphi\!\circ\! o\! =\! f\!\circ\!\varphi$. If $\beta\! =\! o(\alpha )$, then $\varphi (\beta )\! =\! f\big(\varphi (\alpha )\big)$, so that $(\mathcal{C},G_o)\preceq^i_c(X,G_f)$. Similarly, 
$(X,G_f)\preceq^i_c(\mathcal{C},G_o)$, so that $\chi_c(\mathcal{C},G_o)\!\geq\!\chi_c(X,G_f)\!\geq\! 3$ since 
$(X,G_f)\!\in\!\mathfrak{G}_2$. Proposition \ref{cno} then implies that ${\bf d}\!\in\!\mathcal{O}$.\medskip
 
\noindent (b) Lemma \ref{flipo} implies that $\preceq^i_c$ is an equivalence relation on 
$\{ (\mathcal{C},G_o)\mid {\bf d}\!\in\! (\omega\!\setminus\! 2)^\omega\cap\mathcal{O}\}$. Using the axiom of choice, we can pick an element in each equivalence class, which provides the desired $\preceq_c^i$-antichain basis.\hfill{$\square$}\medskip

 We now get a $\preceq_c$-antichain in the style of Theorem \ref{anticha}.

\begin{thm} \label{antichaf} There is a map $\Phi\! :\! 2^\omega\!\rightarrow\!\mathcal{O}$ such that 
$(\mathcal{C}_{\Phi (\alpha )},G_{o_{\Phi (\alpha )}})\not\preceq_c(\mathcal{C}_{\Phi (\beta )},G_{o_{\Phi (\beta )}})$ if 
$\alpha\!\not=\!\beta$.\end{thm}
 
\noindent\emph{Proof.}\ Let $(p_n)_{n\in\omega}$ be the sequence of prime numbers. We define, for each 
$\alpha\!\in\! 2^\omega$, $S_\alpha\!\subseteq\!\omega$ by 
$S_\alpha\! :=\!\{ p_0^{\alpha (0)+1}\ldots p_n^{\alpha (n)+1}\mid n\!\in\!\omega\}$. Note that $S_\alpha$ is infinite, and 
$S_\alpha\cap S_\beta$ is finite if $\alpha\!\not=\!\beta$. In this proof, we consider $(d_\alpha )_j\! =\! 3$ if 
$j\!\notin\! S_\alpha$, $(d_\alpha )_j\! =\! p_{j+1}$ if $j\!\in\! S_\alpha$, so that $\Phi (\alpha )\! :=\! {\bf d}_\alpha\!\in\!\mathcal{O}$ is unbounded, the $(d_\alpha )_j$'s are prime, and $(d_\beta )_l$ is not in $\{ (d_\alpha )_j\mid j\!\in\!\omega\}$ if $\alpha\!\not=\!\beta$, $(d_\beta )_l\!\not=\! 3$ and $l$ is large enough.\medskip

 If $(\mathcal{C}_{\Phi (\alpha )},G_{o_{\Phi (\alpha )}})\preceq_c(\mathcal{C}_{\Phi (\beta )},G_{o_{\Phi (\beta )}})$ with witness 
$\varphi$, then we set $V\! :=\!\varphi [\mathcal{C}_{\Phi (\alpha )}]$ and 
$$E\! :=\! (\varphi\!\times\!\varphi )[G_{o_{\Phi (\alpha )}}]\mbox{,}$$ 
so that $V$ is a compact subset of $\mathcal{C}_{\Phi (\beta )}$ and $E\!\subseteq\! G_{o_{\Phi (\beta )}}$ is a compact graph on $V$ with $\chi_c(V,E)\!\geq\! 3$, by Proposition \ref{cno}. Claim 1 in the proof of Lemma \ref{mino} shows that 
$V\! =\! \mathcal{C}_{\Phi (\beta )}$, so that $\varphi$ is onto, which contradicts Lemma \ref{prepanti}.\hfill{$\square$}\medskip

\noindent\emph{Remark.}\ By Lemma \ref{flipo}, the $o_{\Phi (\alpha )}$'s involved in Theorem \ref{antichaf} are pairwise not  flip-conjugate.\medskip

\noindent\emph{Proof of Theorem \ref{eantichmino}.}\ Theorem \ref{antichaf} provides a map 
$\Phi\! :\! 2^\omega\!\rightarrow\!\mathcal{O}$. We now can set $\mathcal{C}_\alpha\! :=\!\mathcal{C}_{\Phi (\alpha )}$ and 
$f_\alpha\! :=\! o_{\Phi (\alpha )}$, and we are done by Proposition \ref{cno} and Lemma \ref{mino}.\hfill{$\square$}\medskip

 We now turn to the version of Theorem \ref{embed}(b). In fact, we prove something stronger since it is possible to consider always the same space, with restrictions of the graph induced by a fixed odometer to different countable dense subsets.\medskip

\noindent\emph{Notation.}\ Fix ${\bf d}\! :=\! 3^\infty$, so that ${\bf d}\!\in\!\mathcal{O}$ and $\chi_c(3^\omega ,G_o)\! =\! 3$, by Proposition \ref{cno}. We set, for $l\!\in\!\omega$, $i_l\! :=\! 3^{l+2}$, so that $0^{l+2}\!\subseteq\! o^{i_l}(0^\infty )$. We set, for 
$A\!\subseteq\!\omega$ infinite, $S_A\! :=\!\{ 0\}\cup\{ i_l\! +\! i\mid l\!\in\! A\wedge i\! <\! 3^l\}$, so that $S_A\!\subseteq\!\omega$ is infinite and contains arbitrarily large intervals appearing in the definition of $\mathbb{G}_o$.

\begin{lem} \label{S+} The map $A\!\mapsto\! S_A$ is injective. Moreover, $A\!\subseteq\! B$ is equivalent to 
$S_A\!\subseteq\! S_B$.\end{lem}

\noindent\emph{Proof.}\ This comes from the fact that $i_l\! +\! 3^l\! -\! 1\! <\! i_{l+1}$ for each $l\!\in\!\omega$.\hfill{$\square$}\medskip

 We also set, for $S\!\subseteq\!\omega$, $D_S\! :=\!\{ o^i(0^\infty )\mid i\!\in\! S\}$. 

\begin{lem} \label{chromab+} The graphs $(3^\omega ,G_{o_{\vert D_{S_A}}}\! )$ are countable and have CCN three, for each $A\!\subseteq\!\omega$ infinite.\end{lem}
 
\noindent\emph{Proof.}\ Fix $A\!\subseteq\!\omega$ infinite. Let us prove that $D_{S_A}$ is dense in $3^\omega$. Let 
$t\!\in\! 3^{<\omega}$, and $i\! <\! 3^{\vert t\vert}$ with $t\!\subseteq\! o^i(0^\infty )$. We choose $l'\! >\!\vert t\vert$ with 
$l'\!\in\! A$. Then $i_{l'}\! +\! i\!\in\! S_A$ and $t\!\subseteq\! o^{i_{l'}+i}(0^\infty )\!\in\! D_{S_A}$. It remains to prove that 
$\chi_c(3^\omega ,G_{o_{\vert D}})\!\geq\! 3$ if $D$ is dense in $3^\omega$. Towards a contradiction, suppose that there is a clopen subset $C$ of $3^\omega$ with $G_{o_{\vert D}}\cap\big( C^2\cup (3^\omega\!\setminus\! C)^2\big)\! =\!\emptyset$. As 
$\chi_c(3^\omega ,G_o)\! =\! 3$, we may assume that $G_o\cap C^2\!\not=\!\emptyset$, which implies that $C\cap o^{-1}(C)$ is not empty. The density of $D$ gives $\alpha\!\in\! D\cap C\cap o^{-1}(C)$. Then 
$\big(\alpha ,o(\alpha )\big)\!\in\! G_{o_{\vert D}}\cap C^2$ since $o$ is fixed point free, which is the desired contradiction.
\hfill{$\square$}\medskip

 One can prove that if $A,B\!\subseteq\!\omega$ are infinite, then $(3^\omega ,G_{o_{\vert D_{S_A}}})\!\preceq^i_c\! (3^\omega ,G_{o_{\vert D_{S_B}}})$ is equivalent to $A\!\subseteq\! B$. But we will prove a better result.
 
\begin{lem} \label{infinf} We can find a sequence $(S_q)_{q\in\omega}$ of pairwise disjoint infinite subsets of $\omega$ such that, for any $l\!\in\!\omega$, $p\!\not=\! q$, $3l\! <\! r\!\in\! S_p$ and $s\!\in\! S_q$, $\vert r\! -\! s\vert\! >\! l$.\end{lem}
 
\noindent\emph{Proof.}\ Fix a bijection $b\! :\!\omega^2\!\rightarrow\!\omega$, for instance 
$b(q,j)\! :=\!\frac{(q+j)(q+j+1)}{2}\! +\! j$. We set $r^q_j\! :=\! 2^{b(q,j)}$ and $S_q\! :=\!\{ r^q_j\mid j\!\in\!\omega\}$. Then 
$(S_q)_{q\in\omega}$ is a sequence of pairwise disjoint infinite subsets of $\omega$. Fix $l\!\in\!\omega$, $p\!\not=\! q$, and assume that $r\! =\! r^p_i\!\in\! S_p$ and $s\! =\! r^q_j\!\in\! S_q$. Note that
$$r\! -\! s\! =\! 2^{b(p,i)}\! -\! 2^{b(q,j)}\! =\!\left\{\!\!\!\!\!\!\!\!
\begin{array}{ll}
& 2^{b(p,i)}(1\! -\! 2^{b(q,j)-b(p,i)})\mbox{ if }b(p,i)\! >\! b(q,j)\mbox{,}\cr\cr
& 2^{b(q,j)}(2^{b(p,i)-b(q,j)}\! -\! 1)\mbox{ if }b(p,i)\! <\! b(q,j).
\end{array}
\right.$$ 
The first term is at least $2^{b(p,i)-1}$, and is bigger than $l$ if $r\! >\! 2l$. If $2^{b(q,j)-1}\! >\! l$, then the second term is smaller than $-l$ and we are done. If $2^{b(q,j)-1}\!\leq\! l$, then the second term is at least $2^{b(p,i)}\! -\! 2l$, and is bigger than $l$ if $r\! >\! 3l$.\hfill{$\square$}\medskip

\noindent\emph{Proof of Theorem \ref{embed}(b).}\ Let $\psi\! :\!\omega\!\rightarrow\! 3^{<\omega}$ be the bijection defined by the length and $o$: $\emptyset$ for the length $0$, $0$, $1$, $2$ for the length $1$, $0^2$, $10$, $20$, $01$, $1^2$, $21$, $02$, $12$, $2^2$ for the length $2$, $\ldots$ Note that $\vert\psi (k)\vert\!\leq\! k$. Recall the sets 
$S_q\! :=\!\{ r^q_j\mid j\!\in\!\omega\}$ given by Lemma \ref{infinf}. If $s\! =\!\psi (k)\!\in\! 3^{<\omega}$, then there is 
$i'_k\! <\! 3^k$ with $s\!\subseteq\! o^{i'_k}(0^\infty )$, and thus $s\!\subseteq\! o^{i_{r^0_k}+i'_k}(0^\infty )$. Recall that 
${D_S\! :=\!\{ o^i(0^\infty )\mid i\!\in\! S\}}$ if $S\!\subseteq\!\omega$, so that $D_S$ is countable, and $G_{o_{\vert D_S}}$ is a countable graph on $3^\omega$. We proved that $D_S$ is dense in $3^\omega$ if $S$ contains 
$\{ i_{r^0_k}\! +\! i'_k\mid k\!\in\!\omega\}$. In this case, $G_{o_{\vert D_S}}$ is dense in $G_o$, and thus 
$\chi_c(3^\omega ,G_{o_{\vert D_S}})\! =\! 3$ by Proposition \ref{cno} and Lemma \ref{gendense}. We set, for 
$A\!\subseteq\!\omega$, 
$S_A\! :=\!\{ i_{r^0_k}\! +\! i'_k\mid k\!\in\!\omega\}\cup\{ i_{r^{n+1}_k}\mid n\!\in\! A\wedge k\!\in\!\omega\}$, so that 
$\chi_c(3^\omega ,G_{o_{\vert D_{S_A}}})\! =\! 3$ and 
$(3^\omega ,G_{o_{\vert D_{S_A}}})\preceq^i_c(3^\omega ,G_{o_{\vert D_{S_B}}})$ if $A\!\subseteq\! B$.

\vfill\eject

 Assume now that $n\!\in\! A\!\setminus\! B$ and, towards a contradiction, that 
$(3^\omega ,G_{o_{\vert D_{S_A}}})\preceq_c(3^\omega ,G_{o_{\vert D_{S_B}}})$ with witness $\varphi$. If 
$\alpha\!\in\! D_{S_A}$, then $\big(\alpha ,o(\alpha )\big)\!\in\! G_{o_{\vert D_{S_A}}}$, so that 
$\Big(\varphi (\alpha ),\varphi\big( o(\alpha )\big)\Big)\!\in\! G_{o_{\vert D_{S_B}}}$. Thus 
${\varphi\big( o(\alpha )\big)\! =\! o^{\pm 1}\big(\varphi (\alpha )\big)}$. As the set 
$\{\alpha\!\in\! 3^\omega\mid\varphi\big( o(\alpha )\big)\! =\! o^{\pm 1}\big(\varphi (\alpha )\big)\}$ is closed and contains the dense set 
$D_{S_A}$, $\varphi\big( o(\alpha )\big)\! =\! o^{\pm 1}\big(\varphi (\alpha )\big)$ holds for each $\alpha\!\in\! 3^\omega$. In particular, 
$\varphi [\{ o^i(0^\infty )\mid i\!\in\!\omega\} ]\!\subseteq\!\mbox{Orb}_o(0^\infty )$.\medskip

 As $\varphi$ is uniformly continuous, there is $l\!\in\!\omega$ such that $\varphi (\alpha )\vert 1\! =\!\varphi (\beta )\vert 1$ if 
$\alpha\vert l\! =\!\beta\vert l$. As we work with the odometer on $3^\omega$, 
$o\big(\varphi (\alpha )\big)\vert 1\!\not=\! o^{-1}\big(\varphi (\alpha )\big)\vert 1$, and thus 
$\varphi\big( o(\alpha )\big)\! =\! o\big(\varphi (\alpha )\big)$ is equivalent to 
$\varphi\big( o(\alpha )\big)\vert 1\! =\! o\big(\varphi (\alpha )\big)\vert 1$. The previous discussion allows us to define, for each natural number $r\!\leq\! 3^l$, $f(r)\!\in\!\mathbb{Z}$ with $\varphi\big( o^r(0^\infty )\big)\! =\! o^{f(r)}(0^\infty )$. Note that 
$f(r\! +\! 1)\! =\! f(r)\!\pm\! 1$ if $r\! <\! 3^l$, so that 
$$f(0)\! -\! 3^l\! <\! f(r)\! <\! f(0)\! +\! 3^l .$$ 
We set 
$d\! :=\!\mbox{Card}(\{ r\! <\! 3^l\mid f(r\! +\! 1)\! =\! f(r)\! +\! 1\} )\! -\!\mbox{Card}(\{ r\! <\! 3^l\mid f(r\! +\! 1)\! =\! f(r)\! -\! 1\} )$,  
so that $d\!\in\!\mathbb{Z}\!\setminus\!\{ 0\}$ and $-3^l\!\leq\! d\!\leq\! 3^l$. Note that any natural number $i$ has a unique decomposition $3^lq\! +\! r$, where $q\!\in\!\omega$ and $r\! <\! 3^l$. The previous discussion shows that 
$\varphi\big( o^i(0^\infty )\big)\! =\! o^{dq+f(r)}(0^\infty )$. We apply this to ${i_{r^{n+1}_k}\!\in\! S_A}$, where $k$ is large enough to ensure that $r^{n+1}_k\! +\! 2\!\geq\! l$, so that 
${\varphi\big( o^{i_{r^{n+1}_k}}(0^\infty )\big)\! =\! o^{d3^{r^{n+1}_k+2-l}+f(0)}(0^\infty )}$. The previous discussion shows the existence of $i\!\in\! S_B$, say $i_{r^{n_k}_{j_k}}\! +\! i'_{j_k}$, and $\varepsilon\!\in\! 2$ with the property that 
$d3^{r_k^{n+1}+2-l}\! +\! f(0)\! =\! i\! +\!\varepsilon$. In particular, taking $k$ large enough, we see that $d\!\geq\! 1$. Moreover, 
$$\frac{d}{3^l}\! +\!\frac{f(0)}{3^{r_k^{n+1}+2}}\! =\! 
3^{r^{n_k}_{j_k}-r_k^{n+1}}\! +\!\frac{i'_{j_k}+\varepsilon}{3^{r_k^{n+1}+2}}\mbox{,}$$ 
showing that $r^{n_k}_{j_k}\!\leq\! r_k^{n+1}$ if $k$ is large enough. Similarly, 
$d3^{r_k^{n+1}-r^{n_k}_{j_k}-l}\! +\!\frac{f(0)}{3^{r^{n_k}_{j_k}+2}}\! =\! 
1\! +\!\frac{i'_{j_k}+\varepsilon}{3^{r^{n_k}_{j_k}+2}}$, showing that $r_k^{n+1}\!\leq\! r^{n_k}_{j_k}\! +\! l$ if $k$ is large enough. This shows that $0\!\leq\! r_k^{n+1}\! -\! r^{n_k}_{j_k}\!\leq\! l$. As $i\!\in\! S_B$, $n_k\! =\! 0$ or $n_k\! -\! 1\!\in\! B$, showing that $n\! +\! 1\!\not=\! n_k$ since $n\!\notin\! B$. It remains to apply Lemma \ref{infinf} to $p\! :=\! n\! +\! 1$, $q\! :=\! n_k$, 
$r\! :=\! r_k^{n+1}$ with $k$ large enough so that $r_k^{n+1}\! >\! 3l$, and $s\! :=\! r^{n_k}_{j_k}$ to get the desired contradiction. So we proved that $(3^\omega ,G_{o_{\vert D_{S_A}}})\not\preceq_c(3^\omega ,G_{o_{\vert D_{S_B}}})$ if 
$A\!\not\subseteq\! B$.\medskip

 It remains to check that the map $A\!\mapsto\! G_{o_{\vert D_{S_A}}}$ is injective. First, the map $A\!\mapsto\! S_A$ is injective by Lemma \ref{S+}. Then the map $S\!\mapsto\! D_S$ from 
$\mathcal{P}(\omega)\!\equiv\! 2^\omega$ into $2^{\text{Orb}_o(0^\infty )}$ is injective. Finally, the map 
$D\!\mapsto\! G_{o_{\vert D}}$ from $2^{\text{Orb}_o(0^\infty )}$ into $2^{\text{Orb}_o(0^\infty )^2}$ is injective by minimality of $o$.\hfill{$\square$}\medskip

\noindent\emph{Proof of Theorem \ref{absmincomp}.}\ Let $o$ be the odometer on $3^\omega$, and recall $S_\omega$ defined before Lemma \ref{S+}. We set $\mathbb{G}\! :=\! G_{o_{\vert D_{S_\omega}}}$. By Lemma \ref{chromab+}, the graph 
$(3^\omega ,\mathbb{G})$ is countable and has CCN three, and is therefore in $\mathfrak{K}$. Let $(K,G)$ in 
$\mathfrak{K}$ satisfying $(K,G)\preceq^i_c(3^\omega ,\mathbb{G})$, with witness $\varphi$. We set $V\! :=\!\varphi [K]$ and 
$E\! :=\! (\varphi\!\times\!\varphi )[G]$, so that $V$ is a compact subset of $3^\omega$ and $E\!\subseteq\! \mathbb{G}$ is a graph on 
$V$. Also $(K,G)\preceq^i_c(V,E)$ with witness $\varphi$, so that $\chi_c(V,E)\! =\! 3$. Claim 1 in the proof of Lemma \ref{mino} shows that $V\! =\! 3^\omega$. Note that $(3^\omega ,E)\preceq^i_c(K,G)$ with witness $\varphi^{-1}$. So it is enough to find a $\preceq^i_c$-antichain $\big( (3^\omega ,G_\alpha )\big)_{\alpha\in 2^\omega}$ of graphs with CCN three and 
$\preceq^i_c$-below $(3^\omega ,E)$, by Lemma \ref{generalab}. We first inductively construct a sequence $(\alpha_n)_{n\in\omega}$ of points of $D_{S_\omega}$ satisfying the following:
$$\begin{array}{ll}
& (1)~\forall\varepsilon\!\in\! 2\mbox{, }
E_\varepsilon\! :=\!\bigsqcup_{p\in\omega}~\big\{\big(\alpha_{2p+\varepsilon},o(\alpha_{2p+\varepsilon})\big),
\big( o(\alpha_{2p+\varepsilon}),\alpha_{2p+\varepsilon}\big)\big\}\mbox{ is dense in }E\mbox{,}\cr
& (2)~E_0\cap E_1\! =\!\emptyset\mbox{,}\cr
& (3)~\alpha_n\! =\! o^{i_{l_n}+i'_n}(0^\infty )\mbox{,}\cr
& (4)~(l_n)_{n\in\omega}\mbox{ is injective.}
\end{array}$$

\vfill\eject

 Lemma \ref{gendense} implies that $E$ is dense in $G_o$. As $o$ is a homeomorphism of the perfect ($=$ without isolated point) space $3^\omega\!\not=\!\emptyset$, $G_o,E\!\not=\!\emptyset$ are also perfect. Let $(B_q)_{q\in\omega}$ be a basis for the topology of $E$ made up of sets which are not empty. We first choose $(\beta_0,\gamma_0)\!\in\! B_0$. As 
$E\!\subseteq\! G_o$, either $\gamma_0\! =\! o(\beta_0)$, or $\beta_0\! =\! o(\gamma_0)$. The point $\alpha_0$ is $\beta_0$ in the first case, $\gamma_0$ in the second one, so that 
$\big\{\big(\alpha_0,o(\alpha_0)\big),\big( o(\alpha_0),\alpha_0\big)\big\}\!\subseteq\! E$ since $E$ is symmetric, and 
$\alpha_0\!\in\! D_{S_\omega}$ since $o$ is minimal. As $B_0\!\not=\!\emptyset$ is perfect, we can then choose 
$(\beta_1,\gamma_1)\!\in\! B_0\!\setminus\!\big\{\big(\alpha_0,o(\alpha_0)\big),\big( o(\alpha_0),\alpha_0\big)\big\}$. Here again, either 
$\gamma_1\! =\! o(\beta_1)$, or $\beta_1\! =\! o(\gamma_1)$. The point $\alpha_1$ is $\beta_1$ in the first case, $\gamma_1$ in the second one, so that 
$$\big\{\big(\alpha_1,o(\alpha_1)\big),\big( o(\alpha_1),\alpha_1\big)\big\}\!\subseteq\! E\!\setminus\!
\big\{\big(\alpha_0,o(\alpha_0)\big),\big( o(\alpha_0),\alpha_0\big)\big\}$$ 
and $\alpha_1\!\in\! D_{S_\omega}$. It remains to iterate this construction in the other $B_q$'s. At this point, we only ensured (1) and (2). As $i\! <\! 3^l$ in the definition of $S_\omega$, we can also ensure (3) and (4).\medskip

 We ensured that the graph $E_0$ is dense in $G_o$, as well as any $E_0\cup D$ if $D\!\subseteq\! E_1$ is a graph. Lemma \ref{gendense} implies that $\chi_c(3^\omega ,E_0\cup D)\! =\! 3$. Let $(p_n)_{n\in\omega}$ be the sequence of prime numbers. We define, for each $\alpha\!\in\! 2^\omega$, 
$S_\alpha\!\subseteq\!\omega$ by $S_\alpha\! :=\!\{ p_0^{\alpha (0)+1}\ldots p_n^{\alpha (n)+1}\mid n\!\in\!\omega\}$.  
Note that $S_\alpha$ is infinite, and $S_\alpha\cap S_\beta$ is finite if $\alpha\!\not=\!\beta$. We then set 
$D_\alpha\! :=\!
\bigcup_{p\in S_\alpha}~\big\{\big(\alpha_{2p+1},o(\alpha_{2p+1})\big),\big( o(\alpha_{2p+1}),\alpha_{2p+1}\big)\big\}$, so that $D_\alpha\!\subseteq\! E_1$ is a graph. We put $G_\alpha\! :=\! E_0\cup D_\alpha$, so that $G_\alpha$ is a graph with CCN three and $\preceq^i_c$-below $(3^\omega ,E)$.\medskip

 It remains to see that $\big( (3^\omega ,G_\alpha )\big)_{\alpha\in 2^\omega}$ is a $\preceq^i_c$-antichain. So let 
$\alpha\!\not=\!\beta\!\in\! 2^\omega$, and assume, towards a contradiction, that 
$(3^\omega ,G_\alpha )\preceq^i_c(3^\omega ,G_\beta )$ with witness $\psi$. We set
$$D_0\! :=\!\{\alpha_{2k}\mid k\!\in\!\omega\}\cup\{\alpha_{2k+1}\mid k\!\in\! S_\alpha\} .$$ 
\emph{Claim.}\it\ $D_0$ is dense in $3^\omega$.\rm\medskip

 Indeed, let $\emptyset\!\not=\! s\!\in\! 3^{<\omega}$. Note that $G_o$ meets $N_s\!\times\! N_{o(s)}$. As $E_0$ is dense in $G_o$,
$$G_\alpha\! =\! G_{o_{\vert D_0}}\! =\!\textup{Graph}(o_{\vert D_0})\cup\textup{Graph}(o_{\vert D_0})^{-1}$$
also meets this clopen set. As $o$ is the odometer on $3^\omega$, $\textup{Graph}(o_{\vert D_0})$ meets this clopen set, and $D_0$ meets $N_s$.\hfill{$\diamond$}\medskip

 If $x\!\in\! D_0$, then $\big( x,o(x)\big)\!\in\! G_\alpha$, so that 
$\Big(\psi (x),\psi\big( o(x)\big)\Big)\!\in\! G_\beta$, and either $\psi\big (o(x)\big)\! =\! o\big(\psi (x)\big)$, or 
$\psi\big (o(x)\big)\! =\! o^{-1}\big(\psi (x)\big)$. This leads to define 
$P\! :=\!\{ x\!\in\! 3^\omega\mid\psi\big (o(x)\big)\! =\! o\big(\psi (x)\big)\}$ and
$$M\! :=\!\{ x\!\in\! 3^\omega\mid\psi\big (o(x)\big)\! =\! o^{-1}\big(\psi (x)\big)\} .$$
Note that $P,M$ are closed, and disjoint by minimality of $o$. Moreover, $D_0\!\subseteq\! P\cup M$, so that 
$3^\omega\! =\! P\cup M$ by the claim. If $y\!\in\! P$, then $o(y)\!\in\! P$ by injectivity of $\psi$ and minimality of $o$, and similarly with $M$. In particular, the dense set $\{ o^i(0^\infty )\mid i\!\in\!\omega\}$ is contained in either $P$, or $M$. Thus $P\! =\! 3^\omega$ or 
$M\! =\! 3^\omega$, which means that $\psi\!\circ\! o\! =\! o\!\circ\!\psi$ or $\psi\!\circ\! o\! =\! o^{-1}\!\circ\!\psi$.\medskip

 Fix $p\!\in\! S_\alpha\!\setminus\! S_\beta$. If $\psi\!\circ\! o\! =\! o^{-1}\!\circ\!\psi$, then there is $n\!\in\!\omega$ with 
$\psi (\alpha_{2p+1})\! =\! o(\alpha_n)$ by minimality of $o$. In particular, 
$\psi\big( o^{i_{l_{2p+1}}+i'_{2p+1}}(0^\infty )\big)\ =\! o^{i_{l_n}+i'_n+1}(0^\infty )$. If $i\!\in\!\omega$, then
$$\psi\big( o^{i_{l_{2p+1}}+i'_{2p+1}+i}(0^\infty )\big)\ =\! o^{i_{l_n}+i'_n+1-i}(0^\infty ).$$
If we choose $i$ large enough so that 
$o^{i_{l_{2p+1}}+i'_{2p+1}+i}(0^\infty )\!\in\!\{\alpha_{2k}\mid k\!\in\!\omega\}\cup\{\alpha_{2k+1}\mid k\!\in\! S_\alpha\}$ and 
$i_{l_n}\! +\! i'_n\! +\! 1\! -\! i\! <\! 0$, then we get a contradiction with the minimality of $o$. This shows that 
$\psi\!\circ\! o\! =\! o\!\circ\!\psi$. The minimality of $o$ provides $n\!\in\!\omega$ with 
$\psi (\alpha_{2p+1})\! =\!\alpha_n$, and $n\!\not=\! 2p\! +\! 1$ since $p\!\notin\! S_\beta$.

\vfill\eject

 In particular, $\psi\big( o^{i_{l_{2p+1}}+i'_{2p+1}}(0^\infty )\big)\ =\! o^{i_{l_n}+i'_n}(0^\infty )$. If $i\!\in\!\omega$, then 
$$\psi\big( o^{i_{l_{2p+1}}+i'_{2p+1}+i}(0^\infty )\big)\ =\! o^{i_{l_n}+i'_n+i}(0^\infty ).$$ 
Applying this to $i\! =\! i_{l_{2q+1}}\! +\! i'_{2q+1}\! -\! i_{l_{2p+1}}\! -\! i'_{2p+1}\!\geq\! 0$, we get 
$$\psi\big( o^{i_{l_{2q+1}}+i'_{2q+1}}(0^\infty )\big)\ =\! o^{i_{l_n}+i'_n+i_{l_{2q+1}}+i'_{2q+1}-i_{l_{2p+1}}-i'_{2p+1}}(0^\infty )\mbox{,}$$ 
which has to be of the form $o^{i_{l_m}+i'_m}(0^\infty )$ if $q\!\in\! S_\alpha$. In particular, 
$$3^{l_{2q+1}+2}\! +\! i'_{2q+1}\! +\! i_{l_n}\! +\! i'_n\! -\! i_{l_{2p+1}}\! -\! i'_{2p+1}\! =\! 3^{l_m+2}\! +\! i'_m$$ 
and 
$3^{l_{2q+1}-l_m}\! +\!\frac{i'_{2q+1}+i_{l_n}+i'_n-i_{l_{2p+1}}-i'_{2p+1}}{3^{l_m+2}}\! =\! 1\! +\!\frac{i'_m}{3^{l_m+2}}\! <\! 2$, showing that $l_{2q+1}\!\leq\! l_m$ if $q$ is large enough. Also, 
$1\! +\!\frac{i'_{2q+1}+i_{l_n}+i'_n-i_{l_{2p+1}}-i'_{2p+1}}{3^{l_{2q+1}+2}}\! =\! 
3^{l_m-l_{2q+1}}\! +\!\frac{i'_m}{3^{l_{2q+1}+2}}\! <\! 2$, showing that $l_m\!\leq\! l_{2q+1}$ if $q$ is large enough. Thus 
$m\! =\! 2q\! +\! 1$ if $q$ is large enough, by (4), and $i'_m\! =\! i'_{2q+1}$ if $q$ is large enough. This implies that 
$n\! =\! 2p\! +\! 1$, which is the desired contradiction finishing the proof.\hfill{$\square$}\medskip

\noindent\emph{Remark.}\ In fact, $\mathbb{G}$ and the $G_\alpha$'s are induced by a partial homeomorphism with countable domain, so that there is no $\preceq^i_c$-antichain basis for the class of graphs on a 0DMC space induced by a partial homeomorphism with countable domain with CCN at least three.
   
\section{$\!\!\!\!\!\!$ Subshifts}\indent

 We now prove a version of Theorem \ref{infdecrcompact} for graphs induced by a homeomorphism, as announced in the introduction. The proof of Theorem \ref{Ckappa} will provide descending chains of graphs of uncountable CCN, and here we get CCN three. We consider subshifts, which are widely studied particular dynamical systems. We refer to the book [Ku] for basic notions and definitions.

\begin{defi} (a) An \emph{alphabet} is a finite set of cardinality at least two.\smallskip

\noindent (b) Let $A$ be an alphabet, and $\mathbb{X}\!\in\!\{\mathbb{Z},\omega\}$. The \emph{shift map} 
$\sigma\! :\! A^\mathbb{X}\!\rightarrow\! A^\mathbb{X}$ is defined by the formula $\sigma (\alpha )(k)\! :=\!\alpha (k\! +\! 1)$.\end{defi}

 Recall that the sets of the form $[w]_p\! :=\!\{\beta\!\in\! 2^\mathbb{Z}\mid\forall j\! <\!\vert w\vert ~~w(j)\! =\!\beta (p\! +\! j)\}$, where $w\!\in\! 2^{<\omega}$ and $p\!\in\!\mathbb{Z}$, form a basis made up of clopen subsets of the space $2^\mathbb{Z}$, which is therefore homeomorphic to $2^\omega$. If $\mathbb{X}\! =\!\mathbb{Z}$, then the shift map is a homeomorphism, so that $(A^\mathbb{Z},\sigma)$ is a Cantor dynamical system. Corollary \ref{corfp} shows that the fixed points of a homeomorphism $f$ are important in the computation of the CCN of $G_f$.\medskip
  
\noindent\emph{Notation.}\ If $A$ is an alphabet and $\emptyset\!\not=\! w\!\in\! A^{<\omega}$, then 
$w^\mathbb{Z}\!\in\! A^\mathbb{Z}$ is defined by $(w^\mathbb{Z})(k\vert w\vert\! +\! j)\! =\! w(j)$ if $k\!\in\!\mathbb{Z}$ and 
$j\! <\!\vert w\vert$. 

\begin{prop} \label{period} Let $A$ be an alphabet, and $i\! >\! 0$ be a natural number. Then $\sigma^i(\alpha )\! =\!\alpha$ holds exactly when there is $w\!\in\! A^i$ with $\alpha\! =\! w^\mathbb{Z}$ (in this case, we say that $\alpha$ is a \emph{periodic point} of $\sigma$). In particular, the fixed points of $\sigma$ are exactly the constant sequences.\end{prop}

\noindent\emph{Proof.}\ If $k\!\in\!\mathbb{Z}$, then $\sigma^i(\alpha)(k)\! =\!\alpha (k\! +\! i)$. For the left to right implication, we consider $w\! :=\!\alpha\vert i$ defined by $w(j)\! :=\!\alpha (j)$ if $j\! <\! i$.\hfill{$\square$}

\begin{defi} Let $A$ be an alphabet. A \emph{two-sided subshift} is a closed subset $\Sigma$ of $A^\mathbb{Z}$ with the property that $\sigma [\Sigma]\! =\!\Sigma$.\end{defi}

 Note that a two-sided subshift defines a dynamical system, by restriction.\medskip
  
\noindent\emph{Notation.}\ If $A$ is an alphabet, $\alpha\!\in\! A^{\leq\omega}\cup A^\mathbb{Z}$ and 
$w\!\in\! A^{\leq\omega}$, then we write $w\!\sqsubseteq\!\alpha$ when $w$ appears in $\alpha$, i.e., when there is 
$k\!\in\!\mathbb{Z}$ such that $w(j)\! =\!\alpha (k\! +\! j)$ for each $j\! <\!\vert w\vert$. In particular, if 
$\alpha\!\in\! A^\mathbb{Z}$, then $\alpha^+\! :=\!\big(\alpha (0),\alpha (1),\ldots\big)\!\sqsubseteq\!\alpha$.\medskip
 
\noindent\emph{Example.}\ If $A$ is an alphabet and $F\!\subseteq\! A^{<\omega}$, then 
$\Sigma_F\! :=\!\{\alpha\!\in\! A^\mathbb{Z}\mid\forall w\!\sqsubseteq\!\alpha~~w\!\notin\! F\}$ is the set of biinfinite words without subword in $F$. This is a two-sided subshift, and any two-sided subshift is of this form (see [Sa-T\"o, Section 2]).\medskip

 The next notion will be crucial in our study of subshifts.

\begin{defi} Let $A$ be an alphabet. A \emph{substitution} on $A$ is a map $\tau\! :\! A^{<\omega}\!\rightarrow\! A^{<\omega}$ satisfying $\tau (uv)\! =\!\tau (u)\tau (v)$ for all $u,v\!\in\! A^{<\omega}$.\end{defi} 

 A substitution is determined by the images of the letters of the alphabet. Some authors require that 
$\tau^{-1}(\{\emptyset\} )\! =\!\{\emptyset\}$, which will be the case in our examples. We now provide infinite descending chains of graphs induced by a homeomorphism of a 0DMC space with CCN exactly three.

\begin{thm} \label{decfib} There is a $\preceq_c$ and $\preceq^i_c$-descending chain  
$\big( (\Sigma_p,G_{\sigma_{\vert\Sigma_p}})\big)_{p\in\omega}$, where $\Sigma_p$ is a two-sided subshift, 
$(\sigma_{\vert\Sigma_p})^2$ is fixed point free, and $(\Sigma_p,G_{\sigma_{\vert\Sigma_p}})$ has CCN three.\end{thm}
  
\noindent\emph{Proof.}\ We consider the generalized Fibonacci sequence of natural numbers defined by $f_0\! :=\! 2$, 
$f_1\! :=\! 3$, and $f_{p+2}\! :=\! f_p\! +\! f_{p+1}$. Note that $f_p\! >\! 0$, $(f_p)_{p\in\omega}$ is strictly increasing, and $f_p$ is even exactly when 3 divides $p$, by induction. Also, $f_{p+5}\! >\! 8f_p$ since 
$$f_{p+5}\! =\! f_{p+3}\! +\! f_{p+4}\! =\! f_{p+1}\! +\! 2f_{p+2}\! +\! f_{p+3}\! =\! 
f_{p+1}\! +\! 2f_p\! +\! 2f_{p+1}\! +\! f_{p+1}\! +\! f_{p+2}\! =\! 5f_{p+1}\! +\! 3f_p.$$
In particular, $8f_{9p+5}\! <\! f_{9p+14}$.\medskip

 This leads to define a $\subseteq$-increasing sequence $(F_p)_{p\in\omega}$ of subsets of $2^{<\omega}$ by 
$$F_p\! :=\!\{ 0^2,1^3\}\cup\{ w^8\mid w\!\in\! 2^{<\omega}\wedge 0\! <\! 8\vert w\vert\! <\! f_{9p+5}\} .$$ 
This allows us to define the two-sided subshifts $\Sigma_p\! :=\!\Sigma_{F_p}$. Note that $(\Sigma_p)_{p\in\omega}$ is 
$\subseteq$-decreasing, so that $\big( (\Sigma_p,G_{\sigma_{\vert\Sigma_p}})\big)_{p\in\omega}$ is $\preceq^i_c$-decreasing. Also, $w^\mathbb{Z}\!\notin\!\Sigma_p$ if $0\! <\! 8\vert w\vert\! <\! f_{9p+5}$, so that $\sigma^i(\alpha )\!\not=\!\alpha$ if 
$\alpha\!\in\!\bigcup_{0<8i<f_{9p+5}}~\Sigma_p$. In particular, $(\sigma_{\vert\Sigma_p})^2$ is fixed point free for each $p$ since $f_5\! =\! 21$.\medskip

 We finally define a sequence $(w_p)_{p\in\omega}$ of finite binary sequences by $w_0\! :=\! 01$, $w_1\! :=\! 101$ and 
$w_{p+2}\! :=\! w_pw_{p+1}$. Note that $\vert w_p\vert\! =\! f_p$, inductively, so that 
$\sigma^{f_p}(w_p^\mathbb{Z})\! =\! w_p^\mathbb{Z}$. Here is the key fact.\medskip

\noindent\emph{Claim 1.}\it\ Let $p\!\in\!\omega$. Then $w_{9p+5}^\mathbb{Z}\!\in\!\Sigma_p$.\rm\medskip

 Indeed, we consider the subsitution $\tau\! :\! 2^{<\omega}\!\rightarrow\! 2^{<\omega}$ defined by 
$\tau (0)\! :=\! 1$ and $\tau (1)\! :=\! 01$. Note that 
$w_p\! =\!\tau^{p+1}(1)$. Indeed, $\tau^2(1)\! =\!\tau (01)\! =\!\tau (0)\tau (1)\! =\! 101\! =\! w_1$ and 
$$\begin{array}{ll}
\tau^{p+3}(1)\!\!\!\!\!
& =\!\tau^{p+2}\big(\tau (1)\big)\! =\!\tau^{p+2}(01)\! =\!\tau^{p+2}(0)\tau^{p+2}(1)\! =\!\tau^{p+1}\big(\tau (0)\big)\tau^{p+2}(1)
\! =\!\tau^{p+1}(1)\tau^{p+2}(1)\cr
& =\! w_pw_{p+1}\! =\! w_{p+2}.
\end{array}$$

 If $w\! =\!\big( w(0),\!\ldots\! ,w(\vert w\vert\! -\! 1)\big)\!\in\! 2^{<\omega}$, then we set 
$w^{-1}\! =\!\big( w(\vert w\vert\! -\! 1),\!\ldots\! ,w(0)\big)$. The sequence $(w^{-1}_p)_{p\in\omega}$ of \emph{Fibonacci words} is strictly $\subseteq$-increasing, so that its elements are initial segments of the \emph{infinite Fibonacci word} 
$\Phi\!\in\! 2^\omega$. By [Kar, Section 4], $\Phi$ contains no fourth power, i.e., $v^4\!\not\sqsubseteq\!\Phi$ if 
$\emptyset\!\not=\! v\!\in\! 2^{<\omega}$. An induction shows that $w_q^\mathbb{Z}\!\in\!\Sigma_{\{ 0^2,1^3\}}$. We argue by contradiction to prove our claim, which gives $p$ and $w$ with $0\! <\! 8\vert w\vert\! <\! f_{9p+5}$ and 
$w^8\!\sqsubseteq\! w_{9p+5}^\mathbb{Z}$. Note that $w^4\!\sqsubseteq\! w_{9p+5}$ since 
$8\vert w\vert\! <\! f_{9p+5}\! =\!\vert w_{9p+5}\vert$. In particular, if $v\! :=\! w^{-1}$, then 
$v^4\!\sqsubseteq\!\Phi$, which cannot be.\hfill{$\diamond$}\medskip
  
\noindent\emph{Claim 2.}\it\ Let $p\!\in\!\omega$. Then we can find $\alpha\!\in\!\Sigma_p$ and $k\!\in\!\omega$ with 
$2k\! +\! 3\!\leq\! f_{9p+5}$, $\sigma^{2k+3}(\alpha )\! =\!\alpha$ and $\sigma^i(\alpha )\!\not=\!\alpha$ if 
$0\! <\! i\! <\! 2k\! +\! 3$.\rm\medskip

 Indeed, we choose $\alpha\! :=\! w_{9p+5}^\mathbb{Z}$. By Claim 1, $\alpha\!\in\!\Sigma_p$. As 
$\vert w_{9p+5}\vert\! =\! f_{9p+5}$, $\sigma^{f_{9p+5}}(\alpha )\! =\!\alpha$, and $f_{9p+5}$ is odd. Let $n$ be odd and minimal with $\sigma^n(\alpha )\! =\!\alpha$. As $\sigma_{\vert\Sigma_p}$ is fixed point free, $n\!\geq\! 3$, which gives 
$k\!\in\!\omega$ with $n\! =\! 2k\! +\! 3$, so that $2k\! +\! 3\!\leq\! f_{9p+5}$ and $\sigma^{2k+3}(\alpha )\! =\!\alpha$. If 
$0\! <\! i\! <\! 2k\! +\! 3$ and $\sigma^i(\alpha )\! =\!\alpha$, then $i$ has to be even by minimality of $n$. Note then that 
$0\! <\! n\! -\! i\! <\! n$ is odd and $\sigma^{n-i}(\alpha )\! =\!\alpha$, which contradicts the minimality of $n$.\hfill{$\diamond$}\medskip
 
 Claim 2 implies that $\big(\sigma^i(\alpha )\big)_{i\leq 2k+3}$ is a $G_{\sigma_{\vert\Sigma_p}}$-cycle, so that 
$(2k\! +\! 3,C_{2k+3})\preceq^i_c(\Sigma_p,G_{\sigma_{\vert\Sigma_p}})$. In particular, 
$\chi_c(\Sigma_p,G_{\sigma_{\vert\Sigma_p}})\!\geq\! 3$. By Theorem \ref{homeocompact}, 
$\chi_c(\Sigma_p,G_{\sigma_{\vert\Sigma_p}})\! =\! 3$. Assume now, towards a contradiction, that 
$(\Sigma_p,G_{\sigma_{\vert\Sigma_p}})\preceq_c(\Sigma_{p+1},G_{\sigma_{\vert\Sigma_{p+1}}})$ with witness $\varphi$. Let $\alpha$ be given by Claim 2. As $\big(\sigma^i(\alpha )\big)_{i\leq 2k+3}$ is a $G_{\sigma_{\vert\Sigma_p}}$-cycle and odd cycles must map to odd cycles of at most equal length, $(\Sigma_{p+1},G_{\sigma_{\vert\Sigma_{p+1}}})$ contains a cycle of length $2l\! +\! 3\!\leq\! 2k\! +\! 3$. This implies that $\sigma^{2l+3}\big(\varphi (\alpha )\big)\! =\!\varphi (\alpha )$, which cannot be since $0\! <\! 8(2l\! +\! 3)\!\leq\! 8(2k\! +\! 3)\!\leq\! 8f_{9p+5}\! <\! f_{9p+14}$.\hfill{$\square$}\medskip
 
 Lemma \ref{mino} shows that many odometers induce minimal graphs with CCN three. We will now see that it is also the case with subshifts.\medskip
 
\noindent\emph{Notation.}\ Let $r\!\in\! (0,\frac{1}{2})\!\setminus\!\mathbb{Q}$. We consider the irrational rotation 
$R_r\! :\!\mathbb{R}/\mathbb{Z}\!\rightarrow\!\mathbb{R}/\mathbb{Z}$ (well-)defined by $R_r([x])\! :=\! [x\! +\! r]$. We (well-)define $\phi_r\! :\!\mathbb{R}/\mathbb{Z}\!\rightarrow\! 2$ by $\phi_r([x])\! :=\! 0$ if $x\!\in\! [0,r)\mbox{ mod }1$, 
$\phi_r([x])\! :=\! 1$ otherwise, and set 
$\Sigma^2_r\! :=\!\overline{\Big\{\Big(\phi_r\big( R_r^n([x])\big)\Big)_{n\in\mathbb{Z}}\mid [x]\!\in\!\mathbb{R}/\mathbb{Z}\Big\}}$.\medskip

 The following result is mentioned in [MB, Section 4].

\begin{thm} \label{rsub} (Hedlund) Let $r\!\in\! (0,\frac{1}{2})\!\setminus\!\mathbb{Q}$. Then 
$(\Sigma^2_r,\sigma_{\vert\Sigma^2_r})$ is a minimal two-sided subshift.\end{thm}

\noindent\emph{Notation.}\ If $\mathbb{X}\!\in\!\{\mathbb{Z},\omega\}$ and $\Sigma\!\subseteq\! A^\mathbb{X}$, then we denote by $\mathcal{L}(\Sigma )\! :=\!\{ w\!\in\! A^{<\omega}\mid\exists\alpha\!\in\!\Sigma ~~w\!\sqsubseteq\!\alpha\}$ the set of finite words word occurring in $\Sigma$. If moreover $n\!\in\!\omega$, then we set 
$\mathcal{L}_n(\Sigma )\! :=\!\mathcal{L}(\Sigma )\cap A^n$.\medskip

 Let $r\!\in\! (0,\frac{1}{2})\!\setminus\!\mathbb{Q}$. By Theorem \ref{rsub}, $(\Sigma^2_r,\sigma_{\vert\Sigma^2_r})$ is a two-sided subshift. We set $\mathcal{L}^r\! :=\!\mathcal{L}(\Sigma^2_r)$. If moreover $n\!\in\!\omega$, then we set 
$\mathcal{L}^r_n\! :=\!\mathcal{L}^r\cap 2^n$, so that $\mathcal{L}^r_n\! =\!\mathcal{L}_n(\Sigma^2_r)$.

\begin{defi} A subshift $\Sigma$ is \emph{uniformly recurrent} if, for each $w\!\in\!\mathcal{L}(\Sigma )$, there is a natural number $l$ such that, for each $v\!\in\!\bigcup_{n\geq l}~\mathcal{L}_n(\Sigma )$, $w\!\sqsubseteq\! v$.\end{defi}

 By [Sa-T\"o, Section 2], the following result holds.

\begin{thm} \label{minunif} Any minimal two-sided subshift is uniformly recurrent.\end{thm}

\begin{lem} \label{pro} Let $r\!\in\! (0,\frac{1}{2})\!\setminus\!\mathbb{Q}$.\smallskip

\noindent (a) The map $\sigma_{\vert\Sigma^2_r}$ is fixed point free.\smallskip

\noindent (b) $\Sigma^2_r$ is homeomorphic to $2^\omega$.\end{lem}

\vfill\eject

\noindent\emph{Proof.}\ By Theorems \ref{rsub} and \ref{minunif}, $\Sigma^2_r$ is uniformly recurrent. By [MB, Section 4], 
$(\Sigma^2_r,\sigma_{\vert\Sigma^2_r})$ has no periodic point\medskip

\noindent (a) If $\sigma_{\vert\Sigma^2_r}$ has a fixed point $\alpha$, Proposition \ref{period} gives $a\!\in\! 2$ with 
$\alpha\! =\! a^\mathbb{Z}$. In particular, $\alpha$ is a periodic point of $\Sigma^2_r$, which cannot be. Thus 
$\sigma_{\vert\Sigma^2_r}$ is fixed point free.\medskip

\noindent (b) Let $\alpha\!\in\!\Sigma^2_r$. Assume now that $\alpha$ is isolated in $\Sigma^2_r$, which gives 
$\emptyset\!\not=\! w\!\sqsubseteq\!\alpha$ and $p\!\in\!\mathbb{Z}$ with $\Sigma^2_r\cap [w]_p\! =\!\{\alpha\}$. As 
$\Sigma^2_r$ is uniformly recurrent, there is $l$ such that 
$w\!\sqsubseteq\!\alpha\vert [p\! +\!\vert w\vert ,p\! +\!\vert w\vert\! +\! l]$. This gives a natural number $i$ with 
$w\!\subseteq\!\alpha\vert [p\! +\!\vert w\vert\! +\! i,\infty )$, so that $\sigma^{\vert w\vert +i}(\alpha )\!\in\! [w]_p$. Thus 
$\sigma^{\vert w\vert +i}(\alpha )\! =\!\alpha$ and $\alpha$ is periodic. This contradiction shows that $\Sigma^2_r$ is perfect, and thus homeomorphic to $2^\omega$, by [K, 7.4], finishing the proof.\hfill{$\square$}\medskip

 The following result is also mentioned in [MB, Section 4].

\begin{thm} \label{rcomplexity} Let $r\!\not=\! r'\!\in\! (0,\frac{1}{2})\!\setminus\!\mathbb{Q}$. Then the homeomorphisms 
$\sigma_{\vert\Sigma^2_r}$ and $\sigma_{\vert\Sigma^2_{r'}}$ are not conjugate.\end{thm}

 We are now ready to prove that the $\Sigma^2_r$'s induce minimal graphs with CCN exactly three.

\begin{thm} \label{fibgenr} Let $r\!\in\! (0,\frac{1}{2})\!\setminus\!\mathbb{Q}$.\smallskip

\noindent (a) The dynamical system $(\Sigma^2_r,\sigma_{\vert\Sigma^2_r})$ is minimal and the graph 
$(\Sigma^2_r,G_{\sigma_{\vert\Sigma^2_r}})$ has CCN three.\smallskip

\noindent (b) The graph $(\Sigma^2_r,G_{\sigma_{\vert\Sigma^2_r}})$ is $\preceq^i_c$-minimal in $\mathfrak{G}_2$ and in the class of closed graphs on a 0DMC space with CCN at least three.\end{thm}

\noindent\emph{Proof.}\ (a) Note that $\Sigma_r^2$ is not empty. By Theorem \ref{rsub}, 
$(\Sigma^2_r,\sigma_{\vert\Sigma^2_r})$ is minimal. By Lemma \ref{pro}, $\Sigma_r^2$ has cardinality at least two and 
$\sigma_{\vert\Sigma^2_r}$ is fixed point free. By Theorem \ref{homeocompact}, $(\Sigma^2_r,G_{\sigma_{\vert\Sigma^2_r}})$ has CCN two or three. By Theorem \ref{CNG} and minimality, it is enough to find $\alpha\!\in\!\Sigma_r^2$ such that the intersection 
$\overline{\{\sigma^{2p}(\alpha )\mid p\!\in\!\mathbb{Z}\}}\cap\overline{\{\sigma^{2p+1}(\alpha )\mid p\!\in\!\mathbb{Z}\}}$ is not empty.\medskip

 If $[x]\!\in\!\mathbb{R}/\mathbb{Z}$, then 
$\sigma\bigg(\Big(\phi_r\big( R_r^n([x])\big)\Big)_{n\in\mathbb{Z}}\bigg)\! =\!
\Big(\phi_r\big( R_r^{n+1}([x])\big)\Big)_{n\in\mathbb{Z}}$,  
so that 
$$\sigma^{2p+\varepsilon}\bigg(\Big(\phi_r\big( R_r^n([x])\big)\Big)_{n\in\mathbb{Z}}\bigg)\! =\!
\Big(\phi_r\big( R_r^{n+2p+\varepsilon}([x])\big)\Big)_{n\in\mathbb{Z}}.$$ 
Note then that ${\phi_r}_{\vert (\mathbb{R}/\mathbb{Z})\setminus ([0]\cup [r])}$ is continuous. So it is enough to find 
$$[x]\!\in\! (\mathbb{R}/\mathbb{Z})\!\setminus\!\big\{ R_r^l([y])\mid l\!\in\!\mathbb{Z}\wedge y\!\in\!\{ 0,r\}\big\}$$ 
such that $\overline{\big\{\big( R_r^{n+2p+1}([x])\big)_{n\in\mathbb{Z}}\mid p\!\in\!\mathbb{Z}\big\}}$ is not empty. Pick $[x]$ in 
$$(\mathbb{R}/\mathbb{Z})\!\setminus\!\big\{ R_r^l([y])\mid l\!\in\!\mathbb{Z}\wedge y\!\in\!\{ 0,r\}\big\}$$ 
arbitrary. It is enough to check that $\big( R_r^n([x])\big)_{n\in\mathbb{Z}}\!\in\!
\overline{\big\{\big( R_r^{n+2p+1}([x])\big)_{n\in\mathbb{Z}}\mid p\!\in\!\mathbb{Z}\big\}}$. Note that the restriction $b$ of the canonical map from $\pi\! :\!\mathbb{R}\!\rightarrow\!\mathbb{R}/\mathbb{Z}$ to $[0,1)$ is a bijection. The map $b$ is a homeomorphism from $[0,1)$ equipped with $\tau\! :=\!\{ b^{-1}(O)\mid O\mbox{ open in }\mathbb{R}/\mathbb{Z}\}$ onto 
$\mathbb{R}/\mathbb{Z}$. As usual, $[0,1)$ is equipped with the $\tau$-compatible metric defined by 
$d(x,y)\! :=\!\mbox{min}(\vert x\! -\! y\vert ,1\! -\!\vert x\! -\! y\vert )$. The previous identification through $b$ defines a compatible metric $D$ on $\mathbb{R}/\mathbb{Z}$ for which the $R_r$'s are isometries. Let $q$ be a natural number, and, for 
$n\!\in\! [-q,q]\cap\mathbb{Z}$, $0\! <\! a_n\! <\! r_n\! :=\! b^{-1}\big( R_r^n([x])\big)\! <\! b_n\! <\! 1$ mod $1$. We choose a natural number $m\! >\! 0$ such that 
$\frac{1}{m}\! <\!\mbox{min}_{n\in [-q,q]\cap\mathbb{Z}}~\mbox{min}(r_n\! -\! a_n,b_n\! -\! r_n)$.  The next claim is inspired by the proof of [Ku, Proposition 1.32].

\vfill\eject
 
\noindent\emph{Claim.}\it\ Let $r\!\in\! (0,\frac{1}{2})\!\setminus\!\mathbb{Q}$, $x\!\in\!\mathbb{R}$, and $m\! >\! 0$ be a natural number. Then there is $p\!\in\!\mathbb{Z}$ such that, for each $n\!\in\!\mathbb{Z}$, 
$D\big( R_r^{n+2p+1}([x]),R_r^n([x])\big)\! <\!\frac{1}{m}$.\rm\medskip

 Indeed, as $r\!\notin\!\mathbb{Q}$, the $m$ classes $[x]$, $R_r^2([x])$, $\ldots$, $R_r^{2m-2}([x])$ are pairwise different. As $b$ is defined on $[0,1)$, we can find $0\!\leq\! i\! <\! j\! <\! m$ with the property that 
$D\big( R_r^{2i}([x]),R_r^{2j}([x])\big)\! <\!\frac{1}{m}$. As $R_r$ is an isometry, 
$D\big( [x],R_r^{2(j-i)}([x])\big)\! =\! D\big( R_r^{2i}([x]),R_r^{2j}([x])\big)$. Now $R_r^{2(j-i)}\! =\! R_b$ is also a rotation, and either $0\! <\! b\! <\!\frac{1}{m}$, or $1\! -\!\frac{1}{m}\! <\! b\! <\! 1$. In both cases, for any $y\!\in\!\mathbb{R}/\mathbb{Z}$, there is $k\! >\! 0$ with $D\big( R_r^{2k(j-i)}([x]),y\big)\! <\!\frac{1}{m}$. Applying this to $y\! :=\! R_r([x])$, and putting 
$p\! :=\! k(i\! -\! j)$, we get $D\big( R_r^{2p+1}([x]),[x]\big)\! <\!\frac{1}{m}$. It remains to apply again the fact that $R_r$ is an isometry.\hfill{$\diamond$}\medskip

 The claim provides $p\!\in\!\mathbb{Z}$ such that, for each $n\!\in\!\mathbb{Z}$, 
$D\big( R_r^{n+2p+1}([x]),b(r_n)\big)\! <\!\frac{1}{m}$. In particular, $b^{-1}\big( R_r^{n+2p+1}([x])\big)\!\in\! (a_n,b_n)$ mod $1$ if $n\!\in\! [-q,q]\cap\mathbb{Z}$, as desired.\medskip

\noindent (b) We apply (a) and Lemma \ref{mino}.\hfill{$\square$}\medskip

 Theorem \ref{fibgenr} gives a version of Theorem \ref{eantichmino} for subshifts.
 
\begin{cor} \label{controt} There is a $\preceq^i_c$-antichain $(\Sigma_r,G_{\sigma_{\vert\Sigma_r}})_{r\in\mathbb{R}}$, where\smallskip 

\noindent (a) $\Sigma_r$ is a two-sided subshift homeomorphic to $2^\omega$,\smallskip

\noindent (b) $\sigma_{\vert\Sigma_r}$ is a minimal homeomorphism of $\Sigma_r$, and $G_{\sigma_{\vert\Sigma_r}}$ has CCN three,\smallskip

\noindent (c) $(\Sigma_r,G_{\sigma_{\vert\Sigma_r}})_{r\in\mathbb{R}}$ is $\preceq^i_c$-minimal in $\mathfrak{G}_2$ and in the class of closed graphs on a 0DMC space with CCN at least three.\end{cor} 
 
\noindent\emph{Proof.}\ Let $r\!\in\! (0,\frac{1}{2})\!\setminus\!\mathbb{Q}$. By Theorem \ref{rsub} and Lemma \ref{pro}(b), 
$(\Sigma^2_r,G_{\sigma_{\vert\Sigma^2_r}})$ is a minimal two-sided subshift homeomorphic to $2^\omega$. By Lemma \ref{pro}(a), $\sigma_{\vert\Sigma^2_r}$ is fixed point free, so that $\Sigma^2_r$ is a closed graph, with CCN three by Theorem \ref{fibgenr}(a). As $(\Sigma^2_r,G_{\sigma_{\vert\Sigma^2_r}})$ has cardinality at least three and is minimal, the vertices of 
$\Sigma^2_r$ have degree two. By Theorem \ref{fibgenr}(b), $(\Sigma^2_r,G_{\sigma_{\vert\Sigma^2_r}})$ is $\preceq^i_c$-minimal in $\mathfrak{G}_2$ and in the class of closed graphs on a 0DMC space with CCN at least three.\medskip

 So it is enough to find a subfamily of $(\Sigma_r^2,G_{\sigma_{\vert\Sigma_r^2}})_{r\in (0,\frac{1}{2})\setminus\mathbb{Q}}$ which is a $\preceq^i_c$-antichain. By Lemma \ref{flipo}, it is enough to ensure that the homeomorphisms corresponding to the  elements of the subfamily are pairwise not flip conjugate. By Theorem \ref{rcomplexity}, the 
 $\sigma_{\vert\Sigma_r^2}$'s are pairwise not conjugate. Thus $\sigma_{\vert\Sigma_r^2}$ is flip-conjugate to 
 $\sigma_{\vert\Sigma_{r'}^2}$ exactly when $\sigma_{\vert\Sigma_r^2}$ is conjugate to $\sigma_{\vert\Sigma_{r'}^2}^{-1}$. The key remark is that if $\sigma_{\vert\Sigma_r^2}$ is conjugate to $\sigma_{\vert\Sigma_{r'}^2}^{-1}$ and 
$\sigma_{\vert\Sigma_{r''}^2}^{-1}$, then $\sigma_{\vert\Sigma_{r'}^2}$ is conjugate to $\sigma_{\vert\Sigma_{r''}^2}$, which implies that $r'\! =\! r''$. We inductively construct a injective family $(r_\xi )_{\xi <2^{\aleph_0}}$ of elements of 
$(0,\frac{1}{2})\!\setminus\!\mathbb{Q}$ such that the $\sigma_{\vert\Sigma_{r_\xi}^2}$'s are pairwise not flip-conjugate. $r_0$ is an arbitrary element of $(0,\frac{1}{2})\!\setminus\!\mathbb{Q}$. Assume that $1\!\leq\!\eta\! <\! 2^{\aleph_0}$ and 
$(r_\xi )_{\xi <\eta}$ are constructed. The key remark shows that, for each $\xi\! <\!\eta$, there is at most one element $r'_\xi$ of $(0,\frac{1}{2})\!\setminus\!\mathbb{Q}$ such that $\sigma_{\vert\Sigma_{r_\xi}^2}$ is conjugate to 
$\sigma_{\vert\Sigma_{r'_\xi}^2}^{-1}$. If $r'_\xi$ does not exist, then we set $r'_\xi\! :=\! r_\xi$. We choose 
$r_\eta\!\in\! (0,\frac{1}{2})\!\setminus\! (\mathbb{Q}\cup\bigcup_{\xi <\eta}~\{r_\xi ,r'_\xi\} )$,  
so that $(r_\xi )_{\xi\leq\eta}$ is as desired.\hfill{$\square$}\medskip

 Corollary \ref{controt} gives a second proof of the fact, met in Theorem \ref{eantichmino}, that any $\preceq^i_c$-basis for 
$\mathfrak{G}_2$ must have size continuum. We now prove that the basis given by Proposition \ref{me}(a) is far from being a basis for $\mathfrak{G}_2$, as announced in the introduction.

\begin{prop} \label{odoshift} Let ${\bf d}\!\in\! (\omega\!\setminus\! 2)^\omega\cap\mathcal{O}$ and 
$r\!\in\! (0,\frac{1}{2})\!\setminus\!\mathbb{Q}$. Then $(\mathcal{C},G_o)$ and $(\Sigma^2_r,G_{\sigma_{\vert\Sigma^2_r}})$ are $\preceq^i_c$-incompatible in the class of closed graphs on a 0DMC space with CCN at least three.\end{prop} 
 
\noindent\emph{Proof.}\ Towards a contradiction, suppose that we can find a 0DMC space $K$ and a closed graph $G$ on $K$ which is $\preceq^i_c$-below our two graphs. Thus 
$(\Sigma^2_r,G_{\sigma_{\vert\Sigma^2_r}})\preceq^i_c(\mathcal{C},G_o)$, by Corollary \ref{controt}. Theorem \ref{rsub} and Lemma \ref{pro} allow us to apply Lemma \ref{flipo}, so that $\sigma_{\vert\Sigma^2_r}$ and $o$ are flip-conjugate, with witness say $\varphi$. We already saw that $(\mathcal{C},o)$ is equicontinuous, i.e., $(o^n)_{n\in\omega}$ is equicontinuous. In fact, as $o$ is an isometry, $(o^n)_{n\in\mathbb{Z}}$ is equicontinuous. The uniform continuity of $\varphi$ and 
$\varphi^{-1}$ implies that $(\Sigma^2_r,\sigma_{\vert\Sigma^2_r})$ is also equicontinuous. Theorem \ref{rsub} and [Ku, Proposition 3.68(2)] imply that $(\Sigma^2_r,\sigma_{\vert\Sigma^2_r})$ is expansive, which gives $\varepsilon\! >\! 0$ such that, for each $\alpha\!\not=\!\beta\!\in\!\Sigma^2_r$, there is $n\!\in\!\mathbb{Z}$ with 
$d\big(\sigma^n(\alpha ),\sigma^n(\beta )\big)\!\geq\!\varepsilon$. Lemma \ref{pro}(b) gives $\alpha\!\in\!\Sigma^2_r$. The equicontinuity of $(\Sigma^2_r,\sigma_{\vert\Sigma^2_r})$ gives $\delta\! >\! 0$ such that, for each 
$\beta\!\in\! B(\alpha ,\delta )$ and each $k\!\in\!\mathbb{Z}$, $d\big(\sigma^k(\alpha ),\sigma^k(\beta )\big)\! <\!\varepsilon$. Lemma \ref{pro}(b) gives $\beta\!\in\! B(\alpha ,\delta )\!\setminus\!\{\alpha\}$, which is the desired contradiction.\hfill{$\square$}\medskip

\noindent\emph{Remark.}\ Theorem \ref{eantichmino} and Corollary \ref{controt} provide examples 
$\preceq_c^i$-minimal in $\mathfrak{G}_2$. None of them is $\preceq_c$-minimal in $\mathfrak{K}$. Indeed, there is a dense orbit. Let $(X,G_f)$ be one of them. Theorem \ref{corcomp'''''''} provides $\beta\!\in\!\mathcal{J}^c$ such that 
$(\mathbb{K}_\beta ,\mathbb{G}_\beta )\preceq_c(X,G_f)$ and the vertices of $\mathbb{G}_\beta$ have degree at most one. Assume that $(X,G_f)\preceq_c(\mathbb{K}_\beta ,\mathbb{G}_\beta )$, towards a contradiction. Then the dense orbit of $f$ has to be sent to a two-point set because of the degree. So $X$ has to be sent to this closed set, by density. But this contradicts the fact that $(X,G_f)$ has CCN at least three. The examples are not $\preceq_c^i$-minimal in 
$\mathfrak{K}$. Indeed, consider the dense orbit $D\! =\!\mbox{Orb}_f(x)$. Assume that $(X,G_f)\preceq^i_c(X,G_{f_{\vert D}})$ with witness $\varphi$, towards a contradiction. The proof of Lemma \ref{flipo} shows that $\varphi [D]\! =\!\mbox{Orb}_f\big(\varphi (x)\big)\! =\! D$ since $(f_{\vert D})^2$ is fixed point free. Thus $\varphi^2[G_{f_{\vert D}}]\! =\! G_{f_{\vert D}}$. By injectivity of $\varphi$, there is no more room for 
$\varphi [X\!\setminus\! D]$.\medskip

 We now turn to the proof of Theorem \ref{CB1intro}. By [K, 6.C], the countable MC spaces can be analyzed through their Cantor-Bendixson rank. Recall that if $X$ is a topological space, then the 
\emph{Cantor-Bendixson derivative} of $X$ is $X'\! :=\!\{ x\!\in\! X\mid x\mbox{ is a limit point of }X\}$. The \emph{iterated Cantor-Bendixson derivatives} are defined by $X^0\! :=\! X$, $X^{\alpha +1}\! :=\! (X^\alpha )'$, and 
$X^\lambda\! :=\!\bigcap_{\alpha <\lambda}~X^\alpha$ if $\lambda$ is a limit ordinal. Note that if $f$ is a homeomorphism of 
$X$, then all the derivatives are $f$-invariant, i.e., $f[X^\alpha ]\! =\! X^\alpha$ if $\alpha$ is an ordinal. If $X\!\not=\!\emptyset$ is a countable metrizable compact space, then the \emph{Cantor-Bendixson rank} of $X$ is the minimal countable ordinal $\alpha$ with $X^\alpha\! =\!\emptyset$, which is a successor ordinal by compactness. The odd cycles provide examples of graphs induced by a homeomorphism of a countable (0D)MC space with Cantor-Bendixson rank one whose CCN is three, and which are $\preceq^i_c$-minimal in $\mathfrak{G}_2$. We now provide examples for higher ranks, including the example 
$(K_0,h_0)$ mentioned in the introduction.\medskip

\noindent\emph{Proof of Theorem \ref{CB1intro}.}\ (a) For $\xi\! =\! 1$, we can take 
$\Sigma\! :=\!\mbox{Orb}_\sigma\big( (012)^\infty\!\cdot\! (012)^\infty\big)\!\subseteq\! 3^\mathbb{Z}$, which defines a cycle on three points, and we apply Corollary \ref{basisfin}. Assume now that $2\!\leq\!\xi\! =\! n\! +\! 2\! <\!\omega$. We set, for 
$j,m\!\in\!\omega$, $w^0_j\! :=\! 01$, $w^{m+1}_0\! :=\! 1^2$, and 
$w^{m+1}_{j+1}\! :=\! (01)^{j+1}1^2{^\frown}_{k\leq j+1}~w^m_k$. We then set 
${\alpha_0\! :=\! (01)^\infty\!\cdot\! (01)^\infty}$, and, for $m\!\in\!\omega$, 
$\alpha_{m+1}\! :=\! (01)^\infty\!\cdot\! 1^2{^\frown}_{j\in\omega}~w^m_j$ and 
$\beta_m\! :=\! (01)^\infty\!\cdot\! 1{^\frown}_{j\in\omega}~w^m_j$. Finally, 
$\Sigma\! =\!\bigcup_{m\leq n}~\mbox{Orb}_\sigma (\alpha_m)\cup\mbox{Orb}_\sigma (\beta_n)$, so that $\Sigma\! =\! K_0$ and $\sigma_{\vert\Sigma}\! =\! h_0$ if $\xi\! =\! 2$.\medskip

 $\Sigma$ is by definition countable, and $\sigma [\Sigma ]\! =\!\Sigma$. We then set, for 
$1\!\leq\! i\!\leq\! n\! +\! 1$, 
$$\Sigma^{(i)}\! =\!\bigcup_{m\leq n+1-i}~\mbox{Orb}_\sigma (\alpha _m).$$

\vfill\eject

\noindent\emph{Claim 1.}\it\ (a) Let $1\!\leq\! i\!\leq\! n\! +\! 1$. Then 
$\Sigma^{(i)}\! =\!\overline{\mbox{Orb}_\sigma (\alpha_{n+1-i})}$.\smallskip

(b) $\Sigma\! =\!\overline{\mbox{Orb}_\sigma (\beta_n)}$.\rm\medskip

 Indeed, for (a), we argue by induction on $n\! +\! 1\! -\! i$. Note first that 
$$\mbox{Orb}_\sigma (\alpha_0)\! =\!\{ (01)^\infty\!\cdot\! (01)^\infty ,(10)^\infty\!\cdot\! (10)^\infty\}$$ 
is closed, so that $\Sigma^{(n+1)}\! =\!\overline{\mbox{Orb}_\sigma (\alpha_0)}$.\medskip

 Let us prove that $\alpha_m\!\in\!\overline{\mbox{Orb}_\sigma (\alpha_{m+1})}$ if $m\!\in\!\omega$, which holds for $m\! =\! 0$. Note that 
$$\alpha_{m+2}\! =\! (01)^\infty\!\cdot\! 1^2{^\frown}_{j\in\omega}~w^{m+1}_j\! =\! 
(01)^\infty\!\cdot\! 1^21^2{^\frown}_{j\in\omega}~\big( (01)^{j+1}1^2{^\frown}_{k\leq j+1}~w^m_k\big) .$$ 
If $\alpha\!\in\! 2^\mathbb{Z}$ and $a\!\leq\! b$ are integers, then we define 
$\alpha_{[a,b]}\!\in\! 2^{b-a+1}$ by $\alpha_{[a,b]}(l)\! =\!\alpha (a\! +\! l)$ if $l\!\leq\! b\! -\! a$. Note that 
$(01)^{j+1}1^2{^\frown}_{k\leq j+1}~w^m_k\! =\! {\alpha_{m+1}}_{[-2j-2,1+\Sigma_{k\leq j+1}~\vert w^m_k\vert ]}$. An induction shows that $\vert w^m_k\vert\!\geq\! 2$, so that $1+\Sigma_{k\leq j+1}~\vert w^m_k\vert\!\geq\! 2j\! +\! 5$. This implies that 
$\alpha_{m+1}\!\in\!\overline{\mbox{Orb}_\sigma (\alpha_{m+2})}$, as desired.\medskip

 From this we deduce, inductively and by continuity of $\sigma$ and $\sigma^{-1}$, that 
$$\Sigma^{(i)}\! =\!\Sigma^{(i+1)}\cup\mbox{Orb}_\sigma (\alpha_{n+1-i})\!\subseteq\!
\overline{\mbox{Orb}_\sigma (\alpha_{n-i})}\cup\mbox{Orb}_\sigma (\alpha_{n+1-i})
\!\subseteq\!\overline{\mbox{Orb}_\sigma (\alpha_{n+1-i})}.$$
This shows that $\Sigma^{(i)}\!\subseteq\!\overline{\mbox{Orb}_\sigma (\alpha_{n+1-i})}$ if 
$1\!\leq\! i\!\leq\! n\! +\! 1$.\medskip

 Assume then that $1\!\leq\! i\!\leq\! n$, $(k_p)_{p\in\omega}\!\in\!\mathbb{Z}^\omega$ and 
$\big(\sigma^{k_p}(\alpha_{n+1-i})\big)_{p\in\omega}$ converges to $\alpha\!\in\! 2^\mathbb{Z}$. We want to see that 
$\alpha\!\in\!\Sigma^{(i)}$, and we may assume that $i\! <\! n$. If $(k_p)_{p\in\omega}$ has a constant subsequence, then 
$\alpha\!\in\!\mbox{Orb}_\sigma (\alpha_{n+1-i})$ and we are done. So we may assume that $(k_p)_{p\in\omega}$ tends to 
$\pm\infty$. If $(k_p)_{p\in\omega}$ tends to $-\infty$, then $\alpha\!\in\!\mbox{Orb}_\sigma (\alpha_0)$ since 
$\alpha_{n+1-i}\! =\! (01)^\infty\!\cdot\! 1^2{^\frown}_{j\in\omega}~w^{n-i}_j$, and we are done. So we may write 
$k_p\! =\! 2\! +\!\Sigma_{j<l_p}~\vert w_j^{n-i}\vert\! +\! j_p$, where $j_p\! <\!\vert w_{l_p}^{n-i}\vert$, and 
$(l_p)_{p\in\omega}$ tends to $\infty$. As above, $\alpha_{n+1-i}\! =\! (01)^\infty\!\cdot\! 1^21^2
{^\frown}_{j\in\omega}~\big( (01)^{j+1}1^2{^\frown}_{k\leq j+1}~w^{n-i-1}_k\big)$, and 
$$(01)^{j+1}1^2{^\frown}_{k\leq j+1}~w^{n-i-1}_k\! =\! {\alpha_{n-i}}_{[-2j-2,1+\Sigma_{k\leq j+1}~\vert w^{n-i-1}_k\vert ]} .$$ 
If $(j_p)_{p\in\omega}$ and $(\vert w_{l_p}^{n-i}\vert\! -\! j_p)_{p\in\omega}$ also tend to $\infty$, then 
$\alpha\!\in\!\overline{\mbox{Orb}_\sigma (\alpha_{n-i})}\!\subseteq\!\Sigma^{(i+1)}\!\subseteq\!\Sigma^{(i)}$, by induction assumption, as desired. Otherwise, we may assume that $(j_p)_{p\in\omega}$ or 
$(\vert w_{l_p}^{n-i}\vert\! -\! j_p)_{p\in\omega}$ is constant, so that $\alpha\!\in\!\mbox{Orb}_\sigma (\alpha_0)$ since 
$w^{m+1}_j$ starts and ends with $(01)^j$, by induction.\medskip

 The proof of (b) is similar.\hfill{$\diamond$}\medskip
 
 Claim 1 implies that $\Sigma$ is a two-sided subshift.\medskip
 
\noindent\emph{Claim 2.}\it\ $\Sigma$ has Cantor-Bendixson rank $\xi$.\rm\medskip

 Indeed, let us show that $\Sigma^{(i)}$ is the $i$th iterated Cantor-Bendixson derivative of $\Sigma$ if 
${1\!\leq\! i\!\leq\! n\! +\! 1}$. By Claim 1(b), 
$\Sigma\!\setminus\!\mbox{Orb}_\sigma (\beta_n)\!\subseteq\!\Sigma'$. As $\Sigma$ is countable, it has an isolated point, which has therefore to be in $\mbox{Orb}_\sigma (\beta_n)$. As $\sigma$ is a homeomorphism, 
$\mbox{Orb}_\sigma (\beta_n)$ is disjoint from $\Sigma'$, showing that ${\Sigma'\! =\!
\Sigma\!\setminus\!\mbox{Orb}_\sigma (\beta_n)\!\subseteq\!\bigcup_{m\leq n}~\mbox{Orb}_\sigma (\alpha_m)
\! =\!\Sigma^{(1)}}$. It remains to see that $\alpha_m\!\notin\!\mbox{Orb}_\sigma (\beta_n)$ if $m\!\leq\! n$ to get 
$\Sigma^{(1)}\! =\!\Sigma'$. As all the odd coordinates of the $w^m_j$'s and the $\alpha_m$'s are 1, and we can find an even coordinate of $\beta_n$ and an odd coordinate of $\beta_n$ which are 0, we are done.

\vfill\eject
 
 Fix now $1\!\leq\! i\!\leq\! n$. By Claim 1(a), 
$\Sigma^{(i)}\!\setminus\!\mbox{Orb}_\sigma (\alpha_{n+1-i})\!\subseteq\! (\Sigma^{(i)})'$. As 
$\Sigma^{(i)}$ is countable, it has an isolated point, which has therefore to be in 
$\mbox{Orb}_\sigma (\alpha_{n+1-i})$. As $\sigma$ is a homeomorphism, $\mbox{Orb}_\sigma (\alpha_{n+1-i})$ is disjoint from $(\Sigma^{(i)})'$, showing that 
$(\Sigma^{(i)})'\! =\!\Sigma^{(i)}\!\setminus\!\mbox{Orb}_\sigma (\alpha_{n+1-i})\!\subseteq\!
\bigcup_{m\leq n-i}~\mbox{Orb}_\sigma (\alpha_m)\! =\!\Sigma^{(i+1)}$. 
It remains to see that $\alpha_m\!\notin\!\mbox{Orb}_\sigma (\alpha_{m+1+p})$ if $m,p\!\in\!\omega$ to get 
$\Sigma^{(i+1)}\! =\! (\Sigma^{(i)})'$ and, inductively, that $\Sigma^{(i+1)}$ is the $(i\! +\! 1)$th iterated Cantor-Bendixson derivative of $\Sigma$. Note that if $m\!\in\!\omega$, then 
$w^{m+2}_1\! :=\! (01)1^21^2w^{m+1}_1$, so that $(011^21^2)^m011^201\!\subseteq\! w^{m+1}_1$, 
$w^{m+1}_1$ and $w^{m+2+p}_1$ are incompatible if $p\!\in\!\omega$. In particular, 
$\alpha_{m+2}\!\not=\!\alpha_{m+3+p}$, and, because of $(01)^\infty$, 
$\alpha_{m+2}\!\notin\!\mbox{Orb}_\sigma (\alpha_{m+3+p})$. Because of $1^2$, 
$\alpha_0\!\notin\!\mbox{Orb}_\sigma (\alpha_{1+p})$, and because of $1^21^2$, 
$\alpha_1\!\notin\!\mbox{Orb}_\sigma (\alpha_{2+p})$.\medskip
 
 As $\Sigma^{(n+1)}\! =\!\mbox{Orb}_\sigma (\alpha_0)\!\not=\!\emptyset$ is finite, $\Sigma$ has Cantor-Bendixson rank $\xi$.\hfill{$\diamond$}\medskip

 As no sequence in $\Sigma$ is constant, $\sigma_{\vert\Sigma}$ is fixed point free. As moreover 
$\Sigma$ is not empty, $(\Sigma ,G_{\sigma_{\vert\Sigma}})$ has CCN two or three by Theorem 8.1. Note then that 
$$(\alpha_0,\alpha_0)\! =\!\mbox{lim}_{p\rightarrow\infty}~
\big(\sigma^{-2p}(\beta_n),\sigma^{1+\Sigma_{j<p}~\vert w^n_j\vert}(\beta_n)\big)\!\in\!
\overline{\bigcup_{q\in\omega}~G_{\sigma_{\vert\Sigma}}^{2q+1}}\mbox{,}$$ 
so that $\chi_c(\Sigma ,G_{\sigma_{\vert\Sigma}})\! =\! 3$ by Theorem 1.7. For the minimality of 
$(\Sigma ,G_{\sigma_{\vert\Sigma}})$, we first prove the following.\medskip

\noindent\emph{Claim 3.}\it\ Let $V\!\subseteq\!\Sigma$, and 
$E\!\subseteq\! G_{\sigma_{\vert\Sigma}}\cap V^2$ be a graph on $V$ such that $(V,E)$ has CCN three. Then 
$\big(\alpha ,\sigma (\alpha )\big)\!\in\! E$ if $\alpha\!\in\!\mbox{Orb}_\sigma (\beta_n)$.\rm\medskip

 Indeed, we argue by contradiction. Let $k\!\in\!\mathbb{Z}$ with $\alpha\! =\!\sigma^k(\beta_n)$. Recall that the sets of the form $[w]_q\! :=\!\{\beta\!\in\! 2^\mathbb{Z}\mid\forall j\! <\!\vert w\vert ~~w(j)\! =\!\beta (q\! +\! j)\}$, where 
$w\!\in\! 2^{<\omega}$ and $q\!\in\!\mathbb{Z}$, form a basis made up of clopen subsets of the space $2^\mathbb{Z}$.\medskip

 Assume first that $n\! =\! 0$, so that 
$\Sigma_n\! =\!\! =\!\mbox{Orb}_\sigma\big( (01)^\infty\!\cdot\! (01)^\infty\big)\cup
\mbox{Orb}_\sigma\big( (01)^\infty\!\cdot\! 1(01)^\infty\big)$, and 
$$\sigma^{2p+\varepsilon}(\alpha_n)\! =\!\left\{\!\!\!\!\!\!\!
\begin{array}{ll}
& (01)^\infty\!\cdot\! (01)^\infty\mbox{ if }\varepsilon\! =\! 0\mbox{,}\cr
& (10)^\infty\!\cdot\! (10)^\infty\mbox{ if }\varepsilon\! =\! 1\mbox{,}
\end{array}
\right.$$ 
$$\sigma^{2p+\varepsilon}(\beta_n)\! =\!\left\{\!\!\!\!\!\!\!
\begin{array}{ll}
& (01)^\infty 1(01)^{p-1}0\!\cdot\! (10)^\infty\mbox{ if }p\! >\! 0\wedge\varepsilon\! =\! 0\mbox{,}\cr
& (01)^\infty 1(01)^p\!\cdot\! (01)^\infty\mbox{ if }p\!\geq\! 0\wedge\varepsilon\! =\! 1\mbox{,}\cr
& (01)^\infty\!\cdot\! (01)^p1(01)^\infty\mbox{ if }p\!\leq\! 0\wedge\varepsilon\! =\! 0\mbox{,}\cr
& (10)^\infty\!\cdot\! (10)^{-p-1}11(01)^\infty\mbox{ if }p\! <\! 0\wedge\varepsilon\! =\! 1
\end{array}
\right.$$ 
if $p\!\in\!\mathbb{Z}$ and $\varepsilon\!\in\! 2$. We set $C\! :=\! ([0]_{-k}\cup [1^20]_{-k-1})\cap V$, so that $C$ is a clopen subset of $V$ and $E\cap\big( C^2\cup (V\!\setminus\! C)^2\big)\! =\!\emptyset$ since 
$$\begin{array}{ll}
& \ldots,\sigma^{k-2}(\alpha_n),\sigma^k(\alpha_n),\sigma^{k+2}(\alpha_n),\ldots\!\in\! C\mbox{,}\cr
& \ldots,\sigma^{k-3}(\alpha_n),\sigma^{k-1}(\alpha_n),\sigma^{k+1}(\alpha_n),\ldots\!\notin\! C\mbox{,}\cr
\end{array}$$ 
$$\begin{array}{ll}
& \ldots,\sigma^{k-4}(\beta_n),\sigma^{k-2}(\beta_n),\sigma^k(\beta_n),\sigma^{k+1}(\beta_n),
\sigma^{k+3}(\beta_n),\sigma^{k+5}(\beta_n),\ldots\!\in\! C\mbox{,}\cr
& \ldots,\sigma^{k-5}(\beta_n),\sigma^{k-3}(\beta_n),\sigma^{k-1}(\beta_n),\sigma^{k+2}(\beta_n),\sigma^{k+4}(\beta_n),
\sigma^{k+6}(\beta_n),\ldots\!\notin\! C\mbox{,}\cr
\end{array}$$ 
which contradicts the fact that $(V,E)$ has CCN three.\medskip

 In this argument, the case $k\!\not=\! 0$ is similar to the case $k\! =\! 0$, we just have to translate the basic clopen sets of the form $[w]_q$. It will also be the case in the general case $n\!\geq\! 1$ that we now consider, so that we may and will assume that $k\! =\! 0$.
 
\vfill\eject
 
 We set $C\! :=\! ([0]_0\cup [1^20]_0\cup [1^50]_{-1}\cup [01^40]_{-2})\cap V$, so that $C$ is a clopen subset of $V$. We already noticed that the odd coordinates of the $w^m_j$'s and the $\alpha_m$'s are 1, that the $w^m_j$'s have a stricly positive even length, and that $w^{m+1}_j$ starts and ends with $(01)^j$. The definition of the $w^m_j$'s, the 
 $\alpha_m$'s and the $\beta_n$'s then imply that 0 can only be an even coordinate of 
 $\alpha_m$, a negative even coordinate of $\beta_n$ or a positive odd coordinate of $\beta_n$. Moreover, $01^l0$ can be of the form ${\alpha_m}_{[a,b]}$ only if $l\!\in\!\{ 1,3,5\}$ and $a$ is even, and of the form ${\beta_n}_{[a,b]}$ only if $l\! =\! 1$ and $a\!\leq\! -4$ is even, $l\!\in\!\{ 1,3,5\}$ and $a\!\geq\! 3$ is odd, or $l\! =\! 4$ and $a\! =\! -2$. This implies that 
 $E\cap\big( C^2\cup (V\!\setminus\! C)^2\big)\! =\!\emptyset$ since 
$$\begin{array}{ll}
& \ldots,\sigma^{-2}(\alpha_m),\alpha_m,\sigma^2(\alpha_m),\ldots\!\in\! C\mbox{,}\cr
& \ldots,\sigma^{-3}(\alpha_m),\sigma^{-1}(\alpha_m),\sigma(\alpha_m),\ldots\!\notin\! C\cr
\end{array}$$ 
if $m\!\leq\! n$ and 
$$\begin{array}{ll}
& \ldots,\sigma^{-4}(\beta_n),\sigma^{-2}(\beta_n),\beta_n,\sigma (\beta_n),
\sigma^3(\beta_n),\sigma^5(\beta_n),\ldots\!\in\! C\mbox{,}\cr
& \ldots,\sigma^{-5}(\beta_n),\sigma^{-3}(\beta_n),\sigma^{-1}(\beta_n),\sigma^2(\beta_n),
\sigma^4(\beta_n),\sigma^6(\beta_n),\ldots\!\notin\! C\mbox{,}\cr
\end{array}$$ 
which contradicts the fact that $(V,E)$ has CCN three again.\hfill{$\diamond$}\medskip

 Assume now that $(K,G)\!\in\!\mathfrak{G}_2$ and $(K,G)\preceq^i_c(\Sigma ,G_{\sigma_{\vert\Sigma}})$ with witness $\varphi$, which implies that $(K,G)$ has CCN three. As $(K,G)\!\in\!\mathfrak{G}_2$, there is a homeomorphism 
$f\! :\! K\!\rightarrow\! K$ with $G\! =\! G_f$. As $\chi_c(K,G)\! =\! 3$, the set $F_1$ of fixed points of $f$ is a clopen subset of 
$K$, and $\chi_c(K\!\setminus\! F_1,G_f\cap (K\!\setminus\! F_1)^2)\! =\! 3$ by Corollary 7.3. This implies that we may assume that $f$ is fixed point free, so that $G$ is compact. We set ${V\! :=\!\varphi [K]}$ and $E\! :=\! (\varphi\!\times\!\varphi )[G]$, so that $V\!\subseteq\!\Sigma$ is a 0DMC space, $E\!\subseteq\! G_{\sigma_{\vert\Sigma}}$ is a compact graph on $V$, 
${(K,G)\preceq^i_c(V,E)}$ with witness $\varphi$, and $(V,E)\preceq^i_c(K,G)$ with witness 
$\varphi^{-1}$ by compactness.\medskip

 By Claim 3, $\big(\alpha ,\sigma (\alpha )\big)\!\in\! E$ if $\alpha\!\in\!\mbox{Orb}_\sigma (\beta_n)$. The density of $\mbox{Orb}_\sigma (\beta_n)$ in $\Sigma$ given by Claim 1 and the compactness of $E$ then imply that 
$\textup{Graph}(\sigma_{\vert\Sigma})\!\subseteq\! E$. As $E$ is a graph, we get 
$E\! =\! G_{\sigma_{\vert\Sigma}}$ and therefore $V\! =\!\Sigma$. Thus 
$(\Sigma ,G_{\sigma_{\vert\Sigma}})\preceq^i_c(K,G)$ and 
$(\Sigma ,G_{\sigma_{\vert\Sigma}})$ is $\preceq^i_c$-minimal in $\mathfrak{G}_2$ and in the class of closed graphs on a 0DMC space with CCN at least three.\medskip

 Assume now that $\xi\!\geq\!\omega$ is of the form $\eta\! +\! 3$. Using ideas similar to those in [Ce-Da-To-Wy], we now provide a two-sided subshift with Cantor-Bendixson rank of the form $\xi$ having the desired minimality property. The first step of our construction is inspired by [Ce-Da-To-Wy, Theorem 4.6]. Fix an infinite countable ordinal $\eta$, and a closed countable subset $P$ of $2^\omega$ with Cantor-Bendixson rank $\eta\! +\! 1$, which exists by [K, 6.13]. The following fact is known. However, we include a proof for completeness.\medskip

\noindent\emph{Claim 4.}\it\ \label{isoldense} Let $P$ be a countable Polish space. Then $P\!\setminus\! P'$ is dense in $P$.\rm\medskip
 
 Indeed, we set, for $x\!\in\! P'$, $O_x\! :=\! X\!\setminus\!\{ x\}$, so that $O_x$ is a dense open subset of $X$. Moreover, $P\!\setminus\! P'\! =\!\bigcap_{x\in P'}~O_x$ is a $G_\delta$ subset of the countable space $P$. It remains to apply the Baire category theorem (see [K, 8.4]).\hfill{$\diamond$}\medskip
 
 We enumerate $P\!\setminus\! P'\! :=\!\{\gamma_j\mid j\!\in\!\omega\}$ and set, for $j\!\in\!\omega$, 
$$w_j\! :=\! 1^2{^\frown}_{k<j}~\big( (01)^{\gamma_{(j)_0}(k)+\Sigma_{i<k}~(\gamma_{(j)_0}(i)+1)}1^2\big)\mbox{,}$$
and $\delta_\infty\! :=\! (01)^\infty\!\cdot\! 1{^\frown}_{j\in\omega}~(w_j(01)^{j+1})$. Similarly, we define 
$\Phi_0\! :\! 2^\omega\!\rightarrow\! 2^\omega$ by 
$$\Phi_0(\gamma )\! :=\! {^\frown}_{k\in\omega}~\big( (01)^{\gamma (k)+\Sigma_{i<k}~(\gamma (i)+1)}1^2\big) .$$

 We also define 
$\Phi\! :\! 2^\omega\!\rightarrow\! 2^\mathbb{Z}$ by $\Phi (\gamma )\! :=\! (01)^\infty\!\cdot\! 1^2\Phi_0(\gamma )$. We then set 
$$Q\! :=\!\mbox{Orb}_\sigma\big( (01)^\infty\!\cdot\! (01)^\infty\big)\cup
\mbox{Orb}_\sigma\big( (01)^\infty\!\cdot\! 1^2(01)^\infty\big)\cup
\bigcup_{\gamma\in P}~\mbox{Orb}_\sigma\big(\Phi (\gamma )\big)$$ 
and $\Sigma\! :=\! Q\cup\mbox{Orb}_\sigma (\delta_\infty )$. Note that $P$, $Q$ and $\Sigma$ are countable, and 
$\sigma [Q]\! =\! Q$, $\sigma [\Sigma ]\! =\!\Sigma$.\medskip

\noindent\emph{Claim 5.}\it\ $\Sigma$ is a countable two-sided subshift with Cantor-Bendixson rank $\xi$, and 
$\mbox{Orb}_\sigma (\delta_\infty )$ is dense in $\Sigma$.\rm\medskip

 Indeed, as in [Ce-Da-To-Wy], we check that $Q$ is closed and has Cantor-Bendixson rank $\xi$. Note first that 
$\mbox{Orb}_\sigma\big( (01)^\infty\!\cdot\! (01)^\infty\big)\! =\!\{ (01)^\infty\!\cdot\! (01)^\infty ,(10)^\infty\!\cdot\! (10)^\infty\}$ is closed, as well as 
$$\mbox{Orb}_\sigma\big( (01)^\infty\!\cdot\! (01)^\infty\big)\cup\mbox{Orb}_\sigma\big( (01)^\infty\!\cdot\! 1^2(01)^\infty\big) .$$ 
Let $(\delta_n)_{n\in\omega}$ be a sequence of elements of $Q$ converging to $\delta\!\in\! 2^\mathbb{Z}$. By the previous remark, we may assume that the $\delta_n$'s are in 
$\bigcup_{\gamma\in P}~\mbox{Orb}_\sigma\big(\Phi (\gamma )\big)$. This gives  
$(\beta_n)_{n\in\omega}\!\in\! P^\omega$ and $(k_n)_{n\in\omega}\!\in\!\mathbb{Z}^\omega$ with 
${\delta_n\! =\!\sigma^{k_n}\big(\Phi (\beta_n)\big)}$, and we may assume that $(\beta_n)_{n\in\omega}$ converges to 
$\beta\!\in\! P$. If we may assume that $(k_n)_{n\in\omega}$ is constant, then 
$\delta\! =\! \sigma^{k_0}\big(\Phi (\beta )\big)\!\in\! Q$ by continuity of $\Phi_0$ and $\Phi$. Otherwise, we may assume that 
$(k_n)_{n\in\omega}$ is either strictly increasing, or strictly decreasing. In the latter case, 
$\delta\!\in\!\mbox{Orb}_\sigma\big( (01)^\infty\!\cdot\! (01)^\infty\big)\!\subseteq\! Q$. So we may assume that 
$(k_n)_{n\in\omega}$ is strictly increasing and $k_0\!\geq\! 2$. Note that\medskip
 
\leftline{$S\! :=\!\big\{ {^\frown}_{k<l}~\big( (01)^{m_k}1^2\big)(01)^\infty\mid l\!\in\!\omega\wedge
\forall k\! <\! l\! -\! 1~~m_k\! <\! m_{k+1}\big\} ~\cup$}\smallskip

\rightline{$\big\{ {^\frown}_{k\in\omega}~\big( (01)^{m_k}1^2\big)\mid
\forall k\!\in\!\omega~~m_k\! <\! m_{k+1}\big\}$}\medskip

\noindent and $\sigma [S]$ are closed subsets of $2^\omega$, as well as $C\! :=\! S\cup\sigma [S]$. We define, for $\beta\!\in\! 2^\mathbb{Z}$, $\beta^*\!\in\! 2^\omega$ by $\beta^*(i)\! :=\!\beta (i)$ if $i\!\in\!\omega$. As $\Phi_0$ takes values in $S$, $\sigma^2[S]\! =\! S$ and $k_n\!\geq\! 2$, $\delta_n^*\!\in\! C$ for each $n$, and $\delta^*\!\in\! C$. If $\delta^*$ contains at most one $1^3$, then $\delta\!\in\!\mbox{Orb}_\sigma\big( (01)^\infty\!\cdot\! (01)^\infty\big)\cup
\mbox{Orb}_\sigma\big( (01)^\infty\!\cdot\! 1^2(01)^\infty\big)\!\subseteq\! Q$. So we may assume that there are 
$n_0,n_1\!\geq\! 1$ with $1^2(01)^{n_0}1^2(01)^{n_1}1^2\!\subseteq\!\delta^*$, 
$(01)^{n_0}1^2(01)^{n_1}1^2\!\subseteq\!\delta^*$ or $(10)^{n_0-1}11^2(01)^{n_1}1^2\!\subseteq\!\delta^*$. As 
$(\delta_n^*)_{n\in\omega}$ converges to $\delta^*$, we may assume that is also the case for the $\delta_n^*$'s. Note then that, just after this initial segment, $\delta_n^*$ can have at most $n_1\! +\! 2$ blocks $01$ before having a block $1^2$, by definition of $\Phi_0$. This implies that $\delta^*$ is of the form $1^2{^\frown}_{k\in\omega}~\big( (01)^{n_k}1^2\big)$, 
${^\frown}_{k\in\omega}~\big( (01)^{n_k}1^2\big)$ or $(10)^{n_0-1}11^2{^\frown}_{k\geq 1}~\big( (01)^{n_k}1^2\big)$, with $n_{k+1}\! +\! 1\!\leq\! n_{k+2}\!\leq\! n_{k+1}\! +\! 2$. So we may assume that either 
$1^2{^\frown}_{k\leq n}~\big( (01)^{n_k}1^2\big)\!\subseteq\!\delta_n^*$ for each $n$, 
${^\frown}_{k\leq n}~\big( (01)^{n_k}1^2\big)\!\subseteq\!\delta_n^*$ for each $n$, or 
$(10)^{n_0-1}11^2{^\frown}_{1\leq k\leq n}~\big( (01)^{n_k}1^2\big)\!\subseteq\!\delta_n^*$ for each $n$. Note that 
$\Phi_0(\beta_n)$ has an initial segment of the form either 
${^\frown}_{k\leq l}~\big( (01)^{m_k}1^2\big) {^\frown}_{k\leq n}~\big( (01)^{n_k}1^2\big)$, or 
${^\frown}_{k\leq l}~\big( (01)^{m_k}1^2\big) 0(10)^{n_0-1}11^2{^\frown}_{1\leq k\leq n}~\big( (01)^{n_k}1^2\big)$. As 
$m_0\! <\! m_1\! <\!\cdots\! <\! m_l\! <\! n_0$ in both cases, there are only finitely many possible values for the block 
${^\frown}_{k\leq l}~\big( (01)^{m_k}1^2\big)$. So we may assume that this block does not depend on $n$. Note then that 
$\delta_n^*\! =\!\sigma^{k_n-2}\big(\Phi_0(\beta_n)\big)$. This implies that we may assume that $(k_n)_{n\in\omega}$ is constant, which is not the case. This shows that $Q$ is closed.\medskip

 In order to prove that $Q$ has Cantor-Bendixson rank $\xi$, we introduce the notion of the \emph{rank of a point}. If $\mathbb{X}\!\in\!\{\omega ,\mathbb{Z}\}$, $P$ is a countable compact subset of $2^\mathbb{X}$ and $\delta\!\in\! P$, then the rank $rk_P(\delta )$ of $\delta$ in $P$ is the least ordinal $\alpha$ such that $\delta\!\notin\! P^{\alpha +1}$. Under this definition, the Cantor-Bendixson rank of $P$ is $\mbox{sup}\{ rk_P(\delta )\! +\! 1\mid\delta\!\in\! P\}$ (see [Ce-Da-To-Wy, Section 2]).\medskip

 Note that $P_1\! :=\!\big\{\Phi (\gamma )\mid\gamma\!\in\! P\big\}\!\subseteq\! Q$, which implies that 
$rk_{P_1}(\delta )\!\leq\! rk_Q(\delta )$ if $\delta\!\in\! P_1$. By [Ce-Da-To-Wy, Lemma 3.3],  
$rk_Q\big(\sigma^k(\delta )\big)\! =\! rk_Q(\delta )$ if $k\!\in\!\mathbb{Z}$ and $\delta\!\in\! Q$.\medskip

 Thus  
$rk_Q\big(\sigma^k(\delta )\big)\geq\! rk_{P_1}(\delta )$ if $k\!\in\!\mathbb{Z}$ and $\delta\!\in\! P_1$. It follows that 
$rk_Q\big( (01)^\infty\!\cdot\! 1^2(01)^\infty\big)\!\geq\!\eta\! +\! 1$ and hence 
$rk_Q\big( (01)^\infty\!\cdot\! (01)^\infty\big)\!\geq\!\eta\! +\! 2$.\medskip

 For the other direction, note that the map $(k,\gamma )\!\mapsto\!\sigma^k\big(\Phi (\gamma )\big)$ is injective. We now prove by induction on $\rho\! :=\! rk_{P_1}(\delta )$ that $rk_Q(\delta )\! =\!\rho$ if $\delta\!\in\! P_1$. If 
$\rho\! =\! 0$, then $\delta$ is isolated in $P_1$. If $\delta$ is not isolated in $Q$, then there is an injective  sequence 
$(\delta_n)_{n\in\omega}$ of elements of $Q$ converging to $\delta$. The discussion above provides $k_0\!\in\!\mathbb{Z}$ and $(\beta_n)_{n\in\omega}\!\in\! P^\omega$ converging to $\beta\!\in\! P$ with 
$\delta_n\! =\!\sigma^{k_0}\big(\Phi (\beta_n)\big)$. Thus $\delta\! =\!\sigma^{k_0}\big(\Phi (\beta )\big)$. The injectivity property shows that $k_0\! =\! 0$, so that $\delta_n\!\in\! P_1$, contradicting the fact that $\delta$ is isolated in $P_1$. Suppose now that our claim holds for all ordinals strictly below $\rho$, and that $rk_Q(\delta )\! >\!\rho$. Then we can find an injective sequence $(\delta_n)_{n\in\omega}$ of elements of $Q$ converging to $\delta$. The discussion above provides 
$k_0\!\in\!\mathbb{Z}$ and $(\beta_n)_{n\in\omega}\!\in\! P^\omega$ converging to $\beta\!\in\! P$ with 
$\delta_n\! =\!\sigma^{k_0}\big(\Phi (\beta_n)\big)$. The injectivity property shows that $k_0\! =\! 0$, so that 
$\delta_n\!\in\! P_1$. As $rk_{P_1}(\delta )\! =\!\rho$, we may assume that $rk_{P_1}(\delta_n)\! <\!\rho$, so that 
$rk_Q(\delta_n)\! <\!\rho$ by induction assumption. This contradicts the fact that 
$\delta_n\!\in\! Q$. This implies that $rk_Q(\delta )\!\leq\!\eta$ if 
$\delta\!\in\! \bigcup_{\gamma\in P}~\mbox{Orb}_\sigma\big(\Phi (\gamma )\big)$. It follows that 
$rk_Q(\delta )\!\leq\!\eta\! +\! 1$ if 
$\delta\!\in\!\mbox{Orb}_\sigma\big( (01)^\infty\!\cdot\! 1^2(01)^\infty\big)$, and $rk_Q(\delta )\!\leq\!\eta\! +\! 2$ if 
$\delta\!\in\!\mbox{Orb}_\sigma\big( (01)^\infty\!\cdot\! (01)^\infty\big)$. This proves that $Q$ has Cantor-Bendixson rank 
$\xi$.\medskip

 Let us prove that $\Sigma$ is closed. Assume that $(\delta_n)_{n\in\omega}$ is a sequence of elements of $\Sigma$ converging to $\delta\!\in\! 2^\mathbb{Z}$. We may assume that the $\delta_n$'s are in $\mbox{Orb}_\sigma (\delta_\infty )$ since $Q$ is closed. This gives $k_n\!\in\!\mathbb{Z}$ with the property that $\delta_n\! =\!\sigma^{k_n}(\delta_\infty )$. If 
$(k_n)_{n\in\omega}$ has a constant subsequence, then $\delta\!\in\!\mbox{Orb}_\sigma (\delta_\infty )$ and we are done. 
So we may assume that $(k_n)_{n\in\omega}$ tends to $\pm\infty$. If $(k_n)_{n\in\omega}$ tends to $-\infty$, then 
$\delta\!\in\!\mbox{Orb}_\sigma\big( (01)^\infty\!\cdot\! (01)^\infty\big)$, and we are done. So we may write 
${k_n\! =\! 1\! +\!\Sigma_{j<l_n}~(\vert w_j\vert\! +\! 2j\! +\! 2)\! +\! j_n}$, where $j_n\! <\!\vert w_{l_n}\vert\! +\! 2l_n\! +\! 2$, and 
$(l_n)_{n\in\omega}$ tends to $\infty$. In particular, the distance between $\delta_n$ and 
$C\! :=\!\{\beta\!\in\! 2^\mathbb{Z}\mid\exists\varepsilon\!\in\! 2~~\forall i\!\in\!\mathbb{Z}~~\beta (2i\! +\!\varepsilon )\! =\! 1\}$ tends to zero as $n$ tends to infinity, so that $\delta$ is in the closed set $C$. As above we may assume that 
$\delta$ is not in 
$\mbox{Orb}_\sigma\big( (01)^\infty\!\cdot\! (01)^\infty\big)\cup\mbox{Orb}_\sigma\big( (01)^\infty\!\cdot\! 1^2(01)^\infty\big)$, which gives $m$ minimal such that $1^2(01)^m1^2$ is a finite subword of $\delta$, and we may also assume that 
$1^2(01)^m1^2\! =\!\delta_{[0,2m+3]}\! =\! {\delta_n}_{[0,2m+3]}$. This implies that we may assume that $j_n\! =\! 0$. Note then that 
$$\begin{array}{ll}
\delta_n^*\!\!\!\!
& =\! {^\frown}_{j\geq l_n}~(w_j(01)^{j+1})\cr
& =\! {^\frown}_{j\geq l_n}~
\Big( 1^2\Big( {^\frown}_{k<j}~\big( (01)^{\gamma_{(j)_0}(k)+\Sigma_{i<k}~(\gamma_{(j)_0}(i)+1)}1^2\big)\Big) (01)^{j+1}\Big)\cr
& \supseteq\! 1^2\Big( {^\frown}_{k<l_n}~\big( (01)^{\gamma_{(l_n)_0}(k)+\Sigma_{i<k}~(\gamma_{(l_n)_0}(i)+1)}1^2\big)\Big) 
(01)^{l_n+1}
\end{array}$$
converges to $\delta^*$. This implies that $(\gamma_{(l_n)_0})_{n\in\omega}$ converges to some $\gamma\!\in\! P$, and that 
$\delta\! =\!\Phi (\gamma )\!\in\! Q$, showing that $\Sigma$ is closed.\medskip

 This discussion above shows that if $\delta\!\in\! Q$, then $rk_Q(\delta )\! =\! 0$ if and only if there is 
$\gamma\!\in\! P\!\setminus\! P'$ with $\delta\!\in\!\mbox{Orb}_\sigma\big(\Phi (\gamma )\big)$. In particular, 
$rk_\sigma (\delta )\!\geq\! 1$ if $\delta$ is not of this form. If now $j\!\in\!\omega$, then $\Phi (\gamma_j )$ is in 
$\overline{\mbox{Orb}(\delta_\infty )}$, showing that $Q$ is a subset of the Cantor-Bendixson derivative $\Sigma'$ of $\Sigma$. As $\Sigma$ is countable, it has an isolated point which has to be in $\mbox{Orb}_\sigma (\delta_\infty )$, and 
$\Sigma'\! =\! Q$ since $\sigma$ is a homeomorphism. This implies that $\Sigma^{k+1}\! =\! Q^k$ for each natural number $k$, 
and $\Sigma^\theta\! =\! Q^\theta$ if $\theta$ is infinite, so that $\Sigma$ has Cantor-Bendixson rank $\xi$ since $\eta$ is infinite.\medskip
 
 The density assertion comes from the previous discussion and Claim 4.\hfill{$\diamond$}\medskip

 We then partially argue as in the finite case. Note that 
$$\big( (01)^\infty\!\cdot\! (01)^\infty ,(01)^\infty\!\cdot\! (01)^\infty\big)\! =\!\mbox{lim}_{p\rightarrow\infty}~
\big(\sigma^{-2p}(\delta_\infty ),
\sigma^{1+\Sigma_{j<2p}~(\vert w_j\vert +2j+2)+\vert w_{2p}\vert +2p}(\delta_\infty )\big)$$
is in $\overline{\bigcup_{q\in\omega}~G_{\sigma_{\vert\Sigma}}^{2q+1}}$, so that 
$\chi_c(\Sigma ,G_{\sigma_{\vert\Sigma}})\! =\! 3$.

\vfill\eject

\noindent\emph{Claim 6.}\it\ Let $V\!\subseteq\!\Sigma$, and 
$E\!\subseteq\! G_{\sigma_{\vert\Sigma}}\cap V^2$ be a graph on $V$ such that $(V,E)$ has CCN three. Then 
$\big(\alpha ,\sigma (\alpha )\big)\!\in\! E$ if $\alpha\!\in\!\mbox{Orb}_\sigma (\delta_\infty )$.\rm\medskip

 Indeed, we argue by contradiction. Let $k\!\in\!\mathbb{Z}$ with $\alpha\! =\!\sigma^k(\delta_\infty )$. As in the proof of Claim 3, we may assume that $k\! =\! 0$. We set 
$$C\! :=\! ([0]_{0}\cup [1^20]_{0}\cup [01^40]_{-2})\cap V\mbox{,}$$ 
so that $C$ is a clopen subset of $V$. Noted that the odd coordinates of the $w_j$'s are 1, and that the $w_j$'s have a stricly positive even length. The definition of the $w_j$'s and $\delta_\infty$ then imply that 0 can only be an even coordinate of the elements appearing in the definition of $Q$ (which is the union of the orbits of these elements), a negative even coordinate of 
$\delta_\infty$ or a positive odd coordinate of $\delta_\infty$. Moreover, $01^k0$ can be of the form $\alpha_{[a,b]}$ with 
$\alpha$ appearing in the definition of $Q$ only if $k\!\in\!\{ 1,3\}$ and $a$ is even, and of the form ${\delta_\infty}_{[a,b]}$ only if $k\! =\! 1$ and $a\!\leq\! -4$ is even, $k\!\in\!\{ 1,3\}$ and $a\!\geq\! 3$ is odd, or $k\! =\! 4$ and $a\! =\! -2$. In particular, 
$$\begin{array}{ll}
& \ldots,\sigma^{-2}(\alpha ),\alpha ,\sigma^{2}(\alpha ),\ldots\!\in\! C\mbox{,}\cr
& \ldots,\sigma^{-3}(\alpha ),\sigma^{-1}(\alpha ),\sigma (\alpha ),\ldots\!\notin\! C\cr
\end{array}$$ 
if $\alpha$ appears in the definition of $Q$, and   
$$\begin{array}{ll}
& \ldots,\sigma^{-4}(\delta_\infty ),\sigma^{-2}(\delta_\infty ),\delta_\infty ,\sigma (\delta_\infty ),
\sigma^{3}(\delta_\infty ),\sigma^{5}(\delta_\infty ),\ldots\!\in\! C\mbox{,}\cr
& \ldots,\sigma^{-5}(\delta_\infty ),\sigma^{-3}(\delta_\infty ),\sigma^{-1}(\delta_\infty ),
\sigma^{2}(\delta_\infty ),\sigma^{4}(\delta_\infty ),\sigma^{6}(\delta_\infty ),\ldots\!\notin\! C.\cr
\end{array}$$ 
This implies that $E\cap\big( C^2\cup (V\!\setminus\! C)^2\big)\! =\!\emptyset$, which contradicts the fact that $(V,E)$ has CCN three.\hfill{$\diamond$}\medskip

 We now conclude as in the finite case.\medskip
 
\noindent (b) Let $Q\! :=\! (q_j)_{j\in\omega}\!\in\!\omega^\omega$ converging to infinity. We set 
$\alpha_0\! :=\! (01)^\infty\!\cdot\! (01)^\infty$, $\alpha_1\! :=\! (01)^\infty\!\cdot\! 1^2(01)^\infty$ and 
$\beta_Q\! :=\! (01)^\infty\!\cdot\! 1{^\frown}_{j\in\omega}~\big((01)^{q_j}1^2 \big)$. This defines  
$\Sigma_Q\! =\!\bigcup_{m\leq 1}~\mbox{Orb}_\sigma (\alpha_m)\cup\mbox{Orb}_\sigma (\beta_Q)$.\medskip

 We first essentially argue as in the finite case when $n\! =\! 1$ to check the individual properties of $\Sigma_Q$.\medskip

\noindent\emph{Claim 7.}\it\ (a) Let $1\!\leq\! i\!\leq\! 2$. Then 
$\Sigma_Q^{(i)}\! =\!\overline{\mbox{Orb}_\sigma (\alpha_{2-i})}$.\smallskip

(b) $\Sigma_Q\! =\!\overline{\mbox{Orb}_\sigma (\beta_Q)}$.\rm\medskip

 Indeed, let us prove that $\alpha_1\!\in\!\overline{\mbox{Orb}_\sigma (\beta_Q)}$. Fix a natural number $N$. If 
$j$ is large enough, then $q_j\!\geq\! N$, so that $(01)^N1^2(01)^N\! =\! {\alpha_1}_{[-2N,2N+1]}\! =\! 
{\beta_Q}_{[1+(\Sigma_{n<j}~(2q_n+2))+2(q_j-N),(\Sigma_{n\leq j}~(2q_n+2))+2N]}$. This implies that 
$\alpha_1\!\in\!\overline{\mbox{Orb}_\sigma (\beta_Q)}$, as desired.\medskip

 Assume then that $(k_p)_{p\in\omega}\!\in\!\mathbb{Z}^\omega$ and $\big(\sigma^{k_p}(\beta_Q)\big)_{p\in\omega}$ converges to $\alpha\!\in\! 2^\mathbb{Z}$. We want to see that $\alpha\!\in\!\Sigma_Q$. If $(k_p)_{p\in\omega}$ tends to $\infty$, then we may write 
$k_p\! =\! 1\! +\!\big(\Sigma_{j<l_p}~(2q_j\! +\! 2)\big)\! +\! j_p$, where $j_p\! <\! 2q_{l_p}\! +\! 2$, and 
$(l_p)_{p\in\omega}$ tends to $\infty$. If $(j_p)_{p\in\omega}$ and $(2q_{l_p}\! +\! 2\! -\! j_p)_{p\in\omega}$ also tend to 
$\infty$, then $\alpha\!\in\!\overline{\mbox{Orb}_\sigma (\alpha_0)}\!\subseteq\!\Sigma_Q$, as desired. Otherwise, we may assume that $(j_p)_{p\in\omega}$ or $(2q_{l_p}\! +\! 2\! -\! j_p)_{p\in\omega}$ is constant, so that 
$\alpha$ is in $\mbox{Orb}_\sigma (\alpha_1)\!\subseteq\!\Sigma_Q$.\hfill{$\diamond$}\medskip

\noindent\emph{Claim 8.}\it\ $\Sigma_Q$ has Cantor-Bendixson rank 3.\rm\medskip

 Indeed, as all the odd coordinates of the $\alpha_m$'s are 1, and we can find an even coordinate of $\beta_Q$ and an odd coordinate of $\beta_Q$ which are 0, so that $\alpha_m\!\notin\!\mbox{Orb}_\sigma (\beta_Q)$ if $m\!\leq\! 1$. It remains to note that $\alpha_0\!\notin\!\mbox{Orb}_\sigma (\alpha_1)$ to get $\Sigma^{(2)}_Q\! =\! (\Sigma^{(1)}_Q)'$ and, inductively, that $\Sigma^{(2)}_Q$ is the $2$nd iterated Cantor-Bendixson derivative of $\Sigma_Q$.\hfill{$\diamond$}
 
\vfill\eject

 Note then that 
$$(\alpha_0,\alpha_0)\! =\!\mbox{lim}_{p\rightarrow\infty}~
\big(\sigma^{-2p}(\beta_Q),\sigma^{1+(\Sigma_{j<p}~(2q_j+2))+2\lceil\frac{q_p}{2}\rceil}(\beta_Q)\big)\!\in\!
\overline{\bigcup_{q\in\omega}~G_{\sigma_{\vert\Sigma_Q}}^{2q+1}}\mbox{,}$$ 
so that $\chi_c(\Sigma_Q,G_{\sigma_{\vert\Sigma_Q}})\! =\! 3$.\medskip

\noindent\emph{Claim 9.}\it\ Let $V\!\subseteq\!\Sigma_Q$, and 
$E\!\subseteq\! G_{\sigma_{\vert\Sigma_Q}}\cap V^2$ be a graph on $V$ such that $(V,E)$ has CCN three. Then 
$\big(\alpha ,\sigma (\alpha )\big)\!\in\! E$ if $\alpha\!\in\!\mbox{Orb}_\sigma (\beta_Q)$.\rm\medskip

 Indeed, let $k\!\in\!\mathbb{Z}$ with $\alpha\! =\!\sigma^k(\beta_Q)$. As in the proof of Claim 3, we may assume that 
$k\! =\! 0$. Fix $j_0\!\in\!\omega$ such that $q_j\!\geq\! 1$ if $j\!\geq\! j_0$. We set 
$$C\! :=\! ([0]_0\cup\bigcup_{j\leq q_{j_0}}~[01^{2j+1}0]_{-2}\cup\bigcup_{j\leq q_{j_0}}~[01^{2j+2}0]_{-2})\cap V\mbox{,}$$ 
so that $C$ is a clopen subset of $V$. We already noticed that the odd coordinates of the $\alpha_m$'s are 1. The definition of the $\alpha_m$'s and $\beta_Q$ then imply that 0 can only be an even coordinate of $\alpha_m$, a negative even coordinate of $\beta_Q$ or a positive odd coordinate of $\beta_Q$. Moreover, $01^k0$ can be of the form ${\alpha_m}_{[a,b]}$ only if 
$k\!\in\!\{ 1,3\}$ and $a$ is even, and of the form ${\beta_Q}_{[a,b]}$ only if $k\! =\! 1$ and $a\!\leq\! -4$ is even, 
$k\!\in\!\{ 1,3\}$ and $a\!\geq\! 1$ is odd, or $2\!\leq\! k\!\leq\! 2j_0\! +\! 2$ is even and $a\! =\! -2$. This implies that 
 $E\cap\big( C^2\cup (V\!\setminus\! C)^2\big)\! =\!\emptyset$ since 
$$\begin{array}{ll}
& \ldots,\sigma^{-2}(\alpha_m),\alpha_m,\sigma^2(\alpha_m),\ldots\!\in\! C\mbox{,}\cr
& \ldots,\sigma^{-3}(\alpha_m),\sigma^{-1}(\alpha_m),\sigma(\alpha_m),\ldots\!\notin\! C\cr
\end{array}$$ 
if $m\!\leq\! 1$ and 
$$\begin{array}{ll}
& \ldots,\sigma^{-4}(\beta_Q),\sigma^{-2}(\beta_Q),\beta_Q,\sigma (\beta_Q),
\sigma^3(\beta_Q),\sigma^5(\beta_Q),\ldots\!\in\! C\mbox{,}\cr
& \ldots,\sigma^{-5}(\beta_Q),\sigma^{-3}(\beta_Q),\sigma^{-1}(\beta_Q),\sigma^2(\beta_Q),
\sigma^4(\beta_Q),\sigma^6(\beta_Q),\ldots\!\notin\! C\mbox{,}\cr
\end{array}$$ 
which contradicts the fact that $(V,E)$ has CCN three again.\hfill{$\diamond$}\medskip

 We conclude as in the finite case to get the individual properties of $\Sigma_Q$. We now provide a family of size continuum of countable subshifts $(\Sigma_{Q^\nu})_{\nu\in 2^\omega}$. Let $(p_n)_{n\in\omega}$ be the sequence of prime numbers. We set, for $\nu\!\in\! 2^\omega$ and $n\!\in\!\omega$, $q^\nu_0\! :=\! 0$ and 
$q^\nu_{n+1}\! :=\! p_0^{\nu (0)+2}\!\cdots\! p_n^{\nu (n)+2}\! -\! 1$, which defines $Q^\nu\!\in\!\omega^\omega$ converging to infinity.\medskip

 Let us check that the family $\big( (\Sigma_{Q^\nu},G_{\sigma_{\vert\Sigma_{Q^\nu}}})\big)_{\nu\in 2^\omega}$ is a $\preceq^i_c$-antichain. Assume, towards a contradiction, that $\nu\!\not=\!\nu'$ and 
$(\Sigma_{Q^\nu},G_{\sigma_{\vert\Sigma_{Q^\nu}}})\preceq^i_c(\Sigma_{Q^{\nu'}},G_{\sigma_{\vert\Sigma_{Q^{\nu'}}}})$ with witness $\varphi$. Let $m_0$ be minimal with $\nu (m_0)\!\not=\!\nu'(m_0)$. By minimality of 
$(\Sigma_{Q^{\nu'}},G_{\sigma_{\vert\Sigma_{Q^{\nu'}}}})$, we may assume that $\nu (m_0)\! <\!\nu'(m_0)$. If 
$x\!\in\!\Sigma_{Q^{\nu}}$, then 
$\big( x,\sigma (x)\big)\!\in\! G_{\sigma_{\vert\Sigma_{Q^{\nu}}}}$ since $\sigma_{\vert\Sigma_{Q^{\nu}}}$ is fixed point free. Thus $\Big(\varphi (x),\varphi\big(\sigma (x)\big)\Big)\!\in\! G_{\sigma_{\vert\Sigma_{Q^{\nu'}}}}$ and 
$\varphi\big(\sigma (x)\big)\! =\!\sigma^{\pm 1}\big(\varphi (x)\big)$. We choose 
$x\!\in\!\Sigma_{Q^{\nu}}\!\setminus\!\mbox{Orb}_\sigma (\alpha_0)$, so that $\mbox{Orb}_\sigma (x)$ is infinite, as well as 
$\mbox{Orb}_\sigma\big(\varphi (x)\big)\!\supseteq\!\varphi [\mbox{Orb}_\sigma (x)]$, and $\sigma^2(x)\!\not=\! x$. In particular, $\sigma^2_{\vert\mbox{Orb}_\sigma (\varphi (x))}$ is fixed point free. We apply Lemma 5.6 and its proof to 
$V\! =\! X\! =\!\Sigma_{Q^{\nu}}$, $f\! =\!\sigma_{\vert\Sigma_{Q^{\nu}}}$, $I\! =\!\mbox{Orb}_\sigma (x)$, 
$W\! =\! Y\! =\!\Sigma_{Q^{\nu'}}$, $g\! =\!\sigma_{\vert\Sigma_{Q^{\nu'}}}$ and $\varphi$. The proof of Lemma 5.6 shows that $P\cap M\! =\!\emptyset$, either ${\varphi\big(\sigma (z)\big)\! =\!\sigma\big(\varphi (z)\big)}$ for each 
$z\!\in\!\mbox{Orb}_\sigma (x)$ or $\varphi\big(\sigma (z)\big)\! =\!\sigma^{-1}\big(\varphi (z)\big)$ for each 
$z\!\in\!\mbox{Orb}_\sigma (x)$, and $\varphi [\mbox{Orb}_\sigma (x)]\! =\!\mbox{Orb}_\sigma\big(\varphi (x)\big)$. In particular, $\varphi [\mbox{Orb}_\sigma (\alpha_1)],\varphi [\mbox{Orb}_\sigma (\beta_{Q^{\nu}})]$ are disjoint infinite orbits in $\Sigma_{Q^{\nu'}}$, so they are $\mbox{Orb}_\sigma (\alpha_1),\mbox{Orb}_\sigma (\beta_{Q^{\nu'}})$.

\vfill\eject

 As $\mbox{Orb}_\sigma (\beta_{Q^{\nu'}})$ is dense in $\Sigma_{Q^{\nu'}}$, the compact set $\varphi [\Sigma_{Q^{\nu}}]$ is $\Sigma_{Q^{\nu'}}$, so that $\varphi$ is a homeomorphism from 
$\Sigma_{Q^{\nu}}$ onto $\Sigma_{Q^{\nu'}}$. Moreover, $\varphi$ is a witness for the fact that 
$\sigma_{\vert\Sigma_{Q^{\nu}}}$ and $\sigma_{\vert\Sigma_{Q^{\nu'}}}$ are flip-conjugate, by density of 
$\mbox{Orb}_\sigma (\beta_{Q^{\nu}})$ in $\Sigma_{Q^{\nu}}$. In particular, 
$\varphi [\Sigma'_{Q^{\nu}}]\! =\!\Sigma'_{Q^{\nu'}}$ and $\varphi [\Sigma''_{Q^{\nu}}]\! =\!\Sigma''_{Q^{\nu'}}$, so that 
$\varphi [\mbox{Orb}_\sigma (\alpha_0)]\! =\!\mbox{Orb}_\sigma (\alpha_0)$, 
$\varphi [\mbox{Orb}_\sigma (\alpha_1)]\! =\!\mbox{Orb}_\sigma (\alpha_1)$ and 
$\varphi [\mbox{Orb}_\sigma (\beta_{Q^{\nu}})]\! =\!\mbox{Orb}_\sigma (\beta_{Q^{\nu'}})$.  
This gives $n_0,n_1\!\in\!\mathbb{Z}$ with $\varphi (\alpha_1)\! =\!\sigma^{n_1}(\alpha_1)$ and 
$\varphi (\beta_{Q^{\nu}})\! =\!\sigma^{n_0}(\beta_{Q^{\nu'}})$. We then set, for  
$r\!\in\!\omega$, 
$${K^\nu_r\! :=\! 1\! +\! (\Sigma_{j<r}~(2q_j\! +\! 2))\! +\! 2q_r}.$$  
Note that the sequence $\big(\sigma^{K^\nu_r}(\beta_{Q^{\nu}})\big)_{r\in\omega}$ converges to $\alpha_1$, so that 
$\Big(\varphi\big(\sigma^{K^\nu_r}(\beta_{Q^{\nu}})\big)\Big)_{r\in\omega}$ converges to 
$\varphi (\alpha_1)\! =\!\sigma^{n_1}(\alpha_1)$. As 
$\varphi\big(\sigma^{K^\nu_r}(\beta_{Q^{\nu}})\big)\! =\!\sigma^{n_0\pm K^\nu_r}(\beta_{Q^{\nu'}})$, 
this implies that $\big(\sigma^{n_0-n_1\pm K^\nu_r}(\beta_{Q^{\nu'}})\big)_{r\in\omega}$ converges to $\alpha_1$. As 
$(K^\nu_r)_{r\in\omega}$ is strictly increasing, this implies that $\sigma_{\vert\Sigma_{Q^{\nu}}}$ and 
$\sigma_{\vert\Sigma_{Q^{\nu'}}}$ are conjugate and $\big(\sigma^{n_0-n_1+K^\nu_r}(\beta_{Q^{\nu'}})\big)_{r\in\omega}$ converges to $\alpha_1$. In particular, 
$${\sigma^{n_0-n_1+K^\nu_r}(\beta_{Q^{\nu'}})}_{[-2,2]}\! =\!
\big(\beta_{Q^{\nu'}}(n_0\! -\! n_1\! +\! K^\nu_r\! -\! 2),\cdots ,\beta_{Q^{\nu'}}(n_0\! -\! n_1\! +\! K^\nu_r\! +\! 2)\big)
\! =\! {\alpha_1}_{[-2,2]}\! =\! 01^30$$ 
if $r$ is large enough. Using similar notation, this implies that 
$n_0\! -\! n_1\! +\! K^\nu_r\!\in\!\{ K^{\nu'}_m\mid m\!\in\!\omega\}$ if $r$ is large enough. In particular, this gives, for $r$ large enough, $m\! <\! M\!\in\!\omega$ with $n_0\! -\! n_1\! +\! K^\nu_r\! =\! K^{\nu'}_m$ and 
$n_0\! -\! n_1\! +\! K^\nu_{r+1}\! =\! K^{\nu'}_M$. Thus 
$K^\nu_{r+1}\! -\! K^\nu_r\! =\! 2q^\nu_{r+1}\! +\! 2\! =\!\Sigma_{m<j\leq M}~(2q^{\nu'}_j\! +\! 2)$ and 
$$p_0^{\nu (0)+2}\!\cdots\! p_r^{\nu (r)+2}\! =\! q^\nu_{r+1}\! +\! 1\! =\!\Sigma_{m\leq n<M}~(q^{\nu'}_{n+1}\! +\! 1)
\! =\!\Sigma_{m\leq n<M}~(p_0^{\nu'(0)+2}\!\cdots\! p_n^{\nu'(n)+2}).$$ 
We may assume that $r$ is large enough to ensure that $r,m\!\geq\! m_0$, which implies that $p_{m_0}^{\nu'(m_0)+2}$ divides $p_0^{\nu (0)+2}\!\cdots\! p_r^{\nu (r)+2}$, which cannot be since $\nu (m_0)\! <\!\nu'(m_0)$.\hfill{$\square$}\medskip

 Note that Theorem \ref{CB1intro}(b) provides a version of Corollary \ref{controt} (and thus Theorem \ref{eantichmino}) for  countable subshifts (which are not necessary minimal). By minimality and for cardinality reasons, the examples provided by Theorem \ref{CB1intro}(b) are $\preceq^i_c$-incompatible in the class of closed graphs on a 0DMC space with CCN at least three with the examples given by Proposition \ref{odoshift}. 
 
\section{$\!\!\!\!\!\!$ Homeomorphisms of a 0DMS space}\label{MSPolishsub}\indent

 In this section, we prove Theorem \ref{absmincompmain}, among other things.\medskip 

\noindent\emph{Remarks.}\ (a) Let $X$ be a 0DMC space, and $f\! :\! X\!\rightarrow\! X$ be a homeomorphism such that $(X,G_f)$ has CCN at least three. Lemma 7.10 says that if $f$ is minimal, then $(X,G_f)$ is $\preceq^i_c$-minimal in the class of closed graphs on a 0DMC space with CCN at least three. Theorem \ref{CB1intro} shows that the converse is not true since the finite orbit $\mbox{Orb}_\sigma(\alpha_0)$ is not dense in the infinite countable space $K_0$.\medskip

\noindent (b) Theorem \ref{CB1intro} also provides $(K_0,G_{h_0})\!\in\!\mathfrak{G}_2$ for which it is not possible to find $(K,G_f)\!\in\!\mathfrak{G}_2$ with $f$ minimal and $(K,G_f)\preceq^i_c(K_0,G_{h_0})$. In particular, $K_0$ contains no subshift $\Sigma$ such that $(\Sigma ,\sigma_{\vert\Sigma})$ is minimal and has CCN at least three, even if it contains 
$\mbox{Orb}_\sigma(\alpha_0)$ and $(\mbox{Orb}_\sigma(\alpha_0),\sigma_{\vert\text{Orb}_\sigma (\alpha_0)})$ is minimal.\medskip

 We will see that $(K_0,G_{h_0})$ is no more minimal in 0DMS (or 0DP) spaces. Let $\mathcal{T}$ be the set of finer 0DMS topologies $\tau$ on $K_0$ such that $\big( (K_0,\tau ),G_{h_0}\big)$ has CCN at least three and $h_0$ is a homeomorphism of $(K_0,\tau )$.
 
\begin{lem} \label{topol} Let $S$ be a 0DMS (resp., 0DP) space, $f$ be a homeomorphism of $S$ with the properties that $\chi_c(S,G_f)\!\geq\! 3$ and $(S,G_f)\preceq^i_c(K_0,G_{h_0})$. Then there is a finer 0DMS (resp., 0DP) topology $\tau$ in $\mathcal{T}$ with the property that $\big( (K_0,\tau ),G_{h_0}\big)\preceq^i_c(S,G_f)$.\end{lem}
 
\noindent\emph{Proof.}\ By Theorem \ref{CB1intro}, $(K_0,G_{h_0})$ has CCN three, and thus $(S,G_f)$ has CCN three too. Therefore the set $F_1$ of fixed points of $f$ is a clopen subset of $S$, and 
$\chi_c(S\!\setminus\! F_1,G_f\cap (S\!\setminus\! F_1)^2)\! =\! 3$ by Corollary 7.3. So we may assume that $f$ is fixed point free. Let $\varphi$ be a witness for the fact that $(S,G_f)\preceq^i_c(K_0,G_{h_0})$. We define  
$V\! :=\!\varphi [S]$ and ${E\! :=\! (\varphi\!\times\!\varphi )[G_f]}$. The finer topology is 
$$\tau\! :=\!\{ O\!\subseteq\! K_0\mid\varphi^{-1}(O)\!\in\!\boraone (S)\} .$$ 
Note that $\varphi\! :\! S\!\rightarrow\! (V,\tau )$ is a homeomorphism, so that $(V,\tau )$ is a 0DMS (resp., 0DP) space. As 
$\varphi$ is a witness for the fact that $(S,G_f)\preceq^i_c\big( (K_0,\tau ),G_{h_0}\big)$, 
$\chi_c\big( (K_0,\tau ),G_{h_0}\big)\!\geq\! 3$. Also, $\varphi^{-1}$ is a witness for the fact that 
$\big( (V,\tau ),E\big)\preceq^i_c(S,G_f)$.\medskip

 Let us prove that $V\! =\! K_0$ and $E\! =\! G_{h_0}$. As $\varphi$ is a witness for the fact that $(S,G_f)\preceq^i_c(V,E)$, $\chi_c(V,E)\!\geq\! 3$. As $\mbox{Orb}_{h_0}(\beta_0)$ is discrete, there is 
$\varepsilon\!\in\! 2$ with $\sigma^\varepsilon (\alpha_0)\!\in\! V$, which gives 
$x\!\in\! S$ with ${\varphi (x)\! =\!\sigma^\varepsilon (\alpha_0)}$. As $f$ is fixed point free, $f(x)\!\not=\! x$, which implies that 
$\big( x,f(x)\big)\!\in\! G_f$, $\Big(\varphi (x),\varphi\big( f(x)\big)\Big)\!\in\! G_{h_0}$, and 
$\sigma^{1-\varepsilon}(\alpha_0)\! =\!\varphi\big( f(x)\big)\!\in\! V$. This implies that 
$\mbox{Orb}_\sigma (\alpha_0)\!\subseteq\! V$ and 
$\big\{\big(\alpha_0,\sigma (\alpha_0)\big) ,\big(\sigma (\alpha_0),\alpha_0\big)\big\}\!\subseteq\! E$. As 
$\chi_c(V,E)\!\geq\! 3$ again, there is $\alpha\!\in\!\mbox{Orb}_{h_0}(\beta_0)\cap V$, which gives $y\!\in\! S$ with 
$\varphi (y)\! =\!\alpha$.\medskip

 Let us check that $f^2(y)\!\not=\! y$. We argue by contradiction. As just above, 
$\varphi\big( f(y)\big)\! =\! h_0\big(\varphi (y)\big)$ or $\varphi\big( f(y)\big)\! =\! h_0^{-1}\big(\varphi (y)\big)$. Assume that 
$\varphi\big( f(y)\big)\! =\! h_0^{-1}\big(\varphi (y)\big)$, the other case being similar. Then 
$\big(\alpha ,h_0(\alpha )\big)\!\notin\! E$ since $f(y)\! =\! f^{-1}(y)$ is sent to $h_0^{-1}(\alpha )$ by 
$\varphi$ and $h_0^2$ is fixed point free on $\mbox{Orb}_{h_0}(\beta_0)$. This contradicts Claim 3 in the proof of Theorem \ref{CB1intro}.\medskip

 As $f^2(y)\!\not=\! y$, either $\varphi\big( f(z)\big)\! =\! h_0\big(\varphi (z)\big)$ for each $z\!\in\!\mbox{Orb}_f(y)$, or $\varphi\big( f(z)\big)\! =\! h_0^{-1}\big(\varphi (z)\big)$ for each $z\!\in\!\mbox{Orb}_f(y)$, by Lemma 5.6. Lemma 5.2 then implies that $\varphi [\mbox{Orb}_f(y)]\! =\!\mbox{Orb}_{h_0}(\alpha )\! =\!\mbox{Orb}_{h_0}(\beta_0)$. Thus $V\! =\! K_0$ and $E\! =\! G_{h_0}$. In particular, $\big( (K_0,\tau ),G_{h_0}\big)\preceq^i_c(S,G_f)$ and $\tau$ is a finer 0DMS (resp., 0DP) topology on $K_0$.\medskip

 The previous discussion shows that $S\! =\!\mbox{Orb}_f(x)\cup\mbox{Orb}_f(y)$ by injectivity of $\varphi$, and 
$\varphi\!\circ\! f\! =\! h_0^{\pm 1}\!\circ\!\varphi$. Thus $h_0\! =\!\varphi\!\circ\! f^{\pm 1}\!\circ\!\varphi^{-1}$ is a homeomorphism of $(K_0,\tau )$.\hfill{$\square$}\medskip

\noindent\emph{Notation.}\ Lemma \ref{infinf} provides a sequence $(S_q)_{q\in\omega}$ of pairwise disjoint infinite subsets of $\omega$ such that, for any $l\!\in\!\omega$, $p\!\not=\! q$, $3\!\cdot\! 2^{l+1}\! <\! r\!\in\! S_p$ and 
$s\!\in\! S_q$, $\vert r\! -\! s\vert\! >\! 2^{l+1}$. We enumerate, here again, $S_q$ in a strictly increasing way by $\{ r^q_j\mid j\!\in\!\omega\}$.  We set, for $l\!\in\!\omega$, $j_l\! :=\! 3\!\cdot\! 2^{l+1}\! +\! 1$ and, for $A\!\subseteq\!\omega$, 
$$N_A\! :=\!\{ r^q_{j_l}\! +\! r\mid (q\! >\! 0\Rightarrow\! q\! -\! 1\!\in\! A)\wedge l\!\in\!\omega\wedge -l\!\leq\! r\!\leq\! l\} .$$ 
Note that $N_A\!\subseteq\!\omega$. Let $\tau$ be a finer 0DMS topology on $K_0$, and $\mathcal{B}^\tau$ be a countable basis, made up of clopen sets and closed under finite intersections, for $\tau$. We set, for $A\!\subseteq\!\omega$,\medskip

\leftline{$\mathcal{B}^\tau_A\! :=\!\mathcal{B}^\tau\cup\Big\{ C\cap(\{ (01)^\infty\!\cdot\! (01)^\infty\}\cup
\bigcup_{n\in\bigcap_{-p_0\leq r\leq q_0}~(N_A+r)\cap\omega}~\{ (01)^\infty 1(01)^n\!\cdot\! (01)^\infty\} ~\cup$}\smallskip

\rightline{$\bigcup_{n\in\bigcap_{-p_1\leq r\leq q_1}~(N_A+r)\cap\omega}~\{ (01)^\infty\!\cdot\! (01)^{n+1}1(01)^\infty\} )\mid C\!\in\!\mathcal{B}^\tau\wedge p_0,q_0,p_1,q_1\!\in\!\omega\Big\} ~\cup$}\medskip

\rightline{$\Big\{ C\cap(\{ (10)^\infty\!\cdot\! (10)^\infty\}\cup
\bigcup_{n\in\bigcap_{-p_0\leq r\leq q_0}~(N_A+r)\cap\omega}~\{ (01)^\infty 1(01)^n0\!\cdot\! (10)^\infty\} ~\cup$}\smallskip

\rightline{$\bigcup_{n\in\bigcap_{-p_1\leq r\leq q_1}~(N_A+r)\cap\omega}~\{ (10)^\infty\!\cdot\! 1(01)^n1(01)^\infty\} )\mid C\!\in\!\mathcal{B}^\tau\wedge p_0,q_0,p_1,q_1\!\in\!\omega\Big\} .$}

\begin{lem} \label{tauA} Let $A\!\subseteq\!\omega$ and $\tau$ be in $\mathcal{T}$. Then\smallskip

\noindent (a) $\mathcal{B}^\tau_A$ is the basis for a 0DMS (0DP if $\tau$ is) topology $t^\tau_A$ in 
$\mathcal{T}$ finer than $\tau$,\smallskip

\noindent (b) the sequences $\big( (01)^\infty 1(01)^{r^0_{j_l}}\!\cdot\! (01)^\infty\big)_{l\in\omega}$ and 
$\big( (01)^\infty\!\cdot\! (01)^{r^0_{j_l}+1}1(01)^\infty\big)_{l\in\omega}$ are 
$t^\tau_A$-converging to $(01)^\infty\!\cdot\! (01)^\infty$, as well as 
$\big( (01)^\infty 1(01)^{r^{q+1}_{j_l}}\!\cdot\! (01)^\infty\big)_{l\in\omega}$ and 
$\big( (01)^\infty\!\cdot\! (01)^{r^{q+1}_{j_l}+1}1(01)^\infty\big)_{l\in\omega}$ if $q\!\in\! A$.\end{lem}

\noindent\emph{Proof.}\ Note that $\mathcal{B}^\tau_A$ contains $\mathcal{B}^\tau$, is countable, and closed under finite intersections. Thus $\mathcal{B}^\tau_A$ is the basis for a second countable topology $t^\tau_A$ finer than $\tau$. Moreover, 
$\mathcal{B}^\tau_A$ is made up of closed subsets of $K_0$, so that $t^\tau_A$ is zero-dimensional. This shows that 
$(K_0,t^\tau_A)$ is a 0DP space (see [K, 13.2 and 13.3]) if $(K_0,\tau )$ is.\medskip

 We apply the definitions of $N_A$ and $\mathcal{B}^\tau_A$ to see that (b) holds. Note then that 
$$\big( (01)^\infty\!\cdot\! (01)^{r^0_{j_l}+1}1(01)^\infty ,(01)^\infty 1(01)^{r^0_{j_l}}\!\cdot\! (01)^\infty\big)\! =\!\big( h_0^{-2r^0_{j_l}-2}(\beta^0_0),\! h_0^{2r^0_{j_l}+1}(\beta^0_0)\big)
\!\in\! G_{h_0}^{4r^0_{j_l}+3}.$$ 
Thus $\big( (K_0,t^\tau_A),G_{h_0}\big)$ has CCN at least three by (b) and Lemma 3.3.1. It remains to note that 
$\big( (K_0,t^\tau_A),G_{h_0}\big)\!\preceq_c\! (K_0,G_{h_0})$ to see that the CCN is three. Note that\medskip

\leftline{$h_0[\{ (10)^\infty\!\cdot\! (10)^\infty\}\cup\bigcup_{n\in\bigcap_{-p_0\leq r\leq q_0}~(N_A+r)\cap\omega}~
\{ (01)^\infty 1(01)^n0\!\cdot\! (10)^\infty\} ~\cup$}\smallskip

\rightline{$\bigcup_{n\in\bigcap_{-p_1\leq r\leq q_1}~(N_A+r)\cap\omega}~\{ (10)^\infty\!\cdot\! 1(01)^n1(01)^\infty\} ]=$}\smallskip

\rightline{$\{ (01)^\infty\!\cdot\! (01)^\infty\}\cup
\bigcup_{n\in\bigcap_{-p_0\leq r\leq q_0}~(N_A+r)\cap\omega}~\{ (01)^\infty 1(01)^{n+1}\!\cdot\! (01)^\infty\} ~\cup$}\smallskip

\rightline{$\bigcup_{n\in\bigcap_{-p_1\leq r\leq q_1}~(N_A+r)\cap\omega}~\{ (01)^\infty\!\cdot\! (01)^n1(01)^\infty\}$}\medskip

\noindent is equal, up to a finite open discrete set, to\medskip

\rightline{$\{ (01)^\infty\!\cdot\! (01)^\infty\}\cup
\bigcup_{n\in\bigcap_{-p_0+1\leq r\leq q_0+1}~(N_A+r)\cap\omega}~\{ (01)^\infty 1(01)^n\!\cdot\! (01)^\infty\} ~\cup$}\smallskip

\rightline{$\bigcup_{n\in\bigcap_{-p_1-1\leq r\leq q_1-1}~(N_A+r)\cap\omega}~\{ (01)^\infty\!\cdot\! (01)^{n+1}1(01)^\infty\}\mbox{,}$}\medskip

\noindent and thus $t^\tau_A$-open. We argue similarly with $h_0^{-1}$ instead of $h_0$, or with 
$(01)^\infty\!\cdot\! (01)^\infty$ instead of $(10)^\infty\!\cdot\! (10)^\infty$, to see that $h_0$ is a homeomorphism as desired.\hfill{$\square$}

\begin{lem} \label{compari} Let $A,B\!\subseteq\!\omega$ with $A\!\not\subseteq\! B$, $\tau$ be in $\mathcal{T}$, and $\tau'$ be in $\mathcal{T}$ finer than $\tau$. Then 
$$\big( (K_0,t^{\tau'}_A),G_{h_0}\big)\not\preceq^i_c\big( (K_0,t^\tau_B),G_{h_0}\big) .$$
\end{lem}

\noindent\emph{Proof.}\ Towards a contradiction, suppose that there is 
$\varphi\! :\! (K_0,t^{\tau'}_A)\!\rightarrow\! (K_0,t^\tau_B)$. We define  
$V\! :=\!\varphi [K_0]$ and ${E\! :=\! (\varphi\!\times\!\varphi )[G_{h_0}]}$. If $z\!\in\!\mbox{Orb}_{h_0}(\beta_0)$, then 
$\big( z,h_0(z)\big)\!\in\! G_{h_0}$ since $h_0$ is fixed point free, so that 
$\varphi\big( h_0(z)\big)\! =\! h_0^{\pm 1}\big(\varphi (z)\big)$. Thus 
$\varphi [\mbox{Orb}_{h_0}(\beta_0)]\!\subseteq\!\mbox{Orb}_{h_0}(\beta_0)$ or 
$\varphi [\mbox{Orb}_{h_0}(\beta_0)]\!\subseteq\!\mbox{Orb}_{h_0}(\alpha_0)$. By injectivity of 
$\varphi$, the first case holds. As $h_0^2$ is fixed point free on 
$\mbox{Orb}_{h_0}(\beta_0)$, either ${\varphi\big( h_0(z)\big)\! =\! h_0\big(\varphi (z)\big)}$ for each 
$z\!\in\!\mbox{Orb}_{h_0}(\beta_0)$, or ${\varphi\big( h_0(z)\big)\! =\! h_0^{-1}\big(\varphi (z)\big)}$ for each 
$z\!\in\!\mbox{Orb}_{h_0}(\beta_0)$, by Lemma 5.6. Lemma 5.2 then implies that 
$\varphi [\mbox{Orb}_{h_0}(\beta_0)]\! =\!\mbox{Orb}_{h_0}\big(\varphi (\beta_0)\big)$. This provides $n\!\in\!\mathbb{Z}$ with 
$\varphi (\beta_0)\! =\! h_0^n(\beta_0)$. So either 
$\varphi\big( h_0^i(\beta_0)\big)\! =\! h_0^{n+i}(\beta_0)$ for each $i\!\in\!\mathbb{Z}$, or 
$\varphi\big( h_0^i(\beta_0)\big)\! =\! h_0^{n-i}(\beta_0)$ for each $i\!\in\!\mathbb{Z}$, by Lemma 5.2 again. In particular, if 
$q\!\in\! A\!\setminus\! B$ and $l\!\in\!\omega$, then 
$$\varphi\big( (01)^\infty 1(01)^{r^{q+1}_{j_l}}\!\cdot\! (01)^\infty\big)\! =\!
\varphi\big( h_0^{2r^{q+1}_{j_l}+1}(\beta_0)\big)\!\in\!
\{ h_0^{n+2r^{q+1}_{j_l}+1}(\beta_0),h_0^{n-2r^{q+1}_{j_l}-1}(\beta_0)\}$$ 
is $t^\tau_B$-converging to $\varphi\big( (01)^\infty\!\cdot\! (01)^\infty\big)$, in 
$\{ (01)^\infty\!\cdot\! (01)^\infty ,(10)^\infty\!\cdot\! (10)^\infty\}$ by injectivity of $\varphi$ since 
$\varphi [\mbox{Orb}_{h_0}(\beta_0)]\! =\!\mbox{Orb}_{h_0}(\beta_0)$.

\vfill\eject

 Assume, for example, that 
$\varphi\big( (01)^\infty\!\cdot\! (01)^\infty )\! =\! (10)^\infty\!\cdot\! (10)^\infty$, the other case being similar. Then\medskip
 
\leftline{$\{ (10)^\infty\!\cdot\! (10)^\infty\}\cup\bigcup_{n\in N_B}~\{ (01)^\infty 1(01)^n0\!\cdot\! (10)^\infty\}\cup
\bigcup_{n\in N_B}~\{ (10)^\infty\!\cdot\! 1(01)^n1(01)^\infty\}\! =\!$}\smallskip

\rightline{$\{ (10)^\infty\!\cdot\! (10)^\infty\}\cup
\bigcup_{n\in N_B}~\{ h_0^{2n+2}(\beta^0_0),h_0^{-2n-1}(\beta^0_0)\}$}\medskip
 
\noindent is a $t^\tau_B$-neighborhood of $\varphi\big( (01)^\infty\!\cdot\! (01)^\infty\big)$, and thus contains $\varphi\big( (01)^\infty 1(01)^{r^{q+1}_{j_l}}\!\cdot\! (01)^\infty\big)$ if $l$ is large enough.\medskip

 The injectivity of $h_0$ provides, for example, $p\!\in\! B$, 
$k\!\in\!\omega$ and $-k\!\leq\! r\!\leq\! k$ with the property that $n\! -\! 2r^{q+1}_{j_l}\! -\! 1\! =\! -2(r^{p+1}_{j_k}\! +\! r)\! -\! 1$ if $l$ is large enough (the other case is similar). This implies that  
$2\vert r^{q+1}_{j_l}\! -\! r^{p+1}_{j_k}\vert\! =\!\vert n\! +\! 2r\vert\!\leq\!\vert n\vert\! +\! 2k$, 
$k$ goes to $\infty$ as $l$ goes to $\infty$, and $\vert r^{q+1}_{j_l}\! -\! r^{p+1}_{j_k}\vert\!\leq\! 2k$ if $l$ is large enough. As 
$q\! +\! 1\!\not=\! p\! +\! 1$, this contradicts  our choice of $j_k$ since 
$\vert r^{q+1}_{j_l}\! -\! r^{p+1}_{j_k}\vert\! >\! 2^{k+1}\! >\! 2k$.\hfill{$\square$}\medskip
 
\noindent\emph{Proof of Theorem \ref{absmincompmain}.}\ (a) We apply Theorem \ref{CB1intro}.\medskip

\noindent (b) Lemma \ref{topol} provides a finer 0DMS (resp., 0DP) topology $\tau$ on $K_0$ such that $h_0$ is a homeomorphism of $(K_0,\tau )$, $\chi_c\big( (K_0,\tau ),G_{h_0}\big)\!\geq\! 3$ and 
$\big( (K_0,\tau ),G_{h_0}\big)\preceq^i_c(S,G_f)$. We apply Lemmas \ref{tauA} and \ref{compari} to 
$A\! :=\!\{ 2q\mid q\!\in\!\omega\}$, ${B\! :=\!\{ 2q\! +\! 1\mid q\!\in\!\omega\}}$, and $\tau$. As $t^\tau_A,t^\tau_B$ are finer than $\tau$, 
$\big( (K_0,t^\tau_A),G_{h_0}\big) ,\big( (K_0,t^\tau_B),G_{h_0}\big)\preceq^i_c\big( (K_0,\tau ),G_{h_0}\big)$. As 
$A\!\not\subseteq\! B$, 
$$\big( (K_0,t^\tau_A),G_{h_0}\big)\not\preceq^i_c\big( (K_0,t^\tau_B),G_{h_0}\big)$$ 
and thus $\big( (K_0,\tau ),G_{h_0}\big)\not\preceq^i_c\big( (K_0,t^\tau_B),G_{h_0}\big)$. This proves that 
$\big( (K_0,t^\tau_B),G_{h_0}\big)$ is strictly $\preceq^i_c$-below $\big( (K_0,\tau ),G_{h_0}\big)$, and also $(S,G_f)$.\medskip

 Assume now, towards a contradiction, that there is a $\preceq^i_c$-antichain basis $\mathfrak{B}$ for a class in the statement. By Theorem \ref{CB1intro}, $(K_0,G_{h_0})$ is in this class, which gives $(S,G_f)\!\in\!\mathfrak{B}$ with the property that $(S,G_f)\preceq^i_c(K_0,G_{h_0})$. The first part of this theorem provides a finer 0DMS (resp., 0DP) topology 
$\tau'$ on $K_0$ with the properties that $h_0$ is a homeomorphism of $(K_0,\tau')$, $\big( (K_0,\tau'),G_{h_0}\big)$ has CCN at least three, and $\big( (K_0,\tau'),G_{h_0}\big)$ is strictly $\preceq^i_c$-below $(S,G_f)$, which is the desired contradiction.\medskip

 For the size of the basis, towards a contradiction, suppose we can find $\kappa\! <\! 2^{\aleph_0}$ and a basis 
$(B_\gamma )_{\gamma <\kappa}$ for our class. Let $(p_n)_{n\in\omega}$ be the sequence of prime numbers. We define, for each $\alpha\!\in\! 2^\omega$, $S_\alpha\!\subseteq\!\omega$ by 
$$S_\alpha\! :=\!\{ p_0^{\alpha (0)+1}\ldots p_n^{\alpha (n)+1}\mid n\!\in\!\omega\} .$$  
Note that $S_\alpha\!\subseteq\!\omega$ is infinite, and $S_\alpha\cap S_\beta$ is finite if $\alpha\!\not=\!\beta$. By Theorem \ref{CB1intro}, we can apply Lemma \ref{tauA} to any $S_\alpha$ and the topology $\tau_0$ on $K_0$, so that 
$t^{\tau_0}_{S_\alpha}$ is a 0DP topology in $\mathcal{T}$. As $\kappa\! <\! 2^{\aleph_0}$, we can find 
$\gamma\! <\!\kappa$ and $\alpha\!\not=\!\beta$ with 
$B_\gamma\preceq^i_c\big( (K_0,t^{\tau_0}_{S_\alpha}),G_{h_0}\big)$, 
$\big( (K_0,t^{\tau_0}_{S_\beta}),G_{h_0}\big)\preceq^i_c(K_0,G_{h_0})$. Lemma \ref{topol} provides a finer 0DMS (resp., 0DP) topology $\tau$ in $\mathcal{T}$ with $\big( (K_0,\tau ),G_{h_0}\big)\preceq^i_cB_\gamma$. We can apply again Lemma \ref{tauA}, to $A\! :=\!\omega$ and $\tau$, so that $t^\tau_\omega$ is a 0DMS (resp., 0DP) topology in $\mathcal{T}$ finer than $\tau$, so that 
$\big( (K_0,t^\tau_\omega ),G_{h_0}\big)\preceq^i_c\big( (K_0,\tau ),G_{h_0}\big)\preceq^i_c
\big( (K_0,t^{\tau_0}_{S_\alpha}),G_{h_0}\big)$. As $S_\beta$ is infinite, $\omega\not\subseteq\! S_\alpha$. Lemma \ref{compari} then implies that 
$\big( (K_0,t^\tau_\omega ),G_{h_0}\big)\not\preceq^i_c\big( (K_0,t^{\tau_0}_{S_\alpha}),G_{h_0}\big)$, which is the desired contradiction.\hfill{$\square$}\medskip

 We turn to the proof of Theorem \ref{embed}(a). In fact, we prove something stronger since it is possible to consider always the same graph, with different underlying 0DP spaces.\medskip

\noindent\emph{Proof of Theorem \ref{embed}(a).}\ We define, for $A\!\subseteq\!\omega$,\medskip

\leftline{$K_A\! :=\!
\{\varepsilon^{n+1}(\varepsilon\! +\! 1\mbox{ mod }4)1^\infty\mid\varepsilon\!\in\! 4\wedge n\!\in\!\omega\}\cup
\{\varepsilon^{n+1}(\varepsilon\! -\! 1\mbox{ mod }4)0^\infty\mid\varepsilon\!\in\! 4\wedge n\!\in\!\omega\} ~\cup$}\smallskip

\rightline{$\{\varepsilon^\infty\mid\varepsilon\!\in\! 5\}\cup\bigcup_{n\in A}~
\{\varepsilon^{n+2}(\varepsilon\! +\! 1\mbox{ mod }4)s2^\infty\mid\varepsilon\!\in\! 4\wedge s\!\in\! 2^{n+1}\}$.}

\vfill\eject\medskip

 We then enumerate $2^{n+1}\! :=\!\{ s^{n+1}_i\mid i\! <\! 2^{n+1}\}\!\subseteq\! 2^{<\omega}$ and define a function  
$h_A\! :\! K_A\!\rightarrow\! K_A$ by $h_A(\varepsilon^\infty )\! :=\! (\varepsilon\! +\! 1\mbox{ mod }4)^\infty$ if 
$\varepsilon\!\in\! 4$, $h_A(4^\infty )\! :=\! 01^\infty$, ${h_A(320^\infty )\! :=\! 4^\infty}$,
$$h_A(\varepsilon^{n+1}(\varepsilon\! +\! 1\mbox{ mod }4)1^\infty )\! :=\! 
(\varepsilon\! +\! 1\mbox{ mod }4)^{n+1}(\varepsilon\! +\! 2\mbox{ mod }4)1^\infty$$ 
if $\varepsilon\!\not=\! 3$, $h_A(3^{n+1}01^\infty )\! :=\! 0^{n+2}1^\infty$, 
$h_A(\varepsilon^{n+1}(\varepsilon\! -\! 1\mbox{ mod }4)0^\infty )\! :=\! 
(\varepsilon\! +\! 1\mbox{ mod }4)^{n+1}\varepsilon 0^\infty$ if $\varepsilon\!\not=\! 3$, as well as 
$h_A(3^{n+2}20^\infty )\! :=\! 0^{n+1}30^\infty$ on the one hand,
$$h_A(\varepsilon^{n+2}(\varepsilon\! +\! 1\mbox{ mod }4)s^{n+1}_i2^\infty )\! :=\! 
(\varepsilon\! +\! 1\mbox{ mod }4)^{n+2}(\varepsilon\! +\! 2\mbox{ mod }4)s^{n+1}_i2^\infty$$ 
if $\varepsilon\!\not=\! 3$, $h_A(3^{n+2}0s^{n+1}_i2^\infty )\! :=\! 0^{n+2}1s^{n+1}_{i+1\text{ mod }2^{n+1}}2^\infty$ on the other hand. In other words, $K_A$ is the union of the $h_A$-orbit $\{ 0^\infty ,1^\infty ,2^\infty ,3^\infty\}$, the orbit 
$\{ 4^\infty\}\cup\{\varepsilon^{n+1}(\varepsilon\!\pm\! 1\mbox{ mod }4)\eta^\infty\mid\varepsilon\!\in\! 2\wedge n\!\in\!\omega\wedge\eta\!\in\! 2\}$ in the style of the infinite $h_0^0$-orbit of $K_0^0$, and even cycles given by the elements of $A$. The beginning of the proof of Theorem \ref{CB1intro} shows that $K_A$ is a countable 0DMC space, $h_A$ is a homeomorphism of $K_A$, and $(K_A,G_{h_A})$ has CCN three. If $A\!\subseteq\! B$, then $K_A\!\subseteq\! K_B$, 
$G_{h_A}\!\subseteq\! G_{h_B}$, which implies that $(K_A,G_{h_A})\preceq^i_c(K_B,G_{h_B})$. If $A\!\not\subseteq\! B$, then let $n\!\in\! A\!\setminus\! B$. Note that $(K_A,G_{h_A})\not\preceq^i_c(K_B,G_{h_B})$ because $(K_A,G_{h_A})$ has a cycle of length $4\!\cdot\! 2^{n+1}$ and $(K_B,G_{h_B})$ does not.\hfill{$\square$}\medskip

 As announced in the introduction, one can check that the $\sigma_{\vert\Sigma}$'s appearing in the statement of Theorem \ref{CB1intro} are expansive, which leaves the question of infinite Cantor Bendixson ranks uncertain.\medskip

\noindent\emph{Remark.}\ By [K, 33.B], the set $\mathcal{K}_{\aleph_0}(2^\omega )$ of countable compact subsets of 
$2^\omega$ is $\ca$-complete. By [K, 34.18(3)], the Cantor-Bendixson rank $\vert\!\cdot\!\vert_{CB}$ is a co-analytic rank on $\mathcal{K}_{\aleph_0}(2^\omega )$. Thus the map $r\! :\! (X,f)\!\mapsto\!\vert X\vert_{CB}$ defines a co-analytic rank on 
$\mathcal{P}\cap\!\big(\mathcal{K}_{\aleph_0}(2^\omega )\!\times\!\mathcal{H}(2^\omega)\big)$ ($\mathcal{P}$ was defined before Theorem \ref{desc}). By [K, 35.23], $\vert\!\cdot\!\vert_{CB}$ has to be unbounded. Thus $r$ is unbounded, which implies that the co-analytic set $\mathcal{P}\cap\!\big(\mathcal{K}_{\aleph_0}(2^\omega )\!\times\!\mathcal{H}(2^\omega)\big)$ is not Borel, by [K, 35.23] again. By [K, 34.2], $(X,f)\!\mapsto\!\vert X\vert_{CB}$ also defines a co-analytic rank on 
$\mathcal{O}_2^{\aleph_0}\! :=\!
\mathcal{O}_2\cap\!\big(\mathcal{K}_{\aleph_0}(2^\omega )\!\times\!\mathcal{H}(2^\omega)\big)$. Theorem \ref{CB1intro} implies that the co-analytic subset $\mathcal{O}_2^{\aleph_0}$ of $\mathcal{P}$ is not Borel. This set is in fact $\ca$-complete. Indeed, define $\{ s_i\mid i\!\in\! 3\}\! :=\!\{ 0^2,10,1^2\}$, $f_0\!\in\!\mathcal{H}(2^\omega )$ by $f_0(0\alpha )\! :=\! 0\alpha$, 
$f_0(1s\alpha )\! :=\! 1s\alpha$ if $s\!\in\! 2^2\!\setminus\!\{ s_i\mid i\!\in\! 3\}$, and 
$f_0(1s_i\alpha )\! :=\! 1s_{i+1\text{ mod }3}\alpha$. The map 
$X\!\mapsto\! (\{ 0\alpha\mid\alpha\!\in\! X\}\cup\{ 1s_i0^\infty\mid i\!\in\! 3\},f_0)$ is a continuous reduction of 
$\mathcal{K}_{\aleph_0}(2^\omega )$ to $\mathcal{O}_2^{\aleph_0}$, by [K, 4.29] and since $G_{f_0}$ contains the 3-cycle $\{ 1s_i0^\infty\mid i\!\in\! 3\}$.

\section{$\!\!\!\!\!\!$ The classes $\mathfrak{G}_\kappa$}\label{Polishsub}\indent

 We consider, for $\kappa\!\leq\! 3$,\medskip
 
\noindent - the class $\mathfrak{G}_\kappa$ of graphs induced by a homeomorphism of a 0DMC space with CCN strictly bigger than $\kappa$,\medskip

\noindent - the class $\mathfrak{H}_\kappa$ of homeomorphisms of $2^\omega$ whose induced graph has CCN strictly bigger than $\kappa$.\medskip

\noindent\emph{Proof of Theorem \ref{Ckappa}.}\ (a) The CCN is strictly bigger than $0$ if and only if the space is not empty.\medskip

\noindent (b) The CCN is strictly bigger than $1$ if and only if the graph is not empty.\medskip

\noindent (c) By Theorem \ref{eantichmino}, any $\preceq_c^i$-basis for $\mathfrak{G}_2$ must have size continuum.\medskip

\noindent (d) We apply Proposition \ref{LZ1}(b).\medskip

 For the well-foundedness, fix $n\!\in\!\omega$. We enumerate the set of finite binary sequences $2^{n+1}$ by 
$\{ s_i\mid i\! <\! 2^{n+1}\}$, so that $N_{0^n1}\! =\! N_{0^n10}\cup\bigcup_{i<2^{n+1}}~N_{0^n1^2s_i}$.

\vfill\eject

 We consider the map $c_n$ on $N_{0^n1}$ defined by $c_n(0^n10\alpha )\! :=\! 0^n1^2s_0\alpha$, 
$c_n(0^n1^2s_i\alpha )\! :=\! 0^n1^2s_{i+1}\alpha$ if $i\! <\! 2^{n+1}\! -\! 1$, and 
$c_n(0^n1^2s_{2^{n+1}-1}\alpha )\! :=\! 0^n10\alpha$. Note that $c_n$ is a homeomorphism,  
$c_n^{2^{n+1}+1}\! =\!\mbox{Id}_{N_{0^n1}}$, and $c_n^i(\beta )\!\not=\!\beta$ if $i\!\leq\! 2^{n+1}$. Moreover, the function 
$h_p\! :\! 2^\omega\!\rightarrow\! 2^\omega$, defined by $h_p(\beta )\! :=\!\beta$ if $\beta$ is in 
$\{ 0^\infty\}\cup\bigcup_{n<p}~N_{0^n1}$ and $h_p(\beta )\! :=\! c_n(\beta )$ if $n\!\geq\! p$ and $\beta\!\in\! N_{0^n1}$, is a homeomorphism whose set $\{ 0^\infty\}\cup\bigcup_{n<p}~N_{0^n1}$ of fixed points is not open. By Proposition \ref{fp}, 
$\chi_c(2^\omega ,G_{h_p})\! =\! 2^{\aleph_0}$, so that $(2^\omega ,G_{h_p})$ is in all the $\mathfrak{G}_\kappa$'s. As $G_{h_{p+1}}\!\subseteq\! G_{h_p}$, $(2^\omega ,G_{h_{p+1}})\preceq^i_c(2^\omega ,G_{h_p})$. As $G_{h_p}$ contains a cycle of length $2^{p+1}+1$ and all the cycles in $G_{h_{p+1}}$ have length at least $2^{p+2}+1$, 
$(2^\omega ,G_{h_p})\not\preceq_c(2^\omega ,G_{h_{p+1}})$.\medskip

 Theorem \ref{eantichmino} provides $\preceq_c$-antichains of size continuum in $\mathfrak{G}_\kappa$ if 
$\kappa\!\leq\! 2$. For $\mathfrak{G}_3$, we use again the $c_n$'s. Let $(S_\alpha )_{\alpha\in 2^\omega}$ be as in the proof of Theorem \ref{antichaf}. The map $h_\alpha\! :\! 2^\omega\!\rightarrow\! 2^\omega$, defined by 
$h_\alpha (\beta )\! :=\!\beta$ if $\beta\!\in\!\{ 0^\infty\}\cup\bigcup_{n\notin S_\alpha}~N_{0^n1}$ and 
$h_\alpha (\beta )\! :=\! c_n(\beta )$ if $n\!\in\! S_\alpha$ and $\beta\!\in\! N_{0^n1}$, is a homeomorphism whose set of fixed points $\{ 0^\infty\}\cup\bigcup_{n\notin S_\alpha}~N_{0^n1}$ is not open. By Proposition \ref{fp}, 
$\chi_c(2^\omega ,G_{h_\alpha})\! =\! 2^{\aleph_0}$, so that $(2^\omega ,G_{h_\alpha})$ is in all the $\mathfrak{G}_\kappa$'s. If $\alpha\!\not=\!\beta$, then there is $n\!\in\! S_\alpha\!\setminus\! S_\beta$, so that $G_{h_\alpha}$ contains a cycle of length $2^{n+1}+1$, which is not the case of $G_{h_\beta}$. Thus 
$(2^\omega ,G_{h_\alpha})\not\preceq^i_c(2^\omega ,G_{h_\beta})$.\hfill{$\square$}\medskip

 We can also evaluate the descriptive complexity of $\mathfrak{G}_\kappa$ and $\mathfrak{H}_\kappa$. Let $\mathcal{H}(2^\omega )$ be the set of homeomorphisms of $2^\omega$. We equip $\mathcal{H}(2^\omega )$ with the topology whose basic open sets are of the form 
$${O_{U_1,\ldots ,U_n,V_1,\ldots ,V_n}\! :=\!\{ f\!\in\!\mathcal{H}(2^\omega )\mid\forall 1\!\leq\! i\!\leq\! n~~f[U_i]\! =\! V_i\}}\mbox{,}$$ 
where $n$ is a natural number and $U_i,V_i$ are clopen subsets of $2^\omega$. By [I-Me, Section 2], this defines a structure of Polish group on $\mathcal{H}(2^\omega )$. A compatible complete distance is given by
$$d(f,g)\! :=\!\mbox{sup}_{\alpha\in 2^\omega}~d_{2^\omega}\big( f(\alpha ),g(\alpha )\big)\! +\!
\mbox{sup}_{\alpha\in 2^\omega}~d_{2^\omega}\big( f^{-1}(\alpha ),g^{-1}(\alpha )\big) .$$
 
\begin{lem} \label{graphmap} The map $f\!\mapsto\!\textup{Graph}(f)$ from $\mathcal{H}(2^\omega )$ into 
$\mathcal{K}(2^\omega\!\times\! 2^\omega )$ is continuous.\end{lem}

\noindent\emph{Proof.}\ If $O\!\not=\!\emptyset$ is an open subset of $2^\omega\!\times\! 2^\omega$, and 
$(s_n)_{n\in\omega},(t_n)_{n\in\omega}$ are sequences of finite binary sequences with 
$O\! =\!\bigcup_{n\in\omega}~(N_{s_n}\!\times\! N_{t_n})$, then
$$\textup{Graph}(f)\!\subseteq\! O\Leftrightarrow
\exists F\!\subseteq\!\omega\mbox{ finite with }\textup{Graph}(f)\!\subseteq\! U_F\! :=\!\bigcup_{n\in F}~(N_{s_n}\!\times\! N_{t_n}).$$ 

 If $\textup{Graph}(f)\!\subseteq\! U_F$, $l\! :=\!\mbox{max}_{n\in F}~\vert t_n\vert$ and $d(f,g)\! <\! 2^{-l}$, then 
$\textup{Graph}(g)\!\subseteq\! U_F$, which proves that $\{ f\!\in\!\mathcal{H}(2^\omega )\mid\textup{Graph}(f)\!\subseteq\! O\}$ is open. Now
$$\textup{Graph}(f)\cap O\!\not=\!\emptyset\Leftrightarrow
\exists n\!\in\!\omega ~~\textup{Graph}(f)\cap (N_{s_n}\!\times\! N_{t_n})\!\not=\!\emptyset\Leftrightarrow
\exists n\!\in\!\omega ~~\exists\alpha\!\in\! N_{s_n}~~f(\alpha )\!\in\! N_{t_n}\mbox{,}$$ 
so that $\{ f\!\in\!\mathcal{H}(2^\omega )\mid\textup{Graph}(f)\cap O\!\not=\!\emptyset\}$ is open.\hfill{$\square$}
 
\begin{lem} \label{Ghmap} The map $f\!\mapsto\! F^f_1$ from $\mathcal{H}(2^\omega )$ into $\mathcal{K}(2^\omega )$ is Baire class one and not continuous. In fact, $\{ f\!\in\!\mathcal{H}(2^\omega )\mid F^f_1\!\subseteq\! U\}$ is open for each open subset $U$ of $2^\omega$.\end{lem}

\noindent\emph{Proof.}\ If $U$ is an open subset of $2^\omega$, then 
$F^f_1\!\subseteq\! U\Leftrightarrow\textup{Graph}(f)\!\subseteq\!\neg\Delta (2^\omega\!\setminus\! U)$. This implies that 
$\{ f\!\in\!\mathcal{H}(2^\omega )\mid F^f_1\!\subseteq\! U\}$ is open by Lemma \ref{graphmap}. If now 
$U\! =\!\bigcup_{n\in\omega}~N_{s_n}$ is not empty, then 
$F^f_1\cap U\!\not=\!\emptyset\Leftrightarrow\exists n\!\in\!\omega ~~\textup{Graph}(f)\!\not\subseteq\!\neg\Delta (N_{s_n})$, so that $\{ f\!\in\!\mathcal{H}(2^\omega )\mid F^f_1\cap U\!\not=\!\emptyset\}$ is $\boratwo$. This last set is not open if 
$U\! =\! 2^\omega$ since it contains $\mbox{Id}$, which is the limit of $g_n$ defined by 
$g_n(\alpha )(p)\! =\!\alpha (p)\Leftrightarrow p\!\leq\! n$. This finishes the proof.\hfill{$\square$}

\vfill\eject

\noindent\emph{Notation.}\ We define a family $(h_s)_{s\in\omega^{<\omega}}$ of functions from $2^\omega$ into itself as follows. If $s\!\in\!\omega^{<\omega}$, then we set $s(-1)\! :=\! 0$.\medskip

\noindent - If $\vert s\vert$ is even, then we set $h_s(0^\infty )\! :=\! 0^\infty$, 
$$h_s(0^{(\Sigma_{i<j}~(s(2i)+1))+p}1t\varepsilon\alpha )\! :=\! 
0^{(\Sigma_{i<j}~(s(2i)+1))+p}1t(1\! -\!\varepsilon )\alpha$$
if $j\! <\!\frac{\vert s\vert}{2}$, $p\!\leq\! s(2j)$, $t\!\in\! 2^{s(2j-1)}$, $\varepsilon\!\in\! 2$ and $\alpha\!\in\! 2^\omega$, and 
$$h_s(0^{(\Sigma_{i<\frac{\vert s\vert}{2}}~(s(2i)+1))+p}1t\varepsilon\alpha )\! :=\! 
0^{(\Sigma_{i<\frac{\vert s\vert}{2}}~(s(2i)+1))+p}1t(1\! -\!\varepsilon )\alpha$$
if $p\!\in\!\omega$, $t\!\in\! 2^{s(\vert s\vert-1)}$, $\varepsilon\!\in\! 2$ and $\alpha\!\in\! 2^\omega$.\medskip

\noindent - If $\vert s\vert$ is odd, then we set 
$$h_s(0^{(\Sigma_{i<j}~(s(2i)+1))+p}1t\varepsilon\alpha )\! :=\! 
0^{(\Sigma_{i<j}~(s(2i)+1))+p}1t(1\! -\!\varepsilon )\alpha$$
if $j\!\leq\!\frac{\vert s\vert -1}{2}$, $p\!\leq\! s(2j)$, $t\!\in\! 2^{s(2j-1)}$, $\varepsilon\!\in\! 2$ and $\alpha\!\in\! 2^\omega$, and
$$h_s(0^{\Sigma_{i\leq\frac{\vert s\vert -1}{2}}~(s(2i)+1)}\alpha )\! :=\! 0^{\Sigma_{i\leq\frac{\vert s\vert -1}{2}}~(s(2i)+1)}\alpha$$
if $\alpha\!\in\! 2^\omega$.

\begin{lem} \label{rati} The $h_s$'s are continuous involutions, $\chi_c(2^\omega ,G_{h_s})\! =\! 2^{\aleph_0}$ if $\vert s\vert$ is even, $\chi_c(2^\omega ,G_{h_s})\! =\! 2$ if $\vert s\vert$ is odd, and $(h_{sn})_{n\in\omega}$ converges to $h_s$ in 
$\mathcal{H}(2^\omega )$.\end{lem}

\noindent\emph{Proof.}\ Note that $h_s$ is a continuous involution, and thus a homeomorphism. If $\vert s\vert$ is even, then $0^\infty$ is the only fixed point of the map $h_s$, so that $\chi_c(2^\omega ,G_{h_s})\! =\! 2^{\aleph_0}$ by Proposition \ref{fp}. If $\vert s\vert$ is odd, then $F^{h_s}_1\! =\! N_{0^{\Sigma_{i\leq\frac{\vert s\vert -1}{2}}~(s(2i)+1)}}$ is a clopen subset of 
$2^\omega$. By Corollaries \ref{corfp}(a) and \ref{homeocompactcomplete}, $\chi_c(2^\omega ,G_{h_s})\!\in\!\{ 2,3\}$. By Corollary \ref{corfp}(b) and Proposition \ref{invol}, $\chi_c(2^\omega ,G_{h_s})\! =\! 2$. Note that the inequality 
$\mbox{sup}_{\alpha\in 2^\omega}~d_{2^\omega}\big( h_{sn}(\alpha ),h_s(\alpha )\big)\! <\! 2^{-n}$ holds. We are done since the $h_s$'s are involutions.\hfill{$\square$}

\begin{lem} \label{range} The map $(K,f)\!\mapsto\! f[K]$ from $\mathcal{K}(2^\omega )\!\times\!\mathcal{C}(2^\omega ,2^\omega )$ into 
$\mathcal{K}(2^\omega )$ is continuous. This is also the case if we replace $\mathcal{C}(2^\omega ,2^\omega )$ with 
$\mathcal{H}(2^\omega )$, $2^\omega$ with $\mathcal{K}_{2^\infty}$. The map 
$(K,f)\!\mapsto\! (f\!\times\! f)[K]$ from $\mathcal{K}(2^\omega\!\times\! 2^\omega )\!\times\!\mathcal{C}(2^\omega ,2^\omega )$ into $\mathcal{K}(2^\omega\!\times\! 2^\omega )$ is also continuous.\end{lem}
 
\noindent\emph{Proof.}\ Let $O\!\not=\!\emptyset$ be an open subset of $2^\omega$, and $(s_n)_{n\in\omega}$ be a sequence of finite binary sequences with 
$O\! =\!\bigcup_{n\in\omega}~N_{s_n}$. Note that 
$$\begin{array}{ll}
f[K]\!\subseteq\! O\!\!\!\!
& \Leftrightarrow\exists F\!\subseteq\!\omega\mbox{ finite with }f[K]\!\subseteq\!\bigcup_{n\in F}~N_{s_n}\cr
& \Leftrightarrow\exists F\!\subseteq\!\omega\mbox{ finite }\exists C\!\in\!\borone (2^\omega )~
K\!\subseteq\! C\wedge C\!\subseteq\! f^{-1}(\bigcup_{n\in F}~N_{s_n}).
\end{array}$$
Let $l_F\! :=\!\mbox{max}_{n\in F}~\vert s_n\vert$. If 
$\mbox{sup}_{\alpha\in 2^\omega}~d_{2^\omega}\big( f(\alpha ),g(\alpha )\big)\! <\! 2^{-l_F}$ and 
$C\!\subseteq\! f^{-1}(\bigcup_{n\in F}~N_{s_n})$, then 
$$C\!\subseteq\! g^{-1}(\bigcup_{n\in F}~N_{s_n})\mbox{,}$$ 
so that 
$\{ (K,f)\!\in\!\mathcal{K}(2^\omega )\!\times\!\mathcal{C}(2^\omega ,2^\omega )\mid K\!\subseteq\! C\wedge C\!\subseteq\! f^{-1}(\bigcup_{n\in F}~N_{s_n})\}$ is open. This shows that  
$\{ (K,f)\!\in\!\mathcal{K}(2^\omega )\!\times\!\mathcal{C}(2^\omega ,2^\omega )\mid f[K]\!\subseteq\! O\}$ is open (even if $O$ is empty).

\vfill\eject

 Now 
$$\begin{array}{ll}
f[K]\cap O\!\not=\!\emptyset\!\!\!\!
& \Leftrightarrow\exists n\!\in\!\omega~f[K]\cap N_{s_n}\!\not=\!\emptyset\cr
& \Leftrightarrow\exists n\!\in\!\omega~\exists C\!\in\!\borone (2^\omega )~K\cap C\!\not=\!\emptyset\wedge C\!\subseteq\! f^{-1}(N_{s_n})\mbox{,}
\end{array}$$
so that $\{ (K,f)\!\in\!\mathcal{K}(2^\omega )\!\times\!\mathcal{C}(2^\omega ,2^\omega )\mid f[K]\cap O\!\not=\!\emptyset\}$ is open.
\hfill{$\square$}

\begin{thm} \label{descr} $\mathfrak{H}_0\! =\!\mathcal{H}(2^\omega )$, $\mathfrak{H}_1$ is a $\boraone\!\setminus\!\bormone$ subset of $\mathcal{H}(2^\omega )$, while $\mathfrak{H}_2,\mathfrak{H}_3$ are $\bormtwo\!\setminus\!\boratwo$ subsets of 
$\mathcal{H}(2^\omega )$.\end{thm}

\noindent\emph{Proof.}\ We may restrict our attention to $\mathfrak{H}_2,\mathfrak{H}_3$ since 
$\mathfrak{H}_1\! =\!\mathcal{H}(2^\omega )\!\setminus\!\{\mbox{Id}\}$. Fix now $\kappa\! <\!\omega$. Note that 
$\chi_c(2^\omega ,G_f)\!\leq\!\kappa$ holds if and only if 
$$\exists (C_i)_{i<\kappa}\!\in\!\borone (2^\omega )^\kappa ~~(\forall i\!\not=\! j\! <\!\kappa ~~C_i\cap C_j\! =\!\emptyset )\wedge (2^\omega\!\subseteq\!\bigcup_{i<\kappa}~C_i)\wedge (\forall i\! <\!\kappa ~~G_f\cap C_i^2\! =\!\emptyset ).$$
As $\kappa$ is finite and $\borone (2^\omega )$ is countable, $\borone (2^\omega )^\kappa$ is countable. So we can restrict our attention to 
$$\begin{array}{ll}
G_f\cap C_i^2\! =\!\emptyset\!\!\!\!
& \Leftrightarrow 2^\omega\!\subseteq\! (2^\omega\!\setminus\! C_i)\cup f^{-1}(2^\omega\!\setminus\! C_i)\cup F^f_1\cr
& \Leftrightarrow\neg (\exists C\!\in\!\borone (2^\omega )\!\setminus\!\{\emptyset\} ~~(2^\omega\!\setminus\! C_i)\cup 
f^{-1}(2^\omega\!\setminus\! C_i)\cup F^f_1\!\subseteq\! 2^\omega\!\setminus\! C).
\end{array}$$ 
By Lemmas \ref{range} and \ref{Ghmap}, $\mathcal{H}(2^\omega )\!\setminus\!\mathfrak{H}_\kappa$ is a $\boratwo$ subset of 
$\mathcal{H}(2^\omega )$. In particular, $\mathfrak{H}_2,\mathfrak{H}_3$ are $\bormtwo$. It remains to see that 
$\mathfrak{H}_2,\mathfrak{H}_3$ are not $\boratwo$. We will use the family  $(h_s)_{s\in\omega^{<\omega}}$ defined before Lemma \ref{rati}. We set $P\! :=\!\overline{\{ h_s\mid s\in\omega^{<\omega}\}}$, so that $P$ is a Polish space. Note that 
$\mathfrak{H}_3\cap P\!\subseteq\!\mathfrak{H}_2\cap P$ are dense and co-dense $\bormtwo$ subsets of $P$, by Lemma \ref{rati}. By Baire's theorem, $\mathfrak{H}_2,\mathfrak{H}_3$ are not $\boratwo$.\hfill{$\square$}\medskip
 
 We next turn to the $\mathfrak{G}_\kappa$'s.\medskip

\noindent\emph{Proof of Theorem \ref{desc}.}\ Note that $f[X]\! =\! X\Leftrightarrow f[X]\!\subseteq\! X\wedge f^{-1}[X]\!\subseteq\! X$. By Lemma \ref{range} and [K, 4.29], $\mathcal{P}$ is a closed subset of the Polish space $\mathcal{K}(2^\omega )\!\times\!\mathcal{H}(2^\omega )$, and thus a Polish space. Note also that 
${\mathcal{O}_0\! =\!\mathcal{P}\!\setminus\!\big(\{\emptyset\}\!\times\!\mathcal{H}(2^\omega )\big)}$ is a clopen subset of 
$\mathcal{P}$ since $\emptyset$ is an isolated point in $\mathcal{K}(2^\omega )$.\medskip

 Note that $\chi_c(X,G_{f_{\vert X}})\!\leq\! 1\Leftrightarrow X\!\subseteq\! F^f_1$, i.e., $f_{\vert X}\! =\!\mbox{Id}_{\vert X}$. Let 
$\big( (X_n,h_n)\big)_{n\in\omega}$ be a sequence of elements of $\mathcal{P}\!\setminus\!\mathcal{O}_1$ converging to a point 
$(X,h)$ of $\mathcal{P}$. As $\chi_c(X_n,G_{(h_n)_{\vert X_n}})\!\leq\! 1$, $(h_n)_{\vert X_n}\! =\!\mbox{Id}_{\vert X_n}$. Assume, towards a contradiction, that $h_{\vert X}\!\not=\!\mbox{Id}_{\vert X}$. This gives $\alpha\!\in\! X$ with $h(\alpha )\!\not=\!\alpha$, 
$l\!\in\!\omega$ such that $h(\alpha )\vert l\!\not=\!\alpha\vert l$, and $L\!\geq\! l$ such that 
$g(\beta )\vert l\! =\! h(\alpha )\vert l\!\not=\!\alpha\vert l\! =\!\beta\vert l$ if $\beta\!\in\! N_{\alpha\vert L}$ and $d(g,h)\! <\! 2^{-L}$. As 
$\alpha\!\in\! X\cap N_{\alpha\vert L}$, $X_n\cap N_{\alpha\vert L}\!\not=\!\emptyset$ and $d(h_n,h)\! <\! 2^{-L}$ if $n$ is large enough. We pick, for such a $n$, $\beta\!\in\! X_n\cap N_{\alpha\vert L}$, so that $h_n(\beta )\vert l\!\not=\!\beta\vert l$, contradicting 
$(h_n)_{\vert X_n}\! =\!\mbox{Id}_{\vert X_n}$. This shows that $\mathcal{O}_1$ is open. We define, for $n\!\in\!\omega$, a map $g_n\! :\! 2^\omega\!\rightarrow\! 2^\omega$ by $g_n(0^\infty )\! :=\! 0^\infty$, $g_n(0^p1\alpha )\! :=\! 0^p1\alpha$ if 
$p\! <\! 2n$, and $g_n(0^{2p+\varepsilon}1\alpha )\! :=\! 0^{2p+(1-\varepsilon )}1\alpha$ if $p\!\geq\! n$ and 
$\varepsilon\!\in\! 2$, so that $g_n$ is a continuous involution whose set of fixed points 
$\{ 0^\infty\}\cup\bigcup_{p<2n}~N_{0^p1}$ is not open. By Proposition \ref{fp}, $\chi_c(2^\omega ,G_{g_n})\! =\! 2^{\aleph_0}$. As the sequence $(g_n)$ converges in $\mathcal{H}(2^\omega )$ to $\mbox{Id}$, $\mathcal{O}_1$ is not  closed. Thus 
$\mathcal{O}_1$ is $\boraone$-complete, by [K, 22.11].\medskip

 Fix now $\kappa\! <\!\omega$. Note that, by [E, Theorem 2.1(1)],
$$\begin{array}{ll}
\chi_c(X,G_{f_{\vert X}})\!\leq\!\kappa\!\!\!\!
& \Leftrightarrow\exists (C_i)_{i<\kappa}\!\in\!\borone (X)^\kappa ~~(\forall i\!\not=\! j\! <\!\kappa ~~C_i\cap C_j\! =\!\emptyset )\wedge 
(X\!\subseteq\!\bigcup_{i<\kappa}~C_i)\wedge\cr
& \hfill{(\forall i\! <\!\kappa ~~G_f\cap C_i^2\! =\!\emptyset )}\cr
& \Leftrightarrow\exists (C_i)_{i<\kappa}\!\in\!\borone (2^\omega )^\kappa ~~(\forall i\!\not=\! j\! <\!\kappa ~~C_i\cap C_j\! =\!\emptyset )\wedge (2^\omega\!\subseteq\!\bigcup_{i<\kappa}~C_i)\wedge\cr
& \hfill{(\forall i\! <\!\kappa ~~G_f\cap (X\cap C_i)^2\! =\!\emptyset ).}
\end{array}$$

 As $\kappa$ is finite and $\borone (2^\omega )$ is countable, $\borone (2^\omega )^\kappa$ is countable. So we can restrict our attention to 
$$\begin{array}{ll}
G_f\cap (X\cap C_i)^2\! =\!\emptyset\!\!\!\!
& \Leftrightarrow X\!\subseteq\! (2^\omega\!\setminus\! C_i)\cup f^{-1}(2^\omega\!\setminus\! C_i)\cup F^f_1\cr
& \Leftrightarrow\neg (\exists C\!\in\!\borone (2^\omega )~~ X\cap C\!\not=\!\emptyset\wedge (2^\omega\!\setminus\! C_i)\cup 
f^{-1}(2^\omega\!\setminus\! C_i)\cup F^f_1\!\subseteq\! 2^\omega\!\setminus\! C).
\end{array}$$ 
By Lemmas \ref{range} and \ref{Ghmap}, $\mathcal{P}\!\setminus\!\mathcal{O}_\kappa$ is a $\boratwo$ subset of 
$\mathcal{P}$. In particular, the sets $\mathcal{O}_2,\mathcal{O}_3$ are $\bormtwo$. As 
$(\mathcal{O}_\kappa )_{2^\omega}\! =\!\mathfrak{H}_\kappa$, 
$\mathcal{O}_2,\mathcal{O}_3$ are not $\boratwo$ by Theorem \ref{descr}. Thus 
$\mathcal{O}_2,\mathcal{O}_3$ are $\bormtwo$-complete, by [K, 22.11]. The $\ca$-completeness of 
$\mathcal{O}_2^{\aleph_0}$ was proved at the very end of Section \ref{MSPolishsub}.\hfill{$\square$}
 
\section{$\!\!\!\!\!\!$ Equivalence relations}\indent

 Lemma \ref{flipo} and Corollary \ref{corflip} imply that $FCO$ is Borel reducible to different versions of $\equiv^i_c$.\medskip
 
\noindent\emph{Notation.}\ We set $\mathbb{M}\! :=\!\{ f\!\in\!\mathcal{H}(2^\omega )\mid f\mbox{ is minimal}\}$. By [Me, Lemma 4.1], the set $\mathbb{M}$ is a $G_\delta$ subset of $\mathcal{H}(2^\omega )$, and thus a Polish space. If 
$f,g\!\in\!\mathbb{M}$, then $f,g$ are flip-conjugate if and only if there is $\varphi\!\in\!\mathcal{H}(2^\omega )$ with 
$\varphi\!\circ\! f\! =\! g\!\circ\!\varphi$ or $\varphi\!\circ\! f\! =\! g^{-1}\!\circ\!\varphi$, proving that $FCO$ is analytic. Similarly, $CO$ is analytic.\medskip
 
 We first consider the case of graphs induced by a function. As in the introduction, we consider the equivalence relation 
$\equiv^i_c\ :=\ \preceq^i_c\cap~(\preceq^i_c)^{-1}$ on $\mathcal{S}_m$ associated with 
$$(2^\omega ,K)\preceq^i_c(2^\omega ,K')\Leftrightarrow\exists\varphi\! :\! 2^\omega\!\rightarrow\! 2^\omega\mbox{ injective continuous with }K\!\subseteq\! (\varphi\!\times\!\varphi )^{-1}(K').$$ 
We define a map $g\! :\!\mathbb{M}\!\rightarrow\!\mathcal{S}_m$ by $g(f)\! :=\! (2^\omega ,G_f)$ (see Theorem \ref{homeocompact}).
 
\begin{thm} \label{redborp} The equivalence relation $\equiv^i_c$ on the Polish space $\mathcal{S}_m$ is analytic, and $g$ reduces continuously $FCO$ to $\equiv^i_c$. Moreover, the vertices of the graph $g(f)$ have degree two, for each 
$f\!\in\!\mathbb{M}$.\end{thm}
 
\noindent\emph{Proof.}\ By Lemma \ref{flipo}, $g$ reduces $FCO$ to $\equiv^i_c$. Let $O$ be an open subset of 
$2^\omega\!\times\! 2^\omega$, and $(C^0_n)_{n\in\omega}$, $(C^1_n)_{n\in\omega}$ be sequences of clopen subsets of $2^\omega$ with $O\! =\!\bigcup_{n\in\omega}~(C^0_n\!\times\! C^1_n)$. If $f\!\in\!\mathbb{M}$ and $G_f\!\subseteq\! O$, then there is a finite subset $F$ of $\omega$ with 
$G_f\! =\! s\big(\textup{Graph}(f)\big)\!\subseteq\!\bigcup_{n\in F}~(C^0_n\!\times\! C^1_n)$. Note then that 
$\bigcup_{n\in F}~(C^0_n\!\times\! C^1_n)\! =\!
\bigcup_{S\subseteq F}~\big( (\bigcap_{n\in S}~C^0_n\cap\bigcap_{n\in F\setminus S}~2^\omega\!\setminus\! C^0_n)\!\times\! 
(\bigcup_{n\in S}~C^1_n)\big)$. Thus 
$$\begin{array}{ll}
\textup{Graph}(f)\!\subseteq\!\bigcup_{n\in F}~(C^0_n\!\times\! C^1_n)\!\!\!\!
& \Leftrightarrow\forall S\!\subseteq\! F~f[\bigcap_{n\in S}~C^0_n\cap\bigcap_{n\in F\setminus S}~2^\omega\!\setminus\! C^0_n]\!\subseteq\!\bigcup_{n\in S}~C^1_n\cr 
& \Leftrightarrow\forall S\!\subseteq\! F~\exists R_n\!\in\!\borone (2^\omega)\cr
& \hfill{f[\bigcap_{n\in S}~C^0_n\cap\bigcap_{n\in F\setminus S}~2^\omega\!\setminus\! C^0_n]\! =\! R_n
\!\subseteq\!\bigcup_{n\in S}~C^1_n.}
\end{array}$$
This implies that $\{ f\!\in\!\mathbb{M}\mid G_f\!\subseteq\! O\}$ is an open subset of $\mathbb{M}$ since 
$$G_f\!\subseteq\! O\Leftrightarrow
\exists F\!\subseteq\!\omega\mbox{ finite with }\textup{Graph}(f)\!\subseteq\!\bigcap_{\varepsilon\in 2}
\big(\bigcup_{n\in F}~(C^\varepsilon_n\!\times\! C^{1-\varepsilon}_n)\big) .$$ 
Now $G_f\cap O\!\not=\!\emptyset\Leftrightarrow\exists n\!\in\!\omega ~~\exists\varepsilon\!\in\! 2~~
C^\varepsilon_n\cap f^{-1}(C^{1-\varepsilon}_n)\!\not=\!\emptyset\Leftrightarrow
\exists n\!\in\!\omega ~~\exists\varepsilon\!\in\! 2~~\exists\alpha\!\in\! C^\varepsilon_n~~f(\alpha )\!\in\! C^{1-\varepsilon}_n$, so that $\{ f\!\in\!\mathbb{M}\mid G_f\cap O\!\not=\!\emptyset\}$ is an open subset of $\mathbb{M}$. Thus $g$ is continuous.\medskip

 Note that $(2^\omega ,K)\!\in\!\mathcal{S}_m$ if and only if\medskip
 
\leftline{$K\cap\Delta (2^\omega )\! =\!\emptyset\wedge K\!\not=\!\emptyset\wedge
\exists (C_i)_{i<3}\!\in\!\big(\borone (2^\omega )\big)^3~~2^\omega\!\subseteq\!\bigcup_{i<3}~C_i\wedge
\forall i\!\not=\! j\! <\! 3~~C_i\cap C_j\! =\!\emptyset ~\wedge$}\smallskip

\rightline{$K\cap (\bigcup_{i<3}~C_i^2)\! =\!\emptyset\mbox{,}$}\medskip

\noindent so that $\mathcal{S}_m$ is an open subset of $\{ 2^\omega\}\!\times\!\mathcal{K}(2^\omega\!\times\! 2^\omega )$ and thus a Polish space.

\vfill\eject

 Note then that $\varphi\! :\! 2^\omega\!\rightarrow\! 2^\omega$ is injective if and only if 
$\varphi [O\cap U]\! =\!\varphi[O]\cap\varphi [U]$ whenever $O,U$ are clopen subsets of $2^\omega$. By Lemma \ref{range}, and [K, 4.19, 4.29, 27.7], $\preceq^i_c$ and thus $\equiv^i_c$ are analytic.\hfill{$\square$}\medskip
 
\noindent\emph{Notation.}\ We now consider the case of general graphs, and we can ensure that the reduction map associates graphs of continuous chromatic number at least three instead of two or three. Recall that 
$\mathcal{K}_{2^\infty}\! :=\! (2\cup\{ c,a,\overline{a}\} )^\omega$. As $\mathcal{K}_{2^\infty}\!\not=\!\emptyset$ is a perfect  0DMC space, it is homeomorphic to $2^\omega$ via a map $i$, by [K, 7.4]. We equip 
$\mathcal{K}(\mathcal{K}_{2^\infty})$ with the Vietoris topology, so that $\mathcal{K}(\mathcal{K}_{2^\infty})$ is a metrizable compact space, by [K, 4.26]. By [K, 4.29], the map $K\!\mapsto\! i[K]$ defines a homeomorphism from 
$\mathcal{K}(\mathcal{K}_{2^\infty})$ onto $\mathcal{K}(2^\omega )$. We set 
$\mathbb{Q}\! :=\!\big\{ x\!\in\!\mathcal{K}_{2^\infty}\mid
\exists l\!\in\!\omega ~~\exists\varepsilon\!\in\!\{ a,\overline{a}\} ~~\forall k\!\geq\! l~~x(k)\! =\!\varepsilon\big\}$. Note that $\mathbb{Q}$ is countable, as well as $Q\! :=\! i[\mathbb{Q}]$. We set $\mathcal{S}_g\! :=\!\{ (K,R)\!\in\!\mathcal{K}(2^\omega )\!\times\! 2^{Q^2}\mid 
R\!\subseteq\! K^2\!\setminus\!\Delta (K)\wedge\chi_c(K,R)\!\geq\! 3\}$ and equip $2^{Q^2}$ with the product topology of the discrete topology on $2$, so that 
$$\{ (K,R)\!\in\!\mathcal{K}(2^\omega )\!\times\! 2^{Q^2}\mid R\!\subseteq\! K^2\!\setminus\!\Delta (K)\}$$ 
is a metrizable compact space.\medskip

 We consider the equivalence relation $\equiv^i_c$ on $\mathcal{S}_g$ associated with
$$(K,R)\preceq^i_c(K',R')\Leftrightarrow\exists\varphi\! :\! K\!\rightarrow\! K'\mbox{ injective continuous with }
R\!\subseteq\! (\varphi\!\times\!\varphi )^{-1}(R').$$ 
We equip $2^{\mathbb{Q}^2}$ with the product topology of the discrete topology on $2$, so that $2^{\mathbb{Q}^2}$ is homeomorphic to $2^\omega$. The map $R\!\mapsto\! (i\!\times\! i)[R]$ defines a homeomorphism from $2^{\mathbb{Q}^2}$ onto $2^{Q^2}$, and the equality ${i[K]\! =\!\overline{\mbox{proj}\big[(i\!\times\! i)[s(R)]\big]}}$ holds if $K\! =\!\overline{\mbox{proj}[s(R)]}$. Moreover, $(K,R)\equiv^i_c(i[K],(i\!\times\! i)[R])$. We define a map $\mathcal{G}\! :\!\mathbb{M}\!\rightarrow\!\mathcal{S}_g$ by 
$\mathcal{G}(f)\! :=\!\big(\overline{\mbox{proj}\big[ (i\!\times\! i)[\mathbb{G}_f]\big]},(i\!\times\! i)[\mathbb{G}_f]\big)$ (see Lemma \ref{chrom}). 
 
\begin{thm} \label{redborbis} The equivalence relation $\equiv^i_c$ on the Polish space $\mathcal{S}_g$ is analytic, and $\mathcal{G}$ Borel reduces $FCO$ to $\equiv^i_c$. Moreover, the vertices of the graph $\mathcal{G}(f)$ have degree at most one, for each $f\!\in\!\mathbb{M}$.\end{thm}

\noindent\emph{Proof.}\ As $(K,R)\equiv^i_c(i[K],(i\!\times\! i)[R])$, we may replace $2^\omega$ and $Q$ with 
$\mathcal{K}_{2^\infty}$ and $\mathbb{Q}$ respectively. By Corollary \ref{corflip}, $\mathcal{G}$ reduces $FCO$ to 
$\equiv^i_c$ since $\mathcal{C}^+\! =\!\overline{\mbox{proj}[\mathbb{G}_f]}$. Note that, for each $i\!\in\!\mathbb{Z}$, the map 
$f\!\mapsto\! f^i$ defined on $\mathcal{H}(2^\omega )$ is continuous since $\mathcal{H}(2^\omega )$ is a topological group. Note also that the evaluation map $(f,\alpha )\!\mapsto\! f(\alpha )$ is continuous since 
$d_{2^\omega}\big( f(\alpha ),f_0(\alpha_0)\big)\!\leq\! d(f,f_0)\! +\! d_{2^\omega}\big( f_0(\alpha ),f_0(\alpha_0)\big)$. This implies that the map from $\mathbb{M}$ into $(2^\omega )^\mathbb{Z}$ defined by 
$f\!\mapsto\!\big( f^i(0^\infty )\big)_{i\in\mathbb{Z}}$ is continuous. Here we only consider ${\bf d}\! :=\! 2^\infty$. Recall from the notation before Lemma \ref{conv} that $L_{2m}\! :=\! R_{2m+1}\! :=\!\zeta (m)$. The map from $(2^\omega )^\mathbb{Z}$ into $2^{\mathbb{Q}^2}$ associating\medskip

\leftline{$s\big(\{ (c^{l+1}a\overline{a}^\infty ,\gamma_{L_l}\vert (l\! +\! 1)\overline{a}a^\infty)\mid l\!\in\!\omega\} ~\cup$}\smallskip

\rightline{$\{ (\gamma_{L_l+i}\vert (l\! +\! 1)a^{i+1}\overline{a}^\infty ,\gamma_{L_l+i+1}\vert (l\! +\! 1)\overline{a}^{i+2}a^\infty )\mid 
l\!\in\!\omega\wedge i\!\leq\! 2l\} ~\cup$}\smallskip

\rightline{$\{ (\gamma_{R_l}\vert (l\! +\! 1)a^{2l+2}\overline{a}^\infty ,c^{l+1}\overline{a}a^\infty )\mid l\!\in\!\omega\}\big)$}\medskip

\noindent to $(\gamma_i)_{i\in\mathbb{Z}}$ is continuous, as well as $f\!\mapsto\!\mathbb{G}_f$. The map from 
$2^{\mathbb{Q}^2}$ into $2^\mathbb{Q}$ defined by $R\!\mapsto\!\mbox{proj}[s(R)]$ is Baire class one. The map from 
$2^\mathbb{Q}$ into $\mathcal{K}(\mathcal{K}_{2^\infty})$ defined by $S\!\mapsto\!\overline{S}$ is Borel, by [K, 12.11]. Thus $\mathcal{G}$ is Borel.\medskip

 Note that, by [E, Theorem 2.1(1)],
$$\begin{array}{ll}
\chi_c(K,R)\!\leq\! 2\!\!\!\!
& \Leftrightarrow\exists C\!\in\!\borone (K)~~R\cap C^2\! =\! R\cap (K\!\setminus\! C)^2\! =\!\emptyset\cr
& \Leftrightarrow\exists C\!\in\!\borone (2^\omega )~~R\cap C^2\! =\! R\cap (2^\omega\!\setminus\! C)^2\! =\!\emptyset 
\end{array}$$ 
if $R\!\subseteq\! K^2$, so that $\mathcal{S}_g$ is a $G_\delta$ subset of 
$\{ (K,R)\!\in\!\mathcal{K}(2^\omega )\!\times\! 2^{Q^2}\mid R\!\subseteq\! K^2\!\setminus\!\Delta (K)\}$ and thus a Polish space.

\vfill\eject

 If $(K,G),(L,H)\!\in\!\mathcal{S}_g$, then $(K,G)\preceq^i_c(L,H)$ holds if and only if there is $\varphi\! :\! K\!\rightarrow\! L$ injective continuous such that $\big(\varphi (x),\varphi (y)\big)\!\in\! H$ if $(x,y)\!\in\! G$. By [K, 2.8], this holds if and only if there is $\psi\! :\!\mathcal{K}_{2^\infty}\!\rightarrow\!\mathcal{K}_{2^\infty}$ continuous such that 
$\psi [K]\!\subseteq\! L$, $\psi_{\vert K}$ is injective, and 
$\big(\psi (x),\psi (y)\big)\!\in\! H$ if $(x,y)\!\in\! G$. Note that $\psi_{\vert K}$ is injective if and only if 
$\psi [O\cap U\cap K]\! =\!\psi [O\cap K]\cap\psi [U\cap K]$ whenever $O,U$ are clopen subsets of 
$\mathcal{K}_{2^\infty}$. We conclude as in the proof of Theorem \ref{redborp} to see that $\equiv^i_c$ is analytic.
\hfill{$\square$}\medskip
  
\noindent\emph{Remark.}\ As mentioned in the introduction, using oriented graphs instead of graphs, one can prove that $CO$ is Borel reducible to $\equiv^i_c$. In that case, the proof also works in the case of dynamical systems involving continuous maps instead of homeomorphisms, considering forward orbits instead of orbits.

\section{$\!\!\!\!\!\!$ Digraphs and oriented graphs}

\subsection{$\!\!\!\!\!\!$ General digraphs}\indent 
 
 We start with a version of Theorem \ref{corcomp''''''} for digraphs.
 
\begin{them} \label{compdig} We can find a concrete family 
$\big( (\mathbb{K}_\alpha ,\mathbb{D}_\alpha )\big)_{\alpha\in 2^\omega}$, where $\mathbb{K}_\alpha$ is a compact subset of $2^\omega$ and $\mathbb{D}_\alpha$ is a countable digraph on $\mathbb{K}_\alpha$, such that, for any 0DMC space $X$ and any digraph $D$ on $X$, exactly one of the following holds:\smallskip

\noindent (1) $D$ has CCN at most two,\smallskip

\noindent (2) we can find $\alpha\!\in\! 2^\omega$ and $\varphi\! :\!\mathbb{K}_\alpha\!\rightarrow\! X$ injective continuous such that $\mathbb{D}_\alpha\!\subseteq\! (\varphi\!\times\!\varphi )^{-1}(D)$.\smallskip

 In other words, $\big( (\mathbb{K}_\alpha ,\mathbb{D}_\alpha )\big)_{\alpha\in 2^\omega}$ is a $\preceq_c^i$-basis (and thus a $\preceq_c$-basis) for the class of digraphs on a 0DMC space with CCN at least three.\end{them}

\noindent\emph{Proof.}\ We define, for $(\gamma ,\delta )\!\in\!\mathcal{I}\!\times\! 2^{\{ (k,i)\in\omega^2\mid i\leq 2k\}}$, a countable relation $\mathbb{D}_{\gamma ,\delta}$ on $2^\omega$ by 
$$\mathbb{D}_{\gamma ,\delta}\! :=\!
\big\{\big(\gamma^{\delta (k,i)}_k(i),\gamma^{1-\delta (k,i)}_k(i)\big)\mid k\!\in\!\omega\wedge i\!\leq 2k\big\}\mbox{,}$$ 
so that $s(\mathbb{D}_{\gamma ,\delta})\! =\!\mathbb{G}_\gamma$ and $\mathbb{D}_{\gamma ,\delta}$ is a digraph on 
$\mathbb{K}_{\gamma ,\delta}\! :=\!\mathbb{K}_\gamma$. By Proposition \ref{K1''''''}, 
$\chi_c(\mathbb{K}_\gamma ,\mathbb{G}_\gamma )\!\geq\! 3$. As 
$s(\mathbb{D}_{\gamma ,\delta})\! =\!\mathbb{G}_\gamma$, 
$\chi_c(\mathbb{K}_\gamma ,\mathbb{D}_{\gamma ,\delta})\!\geq\! 3$ as well. It will be convenient to replace $2^\omega$ with $\mathcal{I}\!\times\! 2^{\{ (k,i)\in\omega^2\mid i\leq 2k\}}$. We just proved that (1) and (2) cannot hold simultaneously.\medskip

 Assume that (1) does not hold. Then $\chi_c\big( X,s(D)\big)\!\geq\! 3$. Theorem \ref{corcomp''''''} provides 
$\gamma\!\in\!\mathcal{I}$ with $(\mathbb{K}_\gamma ,\mathbb{G}_\gamma )\preceq^i_c\big( X,s(D)\big)$, with witness say $\varphi$. Let $\mathbb{D}\! :=\!\mathbb{G}_\gamma\cap(\varphi\!\times\!\varphi )^{-1}(D)$. Note that 
$s(\mathbb{D})\! =\!\mathbb{G}_\gamma$, which gives $\delta\!\in\! 2^{\{ (k,i)\in\omega^2\mid i\leq 2k\}}$ with 
$\mathbb{D}_{\gamma ,\delta}\!\subseteq\!\mathbb{D}$.\hfill{$\square$}\medskip

 Considering the $\mathbb{D}_\alpha$'s which are oriented graphs, and using the fact that a digraph $\preceq^i_c$-below an oriented graph is also an oriented graph, we get a (less concrete) basis for oriented graphs.
 
\begin{coro} \label{comporg} We can find a $\preceq_c^i$-basis (and thus a $\preceq_c$-basis) of size at most continuum, made up of countable oriented graphs, for the class of oriented graphs on a 0DMC space with CCN at least three.\end{coro}

 By Theorem \ref{eantichmin}, any $\preceq^i_c$-basis for the class of digraphs on a 0DMC space with CCN at least three must have size at least continuum. We will see that the $\preceq^i_c$-basis given by Corollary \ref{comporg} must also have size exactly the continuum later. The second basis given by Theorem \ref{corcomp'''''''} also provides a second basis for digraphs, which is also a basis for oriented graphs, more concrete than the first one we just met.

\begin{them} \label{compdigo} We can find a concrete $\preceq_c$-basis of size continuum, made up of countable oriented graphs, for the class of digraphs on a 0DMC space with CCN at least three.\end{them}

\noindent\emph{Proof.}\ We adapt the proof of Theorem \ref{corcomp'''''''}. We set 
$$\mathcal{J}_o\! :=\!\{ (\beta ,\alpha )\!\in\!\mathcal{J}\!\times\! (2^{<\omega})^\omega\mid
\forall l\!\in\!\omega ~~\vert\alpha (l)\vert\! =\!\lambda_l\! +\! 1\} .$$ 
and $\mathcal{J}^c_o\! :=\!\mathcal{J}_o\cap\big(\mathcal{J}^c\!\times\! (2^{<\omega})^\omega\big)$. We also set, for 
$x,y\!\in\!\mathcal{K}_{\bf d}$ and $\varepsilon\!\in\! 2$,
$$(x,y)_\varepsilon\! :=\!\left\{\!\!\!\!\!\!\!\!
\begin{array}{ll}
& (x,y)\mbox{ if }\varepsilon\! =\! 0\mbox{,}\cr
& (y,x)\mbox{ if }\varepsilon\! =\! 1.
\end{array}
\right.$$

 We then define, for $(\beta ,\alpha )\!\in\!\mathcal{J}_o$, a countable digraph $\mathbb{D}_{\beta ,\alpha}$ on 
$\mathcal{K}_{\bf d}$ by\medskip

\leftline{$\mathbb{D}_{\beta ,\alpha}\! :=\!
\{ (c^{l+1}a\overline{a}^\infty ,s_l(0)\overline{a}a^\infty )_{\alpha (l)(0)}\mid l\!\in\!\omega\} ~\cup$}\smallskip

\rightline{$\{ (s_l(i)a^{i+1}\overline{a}^\infty ,s_l(i\! +\! 1)\overline{a}^{i+2}a^\infty )_{\alpha (l)(i+1)}\mid 
l\!\in\!\omega\wedge i\!\leq\!\lambda_l\! -\! 2\} ~\cup$}\smallskip

\rightline{$\{ (s_l(\lambda_l\! -\! 1)a^{\lambda_l}\overline{a}^\infty ,c^{l+1}\overline{a}a^\infty )_{\alpha (l)(\lambda_l)}\mid 
l\!\in\!\omega\}\mbox{,}$}\medskip

\noindent so that $s(\mathbb{D}_{\beta ,\alpha})\! =\!\mathbb{G}_\beta$ and $\mathbb{D}_{\beta ,\alpha}$ is an oriented graph on $\mathbb{K}_\beta$. By Lemma \ref{chromgen}, $\chi_c(\mathbb{K}_\beta ,\mathbb{G}_\beta )\!\geq\! 3$. As 
$s(\mathbb{D}_{\beta ,\alpha})\! =\!\mathbb{G}_\beta$, $\chi_c(\mathbb{K}_\beta ,\mathbb{D}_{\beta ,\alpha})\!\geq\! 3$ as well. By Theorem \ref{compdig}, it is enough to prove that if $(\gamma ,\delta )$ is in 
$\mathcal{I}\!\times\! 2^{\{ (k,i)\in\omega^2\mid i\leq 2k\}}$, then we can find $(\beta ,\alpha )\!\in\!\mathcal{J}^c_o$ (for 
${\bf d}\! =\! 2^\infty$) such that 
$(\mathbb{K}_\beta ,\mathbb{D}_{\beta ,\alpha})\preceq_c(\mathbb{K}_\gamma ,\mathbb{D}_{\gamma ,\delta})$. We first define $\beta'$ as in the proof of Theorem \ref{corcomp'''''''}, and define $\alpha'\!\in\! (2^{<\omega})^\omega$ by 
${\alpha'(q)(i)\! :=\!\delta (k_q,i)}$ if $q\!\in\!\omega$ and $i\!\leq\!\lambda'_q$. Then $(\beta',\alpha')\!\in\!\mathcal{J}_o$, and 
$(\mathbb{K}_{\beta'},\mathbb{D}_{\beta',\alpha'})\preceq_c(\mathbb{K}_\gamma ,\mathbb{D}_{\gamma ,\delta})$, by the proof of Theorem \ref{corcomp'''''''}. We then define $\beta$ as in the proof of Theorem \ref{corcomp'''''''}, and define 
$\alpha\!\in\! (2^{<\omega})^\omega$ by the formula $\alpha (l)(i)\! :=\!\alpha'(q_0^{2l+2})(i)$ if $l\!\in\!\omega$ and 
$i\!\leq\!\lambda_l$. Then $(\beta ,\alpha )\!\in\!\mathcal{J}^c_o$ and the proof of Theorem \ref{corcomp'''''''} shows that 
$(\mathbb{K}_{\beta},\mathbb{D}_{\beta ,\alpha })\preceq^i_c(\mathbb{K}_\beta',\mathbb{D}_{\beta',\alpha'})$.
\hfill{$\square$}\medskip

 In order to prove that the basis given by Corollary \ref{comporg} has size exactly continuum, we prove the following oriented version of Theorem \ref{Gomin}.

\begin{them} \label{Oomin} Let ${\bf d}\!\in\!\mathfrak{D}$, $V$ be a compact subspace of $\mathcal{C}^+$, and 
$E\!\subseteq\!\mathbb{O}_o\cap V^2$ having CCN three. Then 
$(\mathcal{C}^+,\mathbb{O}_o)\preceq_c^i(V,E)$.\end{them}
 
\noindent\emph{Proof.}\ We essentially copy the proof of Theorem \ref{Gomin}, replacing $E$ with $s(E)$ under the closure symbols, and using the fact that $E\!\subseteq\!\mathbb{O}_o$.\hfill{$\square$}\medskip

 We are now ready to prove a version of Theorem \ref{eantichmin} for oriented graphs.

\begin{them} \label{eantichminor} There is a $\preceq_c$-antichain (and thus $\preceq_c^i$-antichain) 
$\big( (\mathbb{K}_\alpha ,\mathbb{O}_\alpha )\big)_{\alpha\in 2^\omega}$, where\smallskip
 
\noindent (a) $\mathbb{K}_\alpha$ is a 0DMC space,\smallskip

\noindent (b) $\mathbb{O}_\alpha$ is a countable $D_2(\bormone )$ oriented graph on $\mathbb{K}_\alpha$ with CCN three and $\boraone\oplus\bormone$ chromatic number two, and whose vertices have degree at most one,\smallskip

\noindent (c) $(\mathbb{K}_\alpha ,\mathbb{O}_\alpha )$ is $\preceq_c^i$-minimal in the class of digraphs on a 0DMC space with CCN at least three.\smallskip

 In particular, any $\preceq_c^i$-basis for the class of digraphs (or oriented graphs) on a 0DMC space with CCN at least three. must have size at least continuum.\end{them}

\noindent\emph{Proof.}\ Fix ${\bf d}\!\in\!\mathfrak{D}$. By Proposition \ref{belowfin}, $(\mathcal{C}^+,\mathbb{G}_o)$ has CCN three and $\boraone\oplus\bormone$ chromatic number two. As 
$s(\mathbb{O}_o)\! =\!\mathbb{G}_o$, this is also the case of $(\mathcal{C}^+,\mathbb{O}_o)$. We check that 
$\mathbb{O}_\alpha$ is $D_2(\bormone )$ as in the proof of Lemma \ref{D2gen}, so that (a) and (b) hold. For (c), i.e., the minimality of $(\mathcal{C}^+,\mathbb{O}_o)$, we first note that the proof of Lemma \ref{cpctmin} works for digraphs instead of graphs. We then apply Theorem \ref{Oomin}. $(\mathcal{C}^+_{\Phi (\alpha )},
\mathbb{G}_{o_{\Phi (\alpha )}})\not\preceq_c(\mathcal{C}_{\Phi (\beta )}^+,\mathbb{G}_{o_{\Phi (\beta )}})$ if 
$\alpha\!\not=\!\beta$, by Theorem \ref{anticha}. As $s(\mathbb{O}_o)\! =\!\mathbb{G}_o$ again, $(\mathcal{C}^+_{\Phi (\alpha )},
\mathbb{O}_{o_{\Phi (\alpha )}})\not\preceq_c(\mathcal{C}_{\Phi (\beta )}^+,\mathbb{O}_{o_{\Phi (\beta )}})$.\hfill{$\square$}\medskip

 Replacing $\mathbb{G}_p\! =\! s(\mathbb{O}_p)$ with $\mathbb{O}_p$, we get a version of Theorem \ref{infdecrcompact} for oriented graphs in a straightforward way. This kind of argument will be used several times in the sequel, and we will not always repeat it.\medskip
 
  Our version of Theorem \ref{absmin} for digraphs and oriented graphs is as follows.
 
\begin{them} \label{absminor} Let $D$ be a digraph on a 0DMS space $Z$, with CCN at least three and satisfying 
$(Z,D)\preceq^i_c(\mathbb{P},\mathbb{O}_m)$. Then there is a family $\big( (P_\alpha ,O_\alpha )\big)_{\alpha\in 2^\omega}$ of oriented graphs on a 0DP space with CCN three, $\preceq^i_c$-below $(Z,D)$, and pairwise $\preceq^i_c$-incompatible in the class of digraphs on a 0DMS space with CCN at least three.\smallskip

 In particular, there is no $\preceq^i_c$-antichain basis in the class of digraphs (or oriented graphs) on a 0DMS (or 0DP) space with CCN at least three, and any $\preceq^i_c$--basis for one of these classes must have size at least continuum.\end{them}

\noindent\emph{Proof.}\ Note first that we can modify Lemma \ref{below0} as follows. Let $\delta\!\in\! 2^\omega$, and $D$ be a digraph on a 0DMS space $Z$, with CCN at least three and satisfying $(Z,D)\preceq^i_c(\mathbb{P}_\delta ,\mathbb{O}_\delta )$. Then there is $\delta'\!\in\!\mathbb{P}_\infty$ such that 
${\{ k\!\in\!\omega\mid\delta'(k)\! =\! 1\}\!\subseteq\!\{ k\!\in\!\omega\mid\delta (k)\! =\! 1\}}$ and 
$(\mathbb{P}_{\delta'},\mathbb{O}_{\delta'})\preceq^i_c(Z,D)$. We complete the proof of Lemma \ref{below0} as follows. We set $G\! :=\! s(D)$, so that $G$ is a graph on $Z$ with CCN at least three and 
$(Z,G)\preceq^i_c(\mathbb{P}_\delta ,\mathbb{G}_\delta )$. We then set $R\! :=\! (\varphi\!\times\!\varphi )[D]$, so that 
$s(R)\! =\! E$. We can then follow the proof of Lemma \ref{below0}. For the conditions (b)-(d), the couples are not only in 
$s(R)$, but also in $R$ since $R\!\subseteq\! \mathbb{O}_\delta$ and thus $R^{-1}\!\subseteq\! \mathbb{O}_\delta^{-1}$. So we can replace $E$ with $R$ after the first line of the proof, which implies that 
$(\mathbb{P}_{\delta'},\mathbb{O}_{\delta'})\preceq^i_c(Z,D)$.\medskip

 The version of Lemma \ref{less} for oriented graphs is straightforward, and we then follow the proof of Theorem \ref{absmin} to conclude.\hfill{$\square$}\medskip

 Our version of Theorem \ref{absmincomp} for digraphs and oriented graphs is as follows.

\begin{them} \label{absmincompor} There is a countable oriented graph $(3^\omega ,\mathbb{O})$ in the class of digraphs  on a 0DMC space with CCN at least three such that, for each $(K,G)$ in this class satisfying $(K,D)\preceq^i_c(3^\omega ,\mathbb{O})$, there is a $\preceq^i_c$-antichain 
$\big( (3^\omega ,O_\alpha )\big)_{\alpha\in 2^\omega}$ of oriented graphs with CCN three and 
$\preceq^i_c$-below $(K,D)$. In particular, there is no $\preceq^i_c$-antichain basis in this class (or the corresponding one for oriented graphs).\end{them}

\noindent\emph{Proof.}\ We just have to follow the proof of Theorem \ref{absmincomp}. The oriented graphs 
$\mathbb{O}\! :=\!\textup{Graph}(o_{\vert D_{S_\omega}})$ and $O_\alpha\! :=\! G_\alpha\cap\textup{Graph}(o)$ are convenient since $s(\mathbb{O})\! =\!\mathbb{G}$ and $s(O_\alpha )\! =\! G_\alpha$.\hfill{$\square$}\medskip
 
  We can prove a version of Theorem \ref{embed} for oriented graphs in a straightforward way. A straightforward modification of Section \ref{gdn} gives the following version of Corollary \ref{corflip} for oriented graphs. In order to do that, we replace $\mathbb{G}_{f_{\bf d}}$ with $\mathbb{O}_{f_{\bf d}}$ in the definition of continuous tuples.

\begin{them} \label{corflipor} Let ${\bf d},{\bf d}'\!\in\!\mathfrak{C}$, 
$f_{\bf d}\! :\!\mathcal{C}_{\bf d}\!\rightarrow\!\mathcal{C}_{\bf d}$, 
$f_{{\bf d}'}\! :\!\mathcal{C}_{{\bf d}'}\!\rightarrow\!\mathcal{C}_{{\bf d}'}$ be minimal homeomorphisms, and 
$(n_l)_{l\in\omega},(L_l)_{l\in\omega},(R_l)_{l\in\omega}$ defined before Lemma \ref{conv}. Then 
$(\mathcal{C}^+_{\bf d},\mathbb{O}_{f_{\bf d}})\equiv_c^i(\mathcal{C}^+_{{\bf d}'},\mathbb{O}_{f_{{\bf d}'}})$ if and only if 
$f_{\bf d},f_{{\bf d}'}$ are conjugate.\end{them}

\vfill\eject

 A straightforward modification of the proof of Proposition \ref{LZ1}(b) gives the following result. We set 
$\mathbb{O}_1\! :=\!\textup{Graph}(f_0)\!\setminus\!\{ (0^\infty ,0^\infty )\}$. 

\begin{propo} \label{LZ2} 
$\{ (\mathbb{X}_1,\mathbb{O}_1),(\mathbb{X}_1,\mathbb{O}_1^{-1}),(\mathbb{X}_1,\mathbb{R}_1)\}$ is a $\preceq^i_c$-antichain basis for the class of digraphs on a 0DMS space with uncountable CCN.\end{propo}

 Our version of Proposition \ref{LZ1+} for digraphs and oriented graphs is as follows.

\begin{propo} \label{LZ1+or} $(\mathbb{X}_1,\mathbb{R}_1)$ is $\preceq^i_c$-minimal, but not $\preceq_c$-minimal, in the class of digraphs (or oriented graphs) on a 0DMC space with CCN at least three.\end{propo}

\noindent\emph{Proof.}\ We follow the proof of Proposition \ref{LZ1+} since $s(\mathbb{R}_1)\! =\! G_{f_1}$. For $\preceq_c$, we work with the oriented graph $\{ (\varepsilon^{2p+1}(\varepsilon^+)^\infty ,
(\varepsilon^+)^{2p+2}\big( (\varepsilon^+)^+\big)^\infty )\mid\varepsilon\!\in\! 3\wedge p\!\in\!\omega\}$.\hfill{$\square$}

\subsection{$\!\!\!\!\!\!$ Digraphs induced by a partial function}

\emph{Notation.}\ If $f\! :\!\mbox{Domain}(f)\!\subseteq\! X\!\rightarrow\!\mbox{Range}(f)\!\subseteq\! X$ is a partial function, then the digraph \emph{induced} by $f$ is $D_f\! :=\!\textup{Graph}(f)\!\setminus\!\Delta (X)$. Note that $G_f\! =\! s(D_f)$, which gives versions of Proposition \ref{fp} and Corollary \ref{corfp} with $D_f$ instead of $G_f$ in a straightforward way.\medskip
 
 We also have the following versions of Theorem \ref{c2Polish} for digraphs and oriented graphs. 
 
\begin{them} \label{c2Polishor} There is no $\preceq^i_c$-antichain basis for the class of digraphs (or oriented graphs) induced by a partial homeomorphism on a 0DMS (or 0DP) space with CCN at least three. In fact, we can even restrict this class to the case where the spaces are countable Polish and the functions are fixed point free with open domain.\end{them}

\noindent\emph{Proof.}\ We follow the proof of Theorem \ref{c2Polish}. We restrict $f_\delta$ to 
$$D_\delta\! :=\!\big\{ x\!\in\!\mbox{proj}[\mathbb{G}_\delta ]\mid\exists n\!\in\!\omega ~~
x(n)\! =\! a\wedge x\vert n\!\in\! (\omega\cup\{ c\} )^n\big\}\mbox{,}$$ 
so that $\mathbb{O}_\delta\! =\!\textup{Graph}{(f_\delta}_{\vert D_\delta})$ is an oriented graph. We then work with the 
$(\mathbb{P}_\delta ,\mathbb{O}_\delta )$'s  since $\mathbb{G}_\delta\! :=\! s(\mathbb{O}_\delta )$, applying Theorem \ref{absminor}.\hfill{$\square$}

\subsection{$\!\!\!\!\!\!$ Digraphs induced by a total function}\indent 

 Similarly, the versions of Proposition \ref{invol} and Lemma \ref{gendense} for the $D_f$'s are  straightforward. For Lemma \ref{gendense}, we just assume that $E\!\subseteq\! D_f$. We now give a motivating result.
 
\begin{lemm} \label{flipor} Let $X,Y$ be 0DMC spaces, and 
$f\! :\! X\!\rightarrow\! X$, $g\! :\! Y\!\rightarrow\! Y$ be homeomorphisms, $g$ being minimal. Then\smallskip

\noindent (a) $\textup{Graph}(g)$ is an oriented graph on $Y$ if $Y$ has cardinality at least three,\smallskip

\noindent (b) $f,g$ are conjugate with witness $\varphi$ if and only if 
$\big( X,\textup{Graph}(f)\big)\preceq_c^i\big( Y,\textup{Graph}(g)\big)$ with witness $\varphi$.\end{lemm}

\noindent\emph{Proof.}\ (a) As $Y$ has cardinality at least three and $g$ is minimal, $g$ and $g^2$ are fixed point free.\medskip

\noindent (b) The proof is similar to and simpler than the proof of Lemma \ref{flipo}.\hfill{$\square$}\medskip

 The version of Lemma \ref{mino} for the $D_f$'s is straightforward.

\begin{lemm} \label{minor} Let $X$ be a 0DMC space, and $f\! :\! X\!\rightarrow\! X$ be a minimal homeomorphism with $\chi_c\big( X,\textup{Graph}(f)\big)\!\geq\! 3$. Then $\big( X,\textup{Graph}(f)\big)$ is 
$\preceq^i_c$-minimal in the class of closed digraphs (or oriented graphs) on a 0DMC space with CCN at least three. This is the case of $\big(\mathcal{C},\textup{Graph}(o)\big)$ if 
${\bf d}\! =\! (d_j)_{j\in\omega}\!\in\!\mathcal{O}$.\end{lemm}

 This gives a version of Theorem \ref{eantichmino} for digraphs and oriented graphs.
 
\begin{them} \label{eantichminori} There is a $\preceq_c$-antichain (and thus $\preceq_c^i$-antichain) 
$\Big(\big(\mathcal{C}_\alpha ,\textup{Graph}(f_\alpha )\big)\Big)_{\alpha\in 2^\omega}$, where\smallskip
 
\noindent (a) $\mathcal{C}_\alpha$ is homeomorphic to $2^\omega$,\smallskip

\noindent (b) $f_\alpha$ is a minimal homeomorphism of $\mathcal{C}_\alpha$, and $\textup{Graph}(f_\alpha )$ has CCN three,\smallskip

\noindent (c) $\big(\mathcal{C}_\alpha ,\textup{Graph}(f_\alpha )\big)$ is $\preceq_c^i$-minimal in the class of closed digraphs (or oriented graphs) on a 0DMC space with CCN at least three.\smallskip

 In particular, any $\preceq_c^i$-basis for the class of digraphs (or oriented graphs) induced by a homeomorphism of a 0DMC space with CCN at least three must have size continuum.\end{them}

 We consider, for $\kappa\!\leq\! 3$, the class $\mathfrak{G}^o_\kappa$ of digraphs $D_f$ induced by a homeomorphism $f$ of a 0DMC space with CCN strictly bigger than $\kappa$. Replacing $G_o$ with 
$D_o\! =\!\textup{Graph}(o)$, we get the version of Proposition \ref{me} for the $D_f$'s in a straightforward way, using Lemma \ref{flipor}. It is worth noting that the version of Theorem \ref{Ckappa} for the $D_f$'s is different from the one for the graphs.

\begin{them} \label{Ckappaor} (a) $(1,\emptyset )$ is $\preceq^i_c$-minimum in $\mathfrak{G}^o_0$.\smallskip

\noindent (b) Any $\preceq_c^i$-basis for $\mathfrak{G}^o_1$ must have size continuum.\smallskip

\noindent (c) Any $\preceq_c^i$-basis for $\mathfrak{G}^o_2$ must have size continuum.\smallskip

Moreover, the $(\mathfrak{G}^o_\kappa ,\preceq^i_c)$'s and the $(\mathfrak{G}^o_\kappa ,\preceq_c)$'s are not well-founded. They also  contain antichains of size continuum (except maybe for $\preceq_c$ when $\kappa\! =\! 3$).\end{them}

\noindent\emph{Proof.}\ (a) See the proof of Theorem \ref{Ckappa}(a).\medskip

\noindent (b) The situation here is very different from the one for $\mathfrak{G}_1$. Let $K$ be a 0DMC space of cardinality at least two, and $f$ be a minimal homeomorphism of $K$. Note that $f$ is fixed point free, so  that ${D_f\! =\!\textup{Graph}(f)}$ is not empty and thus has CCN strictly bigger than $1$. Assume that $L$ is a 0DMC space, $g$ is a homeomorphism of $L$ such that $D_g$ has CCN strictly bigger than $1$, and ${(L,D_g)\preceq^i_c(K,D_f)}$. We will see that 
$(K,D_f)\preceq^i_c(L,D_g)$, which will prove the $\preceq^i_c$-minimality of $(K,D_f)$ in $\mathfrak{G}^o_1$. If the set of fixed points of $g$ is not open, then $(L,D_g)$ has uncountable CCN by the version of Proposition \ref{fp} for the $D_f$'s, and thus $(K,D_f)$ and $(K,G_f)$ too, which contradicts Theorem \ref{homeocompact}. Thus the set $F_1$ of fixed points of $g$ is open, which by the version of Corollary \ref{corfp} for the $D_f$'s implies that 
$\chi_c(L\!\setminus\! F_1,D_g\cap (L\!\setminus\! F_1)^2)\! =\!\chi_c(L,D_f)\! >\! 1$. This implies that we may assume that $g$ is fixed point free and thus $D_g\! =\!\textup{Graph}(g)$. Lemma \ref{flipor} then implies that $g,f$ are conjugate and thus 
$(K,D_f)\preceq^i_c(L,D_g)$. Theorem \ref{eantichminori} now provides a $\preceq^i_c$-antichain of size continuum made up of minimal elements in $\mathfrak{G}^o_1$, which gives the result.\medskip

\noindent (c) By Theorem \ref{eantichminori}, any $\preceq_c^i$-basis for $\mathfrak{G}^o_2$ must have size continuum.\medskip

 We then argue as in the proof of Theorem \ref{Ckappa}, using the fact that $G_f\! =\! s(D_f)$.\hfill{$\square$}
 
\vfill\eject

\noindent\emph{Remark.}\ If $(K,D_f)$ is in $\mathfrak{G}^o_1$ and the set of fixed points of $f$ is open, then $(K,D_f)$ is 
$\preceq^i_c$-above a similar element $(L,D_g)$ with $g$ fixed point free by the proof of Theorem \ref{Ckappaor}(b). [W, Theorem 5.2] gives a compact subset $M$ of $L$ such that $g[M]\! =\! M\!\not=\!\emptyset$ and the dynamical system 
$(M,g_{\vert M})$ is minimal. As $g$ is fixed point free, $M$ has cardinality at least two. The proof of Theorem \ref{Ckappaor}(b) shows that $(M,D_{g_{\vert M}})$ is $\preceq^i_c$-minimal in $\mathfrak{G}^o_1$, and is $\preceq^i_c$-below $(K,D_f)$. So the elements of $\mathfrak{G}^o_1$ induced by a minimal homeomorphism form a 
$\preceq^i_c$-basis for the subclass of $\mathfrak{G}^o_1$ whose elements are induced by a homeomorphism with an open set of fixed points. So in order to get an interesting basis, we need to understand the elements $(K,D_f)$ of $\mathfrak{G}^o_1$ whose set of fixed points is not open. In such a case the CCN is $2^{\aleph_0}$, by the version of Proposition \ref{fp} for the $D_f$'s. Proposition \ref{LZ2} implies that $(\mathbb{X}_1,\mathbb{R}_1)\preceq^i_c(K,D_f)$. The problem is that $\mathbb{R}_1$ is not of the form $D_g$.\medskip

 The versions of Theorems \ref{desc}, \ref{decfib}, \ref{fibgenr} and Proposition \ref{odoshift} for the $D_f$'s are direct.

\begin{coro} \label{controtor} There is a $\preceq^i_c$-antichain 
$\big(\Sigma_r,\textup{Graph}(\sigma_{\vert\Sigma_r})\big)_{r\in\mathbb{R}}$ of size continuum, where\smallskip 

\noindent (a) $\Sigma_r$ is a two-sided subshift homeomorphic to $2^\omega$,\smallskip

\noindent (b) $\sigma_{\vert\Sigma_r}$ is a minimal homeomorphism of $\Sigma_r$, $\textup{Graph}(\sigma_{\vert\Sigma_r})$ is an oriented graph with CCN three,\smallskip

\noindent (c) $\big(\Sigma_r,\textup{Graph}(\sigma_{\vert\Sigma_r})\big)_{r\in\mathbb{R}}$ is $\preceq^i_c$-minimal in the class of closed digraphs (or oriented graphs) on a 0DMC space with CCN at least three.\end{coro} 

 The version of Theorem \ref{absmincompmain} for digraphs, $\mathfrak{G}^o_2$ and the $D_f$'s is straightforward. For oriented graphs, we modify $K_0$ and $h_0$  (note that 
$\big(\sigma^\varepsilon (\alpha_0),\sigma^{1-\varepsilon}(\alpha_0)\big)\!\in\! D_{h_0}$ for each $\varepsilon\!\in\! 2$, so that $D_{h_0}$ is not an oriented graph). In order to get a version of Theorem \ref{absmincompmain} for oriented graphs, we set 
$\alpha^o_0\! :=\! (0123)^\infty\!\cdot\! (0123)^\infty$, $\beta^o_0\! :=\! (0123)^\infty\!\cdot\! 4(0123)^\infty$, 
$K^o_0\! :=\!\mbox{Orb}_\sigma (\alpha^o_0)\cup\mbox{Orb}_\sigma (\beta^o_0)$ and $h^o_0\! :=\!\sigma_{\vert K^o_0}$. The version of Lemma \ref{topol} is as follows.
 
\begin{lemm} \label{topol+} Let $S$ be a 0DMS (resp., 0DP) space, $f$ be a homeomorphism of $S$ with the properties that 
$\chi_c(S,D_f)\!\geq\! 3$ and $(S,D_f)\preceq^i_c(K^o_0,D_{h^o_0})$. Then there is a finer 0DMS (resp., 0DP) topology $\tau$ in $\mathcal{T}$ with the property that $\big( (K^o_0,\tau ),D_{h^o_0}\big)\preceq^i_c(S,D_f)$.\end{lemm}
 
\noindent\emph{Proof.}\ The argument is a slight variation of that in the proof of Lemma \ref{topol}. For instance, as 
$\mbox{Orb}_{h^o_0} (\beta^o_0)$ is discrete, there is $\varepsilon\!\in\! 4$ with $\sigma^\varepsilon (\alpha^o_0)\!\in\! V$, which gives $x\!\in\! S$ with ${\varphi (x)\! =\!\sigma^\varepsilon (\alpha^o_0)}$. As $f$ is fixed point free, $f(x)\!\not=\! x$, which implies that $\big( x,f(x)\big)\!\in\! D_f$, $\Big(\varphi (x),\varphi\big( f(x)\big)\Big)\!\in\! D_{h^o_0}$, and 
$\sigma^{\varepsilon +1\text{ mod }4} (\alpha^o_0)\! =\!\varphi\big( f(x)\big)\!\in\! V$. Iterating this argument, we see that 
$\mbox{Orb}_\sigma (\alpha^o_0)\!\subseteq\! V$ and 
$\big\{\big(\sigma^\varepsilon (\alpha^o_0),\sigma^{\varepsilon +1\text{ mod }4}(\alpha^o_0)\big)\mid
\varepsilon\!\in\! 4\big\}\!\subseteq\! E$.\hfill{$\square$}\medskip

\noindent\emph{Notation.}\ We set, for $A\!\subseteq\!\omega$,\medskip

\leftline{$\mathcal{B}^\tau_A\! :=\!\mathcal{B}^\tau\cup\bigcup_{\varepsilon\in 4}~\Big\{ C\cap(\{\sigma^\varepsilon (\alpha_0)\}\cup\bigcup_{n\in\bigcap_{-p_0\leq r\leq q_0}~(N_A+r)\cap\omega}~\{\sigma^{\varepsilon +4n+1}(\beta_0)\} ~\cup$}\smallskip

\rightline{$\bigcup_{n\in\bigcap_{-p_1\leq r\leq q_1}~(N_A+r)\cap\omega}~
\{\sigma^{\varepsilon -4n-2}(\beta_0)\} )\mid C\!\in\!\mathcal{B}^\tau\wedge p_0,q_0,p_1,q_1\!\in\!\omega\Big\} .$}\medskip

\noindent The version of Lemma \ref{tauA}, that of Lemma \ref{compari}, as well as the rest of the proof of Theorem \ref{absmincompmain} for oriented graphs, are then straightforward. The version of Theorem \ref{CB1intro} for digraphs and 
$\mathfrak{G}^o_2$ is straightforward.

\subsection{$\!\!\!\!\!\!$ Equivalence relations}\indent
   
 We define a map $d\! :\!\mathbb{M}\!\rightarrow\!\mathcal{S}_m$ by $d(f)\! :=\!\big( 2^\omega ,\textup{Graph}(f)\big)$. Applying Lemma \ref{flipor}, we get the following result.
 
\begin{them} \label{redborpor} The map $d$ reduces continuously $CO$ to $\equiv^i_c$. Moreover, the vertices of the digraph $d(f)$ have degree one, for each $f\!\in\!\mathbb{M}$.\end{them}
 
 We then define a map $\mathcal{D}\! :\!\mathbb{M}\!\rightarrow\!\mathcal{S}_g$ by 
$\mathcal{D}(f)\! :=\!\big(\overline{\mbox{proj}\big[ (i\!\times\! i)[\mathbb{G}_f]\big]},(i\!\times\! i)[\mathbb{O}_f]\big)$. Applying Theorem \ref{corflipor}, we get the following result.
 
\begin{them} \label{redborbisor} The map $\mathcal{D}$ Borel reduces $CO$ to $\equiv^i_c$. Moreover, the vertices of the digraph $\mathcal{D}(f)$ have degree at most one, for each $f\!\in\!\mathbb{M}$.\end{them}

\section{$\!\!\!\!\!\!$ A summary for future work}\indent

 We summarize a number of our results in the following table, which leaves open questions about graphs on a 0DMS space with CCN at least three.\bigskip
 
\centerline{\scalebox{0.66}{$$\begin{tabular}{|c|c|c|c|}
  \hline
  & finite & metrizable compact & Polish or metrizable separable \\
  \hline
  $\preceq^i_c$ 
  & $\begin{array}{ll} 
  & \mbox{(1) concrete antichain basis of size }\aleph_0\cr 
  & \mbox{(2) concrete basis of size }\aleph_0\cr 
  & \mbox{(3) any basis is infinite}\cr 
  & \mbox{(4) antichain of size }\aleph_0\mbox{ made up of minimals}\cr 
  & \mbox{(5) no infinite descending chain}\cr
  & \mbox{(6) minimal elements}\cr
  & \mbox{(7) embed }\subseteq\mbox{ on }\mathcal{P}_{<\infty}(\omega )
  \end{array}$ 
  & $\begin{array}{ll} 
  & \mbox{(1) no antichain basis}\cr 
  & \mbox{(2) concrete basis of size }2^{\aleph_0}\cr 
  & \mbox{(3) any basis has size at least }2^{\aleph_0}\cr 
  & \mbox{(4) antichain of size }2^{\aleph_0}\mbox{ made up of minimals}\cr 
  & \mbox{(5) infinite descending chain}\cr 
  & \mbox{(6) minimal elements}\cr
  & \mbox{(7) embed }\subseteq\mbox{ on }\mathcal{P}(\omega )\end{array}$ 
  & $\begin{array}{ll} 
  & \mbox{(1) no antichain basis}\cr 
  & \mbox{(3) any basis has size at least }2^{\aleph_0}\cr
  & \mbox{(4) antichain of size }2^{\aleph_0}\cr 
  & \mbox{(5) infinite descending chain}\cr 
  & \mbox{(6) minimal elements}\cr
  & \mbox{(7) embed }\subseteq\mbox{ on }\mathcal{P}(\omega )
  \end{array}$\\
  \hline
  $\preceq_c$ 
  & $\begin{array}{ll} 
  & \mbox{(1) no antichain basis}\cr 
  & \mbox{(2) concrete basis of size }\aleph_0\cr 
  & \mbox{(3) any basis is infinite}\cr 
  & \mbox{(4) antichain of size }\aleph_0\cr 
  & \mbox{(5) infinite descending chain}\cr 
  & \mbox{(6) no minimal element}
  \end{array}$ 
  & $\begin{array}{ll} 
  & \mbox{(2) concrete basis of size }2^{\aleph_0}\cr 
  & \mbox{(4) antichain of size }2^{\aleph_0}\cr 
  & \mbox{(5) infinite descending chain}\cr 
  & \mbox{(7) embed }\subseteq\mbox{ on }\mathcal{P}(\omega )
  \end{array}$ 
  & $\begin{array}{ll} 
  & \mbox{(4) antichain of size }2^{\aleph_0}\cr 
  & \mbox{(5) infinite descending chain}\cr 
  & \mbox{(7) embed }\subseteq\mbox{ on }\mathcal{P}(\omega )
  \end{array}$\\
  \hline
\end{tabular}$$}}\bigskip

 It is remarkable that the properties in the last two columns are the same for graphs induced by a partial homeomorphism with countable domain, (possibly) up to (2)-$\preceq^i_c$ in the compact case. For graphs induced by a total homeomorphism, (4) and (5) hold, as well as (3), (6) and (7)-$\preceq^i_c$, and (1)-$\preceq^i_c$ in the case of spaces which are not compact. All these results admit versions for digraphs and oriented graphs. 

\vfill\eject

\section{$\!\!\!\!\!\!$ References}

\noindent [A-D-H]\ \ A. S. Asratian, T. M. J. Denley and R. H\"aggkvist,~\it Bipartite graphs and their applications,~\rm Cambridge Tracts in Mathematics, 131, Cambridge University Press, Cambridge, 1998, xii+259 pp

\noindent [Ba]\ \ S. Banach, Un th\'eor\`eme sur les transformations biunivoques,\ \it Fund. Math.\rm\ 6 (1924), 236-239

\noindent [B0]\ \ A. Bernshteyn, Distributed Algorithms, the Lov\'asz Local Lemma, and Descriptive Combinatorics,\ \it arXiv:2004.04905\ \rm

\noindent [B1]\ \ A. Bernshteyn, Probabilistic constructions in continuous combinatorics and a bridge to distributed algorithms,\ \it arXiv:2102.08797\ \rm

\noindent [Bo]\ \ N. Bourbaki,~\it Topologie g\'en\'erale, ch. 1-4,~\rm Hermann, 1974

\noindent [C]\ \ R. Carroy, A quasi-order on continuous functions,,\ \it J. Symbolic Logic\ \rm 78, 2 (2013), 633-648

\noindent [C-M-Sc-V1]\ \ R. Carroy, B. D. Miller, D. Schrittesser and Z. Vidny\'anszky, Minimal definable graphs of definable chromatic number at least three,\ \it Forum of Mathematics Sigma\ \rm 9 (2021), e7, 1-16

\noindent [C-M-Sc-V2]\ \ R. Carroy, B. D. Miller, D. Schrittesser and Z. Vidny\'anszky, Minimal definable graphs with no definable two-colorings,\ \it http://www.logic.univie.ac.at/zoltan.vidnyanszky/summary2.pdf\rm\ \\ (2018)

\noindent [C-M-So]\ \ R. Carroy, B. D. Miller and D. T. Soukup, The open dihypergraph dichotomy and the second level of the Borel hierarchy,\ \it Contemporary Mathematics\ \rm 752 (2020), 1-19

\noindent [Ce-Da-To-Wy]\ \ D. Cenzer, A. Dashti, F. Toska and S. Wyman, Computability of countable subshifts in one dimension,\ \it Theory Comput. Syst.\rm\ 51 (2012), 352-371

\noindent [Co-M]\ \ C. T. Conley and B. D. Miller, An antibasis result for graphs of infinite Borel chromatic number,\ \it Proc. Amer. Math. Soc.\rm 142, 6 (2014) 2123-2133

\noindent [E]\ \ R. L. Ellis, Extending continuous functions on zero-dimensional spaces,\ \it Math. Ann.\ \rm 186 (1970), 114-122

\noindent [G]\ \ S. Gao,~\it Invariant Descriptive Set Theory,~\rm Pure and Applied Mathematics, A Series of Monographs and Textbooks, 293, Taylor and Francis Group, 2009

\noindent [G-J-Kr-Se]\ \ S. Gao, S. C. Jackson, E. W. Krohne and B. Seward, Continuous Combinatorics of Abelian Group Actions,\ \it arXiv:1803.03872\ \rm

\noindent [He-N]\ \ P. Hell and J. Ne\v set\v ril,~\it Graphs and homomorphisms,~\rm Oxford Lecture Series in Mathematics and its Applications, 28, Oxford University Press, Oxford, 2004, xii+244 pp

\noindent [I-Me]\ \ T. Ibarluc\' ia and J. Melleray, Full groups of minimal homeomorphisms and Baire category methods,\ \it Ergodic Theory Dynam. Systems\rm\ 36, 2 (2016), 550-573

\noindent [Kar]\ \ J. Karhum\"aki, On cube-free $\omega$-words generated by binary morphisms,\ \it Discrete Appl. Math. \rm 5, 3 (1983), 279-297

\noindent [Ka]\ \ B. Kaya, The complexity of topological conjugacy of pointed Cantor minimal systems,\ \it  Arch. Math. Logic\rm\ 56, 3-4 (2017), 215-235

\noindent [K]\ \ A. S. Kechris,~\it Classical Descriptive Set Theory,~\rm Springer-Verlag, 1995

\noindent [K-Ma]\ \ A. S. Kechris and A. S. Marks,~\it Descriptive Graph Combinatorics,~\rm preprint, 2020 (see the first author's web page at http://www.math.caltech.edu/~kechris/papers/combinatorics20book.pdf)

\noindent [K-S-T]\ \ A. S. Kechris, S. Solecki and S. Todor\v cevi\' c, Borel chromatic numbers,\ \it Adv. Math.\rm\ 141, 1 (1999), 1-44

\noindent [Ki-Kat-Pa]\ \ I. S. Kim, H. Kato and J. J. Park, On the countable compacta and expansive homeomorphisms,\ \it Bull. Korean Math. Soc.\rm\ 36, 2 (1999), 403-409

\noindent [Kn-R]\ \ B. Knaster and M. Reichbach, Notion d'homog\'en\'eit\'e et prolongements des hom\'eomorphies,\ \it Fund. Math.\rm\ 40 (1953), 180-193

\noindent [Kra-St]\ \ A. Krawczyk and J. Steprans, Continuous colorings of closed graphs,\ \it Topology Appl.~\rm 51 (1993), 13-26

\noindent [Ku]\ \ P. K\r{u}rka,~\it Topological and symbolic dynamics,~\rm Cours Sp\'ecialis\'es (Specialized Courses), 11, 
Soci\'et\'e Math\'ematique de France, Paris, 2003

\noindent [L]\ \ D. Lecomte, On minimal non potentially closed subsets of the plane,\ \it Topology Appl.\rm\ 154, 1 (2007), 241-262

\noindent [L-Za]\ \ D. Lecomte and R. Zamora, Injective tests of low complexity in the plane,\ \it Math. Logic Quart.\rm\ 65, 2 (2019) 134-169

\noindent [L-Z1]\ \ D. Lecomte and M. Zelen\'y, Baire-class $\xi$ colorings: the first three levels,\ \it Trans. Amer. Math. Soc.\rm\ 366, 5 (2014), 2345-2373

\noindent [L-Z2]\ \ D. Lecomte and M. Zelen\'y, Analytic digraphs of uncountable Borel chromatic number under injective definable homomorphism,\ \it arXiv:1811.04738\rm\ 

\noindent [Lo]\ \ M. Lothaire,~\it Algebraic combinatorics on words,~\rm Cambridge University Press, 2002

\noindent [MB]\ \ N. Matte Bon, Subshifts with slow complexity and simple groups with the Liouville property,\ \it Geom. Funct. Anal.\rm\ 24, 5 (2014), 1637-1659

\noindent [Me]\ \ J. Melleray, Dynamical simplicies and Borel complexity of orbit equivalence,\ \it Isr. J. Math.\rm\ 236 (2020), 317-344

\noindent [Pe]\ \ Y. Pequignot, Finite versus infinite: an insufficient shift,\ \it Adv. Math.\rm\ 320, 7 (2017), 244-249
 
\noindent [P]\ \ I. F. Putnam,~\it Cantor minimal systems,~\rm University Lecture Series, 70. American Mathematical Society, Providence, RI, 2018. xiii+149 pp

\noindent [Sa-T\"o]\ \ V. Salo and I. T\"orm\"a, Block maps between primitive uniform and Pisot substitutions,\ \it Ergodic Theory Dynam. Systems\rm\ 35, 7 (2015), 2292-2310

\noindent [T-V] S. Todor\v cevi\' c and Z. Vidny\'anszky, A complexity problem for Borel graphs,\ \it Invent. Math.\ \rm 226 (2021), 225-249 

\noindent [W]\ \ P. Walters,~\it Ergodic theory-introductory lectures,~\rm Lecture Notes in Mathematics, Vol. 458, Springer-Verlag, Berlin-New York, 1975, vi+198 pp

\vfill\eject

\centerline{\bf List of symbols}

\noindent $\mathcal{G}_0$~~~~~ 2

\noindent $\Delta (X)$~~~~~ 2

\noindent $(X,R)$, $\chi_B(X,R)$~~~~~ 2

\noindent $\preceq_c$~~~~~ 2

\noindent $\mathbb{L}_0$~~~~~ 2

\noindent $(\mathbb{A}_\xi ,\mathbb{G}_\xi )$~~~~~ 2

\noindent $\preceq^i_c$, $\preceq^i_B$~~~~~ 3

\noindent CCN, $\chi_c(X,R)$~~~~~ 3

\noindent 0DMS, 0DP, 0DMC~~~~~ 3

\noindent $\mathfrak{K}$~~~~~ 3

\noindent $D_2(\bormone )$, $\boraone\oplus\bormone$~~~~~ 4

\noindent $(X,f)$~~~~~ 4

\noindent $\mbox{Orb}_f(x)$~~~~~ 4

\noindent $R^l$, $R^{-1}$, $s(R)$~~~~~ 5

\noindent $\mathcal{N}$~~~~~ 5

\noindent $\mathbb{O}_m$, $\mathbb{G}_m$, $\mathbb{P}$~~~~~ 5

\noindent $G_f$~~~~~ 6

\noindent $\mathfrak{G}_\kappa$~~~~~ 6

\noindent $\mathbb{X}_1$, $\mathbb{R}_1$, $f_1$~~~~~ 6

\noindent $\mathcal{K}(X)$~~~~~ 7

\noindent $\mathcal{H}(2^\omega )$, $\mathcal{P}$, $\mathcal{O}_\kappa$, $\mathcal{O}^{\aleph_0}_2$~~~~~ 7

\noindent $\sigma$, $(01)^\infty\!\cdot\! (01)^\infty$, $K_0$, $h_0$~~~~~ 7

\noindent $\mathbb{M}$, $FCO$, $CO$~~~~~ 9

\noindent $\mathcal{S}_m$~~~~~ 9

\noindent $\leq_B$~~~~~ 9

\noindent $=^+$~~~~~ 9

\noindent $(2p\! +\! 3,C_{2p+3})$~~~~~ 10

\noindent $\mathcal{S}$, $(\lambda_l)_{l\in\omega}$, $\big(s_l(i)\big)_{i<\lambda_l}$, $\mathcal{I}$, 
$\mathbb{G}_\gamma$, $\mathbb{K}_\gamma$~~~~~ 11

\noindent $\mathfrak{C}$, $\bf d$, $\pi_{j\in S}~d_j$, $\prod_{j\in S}~d_j$, $\prod_l$, $\mathcal{C}$, $\mathcal{C}_{\bf d}$, 
$N_s$~~~~~ 11

\noindent $R_n$, ${}^nR$~~~~~ 11

\noindent $\mathcal{J}$, $\mathcal{J}^c$, $c$, $a$, $\overline{a}$, $\mathcal{K}_{\bf d}$, $\mathbb{O}_\beta$~~~~~ 14
 
\noindent $\mathbb{G}_\beta$, $\mathbb{K}_\beta$~~~~~ 15

\noindent $\mathbb{O}_\delta$, $\mathbb{G}_\delta$, $\mathbb{P}_\delta$, $\mathbb{P}_\infty$~~~~~ 20

\noindent $\mathbb{T}$~~~~~ 22

\noindent $(n_l)_{l\in\omega}$, $(n^{\bf d}_l)_{l\in\omega}$, $(L_l)_{l\in\omega}$, $(R_l)_{l\in\omega}$, $f_{l,i}$, $f^{\bf d}_{l,i}$, 
$\mathbb{O}_f$, $\mathbb{G}_f$, $\mathcal{C}^+$~~~~~ 26

\noindent $\zeta$~~~~~ 28

\noindent $\equiv^i_c$~~~~~ 29

\noindent $o$, $o_{\bf d}$, $\mbox{Orb}^+_f(x)$~~~~~ 29

\noindent $\mathfrak{D}$, $\mu$, $\mathbb{O}_o$, $\mathbb{G}_o$~~~~~ 30

\noindent $F_1$, $F_1^f$~~~~~ 36

\noindent $f_0$~~~~~ 37

\noindent $r_C$, $r'_C$~~~~~ 39

\noindent $\mathcal{O}$~~~~~ 45

\noindent $\mathfrak{G}_2^e$~~~~~ 46

\noindent $i_l$, $S_A$, $D_S$~~~~~ 47

\noindent $\mathbb{G}$~~~~~ 48

\noindent $\mathcal{G}_\alpha$~~~~~ 49

\noindent $[w]_p, w^\mathbb{Z}$~~~~~ 50

\noindent $\Sigma$, $\Sigma_F$~~~~~ 51

\noindent $\sqsubseteq$, $\alpha^+$~~~~~ 51

\noindent $R_r$, $\phi_r$, $\Sigma^2_r$~~~~~ 52

\noindent $\mathcal{L}(\Sigma )$, $\mathcal{L}_n(\Sigma )$, $\mathcal{L}^r$, $\mathcal{L}^r_n$~~~~~ 52

\noindent $X'$, $X^\alpha$~~~~~ 55

\noindent $\alpha_0$~~~~~ 55

\noindent $\mathcal{T}$~~~~~ 63

\noindent $(S_q)_{q\in\omega}$, $(r^q_j)_{j\in\omega}$, $j_l$, $N_A$, $\mathcal{B}^\tau$, $\mathcal{B}^\tau_A$~~~~~ 64

\noindent $t^\tau_A$~~~~~ 65

\noindent $\mathfrak{H}_\kappa$~~~~~ 67

\noindent $(h_s)_{s\in\omega^{<\omega}}$~~~~~ 69

\noindent $\mathbb{Q}$, $Q$, $\mathcal{S}_g$~~~~~ 72

\noindent $D_f$~~~~~ 76

\noindent $\mathfrak{G}^o_\kappa$~~~~~ 77

\noindent $K^o_0$, $h^o_0$~~~~~ 78

\vfill\eject

\centerline{\bf Index}

\noindent alphabet~~~~~ 50

\noindent antichain~~~~~ 3

\noindent basis~~~~~ 3

\noindent Borel chromatic number~~~~~ 2

\noindent Cantor-Bendixson derivative, iterated Cantor-Bendixson derivative, Cantor-Bendixson rank~~~~~ 55

\noindent Cantor dynamical system~~~~~ 4

\noindent coloring~~~~~ 2

\noindent conjugate~~~~~ 4

\noindent continuous chromatic number~~~~~ 3

\noindent continuous tuple~~~~~ 26

\noindent diagonal~~~~~ 2

\noindent digraph~~~~~ 2

\noindent digraph induced by a function~~~~~ 76

\noindent dynamical system~~~~~4

\noindent equicontinuous dynamical system~~~~~ 46

\noindent expansive dynamical system~~~~~ 9

\noindent flip-conjugate~~~~~ 4

\noindent graph~~~~~ 2

\noindent graph induced by a function~~~~~ 6

\noindent homomorphism~~~~~ 2

\noindent minimal dynamical system~~~~~ 4

\noindent minimum element for a quasi-order~~~~~ 3

\noindent odometer~~~~~ 29

\noindent orbit-equivalent~~~~~ 4

\noindent oriented graph~~~~~ 3

\noindent periodic point~~~~~ 50

\noindent quasi-order~~~~~ 3

\noindent rank of a point~~~~~ 59

\noindent shift map~~~~~ 50

\noindent two-sided subshift~~~~~ 51

\noindent substitution~~~~~ 51

\noindent uniformly recurrent subshift~~~~~ 52

\noindent walk, odd walk, closed walk, cycle~~~~~ 10

\noindent well-quasi-order~~~~~ 4

\end{document}